\documentclass[12pt]{amsart}
\usepackage{amsmath,amsthm,amscd,amssymb,amsfonts}
\usepackage{fullpage,colonequals}
\usepackage{mathrsfs, upgreek, mathabx}
\usepackage{stmaryrd}
\usepackage{tikz}
\usepackage[all,cmtip]{xy}
\usepackage{xcolor}
\usepackage{graphicx, color}
\usepackage{mathtools}
\usepackage[hidelinks]{hyperref}
\hypersetup{
  colorlinks,
  linkcolor={red!85!black},
  citecolor={blue!85!black},
  urlcolor={blue!95!black}
}

\usepackage{mathrsfs,upgreek}

\numberwithin{equation}{section} 
\numberwithin{figure}{section} 
\mathtoolsset{showonlyrefs} 
\usepackage{enumitem} 
\setlist[enumerate]{label=$(\roman*)$, ref=$(\roman*)$}

\theoremstyle{plain}
\newtheorem{theoalph}{Theorem}

\newtheorem{theorem}{Theorem}[section]
\newtheorem{coro}[theorem]{Corollary}
\newtheorem{proposition}[theorem]{Proposition}
\newtheorem{lemma}[theorem]{Lemma}

\newtheorem{innercustomthm}{Theorem}
\newenvironment{custtheo}[1]
{\renewcommand\theinnercustomthm{#1}\innercustomthm}
{\endinnercustomthm}

\theoremstyle{definition}

\theoremstyle{remark}
\newtheorem{remark}[theorem]{Remark}

\newtheoremstyle{citing}
{3pt}
{3pt}
{\itshape}
{}
{\bfseries}
{.}
{.5em}
{\thmnote{#3}}

\theoremstyle{citing}

\newcommand{\A}{\mathbb{A}}

\newcommand{\C}{\mathbb{C}}

\newcommand{\F}{\mathbb{F}}
\newcommand{\K}{\mathbb{K}}
\newcommand{\Q}{\mathbb{Q}}
\newcommand{\R}{\mathbb{R}}

\newcommand{\Z}{\mathbb{Z}}

\renewcommand{\H}{\mathbb{H}}

\newcommand{\cB}{\mathcal{B}}
\newcommand{\cC}{\mathcal{C}}

\newcommand{\cE}{\mathcal{E}}
\newcommand{\cF}{\mathcal{F}}

\newcommand{\cK}{\mathcal{K}}

\newcommand{\cM}{\mathcal{M}}

\newcommand{\cO}{\mathcal{O}}

\newcommand{\cQ}{\mathcal{Q}}

\newcommand{\fD}{\mathfrak{D}}
\newcommand{\fE}{\mathfrak{E}}

\newcommand{\sC}{\mathscr{C}}

\newcommand{\sQ}{\mathscr{Q}}

\newcommand{\sT}{\mathscr{T}}

\newcommand{\hC}{\widehat{C}}

\newcommand{\hF}{\widehat{F}}

\newcommand{\halpha}{\widehat{\alpha}}

\newcommand{\hmu}{\widehat{\mu}}
\newcommand{\hnu}{\widehat{\nu}}
\newcommand{\hxi}{\widehat{\xi}}

\newcommand{\hrho}{\widehat{\rho}}

\newcommand{\hsigma}{\widehat{\sigma}}

\newcommand{\hphi}{\widehat{\phi}}

\newcommand{\tE}{\widetilde{E}}

\newcommand{\tgamma}{\widetilde{\gamma}}

\newcommand{\tdelta}{\widetilde{\delta}}

\newcommand{\tphi}{\widetilde{\phi}}
\newcommand{\tvarphi}{\widetilde{\varphi}}
\newcommand{\tPhi}{\widetilde{\Phi}}

\newcommand{\tpsi}{\widetilde{\psi}}

\renewcommand{\=}{\coloneqq}
\renewcommand{\:}{\colon}
\newcommand{\dd}{\hspace{1pt}\operatorname{d}\hspace{-1pt}}

\DeclareMathOperator{\Id}{Id}

\DeclareMathOperator{\dist}{dist}
\DeclareMathOperator{\supp}{supp} 

\newcommand{\bfB}{\mathbf{B}}
\newcommand{\bfD}{\mathbf{D}}
\newcommand{\bfG}{\mathbf{G}}

\newcommand{\bfL}{\mathbf{L}}

\newcommand{\bfR}{\mathbf{R}}

\newcommand{\bfX}{\mathbf{X}}
\newcommand{\bfY}{\mathbf{Y}}

\DeclareMathOperator{\bfx}{\mathbf{x}}
\DeclareMathOperator{\bfzero}{\mathbf{0}}

\DeclareMathOperator{\oO}{O}
\newcommand{\bfone}{\mathbf{1}}

\newcommand{\wtf}{\widetilde{f}}

\newcommand{\wtg}{\widetilde{g}}

\newcommand{\tcF}{\widetilde{\cF}}

\newcommand{\fra}{\mathfrak{a}}
\newcommand{\frb}{\mathfrak{b}}
\newcommand{\frp}{\mathfrak{p}}
\newcommand{\frq}{\mathfrak{q}}

\newcommand{\ssetminus}{\smallsetminus}

\DeclareMathOperator{\Aut}{Aut}
\DeclareMathOperator{\End}{End}
\DeclareMathOperator{\Hom}{Hom}
\DeclareMathOperator{\Iso}{Iso}
\DeclareMathOperator{\Ker}{Ker}

\DeclareMathOperator{\Div}{Div}
\DeclareMathOperator{\Frob}{Frob}

\DeclareMathOperator{\lcm}{lcm}

\DeclareMathOperator{\GL}{GL}
\DeclareMathOperator{\SL}{SL}

\newcommand{\Zp}{\Z_{p}}
\newcommand{\Fp}{\F_p}
\newcommand{\Fpalg}{\overline{\F}_p}

\newcommand{\Qp}{\Q_p}
\newcommand{\Qpalg}{\overline{\Q}_p}
\newcommand{\MQpalg}{\cM_{\Qpalg}}
\newcommand{\OQpalg}{\cO_{\Qpalg}}

\newcommand{\Cp}{\C_p}
\newcommand{\Op}{\cO_p}
\newcommand{\Mp}{\cM_p}
\newcommand{\MK}{\cM_{\cK}}
\newcommand{\OK}{\cO_{\cK}}
\newcommand{\OQ}{\cO_{\cQ}}

\DeclareMathOperator{\ord}{ord}

\DeclareMathOperator{\sups}{sups}
\newcommand{\Ell}{Y}

\newcommand{\Sups}{Y_{\sups}(\Cp)}
\newcommand{\tSups}{Y_{\sups}(\Fpalg)}
\newcommand{\SupsQpalg}{Y_{\sups}(\Qpalg)}

\newcommand{\CM}{CM}
\newcommand{\odelta}{\overline{\delta}}

\DeclareMathOperator{\Berk}{Berk}
\newcommand{\AKber}{\A^1_{\Berk}}
\DeclareMathOperator{\can}{can}
\newcommand{\xcan}{x_{\can}}

\renewcommand{\t}{\mathbf{t}}

\newcommand{\kval}{v_p}

\newcommand{\kproj}{\widehat{v}_p}

\newcommand{\FE}{\cF_{E}}
\newcommand{\FEpr}{\cF_{E'}}
\newcommand{\tFE}{\tcF_{E}}
\newcommand{\tFEpr}{\tcF_{E'}}

\newcommand{\Qua}{\sQ}
\newcommand{\Quap}{{\Qua}_{p}} 

\newcommand{\disc}{\Delta} 
\newcommand{\pd}{\mathfrak{D}} 
\newcommand{\pfd}{\mathfrak{d}} 
\newcommand{\Qpd}{\Qp(\sqrt{\pd})}
\newcommand{\Qpfd}{\Qp(\sqrt{\pfd})}

\newcommand{\dB}{\dist_{\bfB}} 

\newcommand{\Bp}{\bfB_p}
\newcommand{\qval}{\ord_{\bfB}}

\newcommand{\qvalss}{\ord_{\Bss}}
\DeclareMathOperator{\nr}{nr}
\DeclareMathOperator{\tr}{tr}

\newcommand{\rf}{\Bbbk}
\newcommand{\rfk}{\rf_0}
\newcommand{\piem}{\pi_0}

\renewcommand{\ss}{e}
\newcommand{\sspr}{\ss'}
\newcommand{\Bss}{\bfB_{\ss}}

\newcommand{\Dss}{\bfD_{\ss}}
\newcommand{\Dsspr}{\bfD_{\sspr}}
\newcommand{\hDss}{\widehat{\bfD}_{\ss}}
\newcommand{\hDsspr}{\widehat{\bfD}_{\sspr}}
\newcommand{\Fss}{\cF_{\ss}}
\newcommand{\Fsspr}{\cF_{\sspr}}
\newcommand{\Gss}{\bfG_{\ss}}
\newcommand{\Gsspr}{\bfG_{\sspr}}
\newcommand{\Piss}{\Pi_{\ss}}
\newcommand{\Pisspr}{\Pi_{\sspr}}
\newcommand{\Rss}{\bfR_{\ss}}
\newcommand{\Rsspr}{\bfR_{\sspr}}
\newcommand{\Xss}{\bfX_{\ss}}

\DeclareMathOperator{\Var}{Var}
\newcommand{\red}{\mathbf{red}}

\newcommand{\Lss}{\bfL_{\ss}}

\newcommand{\Lfss}{\bfL_{\ss, \operatorname{f}}}

\DeclareMathOperator{\Fix}{Fix}
\newcommand{\Fixss}{\Fix_{\ss}}

\DeclareMathOperator{\Tr}{Tr}
\newcommand{\Trss}{\Tr_{\ss}}
\newcommand{\Trsspr}{\Tr_{\sspr}}
\newcommand{\Uss}{U_{\ss}}

\newcommand{\Nr}{\mathbf{Nr}}
\newcommand{\NE}{\Nr_{E}}
\newcommand{\NEpr}{\Nr_{E'}}
\newcommand{\Npfd}{\Nr_{\pfd}}
\newcommand{\coset}{\mathfrak{N}}
\DeclareMathOperator{\Orb}{Orb}
\newcommand{\corbit}{\Orb_{\coset}(E)}
\newcommand{\corbitc}{\overline{\corbit}}
\DeclareMathOperator{\Ev}{Ev}

\newcommand{\Rsups}{\R_{\sups}}
\newcommand{\vsups}{v^{\sups}}

\DeclareMathOperator{\Cl}{Cl}

\newcommand{\sm}{\mathfrak{j}}

\title[Linnik equidistribution on the supersingular locus]{$p$-Adic distribution of \CM{} points and Hecke orbits.
  \\ \small{II: Linnik equidistribution on the supersingular locus}}

\author{Sebasti\'an Herrero}
\address{Instituto de Matem\'aticas, Pontificia Universidad Cat\'olica de Valpara\'iso, Blanco Viel 596, Cerro Bar\'on, Valpara\'iso,
  Chile.}
\email{sebastian.herrero.m@gmail.com}

\author{Ricardo Menares}

\address{
  Facultad de Matem\'aticas, Pontificia Universidad Cat\'olica de Chile, Vicu\~na Mackenna 4860, Santiago, Chile.}

\email{rmenares@mat.uc.cl}

\author{Juan Rivera-Letelier}
\address{Department of Mathematics, University of Rochester.
  Hylan Building, Rochester, NY~14627, U.S.A.}
\email{riveraletelier@gmail.com}
\urladdr{\url{http://rivera-letelier.org/}}

\usepackage{microtype}
\begin{document}

\begin{abstract}
  For every prime number~$p$, this paper determines the asymptotic distribution of \CM{} points on the moduli space of elliptic curves over~$\C_p$.
  In stark contrast to the complex case, the $p$\nobreakdash-adic setting admits infinitely many distinct accumulation measures.
  The paper identifies all such measures by translating this problem as a $p$\nobreakdash-adic analog of Linnik’s classical problem on the asymptotic distribution of integer points on spheres.
  This translation draws on the relationship between the deformation theories of elliptic curves and formal modules, together with results of Gross and Hopkins.
  The bounds for Fourier coefficients of cuspidal modular forms of Deligne, Iwaniec, and Duke yield a deviation estimate that solves this $p$\nobreakdash-adic Linnik problem.
  This paper also identifies all the accumulation measures of an arbitrary Hecke orbit.
\end{abstract}

\maketitle

\setcounter{tocdepth}{1}
\tableofcontents

\section{Introduction}

For every prime number~$p$, this paper determines the asymptotic distribution of \CM{} points on the moduli space of elliptic curves over~$\Cp$, thereby completing the work initiated in~\cite{HerMenRiv20}.
This is motivated by arithmetic applications in the companion paper~\cite{HerMenRiv24}, and by the results of Linnik, Duke, and Clozel and Ullmo in the complex setting, see~\cite{Lin68,Duk88,CloUll04} and the proceedings article~\cite{MicVen06}.
The paper~\cite{0HerMenRiv2607} applies the results of~\cite{HerMenRiv20,HerMenRiv24} and the current paper to study statistical properties of Hecke correspondences.

The precise statements of the main results rely on the following notation and terminology.
Throughout the rest of this paper, fix a prime number~$p$ and a completion ${(\Cp, | \cdot |_p)}$ of an algebraic closure of the field of $p$\nobreakdash-adic numbers~$\Qp$.
The endomorphism ring of an elliptic curve over~$\Cp$ is isomorphic to~$\Z$ or to an order in a quadratic imaginary extension of~$\Q$.
In the latter case, the order only depends on the class~$E$ in the moduli space~$\Ell(\Cp)$ of elliptic curves over~$\Cp$ and the class~$E$ has \emph{complex multiplication} or is a \emph{\CM{} point}.
The \emph{discriminant of a \CM{} point} is the discriminant of the endomorphism ring of a representative elliptic curve.
A \emph{discriminant} is the discriminant of an order in a quadratic imaginary extension of~$\Q$.
For every discriminant~$D$, the set~$\Lambda_D$ defined by
\begin{displaymath}
  \Lambda_D
  \=
  \{ E \in \Ell(\Cp) \: \text{\CM{} point of discriminant~$D$} \},
\end{displaymath}
is finite and nonempty.
For each~$x$ in~$\Ell(\Cp)$, denote by~$\delta_x$ the Dirac measure on~$\Ell(\Cp)$ at~$x$ and put
\begin{displaymath}
  \odelta_D
  \=
  \frac{1}{\# \Lambda_D} \sum_{E \in \Lambda_D} \delta_E.
\end{displaymath}
Then, ${\odelta_D}$ is a Borel probability measure on~$\Ell(\Cp)$.

This paper identifies all accumulation measures of
\begin{equation}
  \label{eq:1}
  \left\{ \odelta_D \: D \text{ discriminant} \right\}
\end{equation}
with respect to the weak topology on the space of Borel measures on the Berkovich space associated with~$\Ell(\Cp)$.
In stark contrast to the complex case where the limit exists \cite{Duk88,CloUll04}, there are infinitely many different accumulation measures of~\eqref{eq:1}.
The companion paper~\cite{HerMenRiv20} identifies all subsequences of~\eqref{eq:1} converging to the Dirac measure at the ``Gauss'' or ``canonical'' point.
They correspond to the sequences of \CM{} points that are either in the ordinary reduction locus, or that are in the supersingular reduction locus and have discriminants with $p$\nobreakdash-adic norms tending to zero \cite[Theorem~A]{HerMenRiv20}.
This paper treats the remaining case, of sequences of \CM{} points in the supersingular locus whose discriminants have $p$\nobreakdash-adic norm bounded from below by a strictly positive constant.

A key special case is that of a sequence of discriminants~$(D_n)_{n = 1}^{\infty}$ tending to~$-\infty$ such that, for every~$n$, the conductor of~$D_n$ is a $p$\nobreakdash-adic unit and~$\Q(\sqrt{D_n})$ embeds inside a fixed quadratic field extension of~$\Qp$.
The corresponding \CM{} points are naturally related to points in certain quadratic ${\Z}$-lattices embedded in the trace-zero subspace of a $p$\nobreakdash-adic quaternion algebra.
Figuratively, for each~$n$, the set of \CM{} points~$\Lambda_{D_n}$ corresponds to the integer points in the sphere of radius~$|D_n|$ of a three dimensional $p$\nobreakdash-adic space.
Thus, the problem of determining the accumulation measures of~$(\odelta_{D_n})_{n = 1}^{\infty}$ translates to a $p$\nobreakdash-adic version of \mbox{Linnik's} classical problem on the asymptotic distribution of integer points on spheres.
This paper solves this $p$\nobreakdash-adic Linnik problem using a deviation estimate extracted from the bounds for Fourier coefficients of cuspidal modular forms of Deligne~\cite{Del74}, Iwaniec~\cite{Iwa87}, and Duke~\cite{Duk88} (Theorem~\ref{t:p-adic-Linnik} in Section~\ref{s:p-adic-linnik}).
In the key special case, the sequence of measures~$(\odelta_{D_n})_{n = 1}^{\infty}$ converges, except in a peculiar case where there are precisely two accumulation measures (Theorems~\ref{t:CM-symmetric} and~\ref{t:CM-broken-symmetry} in Section~\ref{ss:CM}).
The peculiar case is that of a fixed fundamental discriminant such that~$p$ is the only prime number dividing it and a varying conductor tending to~$\infty$.
Genus theory elucidates the phenomenon, reminiscent of symmetry breaking in physics, behind the emergence of two accumulation measures in this case.
To pass from the key special case to the general case, this paper deduces an analogous equidistribution result for Hecke orbits (Theorem~\ref{t:Hecke-orbits} in Section~\ref{ss:Hecke-orbits}) from the $p$\nobreakdash-adic Linnik equidistribution result mentioned earlier.

The companion paper~\cite{HerMenRiv24} relies on this paper and~\cite{HerMenRiv20} to prove that, for every finite set of prime numbers~$S$, there are at most finitely many singular moduli that are $S$-units.
This extends a result of Habegger in~\cite{Hab15} when ${S = \emptyset}$.

The precise statements of the main results follow.

\subsection{Equidistribution of \CM{} points along a $p$-adic discriminant}
\label{ss:CM}
A \emph{fundamental discriminant} is the discriminant of the ring of integers of a quadratic imaginary extension of~$\Q$.
The \emph{fundamental discriminant} of a discriminant~$D$ is the discriminant~$d$ of~$\Q(\sqrt{D})$.
It divides~$D$ and the quotient~$\frac{D}{d}$ is the square of an integer in~$\Z_{>0}$ that is the \emph{conductor of~$D$}.
A discriminant is \emph{prime}, if it is fundamental and divisible by a unique prime number.
If~$d$ is a prime discriminant divisible by~$p$, then
\begin{displaymath}
  p \equiv -1 \mod 4
  \text{ and }
  d = -p,
  \text{ or }
  p = 2
  \text{ and~$d = -4$ or~$d = -8$.}
\end{displaymath}

A \emph{$p$\nobreakdash-adic quadratic order} is a $\Zp$-order in a quadratic field extension of~$\Qp$, and a \emph{$p$\nobreakdash-adic discriminant} is a set formed by the discriminants of all $\Zp$-bases of a $p$\nobreakdash-adic quadratic order.
Every $p$\nobreakdash-adic discriminant is thus a coset in~$\Qp^{\times} / (\Zp^{\times})^2$ contained in~$\Zp$.
The $p$\nobreakdash-adic discriminant is a complete isomorphism invariant for $p$\nobreakdash-adic quadratic orders (Lemma~\ref{l:p-adic-discriminants-Appendix}$(ii)$ in Appendix~\ref{ss:Apendix-A}).\footnote{
  The notion of $p$\nobreakdash-adic discriminant of a $p$\nobreakdash-adic quadratic order~$\cO$ used here differs from the classical discriminant of a $\Z_p$-lattice~$I$, as defined, \emph{e.g.}, in \cite[Definition~15.2.2]{Voi21}.
  For ${I = \cO}$, the classical discriminant is the~$\Zp$-submodule of~$\Qp$ generated by the $p$\nobreakdash-adic discriminant.
  In contrast to the $p$\nobreakdash-adic discriminant, the classical discriminant is not a complete isomorphism invariant for $p$\nobreakdash-adic quadratic orders.
}

Denote by~$\Qpalg$ the algebraic closure of~$\Qp$ inside~$\Cp$, and by~$\Op$ and~$\OQpalg$ the ring of integers of~$\Cp$ and~$\Qpalg$, respectively.
For each~$E$ in~$\Ell(\Cp)$ represented by a Weierstrass equation with coefficients in~$\OQpalg$ having smooth reduction, denote by~$\FE$ its formal group over~$\Qpalg$ and by~$\End(\FE)$ the ring of endomorphisms of~$\FE$ that are defined over the ring of integers of a finite field extension of~$\Qp$.
Then, $\End(\FE)$ is either isomorphic to~$\Zp$ or to a $p$\nobreakdash-adic quadratic order, see, \emph{e.g.}, \cite[Chapter~IV, Section~1, Theorem~1$(iii)$]{Fro68}.
In the latter case, $E$ has \emph{formal complex multiplication} or is a \emph{formal \CM{} point}.
Moreover, the \emph{$p$\nobreakdash-adic discriminant of~$E$} is the $p$\nobreakdash-adic discriminant of the $p$\nobreakdash-adic quadratic order~$\End(\FE)$.
For each $p$\nobreakdash-adic discriminant~$\pd$, put
\begin{displaymath}
  \Lambda_{\pd}
  \=
  \{ E \in \Ell(\Cp) \: \text{formal \CM{} point of $p$\nobreakdash-adic discriminant~$\pd$} \}.
\end{displaymath}

\begin{theoalph}
  \label{t:CM-symmetric}
  For every $p$\nobreakdash-adic discriminant~$\pd$, the subset~$\Lambda_{\pd}$ of~$\Ell(\Cp)$ is compact and there is a Borel probability measure~$\nu_{\pd}$ on~$\Ell(\Cp)$ such that ${\supp(\nu_{\pd}) = \Lambda_{\pd}}$ and the following equidistribution property holds.
  Let~$(D_n)_{n = 1}^{\infty}$ be a sequence of discriminants in~$\pd$ tending to~$-\infty$, such that, for every~$n$, the fundamental discriminant of~$D_n$ is either not divisible by~$p$ or not a prime discriminant.
  Then, the following weak convergence of measures holds,
  \begin{equation}
    \label{eq:2}
    \odelta_{D_n} \to \nu_{\pd}
    \text{ as }
    n \to \infty.
  \end{equation}
\end{theoalph}

The next result addresses the case left out in the theorem above.
Namely, that for some~$n_0$ in~$\Z_{>0}$ the fundamental discriminant~$d$ of~$D_{n_0}$ is a prime discriminant divisible by~$p$.
Passing to a subsequence if necessary, it suffices to consider the case where, for every~$n$, the fundamental discriminant of~$D_n$ equals~$d$.
In the following result, $\left( \frac{\cdot}{\cdot} \right)$ denotes the Kronecker symbol.

\begin{theoalph}
  \label{t:CM-broken-symmetry}
  Let~$d$ be a prime discriminant divisible by~$p$.
  Moreover, let~$m$ be in~$\Z_{\ge 0}$, put ${D \= d p^{2m}}$, and let~$\pd$ be the $p$\nobreakdash-adic discriminant containing~$D$.
  Then, there is a partition of~$\Lambda_{\pd}$ into disjoint compact sets,
  \begin{displaymath}
    \Lambda_{\pd}
    =
    \Lambda_{\pd}^+ \sqcup \Lambda_{\pd}^-,
  \end{displaymath}
  such that the measures~$\nu_{\pd}^+$ and~$\nu_{\pd}^-$, defined by
  \begin{displaymath}
    \nu_{\pd}^+
    \=
    2\nu_{\pd}|_{\Lambda_{\pd}^+}
    \text{ and }
    \nu_{\pd}^-
    \=
    2\nu_{\pd}|_{\Lambda_{\pd}^-},
  \end{displaymath}
  are both probability measures and the following equidistribution property holds.
  For every sequence~$(f_n)_{n = 0}^{\infty}$ in~$\Z_{>0}$ tending to~$\infty$ and satisfying, for every~$n$, the equality~$\left( \frac{d}{f_n} \right) = 1$ (resp. ${\left( \frac{d}{f_n} \right) = - 1}$), the following weak convergence of measures holds,
  \begin{displaymath}
    \odelta_{D f_n^2} \to \nu_{\pd}^+
    \text{ (resp. $\odelta_{D f_n^2} \to \nu_{\pd}^-$)}
    \text{ as }
    n \to \infty.
  \end{displaymath}
\end{theoalph}

An elliptic curve class in~$\Ell(\Cp)$ has \emph{supersingular} (resp. \emph{ordinary}) \emph{reduction}, if there is a representative Weierstrass equation with coefficients in~$\Op$ whose reduction is smooth and supersingular (resp. ordinary).
Every formal \CM{} point has supersingular reduction (Theorem~\ref{t:formal-CM-formulae}).
So, for every $p$\nobreakdash-adic discriminant~$\pd$, the set~$\Lambda_{\pd}$ is contained in the set~$\Sups$ of elliptic curve classes in~$\Ell(\Cp)$ with supersingular reduction.

Every $p$\nobreakdash-adic discriminant~$\pd$ contains a dense subset of discriminants, so there are plenty of sequences~$(D_n)_{n = 1}^{\infty}$ as in Theorem~\ref{t:CM-symmetric}.
A \CM{} point of discriminant~$D$ belongs to~$\Lambda_{\pd}$ if and only if~$D$ belongs to~$\pd$ (Corollary~\ref{c:formalization}$(ii)$).
In particular, for every discriminant~$D$ in~$\pd$, the inclusion ${\Lambda_D \subseteq \Lambda_{\pd}}$ holds and all \CM{} points appearing in Theorems~\ref{t:CM-symmetric} and~\ref{t:CM-broken-symmetry} belong to~$\Sups$.

The construction of the limit measure~$\nu_{\pd}$ in Theorem~\ref{t:CM-symmetric} relies on the following terminology.
A $p$\nobreakdash-adic discriminant is \emph{fundamental} if it is the $p$\nobreakdash-adic discriminant of the ring of integers of a quadratic field extension of~$\Qp$.
For a fundamental $p$\nobreakdash-adic discriminant~$\pfd$, the measure~$\nu_{\pfd}$ is a finite convex combination of projections of Haar measures on certain $p$\nobreakdash-adic compact Lie groups attached to supersingular elliptic curves over~$\Fpalg$.
Proposition~\ref{p:trace-zero-spheres} defines the Haar measures, Proposition~\ref{p:trace-zero-spheres-to-CM} and~\eqref{eq:79} their projections, and~\eqref{eq:80} defines~$\nu_{\pfd}$.
For a non-fundamental $p$\nobreakdash-adic discriminant~$\pd$, let~$\pfd$ be the fundamental $p$\nobreakdash-adic discriminant and~$m$ in~$\Z_{> 0}$ such that ${\pd = \pfd p^{2m}}$ (Lemma~\ref{l:p-adic-discriminants-Appendix}$(i)$).
The formula~\eqref{eq:117} defines~$\nu_{\pfd p^{2m}}$ as a suitable pull-back of~$\nu_{\pfd}$ under the $m$-th iterate of the canonical branch~$\t$ of the Hecke correspondence~$T_p$.

The special r{\^o}le of prime discriminants in the symmetry breaking of Theorem~\ref{t:CM-broken-symmetry} is related to genus theory, as follows.
Let~$\cQ$ be a ramified quadratic field extension of~$\Qp$ inside~$\Cp$, ${\cO_{\cQ}}$ its ring of integers, and~$\pfd$ the $p$\nobreakdash-adic discriminant of~$\cO_{\cQ}$.
Denoting by ${\nr \: \cQ \to \Qp}$ the norm map, $\nr(\cO_{\cQ}^{\times})$ is a subgroup of index two of~$\Z_p^{\times}$ containing~$(\Z_p^{\times})^2$.
Then, ${\Lambda_{\pfd}}$ partitions into two compact subsets, $\Lambda_{\pfd}^+$ and~$\Lambda_{\pfd}^-$, where two elliptic curve classes~$E$ and~$E'$ belong to the same part if~$E'$ is accumulated in~$\Ell(\C_p)$ by classes~$E''$ for which there exists an isogeny ${E \to E''}$ whose degree belongs to~$\nr(\cO_{\cQ}^{\times})$.
Let~$d$ in~$\pfd$ be a fundamental discriminant and put ${K \= \Q(\sqrt{d})}$.
Then, ${p \mid d}$ because~$\cQ$ is ramified over~$\Qp$.
The ideal class group~$\Cl(K)$ of~$K$ acts simply transitively on~$\Lambda_d$.
Denoting by~$t$ the number of prime divisors of~$d$, genus theory yields that the 2-torsion subgroup of~$\Cl(K)$ has~$2^{t - 1}$ elements.
When the discriminant~$d$ is prime, the group~$\Cl(K)$ has no element of order two, hence every element of~$\Cl(K)$ is a square.
This implies that~$\Lambda_d$ is completely inside either~$\Lambda_{\pfd}^+$ or~$\Lambda_{\pfd}^-$.
Similarly, for every~$f$ in ${\Z_{>0} \ssetminus p \Z}$, the set~$\Lambda_{df^2}$ is completely inside either~$\Lambda_{\pfd}^+$ or~$\Lambda_{\pfd}^-$ (Proposition~\ref{p:symmetry-breaking}$(ii)$).
In contrast, when the discriminant~$d$ is not prime, there is a nontrivial genus character ${\chi \: \Cl(K) \to \{1, -1\}}$ such that the partition of~$\Lambda_d$ induced by the partition ${\Cl(K) = \chi^{-1}(1) \sqcup \chi^{-1}(-1)}$ coincides with
\begin{equation}
  \label{eq:3}
  \Lambda_d
  =
  (\Lambda_d \cap \Lambda_{\pfd}^+) \sqcup (\Lambda_d \cap \Lambda_{\pfd}^-).
\end{equation}
In particular, ${\# (\Lambda_d \cap \Lambda_{\pfd}^+) = \# (\Lambda_d \cap \Lambda_{\pfd}^-)}$ (Proposition~\ref{p:symmetry-breaking}$(i)$).

In what follows, consider~$\Ell(\Cp)$ as a subspace of the Berkovich affine line~$\AKber$ over~$\Cp$, using the $j$-invariant to identify~$\Ell(\Cp)$ with the subspace~$\Cp$ of~$\AKber$.
Moreover, denote by~$\xcan$ the ``canonical'' or ``Gauss point'' of~$\AKber$.

The following corollary is a consequence of Theorem~\ref{t:CM-symmetric} and \cite[Theorem~A]{HerMenRiv20}.

\begin{coro}
  \label{c:CM-fundamental}
  The set of accumulation measures of
  \begin{equation}
    \label{eq:4}
    \left\{ \odelta_d \: d \text{ fundamental discriminant} \right\}
  \end{equation}
  in the space of Borel measures on~$\AKber$ equals
  \begin{displaymath}
    \left\{ \nu_{\pfd} \: \pfd \text{ fundamental $p$\nobreakdash-adic discriminant} \right\} \cup \{ \delta_{\xcan} \}.
  \end{displaymath}
\end{coro}

For all distinct $p$\nobreakdash-adic discriminants~$\pd$ and~$\pd'$, the compact sets~$\Lambda_{\pd}$ and~$\Lambda_{\pd'}$ are disjoint by definition, so the measures~$\nu_{\pd}$ and~$\nu_{\pd'}$ are different.
Since there are three fundamental $p$\nobreakdash-adic discriminants if~$p$ is odd and seven if~$p = 2$, see, \emph{e.g.}, Lemma~\ref{l:p-adic-discriminants-Appendix}$(iii)$ in Appendix~\ref{ss:Apendix-A},
Corollary~\ref{c:CM-fundamental} implies that~\eqref{eq:4} has precisely four accumulation measures if~$p$ is odd and eight if~$p = 2$.
This is in contrast to Duke's result that in the complex setting the limit exists~\cite{Duk88}.

To deduce Corollary~\ref{c:CM-fundamental} from Theorem~\ref{t:CM-symmetric} and \cite[Theorem~A]{HerMenRiv20}, consider the following consequence of the latter.
For each sequence of discriminants~$(D_j)_{j = 1}^{\infty}$ tending to~$-\infty$, \cite[Theorem~A]{HerMenRiv20} implies the convergence of measures ${\odelta_{D_j} \to \delta_{x_{\can}}}$ as ${j \to \infty}$ in each of the following situations:
\begin{itemize}
\item [$(i)$]
  For every~$j$, the set~$\Lambda_{D_j}$ is disjoint from~$\Sups$.
\item [$(ii)$]
  For every~$j$, the inclusion ${\Lambda_{D_j} \subseteq \Sups}$ holds and ${|D_j|_p \to 0}$ as ${j \to \infty}$.
\end{itemize}
Corollary~\ref{c:CM-fundamental} is a direct consequence of this property, Theorem~\ref{t:CM-symmetric}, and the fact that a \CM{} point of fundamental discriminant~$d$ belongs to~$\Sups$ if and only if~$d$ belongs to a fundamental $p$\nobreakdash-adic discriminant (Lemma~\ref{l:p-adic-discriminants}).

The following corollary is an immediate of Theorems~\ref{t:CM-symmetric} and~\ref{t:CM-broken-symmetry}, the consequence of \cite[Theorem~A]{HerMenRiv20} stated above, and the fact that a \CM{} point belongs to~$\Sups$ if and only if its discriminant belongs to a $p$\nobreakdash-adic discriminant (Lemma~\ref{l:p-adic-discriminants}).

\begin{coro}
  \label{c:CM}
  In the case where ${p \equiv -1 \mod 4}$, denote by~$\widehat{\pfd}$ the $p$\nobreakdash-adic discriminant containing~$-p$.
  In the case where ${p = 2}$, denote by~$\widehat{\pfd}$ (resp. $\widehat{\pfd}'$) the $p$\nobreakdash-adic discriminant containing~$-4$ (resp.~$-8$).
  Then, the set of accumulation measures of~\eqref{eq:1} in the space of Borel measures on~$\AKber$ equals
  \begin{displaymath}
    \left\{ \nu_{\pd} \: \pd \text{ $p$\nobreakdash-adic discriminant} \right\} \cup \{ \delta_{\xcan} \},
  \end{displaymath}
  \begin{displaymath}
    \left\{ \nu_{\pd} \: \pd \text{ $p$\nobreakdash-adic discriminant} \right\}
    \cup
    \left\{ \nu_{\widehat{\pfd} p^{2m}}^+, \nu_{\widehat{\pfd} p^{2m}}^- \: m \in \Z_{\ge 0} \right\}
    \cup
    \{ \delta_{\xcan} \},
  \end{displaymath}
  or
  \begin{displaymath}
    \left\{ \nu_{\pd} \: \pd \text{ $p$\nobreakdash-adic discriminant} \right\}
    \cup
    \left\{ \nu_{\widehat{\pfd} p^{2m}}^+, \nu_{\widehat{\pfd} p^{2m}}^-, \nu_{\widehat{\pfd}' p^{2m}}^+, \nu_{\widehat{\pfd}' p^{2m}}^- \: m \in \Z_{\ge 0} \right\}
    \cup
    \{ \delta_{\xcan} \},
  \end{displaymath}
  depending on whether ${p \equiv 1 \mod 4}$, ${p \equiv -1 \mod 4}$, or~${p = 2}$, respectively.

  In particular, in all cases the set of accumulation measures of~\eqref{eq:1} is countably infinite.
\end{coro}

This is in contrast to the complex setting where the limit exists~\cite{Duk88,CloUll04}.

The companion paper~\cite{HerMenRiv24} proves that every accumulation measure of~\eqref{eq:1} different from~$\delta_{\xcan}$ is nonatomic \cite[Theorem~B]{HerMenRiv24}.
This is one of the main ingredients in showing that, for every finite set of prime numbers~$S$, there are at most finitely many singular moduli that are $S$-units \cite[Theorem~A]{HerMenRiv24}.

The following corollary is an immediate consequence of Theorem~\ref{t:CM-symmetric} and the fact that every formal \CM{} point belongs to~$\Sups$ (Theorem~\ref{t:formal-CM-formulae}).

\begin{coro}
  \label{c:formal-density}
  Every formal \CM{} point is accumulated by \CM{} points in~$\Ell(\Cp)$.
  In particular, \CM{} points in~$\Sups$ form a dense subset of the set of formal \CM{} points.
\end{coro}

Coleman and McMurdy proved the first result of this type \cite[Theorem~4.1]{ColMcM06}, for~$p \ge 5$ and in the context of certain modular curves of level bigger than one.
Their result implies that every formal \CM{} point~$E$ such that~$\End(\FE)$ is the ring of integers of a ramified quadratic field extension of~$\Qp$ is approximated by \CM{} points.

The following corollary is a direct consequence of Corollary~\ref{c:formal-density}, together with the isolation of \CM{} points with ordinary reduction \cite[Corollary~B]{HerMenRiv20}.

\begin{coro}
  \label{c:ordinary-isolation}
  A \CM{} point is isolated in~$\Ell(\Cp)$ if and only if it has ordinary reduction.
\end{coro}

\subsection{Equidistribution of partial Hecke orbits}
\label{ss:Hecke-orbits}
The statement of the next main result relies on the following definition of Hecke correspondences.
See Section~\ref{ss:Hecke-correspondences} for background.

A \emph{divisor on~$\Ell(\Cp)$} is an element of the free Abelian group
\begin{displaymath}
  \Div(\Ell(\Cp))
  \=
  \bigoplus_{E\in \Ell(\Cp)} \Z E.
\end{displaymath}
For each divisor~$\fD$ in~$\Div(\Ell(\Cp))$, given by ${\fD = \sum_{E \in \Ell(\Cp)} n_E E}$, the \emph{degree} and \emph{support} of~$\fD$ are defined by
\begin{displaymath}
  \deg(\fD)
  \=
  \sum_{E \in \Ell(\Cp)} n_E
  \text{ and }
  \supp(\fD)
  \=
  \{ E \in \Ell(\Cp) \: n_E \neq 0 \},
\end{displaymath}
respectively.
If in addition ${\deg(\fD) \ge 1}$ and, for every~$E$ in~$\Ell(\Cp)$, the inequality ${n_E \ge 0}$ holds, then
\begin{displaymath}
  \odelta_{\fD}
  \=
  \frac{1}{\deg(\fD)}\sum_{E \in \Ell(\Cp)} n_E \delta_E
\end{displaymath}
is a Borel probability measure on~$\Ell(\Cp)$.

For each~$n$ in~$\Z_{>0}$, the $n$-th \emph{Hecke correspondence} is the linear map
\begin{displaymath}
  T_n \: \Div(\Ell(\Cp)) \to \Div(\Ell(\Cp))
\end{displaymath}
defined, for~$E$ in~$\Ell(\Cp)$, by
\begin{displaymath}
  T_n(E)
  \=
  \sum_{C\leq E \text{ of order } n}E/C,
\end{displaymath}
where the sum runs over all subgroups~$C$ of~$E$ of order~$n$.
So, $\supp(T_n(E))$ is the set of~$E'$ in~$\Ell(\Cp)$ for which there is an isogeny ${E \to E'}$ of degree~$n$.

For each~$E$ in ${\Ell(\Cp) \ssetminus \Sups}$, the asymptotic distribution of the Hecke orbit~$(T_n(E))_{n = 1}^{\infty}$ is described by \cite[Theorem~C]{HerMenRiv20}.
The following theorem addresses the more difficult case where~$E$ belongs to~$\Sups$.
The description depends on a subgroup~$\NE$ of~$\Zp^{\times}$ defined as follows.
In the case where~$E$ is not a formal \CM{} point, ${\NE \= (\Zp^{\times})^2}$.
In the case where~$E$ is a formal \CM{} point, denote by~$\Aut(\FE)$ the group of isomorphisms of~$\FE$ defined over~$\OQpalg$, and by~$\nr$ the norm map from the field of fractions of~$\End(\FE)$ to~$\Qp$.
Then,
\begin{displaymath}
  \NE
  \=
  \left\{ \nr \left( \varphi \right) \: \varphi \in \Aut(\FE) \right\}.
\end{displaymath}
In all cases~$\NE$ is a subgroup of~$\Zp^{\times}$ containing~$(\Zp^{\times})^2$.
In particular, the index of~$\NE$ in~$\Zp^{\times}$ is at most two if~$p$ is odd, and at most four if~$p=2$.

\begin{theoalph}[Equidistribution of partial Hecke orbits]
  \label{t:Hecke-orbits}
  Let~$E$ be in~$\Sups$, let~$\coset$ be a coset in~$\Qp^{\times}/ \NE$ contained in~$\Zp$, and define the \emph{partial Hecke orbit} by
  \begin{equation}
    \label{eq:5}
    \corbit
    \=
    \bigcup_{n \in \coset \cap \Z_{>0}} \supp(T_n(E)).
  \end{equation}
  Then, the closure~$\corbitc$ in~$\Sups$ of this set is compact.
  Moreover, there is a Borel probability measure~$\mu_{\coset}^E$ on~$\Ell(\Cp)$ such that ${\supp(\mu_{\coset}^E) = \corbitc}$ and, for every sequence~$(n_j)_{j = 1}^{\infty}$ in~$\coset \cap \Z_{>0}$ tending to~$\infty$, the following weak convergence of measures holds,
  \begin{displaymath}
    \odelta_{T_{n_j}(E)} \to \mu_{\coset}^E
    \text{ as }
    j \to \infty.
  \end{displaymath}
\end{theoalph}

See Theorem~\ref{t:Hecke-orbits-pr} in Section~\ref{s:Hecke-orbits} for a quantitative version of this result.

Together with \cite[Theorem~C]{HerMenRiv20}, Theorem~\ref{t:Hecke-orbits} identifies all limits of Hecke orbits in~$\Ell(\Cp)$.
Indeed, \cite[Theorem~C]{HerMenRiv20} implies that, for all~$E$ in~$\Ell(\Cp)$ and sequence~$(n_j)_{j = 1}^{\infty}$ in~$\Z_{>0}$ tending to~$\infty$, the weak convergence of measures ${\odelta_{T_{n_j}(E)} \to \delta_{x_{\can}}}$ as ${j \to \infty}$ holds in each of the following situations:
\begin{enumerate}
\item [$(i)$]
  $E$ is outside~$\Sups$.
\item [$(ii)$]
  $E$ belongs to~$\Sups$ and ${|n_j|_p \to 0}$ as ${j \to \infty}$.
\end{enumerate}
Combined with Theorem~\ref{t:Hecke-orbits}, this implies the following as an immediate consequence.

\begin{coro}
  \label{c:Hecke-orbits}
  For every~$E$ in~$\Ell(\Cp)$, the set of accumulation measures of~$(\odelta_{T_n(E)})_{n = 1}^{\infty}$ in the space of Borel probability measures on~$\AKber$ equals
  \begin{equation}
    \label{eq:6}
    \left\{ \mu_{\coset}^E \: \coset \in \Qp^{\times} / \NE, \coset \subset \Zp \right\} \cup \{ \delta_{\xcan} \}.
  \end{equation}
\end{coro}

For all distinct cosets~$\coset$ and~$\coset'$ in ${\Qp^{\times} / \NE}$ contained in~$\Zp$, the measures~$\mu_{\coset}^E$ and~$\mu_{\coset'}^E$ are different (Proposition~\ref{p:orbit-measures}$(ii)$).
In particular, the set of accumulation measures~\eqref{eq:6} is countably infinite.
This is in contrast to the complex setting where the limit exists, see~\cite{CloUll04,CloOhUll01,EskOh06}.
The companion paper \cite{HerMenRiv24} proves that~$\mu_{\coset}^E$ is nonatomic.

\subsection{Asymptotic distribution of integer points on $p$-adic spheres}
\label{s:p-adic-linnik}
The proofs of Theorems~\ref{t:CM-symmetric}, \ref{t:CM-broken-symmetry}, and~\ref{t:Hecke-orbits} rely on the $p$\nobreakdash-adic equidistribution result stated below, which is inspired by Linnik's classical problem on the asymptotic distribution of integer points on spheres.
See~\cite{Duk88,EllMicVen13} for refinements and a historical perspective.

Fix an integer~$n$ satisfying ${n \ge 3}$ and a positive definite quadratic form~$Q$ in~$\Z[X_1, \ldots, X_n]$.
For each~$m$ in~$\Z_{>0}$, put
\begin{align*}
  V_m(Q)
  & \=
    \{ \bfx \in \Z^n \: Q(\bfx)=m\}
    \intertext{and, for each~$\ell$ in~$\Zp$, define the sphere}
    S_{\ell}(Q)
  & \=
    \{x\in \Zp^n \: Q(x) = \ell\}.
\end{align*}
The orthogonal group of~$Q$ with coefficients in~$\Zp$, defined as
\begin{displaymath}
  \oO_Q(\Zp)
  \=
  \{ T \in \GL_n(\Zp) \: Q(T \cdot X) = Q(X) \},
\end{displaymath}
is compact, acts on~$\Zp^n$, and, for every~$\ell$ in~$\Zp$, preserves~$S_{\ell}(Q)$.

For the following result, let~$\ell$ in~$\Zp \ssetminus \{ 0 \}$ be such that ${S_{\ell}(Q) \neq \emptyset}$ and the compact group~$\oO_Q(\Zp)$ acts transitively on~$S_{\ell}(Q)$.
Then, there is a unique Borel probability measure on~$S_{\ell}(Q)$ that is invariant under the action of~$\oO_Q(\Zp)$, see, \emph{e.g.}, Lemma~\ref{l:existence-invariant-measure}.
For every~$u$ in~$\Zp^{\times}$, denote by~$M_u$ the element of~$\GL_n(\Zp)$ defined by
\begin{displaymath}
  M_{u} (X_1, \ldots, X_n)
  \=
  (u X_1, \ldots, u X_n).
\end{displaymath}
Then, ${M_u(S_\ell(Q)) = S_{\ell u^2}(Q)}$.

\begin{theoalph}[$p$-Adic Linnik equidistribution]
  \label{t:p-adic-Linnik}
  Put ${\kappa_n \= \frac{1}{2}}$ if~$n$ is even and ${\kappa_n \= \frac{2}{7}}$ if~$n$ is odd, and fix a real number~$c$ satisfying ${c > \frac{n}{4}-\kappa_n}$.
  Let~$\ell$ in~$\Zp \ssetminus \{ 0 \}$ be such that ${S_\ell(Q) \neq \emptyset}$ and~$\oO_Q(\Zp)$ acts transitively on~$S_{\ell}(Q)$, and~$\mu_{\ell}$ the unique Borel probability measure on~$S_{\ell}(Q)$ that is invariant under the action of~$\oO_Q(\Zp)$.
  Moreover, let~$(m_j)_{j = 1}^{\infty}$ be a sequence in~$\Z_{>0}$ tending to~$\infty$ that is contained in the coset~$\ell (\Zp^{\times})^2$ of~$\Qp^{\times} / (\Zp^{\times})^2$ and such that, for every sufficiently large~$j$, the inequality ${\# V_{m_j}(Q) \ge m_j^c}$ holds.
  For each~$j$ in~$\Z_{>0}$, let~$u_j$ in~$\Zp^{\times}$ be such that ${m_j = \ell u_j^2}$.
  If~$n = 3$, then suppose in addition that there is~$S$ in~$\R$ such that ${S \ge 1}$ and, for each~$j$, the largest square diving~$m_j$ is less than or equal to~$S$.
  Then, the following weak convergence of measures on~$\Zp^n$ holds,
  \begin{displaymath}
    \frac{1}{\# V_{m_j}(Q)} \sum_{\bfx \in V_{m_j}(Q)} \delta_{M_{u_j}^{-1}(\bfx)} \to \mu_\ell
    \text{ as }
    j \to \infty.
  \end{displaymath}
\end{theoalph}

See also Theorem~\ref{t:finitary-deviation-estimate} and Corollary~\ref{c:functional-deviation} for quantitative variants of this result.

The circle method yields that~$\# V_m(Q)$ grows at least like~$m^{\frac{n}{2} - 1}$, provided ${n \ge 5}$ and, for every prime number~$q$, the equation ${Q(x) = m}$ is solvable in~$\Z_q^n$.
For ${n = 4}$, the circle method yields that, for every~$\varepsilon$ in~$\R_{>0}$, the number~$\# V_m(Q)$ grows at least like~$m^{1 - \varepsilon}$, provided that, for every prime number~$q$, the equation ${Q(x) = m}$ has a solution~$x$ in~$\Z_q^n$ for which~$\nabla Q(x)$ is a unit in~$\Z_q$, see, \emph{e.g.}, \cite[Theorem~4 and Corollary~1]{Hea96}.
For~$n \ge 3$, the quantity~$\# V_m(Q)$ can also be estimated in some situations using the theory of modular forms and theta series, see, \emph{e.g.}, \cite{Sie67}, \cite[Chapter~VI]{Ogg69}, or \cite[Chapter~11]{Iwa97}.
This paper does not rely on these general results because the growth of~$\# V_{m}(Q)$ in~$m$ is well understood in the applications considered here.

Besides the aforementioned results on the asymptotic distribution of integer points on spheres, previous works on the distribution of lattice points on spheres have focused primarily on obtaining asymptotic formulae for the number of such points in convex regions of~$\R^n$ and in fixed residue classes in~$\Q^n/\Z^n$.
These asymptotic formulae are usually expressed in terms of products of local densities and/or appropriate singular series.
See, \emph{e.g.}, \cite{DukSch90} and the references therein.
Siegel introduced local densities in~\cite{Sie35a} to measure the density of representations of an integer by a quadratic form over each local completion of~$\Q$.
See also \cite[p.~184]{Iwa97}.
Weil later extended Siegel's work to the adelic setting, leading to what is now known as the Siegel--Weil formula.
Singular series trace their origins back to the circle method of Hardy, Ramanujan, and Littlewood, and were later refined by Kloosterman.
Theorem~\ref{t:p-adic-Linnik} focuses on the $p$\nobreakdash-adic asymptotic distribution of lattice points and is local in nature.

\subsection{Notes and references}
\label{ss:notes-references}

In~\cite{Dis22}, Disegni investigates the $p$\nobreakdash-adic asymptotic distribution of \CM{} points on modular curves with level structure as the conductor tends to zero $p$\nobreakdash-adically.
As explained in~\cite[p.~640]{Dis22}, these results do not overlap with those in the present paper.
Indeed, when restricted to the modular curve of level one, the main theorem of~\cite{Dis22} is subsumed by \cite[Theorem~A]{HerMenRiv20}.

Since the completion of the present paper, other works have appeared on the distribution of \CM{} points.
P{\'e}rez-Pi{\~n}a studies the $p$\nobreakdash-adic distribution of \CM{} points on a Shimura curve in \cite{Per25} when~$p$ divides the discriminant of the underlying quaternion algebra.
In contrast to the case studied here, the $p$\nobreakdash-adic space in the latter is uniformized by the $p$\nobreakdash-adic upper-half plane, and a unique asymptotic distribution measure arises.
Saettone established analogous results for function fields and \mbox{Drinfel'd} modules in \cite{0Sae2606}.

\subsection{Strategy and organization}
\label{ss:organization}
This section explains the proof strategy of the main results and simultaneously describes the organization of the paper.

After some preliminaries in Section~\ref{s:preliminaries}, Section~\ref{s:p-adic-Linnik} proves Theorem~\ref{t:p-adic-Linnik} on the asymptotic distribution of integer points on $p$\nobreakdash-adic spheres.
The proof relies on a deviation estimate modulo large powers of~$p$ (Theorem~\ref{t:finitary-deviation-estimate}).
Section~\ref{ss:modular form} provides the main ingredient: The construction of an auxiliary cuspidal modular form (Proposition~\ref{p:modular form}).
The formula~\eqref{eq:21} expresses this modular form in terms of the ``special congruent theta functions'' of \cite[Chapter 10]{Iwa97}, with spherical function ${P = 1}$.
See also \cite[Chapter~VI]{Ogg69}.
The proof of the deviation estimate relies on the bounds for Fourier coefficients of cuspidal modular forms of Deligne~\cite{Del74}, Iwaniec~\cite{Iwa87}, and Duke~\cite{Duk88}.
Section~\ref{ss:proof of finitary deviation estimate} completes the proof of the deviation estimate and derives Theorem~\ref{t:p-adic-Linnik} from it.

Section~\ref{s:CM-formulae} provides several formulae for (formal) \CM{} points having supersingular reduction.
Section~\ref{ss:fixed-points-formula} provides the first formula, which relies on the Gross--Hopkins group action on the Lubin--Tate deformation space~\cite{HopGro94b}, recalled in Section~\ref{ss:from-elliptic-curves}.
It interprets each \CM{} point of fundamental discriminant as projections of fixed points of certain elements of this action (Theorem~\ref{t:fixed-points-formula}).
Sections~\ref{ss:CM-from-canonical} and~\ref{ss:formal-CM-formulae} provide the remaining formulae, which rely on the canonical branch~$\t$ of~$T_p$ to relate (formal) \CM{} points whose conductors differ by a power of~$p$ (Theorems~\ref{t:CM-from-canonical} and~\ref{t:formal-CM-formulae}).

Section~\ref{s:CM-fundamental} describes the asymptotic distribution of \CM{} points of fundamental discriminant in a quantitative form (Theorem~\ref{t:CM-fundamental}).
It is one of the main ingredients in the proof of Theorem~\ref{t:CM-symmetric}.
To explain the strategy of proof, fix a supersingular elliptic curve class~$\ss$ in the moduli space~$\Ell(\Fpalg)$ of elliptic curves over~$\Fpalg$ and denote by~$\Rss$ the $p$\nobreakdash-adic space of endomorphisms of the formal $\Zp$\nobreakdash-module of~$\ss$.
The proof begins in Section~\ref{ss:trace-zero-spheres} defining certain spheres in the zero-trace subspace of~$\Rss$, called \emph{zero-trace spheres}, and showing that each of these sets carries a natural homogeneous measure (Proposition~\ref{p:trace-zero-spheres}).
Sections~\ref{ss:CM-parametrization} and~\ref{ss:trace-zero-spheres-to-CM} provide a key step in the proof of Theorem~\ref{ss:trace-zero-spheres}: For every fundamental $p$\nobreakdash-adic discriminant~$\pfd$, the set of formal \CM{} points in~$\Lambda_{\pfd}$ in the residue disc associated with~$\ss$ is naturally parametrized by a zero-trace sphere (Propositions~\ref{p:CM-parametrization} and~\ref{p:trace-zero-spheres-to-CM}).
Using this parametrization, Section~\ref{ss:proof-CM-fundamental} deduces Theorem~\ref{t:CM-fundamental} from the asymptotic distribution of integer points on $p$\nobreakdash-adic spheres in Section~\ref{s:p-adic-Linnik}, and the equidistribution of \CM{} points on supersingular residue discs in Section~\ref{ss:residual} (Theorem~\ref{t:residual}).

Section~\ref{s:Hecke-orbits} proves the results on the asymptotic distribution of Hecke orbits (Theorem~\ref{t:Hecke-orbits} in Section~\ref{ss:Hecke-orbits}) and that distinct partial Hecke orbits have different limit measures (Proposition~\ref{p:orbit-measures}).
The proof of the former occupies Sections~\ref{ss:supersingular-spheres}, \ref{ss:supersingular-sphere-to-orbit}, and~\ref{ss:Hecke-orbits-proof}, and the latter Section~\ref{ss:applications}.
Section~\ref{s:Hecke-orbits} begins stating a quantitative version of Theorem~\ref{t:Hecke-orbits} with a convergence rate that is uniform on the initial point (Theorem~\ref{t:Hecke-orbits-pr}).
The proof of Theorems~\ref{t:CM-symmetric} and~\ref{t:CM-broken-symmetry} rely on Theorem~\ref{t:Hecke-orbits-pr}.
To explain the strategy of proof of Theorem~\ref{t:Hecke-orbits-pr}, fix supersingular elliptic curve classes~$\ss$ and~$\ss'$ in~$\Ell(\Fpalg)$ and denote by~$\bfR_{\ss, \ss'}$ the $p$\nobreakdash-adic space of morphisms from the formal $\Zp$\nobreakdash-module of~$\ss$ to that of~$\ss'$.
The proof begins in Section~\ref{ss:supersingular-spheres} by introducing certain spheres in~$\bfR_{\ss, \ss'}$, called \emph{supersingular spheres},
and by showing that each supersingular sphere carries a natural homogeneous measure (Proposition~\ref{p:supersingular-spheres}).
Section~\ref{ss:supersingular-sphere-to-orbit} provides a key step: Each partial Hecke orbit restricted to the residue disc associated with~$\ss'$ is parametrized by a supersingular sphere (Proposition~\ref{p:supersingular-sphere-to-orbit}).
After these considerations, Section~\ref{ss:Hecke-orbits-proof} deduces Theorem~\ref{t:Hecke-orbits-pr} from the asymptotic distribution of integer points on $p$\nobreakdash-adic spheres in Section~\ref{s:p-adic-Linnik}.

Section~\ref{s:CM} proves the results on the asymptotic distribution of \CM{} points (Theorems~\ref{t:CM-symmetric} and~\ref{t:CM-broken-symmetry} in Section~\ref{ss:Hecke-orbits}).
The (formal) \CM{} points formulae in Sections~\ref{ss:CM-from-canonical} and~\ref{ss:formal-CM-formulae} reduce the analysis to the case of fundamental $p$\nobreakdash-adic discriminants.
For each fundamental $p$\nobreakdash-adic discriminant~$\pfd$, the proof begins in Section~\ref{ss:CM-orbits} by studying how~$\Lambda_{\pfd}$ partitions into closures of partial Hecke orbits (Proposition~\ref{p:CM-orbits}).
The set~$\Lambda_{\pfd}$ coincides with the closure of a partial Hecke orbit if~$\Qpfd$ is unramified over~$\Qp$.
If~$\Qpfd$ is ramified over~$\Qp$, then~$\Lambda_{\pfd}$ partitions into precisely two closures of partial Hecke orbits.
In the latter case, Section~\ref{ss:symmetry-breaking} relies on genus theory to determine, for each discriminant~$D$ in~$\pfd$, how~$\Lambda_D$ distributes between these parts (Proposition~\ref{p:symmetry-breaking}).
Here is where prime discriminants divisible by~$p$ play a special role.
Section~\ref{s:proof-of-CM} relies on these results to deduce Theorems~\ref{t:CM-symmetric} and~\ref{t:CM-broken-symmetry} from Theorems~\ref{t:CM-fundamental} and~\ref{t:Hecke-orbits-pr}.

For the reader's convenience, Appendix~\ref{ss:Apendix-A} gathers some basic facts about quadratic field extensions of~$\Qp$ and $p$\nobreakdash-adic discriminants.

\subsection*{Acknowledgments}
The authors thank the referees for their numerous and varied comments that helped us improve the exposition.

The first named author was supported by ANID/CONICYT, FONDECYT Postdoctorado Nacional grant~3190086.
The second named author was partially supported by FONDECYT grant 1171329.
The first and second named author were partially supported by ANID FONDECYT Regular grant 1250734.
The third named author acknowledges partial support from NSF grant DMS-1700291.
The authors would like to thank Pontificia U.~Cat{\'o}lica de Valpara{\'{\i}}so, U.~of Rochester, and U.~de Barcelona for hospitality during the preparation of this work.

\section{Preliminaries}
\label{s:preliminaries}

For each~$n$ in~$\Z_{>0}$, put
\begin{displaymath}
  d(n)
  \=
  \sum_{d \in \Z_{>0}, d \mid n} 1
  \text{ and }
  \sigma_1(n)
  \=
  \sum_{d \in \Z_{>0}, d \mid n} d.
\end{displaymath}
Then,
\begin{equation}
  \label{eq:lowb-sigma1}
  \sigma_1(n)
  \ge
  n,
\end{equation}
and, for every~$\varepsilon$ in~$\R_{>0}$,
\begin{equation}
  \label{eq:7}
  d(n)
  =
  o(n^{\varepsilon}).
\end{equation}

Given an algebraically closed field~$\K$, denote by~$\Ell(\K)$ the set of ${\K}$-isomorphism classes of elliptic curves over~$\K$.
For each class~$E$ in~$\Ell(\K)$, the $j$-invariant~$j(E)$ of~$E$ is an element of~$\K$ determining~$E$ completely and the map ${j \: \Ell(\K) \to \K}$ is a bijection.

Given a field extension~$\cK$ of~$\Qp$, denote by~$\OK$ its ring of integers and by~$\MK$ the maximal ideal of~$\OK$.
When ${\cK = \Cp}$, denote~$\OK$ and~$\MK$ by~$\Op$ and~$\Mp$, respectively.
Moreover, identify the residue field of~$\Cp$ with an algebraic closure~$\Fpalg$ of the field with~$p$ elements~$\Fp$ and denote by ${\pi \: \Op\to \Fpalg}$ the reduction map.
For every finite field extension~$\cK$ of~$\Qp$ inside~$\Cp$, the equality ${\cM_{\cK} = \Mp \cap \OK}$ holds.

For each quadratic field extension~$\cK$ of~$\Qp$, denote by~$x \mapsto \overline{x}$ the unique field automorphism of~$\cK$ over~$\Qp$ different from the identity.
Moreover, for every~$x$ in~$\cK$, put
\begin{displaymath}
  \tr(x)
  \=
  x + \overline{x},
  \nr(x)
  \=
  x \overline{x},
  \text{ and }
  \disc(x)
  \=
  (x-\overline{x})^2.
\end{displaymath}
Note that~$\tr(x)$, ${\nr(x)}$, and~$\Delta(x)$ belong to~$\Qp$ and ${\Delta(x) = \tr(x)^2 - 4 \nr(x)}$.

Denote by~$\Q_{p^2}$ the unique unramified quadratic field extension of~$\Qp$ inside~$\Cp$, $\Z_{p^2}$ its ring of integers, and~$\F_{p^2}$ its residue field.
For each~$\Delta$ in~$\Qp$, denote by~$\Qp(\sqrt{\Delta})$ the smallest field extension of~$\Qp$ inside~$\Cp$ containing a root of~$X^2 - \Delta$.
Lemma~\ref{l:quadratic-extensions}$(i)$ in Appendix~\ref{ss:Apendix-A} identifies all quadratic field extensions of~$\Qp$ inside~$\Cp$.

The endomorphism ring of an elliptic curve over~$\Fpalg$ is isomorphic to an order in either a quadratic imaginary extension of~$\Q$ or a quaternion algebra over~$\Q$.
In the latter case the corresponding elliptic curve class is \emph{supersingular}.

An elliptic curve class~$E$ has \emph{good reduction}, if it is represented by a Weierstrass equation with coefficients in~$\Op$ whose reduction is smooth.
In this case the reduction is an elliptic curve over~$\Fpalg$, whose class~$\tE$ only depends on~$E$ and is the \emph{reduction of~$E$}.

A \emph{divisor} on a set~$X$\footnote{
  This paper only uses this definition when~$X$ is one of several types of one-dimensional objects.
  For such~$X$, the notion of divisor introduced here can be seen as a natural extension of the usual notion of Weil divisor.
}
is a formal finite sum~$\sum_{x\in X}n_xx$ in~$\bigoplus_{x \in X} \Z x$.
The \emph{degree} and \emph{support} of a divisor~$\fD$, given by ${\fD = \sum_{x \in X}n_x x}$, are defined by
\begin{displaymath}
  \deg(\fD)
  \=
  \sum_{x \in X} n_x
  \text{ and }
  \supp(\fD)
  \=
  \{ x \in X \: n_x \neq 0 \},
\end{displaymath}
respectively.
For all set~$X'$ and map ${f \: X \to X'}$, the \emph{push-forward action of~$f$ on divisors~$f_* \: \Div(X) \to \Div(X')$} is the linear extension of the action of~$f$ on points.

\subsection{Discriminants and their $p$-adic counterparts}
\label{ss:discriminants}

Recall that a fundamental discriminant~$d$ is the discriminant of the ring of integers of a quadratic imaginary extension~$K$ of~$\Q$.
If~$d_0$ is the unique square-free integer such that ${K = \Q (\sqrt{d_0})}$, then
\begin{equation}
  \label{eq:8}
  d
  =
  \begin{cases}
    d_0 & \text{if } d_0 \equiv 1 \mod 4;
    \\
    4d_0 & \text{if } d_0 \equiv -1, 2 \mod 4.
  \end{cases}
\end{equation}

Recall that a discriminant~$D$ is the discriminant of an order in a quadratic imaginary extension of~$\Q$.
Moreover, the fundamental discriminant of~$D$ is the discriminant~$d$ of~$\Q(\sqrt{D})$, it divides~$D$, and the quotient~$\frac{D}{d}$ is the square of the conductor of~$D$.
Conversely, for all fundamental discriminant~$d$ and~$f$ in~$\Z_{>0}$, the integer~$d f^2$ is the unique discriminant of fundamental discriminant~$d$ and conductor~$f$.
Moreover, putting ${D \= df^2}$, there is a unique order~$\cO_{d, f}$ of discriminant~$D$ in the quadratic imaginary extension~$\Q(\sqrt{d})$ of~$\Q$, and it is given by
\begin{equation}
  \label{eq:9}
  \cO_{d, f} \= \Z\left[ \tfrac{D + \sqrt{D}}{2} \right].
\end{equation}
Conversely, every order in~$\Q(\sqrt{d})$ is of this form, see, \emph{e.g.}, \cite[Chapter~8, Section~1, Theorem~3]{Lan87}.
Note that ${\cO_{d, f} = \Z + f \cO_{d, 1}}$ and the index of~$\cO_{d, f}$ in~$\cO_{d, 1}$ equals~$f$.
Moreover, $\cO_{d, 1}$ is the unique maximal order in~$\Q(\sqrt{d})$ and it coincides with the ring of integers of~$\Q(\sqrt{d})$.

A discriminant~$D$ is \emph{$p$-supersingular}, if the reduction of some \CM{} point of discriminant~$D$ is supersingular.
In this case, the reduction of every \CM{} point of discriminant~$D$ is supersingular.
Equivalently, a discriminant~$D$ is $p$\nobreakdash-supersingular if~$p$ is ramified or inert in~$\Q(\sqrt{D})$, see~\cite{Deu41} or \cite[Chapter~13, Section~4, Theorem~12]{Lan87}.
A discriminant is $p$\nobreakdash-supersingular if and only if its fundamental discriminant is.
A fundamental discriminant~$d$ is $p$\nobreakdash-supersingular if and only if~$\left( \frac{d}{p} \right) \neq 1$ if~$p$ is odd and~$d \not\equiv 1 \mod 8$ if~$p = 2$.

Recall that a $p$\nobreakdash-adic quadratic order is a $\Zp$-order in a quadratic field extension of~$\Qp$.
For each quadratic field extension~$\cK$ of~$\Qp$, the ring of integers~$\OK$ is the unique maximal $\Zp$-order in~$\cK$.
Moreover, for every~$m$ in~$\Z_{\ge 0}$, the set~$\Zp + p^m \OK$ is a $\Zp$-order in~$\cK$ and every $\Zp$-order in~$\cK$ is of this form.

Recall that a $p$\nobreakdash-adic discriminant is a coset in~$\Qp^{\times} / (\Zp^{\times})^2$ formed by the discriminants of all $\Zp$-bases of a $p$\nobreakdash-adic quadratic order.
Furthermore, a $p$\nobreakdash-adic discriminant is \emph{fundamental}, if it is the $p$\nobreakdash-adic discriminant of the ring of integers of a quadratic field extension of~$\Qp$.
The $p$\nobreakdash-adic discriminant is an isomorphism invariant of $p$\nobreakdash-adic quadratic orders.
Lemma~\ref{l:p-adic-discriminants-Appendix} in Appendix~\ref{ss:Apendix-A} identifies all $p$\nobreakdash-adic quadratic orders and $p$\nobreakdash-adic discriminants.
For all $p$\nobreakdash-adic discriminant~$\pd$ and~$\Delta$ in~$\pd$, the field~$\Qp(\sqrt{\Delta})$ is a quadratic field extension of~$\Qp$ inside~$\Cp$ that depends only on~$\pd$, but not on~$\Delta$.
Denote it by~$\Qpd$.

The following lemma gathers basic facts that are important in what follows.
For the reader's convenience, Appendix~\ref{ss:Apendix-A} includes a proof.

\begin{lemma}
  \label{l:p-adic-discriminants}
  A discriminant (resp.~a discriminant whose conductor is not divisible by~$p$) belongs to a $p$\nobreakdash-adic discriminant (resp. fundamental $p$\nobreakdash-adic discriminant) if and only if it is $p$-supersingular.
  Moreover, for each $p$\nobreakdash-adic discriminant (resp. fundamental $p$\nobreakdash-adic discriminant)~$\pd$, the set of discriminants (resp. fundamental discriminants) belonging to~$\pd$ is dense in~$\pd$.
\end{lemma}

\subsection{$p$-Adic division quaternion algebras}
\label{ss:p-adic-quaternion}

Recall that there is a unique division quaternion algebra over~$\Qp$ up to isomorphism.
For the rest of this paper, fix such an algebra~$\Bp$.
Refer to~\cite{Vig80} and~\cite{Voi21} for background on quaternion algebras.

Let~$\bfB$ be an algebra over~$\Qp$ isomorphic to~$\bfB_p$.
Denote by~$1_{\bfB}$ its multiplicative identity, and identify~$\Qp$ with its image in~$\bfB$ by the map~$\ell \mapsto \ell \cdot 1_{\bfB}$.
Moreover, denote by~$g \mapsto \overline{g}$ the involution of~$\bfB$, and, for each~$g$ in~$\bfB$, denote by
\begin{displaymath}
  \tr(g) \= g + \overline{g},
  \nr(g) \= g \overline{g},
  \text{ and }
  \disc(g) \= \tr(g)^2 - 4 \nr(g),
\end{displaymath}
the \emph{reduced trace}, the \emph{reduced norm}, and the \emph{discriminant} of~$g$, respectively.
Each of these functions takes images in~$\Qp$.
Furthermore, the function ${\qval \: \bfB \to \Z \cup \{\infty\}}$ defined for~$g$ in~$\bfB$ by~$\qval(g) \= \ord_p(\nr(g))$, is the unique valuation extending the valuation~$2 \ord_p$ on~$\Qp$.
The valuation ring of~$\bfB$, defined by
\begin{displaymath}
  \bfR
  \=
  \{ g \in \bfB \: \qval(g) \ge 0 \},
\end{displaymath}
is the unique maximal $\Zp$-order in~$\bfB$ and coincides with the set of elements of~$\bfB$ that are integral over~$\Zp$.
The function ${\dB \: \bfB \times \bfB \to \R}$, defined for~$g$ and~$g'$ in~$\bfB$ by
\begin{displaymath}
  \dB(g, g')
  \=
  p^{- \frac{1}{2} \qval(g - g')},
\end{displaymath}
is an ultrametric distance on~$\bfB$ that makes~$\bfB$ into a topological algebra over~$\Qp$.
Put
\begin{displaymath}
  \bfG
  \=
  \{ g \in \bfB \: \qval(g) = 0 \}.
\end{displaymath}
Then, ${\bfG}$ is the group of units of~$\bfR$ and each right (resp.~left) multiplication map on~$\bfB$ by an element of~$\bfG$ is an isometry.

Section~\ref{s:CM} relies on the following consequence of the Skolem--Noether theorem.

\begin{lemma}
  \label{l:skew-generator}
  Let~$\bfB$ be an algebra over~$\Qp$ isomorphic to~$\bfB_p$, $\varphi$ in ${\bfB \ssetminus \Qp}$, and~$\theta$ in ${\Zp^{\times} \ssetminus \nr(\cO_{\Qp(\varphi)}^{\times})}$.
  Then, there is~$\gamma$ in~$\bfG$ satisfying
  \begin{equation}
    \label{eq:10}
    \gamma \varphi \gamma^{-1}
    =
    \overline{\varphi}
    \text{ and }
    \gamma^2 = \theta.
  \end{equation}
\end{lemma}

The proof of Lemma~\ref{l:skew-generator} relies on the following basic lemma.
For the reader's convenience, Appendix~\ref{ss:Apendix-A} includes a proof of (a more detailed version of) this lemma (Lemma~\ref{l:quadratic-extensions}$(ii)$).

\begin{lemma}
  \label{l:integer-norms}
  For every quadratic field extension~$\cK$ of~$\Qp$, the subgroup~$\nr(\cO_{\cK}^{\times})$ of~$\Zp^{\times}$ equals~$\Zp^{\times}$ if~$\cK$ is unramified over~$\Qp$ and has index two in~$\Zp^{\times}$ if~$\cK$ is ramified over~$\Qp$.
\end{lemma}

\begin{proof}[Proof of Lemma~\ref{l:skew-generator}]
  If~$\Qp(\varphi)$ is unramified over~$\Qp$, then~$\nr(\cO_{\Qp(\varphi)}^{\times}) = \Zp^{\times}$ by Lemma~\ref{l:integer-norms} and there is nothing to prove.
  Suppose that~$\Qp(\varphi)$ is ramified over~$\Qp$ and let~$\varpi$ be a uniformizer of~$\cO_{\Qp(\varphi)}$.
  Then, ${\ord_p(\nr(\varpi)) = 1}$ and~$\nr(\cO_{\Qp(\varphi)}^{\times})$ has index two in~$\Zp^{\times}$ by Lemma~\ref{l:integer-norms}.

  By \cite[\emph{Chapitre}~I, \emph{Corollaire}~2.2 and \emph{Corollaire}~2.4]{Vig80}, there exists a nonzero element~$\gamma_0$ of~$\bfB$ satisfying
  \begin{displaymath}
    \gamma_0 \varphi \gamma_0^{-1}
    =
    \overline{\varphi}
    \text{ and }
    \gamma_0^2 \in \Qp^{\times} \ssetminus \nr(\Qp(\varphi)^{\times}).
  \end{displaymath}
  In particular, ${\tr(\gamma_0) = 0}$ and ${\nr(\gamma_0) = - \gamma_0^2}$.
  Let~$\theta_0$ be in~$\Zp^{\times}$ and~$n$ in~$\Z$ such that ${\gamma_0^2 = \nr(\varpi)^n \theta_0}$.
  Then, $\theta_0$ is outside~$\nr(\cO_{\Qp(\varphi)}^{\times})$ and, since~$\nr(\cO_{\Qp(\varphi)}^{\times})$ has index two in~$\Zp^{\times}$, it follows that the quotient~$\theta / \theta_0$ belongs to~$\nr(\cO_{\Qp(\varphi)}^{\times})$.
  Let~$\rho$ be in~$\cO_{\Qp(\varphi)}$ such that ${\nr(\rho) = \theta / \theta_0}$ and put ${\gamma \= \gamma_0 \rho \varpi^{-n}}$.
  Then, for every~$\varphi'$ in~$\Qp(\varphi)$,
  \begin{displaymath}
    \gamma \varphi' \gamma^{-1}
    =
    \gamma_0 (\rho \varpi^{-n} \varphi' (\rho \varpi^{-n})^{-1}) \gamma_0^{-1}
    =
    \gamma_0 \varphi' \gamma_0^{-1}
    =
    \overline{\varphi'}.
  \end{displaymath}
  This applies to ${\varphi' = \rho \varpi^{-n}}$, so
  \begin{displaymath}
    \tr(\gamma)
    =
    \gamma_0 \rho \varpi^{-n} + \overline{\rho \varpi^{-n}} \overline{\gamma_0}
    =
    \gamma_0 \rho \varpi^{-n} - (\gamma_0 \rho \varpi^{-n} \gamma_0^{-1}) \gamma_0
    =
    0,
  \end{displaymath}
  and thus
  \begin{displaymath}
    \gamma^2
    =
    - \nr(\gamma)
    =
    - \nr(\gamma_0) \nr(\rho) \nr(\varpi)^{- n}
    =
    \gamma_0^2 (\theta / \theta_0) \nr(\varpi)^{- n}
    =
    \theta.
  \end{displaymath}
  In particular, ${\nr(\gamma) = - \theta}$, so~$\nr(\gamma)$ belongs to~$\Zp^{\times}$ and thus~$\gamma$ belongs to~$\bfG$.
\end{proof}

\subsection{Supersingular elliptic curves}
\label{ss:ss}
Denote by~$\tSups$ the finite subset of~$\Ell(\Fpalg)$ of supersingular elliptic curve classes.
Using ${j \: \Ell(\Fpalg) \to \Fpalg}$ to identify~$\Ell(\Fpalg)$ with~$\Fpalg$, the inclusion ${\tSups \subseteq \F_{p^2}}$ holds and~$\tSups$ is the zero set of a polynomial with coefficients in~$\Fp$, see, \emph{e.g.}, \cite{Deu41} and \cite[Chapter~V, Theorems~3.1 and~4.1]{Sil09}.
In particular, the Frobenius map ${\Frob \: \Fpalg \to \Fpalg}$ maps~$\tSups$ onto itself and induces an involution on this set.

For each~$\ss$ in~$\tSups$, denote by~$\End(\ss)$ and~$\Aut(\ss)$ the ring of endomorphisms and the group of automorphisms of~$\ss$ defined over~$\Fpalg$, respectively.
This paper relies several times on the mass formula of Deuring and Eichler,
\begin{equation}
  \label{eq:Deuring-Eichler}
  \sum_{\ss \in \tSups}\frac{1}{\#\Aut(\ss)}=\frac{p-1}{24},
\end{equation}
see, \emph{e.g.}, \cite{Eic55b} or \cite[Exercise~5.9]{Sil09}.

For all~$\ss$ and~$\ss'$ in~$\tSups$ and~$m$ in~$\Z_{>0}$, denote by~$\Hom_m(\ss, \ss')$ the set of isogenies from~$\ss$ to~$\ss'$ of degree~$m$.
If~$\ss$ is supersingular, then the ring ${\End(\ss) \otimes \Qp}$ is isomorphic to~$\Bp$.
For each~$g$ in~$\End(\ss)$, the discriminant~$\disc(g)$, viewed as an element of~$\End(\ss) \otimes \Zp$, belongs to~$\Zp$.

For each~$\ss$ in~$\tSups$, denote by~$\Dss$ the set of~$E$ in~$\Ell(\Cp)$ having good reduction, and such that the reduced elliptic curve is isomorphic to~$\ss$.
The set~$\Dss$ is a residue disc in~$\Ell(\Cp)$.

\subsection{Formal $\Zp$\nobreakdash-modules}
\label{ss:formal-modules}
This section reviews formal $\Zp$\nobreakdash-modules.
Refer to~\cite{Fro68,Haz78} for background.

Fix a complete, local, and Noetherian $\Zp$-algebra~$R_0$ with structural map ${\piem \: \Zp\to R_0}$, maximal ideal~$\cM_0$, and residue field isomorphic to a subfield~$\rfk$ of~$\Fpalg$.
Endow~$R_0$ with its natural $\cM_0$\nobreakdash-adic topology, fix a reduction morphism~$R_0\to \rfk$, and denote it by ${z \mapsto \widetilde{z}}$.
In most cases, $R_0$ is a subfield of~$\Fpalg$ or the ring of integers of a finite field extension of~$\Qp$ inside~$\Cp$, and the corresponding reduction morphism is the inclusion map or the restriction of~$\pi$.

For all formal groups~$\cF$ and~$\cF'$ defined over a ring~$R_0$, denote by~$\Hom_{R_0}(\cF,\cF')$ the set of morphisms ${\cF \to \cF'}$ defined over~$R_0$ and put ${\End_{R_0}(\cF) \= \Hom_{R_0}(\cF,\cF)}$.
Moreover, denote by~$\Iso_{R_0}(\cF, \cF')$ the set of isomorphisms ${\cF \to \cF'}$ defined over~$R_{0}$ and put ${\Aut_{R_0}(\cF)\=\Iso_{R_0}(\cF,\cF)}$.

For all rings~$R$ and~$R'$, ring morphism ${\sigma \: R \to R'}$, and formal power series~$f$ with coefficients in~$R$, the \emph{base change~$\sigma f$ of~$f$ under~$\sigma$} is the power series with coefficients in~$R'$ obtained by applying~$\sigma$ to the coefficients of~$f$.

For each formal group~$\cF$ over~$R_0$, denote by~$\tcF$ its reduction, which is the formal group over~$\rfk$ obtained as base change of~$\cF$ under the reduction map~$R_0\to \rfk$.
A \emph{formal $\Zp$\nobreakdash-module over~$R_0$} (resp.~$\rfk$) is a formal group~$\cF$ over~$R_0$ (resp.~$\rfk$) of dimension~1, together with a ring homomorphism ${\theta \: \Zp \to \End_{R_0}(\cF)}$ (resp. ${\theta \: \Zp \to \End_{\rfk}(\cF)}$) such that, in coordinates, for every~$\ell$ in~$\Zp$,
\begin{displaymath}
  \theta(\ell)(X) \equiv \piem(\ell) X \mod X^2
  \text{ (resp. ${\widetilde{\piem(\ell)} X \mod X^2}$)}.
\end{displaymath}

Every formal group~$\cF$ over~$R_0$ admits a unique structure of formal $\Zp$\nobreakdash-module over~$R_0$, such that the structural ring homomorphism~$\theta$ is continuous with respect to the $p$\nobreakdash-adic filtration on~$\Zp$ and the height filtration on $\End_{R_0}(\cF)$, see \cite[Chapter~IV, Section~1, proof of Theorem~1 and Chapter~III, Section~2, Corollary of Proposition~2]{Fro68}.

If~$R_1$ is another complete, local, and Noetherian $\Zp$\nobreakdash-algebra, $\cF$ is a formal~$\Zp$\nobreakdash-module over~$R_0$, and ${\sigma_0 \: R_0\to R_1}$ is a morphism of~$\Zp$-algebras, then~$\sigma_0\cF$ has a canonical structure of formal~$\Zp$\nobreakdash-module over~$R_1$.

\subsection{Deformation spaces of formal $\Zp$\nobreakdash-modules}
\label{ss:deformation-spaces}
This section reviews deformation theory of formal $\Zp$\nobreakdash-modules.
Refer to~\cite{Dri74,Haz78,HopGro94b} for background.

Let~$R_0$, $\piem$, $\cM_0$, and~$\rfk$ be as in Section~\ref{ss:formal-modules}, $\Bbbk$ a subfield of~$\rfk$, and~$\cF_0$ a formal~$\Zp$\nobreakdash-module over~$\Bbbk$.
A \emph{deformation of~$\cF_0$ over~$R_0$} is a pair~$(\cF,\alpha)$, where~$\cF$ is a formal $\Zp$\nobreakdash-module over~$R_0$ and~$\alpha \: \tcF\to \cF_0$ is an isomorphism of formal~$\Zp$\nobreakdash-modules defined over~$\rfk$.
Two such deformations~$(\cF,\alpha)$ and~$(\cF',\alpha')$ are \emph{isomorphic}, if there exists an isomorphism~$\varphi$ in $\Iso_{R_0} ( \cF , \cF')$ such that $\alpha'\circ \tvarphi = \alpha$.
Denote by~$\bfX(\cF_0,R_0)$ the set of isomorphism classes of deformations of~$\cF_0$ over~$R_0$.

The work of Lubin and Tate~\cite{LubTat66} and \mbox{Drinfel'd}~\cite{Dri74}, see also \cite[Section~12]{HopGro94b} or \cite[Section~21.5]{Haz78}, yields
a formal $\Zp$\nobreakdash-module~$\cF^{\star}(t)$ over~$\Zp[\![t]\!]$, called a \emph{universal} formal $\Zp$\nobreakdash-module of height two, satisfying the following properties:
\begin{enumerate}
\item[$(i)$]
  The reduction~$\widetilde{\cF^{\star}(0)}$ is a formal~$\Zp$\nobreakdash-module over~$\Fp$ of height two.
\item[$(ii)$]
  The~$p$-th power Frobenius endomorphism~$\varphi$ on~$\widetilde{\cF^{\star}(0)}$, given in coordinates by ${\varphi(X) = X^p}$, satisfies the relation~$\varphi^2=-p$ in~$\End_{\Fpalg}(\widetilde{\cF^{\star}(0)})$.
\item[$(iii)$]
  Denoting by~$\Id$ the identity automorphism of~$\widetilde{\cF^{\star}(x)}$, which is equal to~$\widetilde{\cF^{\star}(0)}$, the map
  \begin{equation}
    \label{eq:parametrization-X-e}
    \begin{array}{rcl}
      \cM_{0} & \to & \bfX(\widetilde{\cF^{\star}(0)},R_0)
      \\
      x & \mapsto & (\cF^{\star}(x),\Id),
    \end{array}
  \end{equation}
  is a bijection.
\end{enumerate}
Moreover, the bijection~\eqref{eq:parametrization-X-e} is functorial on~$R_0$, see \cite[Theorem~21.5.6]{Haz78}.
As a consequence of~$(iii)$, the set~$\bfX(\widetilde{\cF^{\star}(0)},R_0)$ is parametrized by the open disc~$\cM_{0}$ of~$R_0$.
Use this parametrization to endow~$\bfX(\widetilde{\cF^{\star}(0)},R_0)$ with the topology coming from the topology on~$\cM_{0}$ inherited from~$R_0$.

For all~$\Zp$-formal modules~$\cF_0$ and~$\cF_0'$ defined over a subfield of~$\rfk$, consider the natural map
\begin{displaymath}
  \begin{array}{rcl}
    \Iso_{\rfk}(\cF_0,\cF_0') \times \bfX(\cF_0,R_0) & \to & \bfX(\cF_0',R_0)
    \\
    (\beta, (\cF,\alpha)) & \mapsto & (\cF, \beta \circ \alpha)
  \end{array}
\end{displaymath}
and define ${\beta \cdot (\cF,\alpha) \= (\cF, \beta \circ \alpha)}$.
In particular, the group~$\Aut_{\rfk}(\cF_0)$ acts on~$\bfX(\cF_0,R_0)$.
Fixing the parametrization~\eqref{eq:parametrization-X-e} yields an action of the group~$\Aut_{\rfk}(\widetilde{\cF^{\star}(0)})$ on~$\cM_0$.

The rest of the paper relies several times on the following lemma.

\begin{lemma}[\textcolor{black}{\cite[Proposition~14.13]{HopGro94b}}]
  \label{l:fixed}
  An element~$g$ of~$\Aut_{\rfk}(\cF_0)$ fixes the point~$(\cF,\alpha)$ in~$\bfX(\cF_0,R_0)$ if and only if~$g$ belongs to the image of the injective group homomorphism ${\Aut_{R_0}(\cF) \to \Aut_{\rfk}(\cF_0)}$ given by ${\varphi \mapsto \alpha \circ \tvarphi \circ \alpha^{-1}}$.
\end{lemma}

\subsection{From elliptic curves to formal~$\Zp$\nobreakdash-modules}
\label{ss:from-elliptic-curves}
Let~$R$ be either a subfield of~$\Fpalg$ or the ring of integers of a finite field extension of~$\Qp$ inside~$\Cp$.
For each elliptic curve~$\ss$ given by a Weierstrass equation with coefficients in~$R$ and having smooth reduction if~$R$ has characteristic zero, denote by~$\Fss$ the formal group of~$\ss$ considered as a formal $\Zp$\nobreakdash-module, see, \emph{e.g.}, \cite[Section~4]{Blu98}.
Given such~$\ss$ and~$\ss'$, denote by ${\phi \mapsto \hphi}$ the natural morphism ${\Hom_{R}(\ss, \ss') \to \Hom_{R}(\Fss,\Fsspr)}$.
This morphism is injective and compatible with addition and composition, see \cite[Proposition~5.1]{Blu98}.
When ${\ss' = \ss}$, it is a ring homomorphism ${\End_{R}(\ss) \to \End_{R}(\Fss)}$.

Fix a universal formal $\Zp$\nobreakdash-module~$\cF^{\star}(t)$ as in Section~\ref{ss:deformation-spaces}.
Recall that if~$\ss$ is a supersingular elliptic curve over~$\Fpalg$, then the height of~$\Fss$ is two, see, \emph{e.g.}, \cite[Chapter~V, Theorem~3.1]{Sil09}.
Moreover, $\Q_{p^2}$ denotes the unique unramified quadratic field extension of~$\Qp$ inside~$\Cp$, and~$\Z_{p^2}$ and~$\F_{p^2}$ the ring of integers and the residue field of~$\Q_{p^2}$, respectively.

\begin{lemma}
  \label{l:supersingular-fixation}
  For every~$\ss$ in~$\tSups$, there is an elliptic curve~$\ss_0$ defined over~$\F_{p^2}$ representing~$\ss$ and such that there exists an isomorphism ${\phi_0 \: \widetilde{\cF^{\star}(0)}\to \cF_{\ss_0}}$ of formal~$\Zp$\nobreakdash-modules defined over~$\F_{p^2}$.
\end{lemma}

\begin{proof}
  Let~$\ss_0$ be an elliptic curve representing~$\ss$ that is given by a Weierstrass equation with coefficients in~$\F_{p^2}$ such that the~$p^2$-th power Frobenius endomorphism~$\Frob^2$ of~$\ss_0$, defined in affine coordinates by ${\Frob^2(x,y) = (x^{p^2},y^{p^2})}$, satisfies ${\Frob^2 = -p}$ in~$\End(\ss_0)$, see, \emph{e.g.}, \cite[Lemma~3.21]{BakGonGonPon05}.
  Thus, the minimal polynomial of~$\Frob^2$ over~$\Zp$ is~$X + p$.
  Since the minimal polynomial of the~$p^2$-th Frobenius endomorphism of~$\widetilde{\cF^{\star}(0)}$ over~$\Zp$ is the same, the existence of~$\phi$ follows from \cite[Proposition~24.2.9]{Haz78}.
\end{proof}

For each~$\ss$ in~$\tSups$, fix~$\ss_0$ and~$\phi_0$ as in Lemma~\ref{l:supersingular-fixation} above and put
\begin{equation}
  \label{eq:11}
  \Fss
  \=
  \cF_{\ss_0},
  \phi_{\ss}
  \=
  \phi_0,
  \Rss
  \=
  \End_{\Fpalg}(\Fss),
  \text{ and }
  \Gss
  \=
  \Aut_{\Fpalg}(\Fss).
\end{equation}
The algebra~$\Bss$ over~$\Qp$, defined by
\begin{displaymath}
  \Bss
  \=
  \End_{\Fpalg}(\Fss)\otimes \Qp,
\end{displaymath}
is isomorphic to~$\Bp$.
Endow~$\Bss$ with its natural metric, as in Section~\ref{ss:p-adic-quaternion}, and identify~$\Rss$ and~$\Gss$ with the unique maximal $\Zp$-order in~$\Bss$ and with the group of units of this order, respectively, see, \emph{e.g.}, \cite[Chapter~III, Section~2, Theorem~3]{Fro68}.
In particular, both of these sets are metric subspaces of~$\Bss$ and thus each right (resp.~left) multiplication map on~$\Rss$ by an element of~$\Gss$ is an isometry.

Since~$\End(\ss)$ is a maximal order in~$\End(\ss)\otimes \Q$, and being a maximal order is a local property, see, \emph{e.g.}, \cite[\emph{Chapitre}~III, Section~5.A]{Vig80}, it follows that the natural map ${\End(\ss) \otimes \Zp\to \Rss}$ is an isomorphism.
This natural map extends to an isomorphism ${\End(\ss) \otimes \Qp\to \Bss}$.

For each finite field extension~$\cK$ of~$\Q_{p^2}$ inside~$\Cp$, put ${\Xss(\OK) \= \bfX(\Fss,\OK)}$.
The map
\begin{displaymath}
  \begin{array}{rcl}
    \bfX(\widetilde{\cF^{\star}(0)},\OK) & \to & \Xss(\OK)
    \\
    (\cF,\alpha) & \mapsto & \phi_{\ss} \cdot(\cF,\alpha) = (\cF,\phi_{\ss}\circ \alpha).
  \end{array}
\end{displaymath}
is a bijection.
If~$\cK'$ is a finite field extension of~$\cK$ inside~$\Cp$, then each deformation of~$\Fss$ over~$\OK$ can be considered as a deformation of~$\Fss$ over~$\cO_{\cK'}$, and this induces a natural map ${\Xss(\OK)\to \Xss(\cO_{\cK'})}$ that is injective \cite[Proposition~12.10]{HopGro94b}.
Consider the direct limit
\begin{displaymath}
  \Xss (\OQpalg)
  \=
  \varinjlim \Xss(\OK),
\end{displaymath}
over the directed set of finite field extensions~$\cK$ of~$\Q_{p^2}$ inside~$\Cp$ ordered by inclusion.
As~$\cK$ runs through the finite field extensions of~$\Q_{p^2}$ inside~$\Cp$, the parametrization of~$\Xss(\OK)$ by~$\MK$, given by~\eqref{eq:parametrization-X-e}, and the action of~$\Aut_{\pi(\OK)}(\cF_{\ss})$ on~$\Xss(\OK)$, defined in Section~\ref{ss:deformation-spaces}, induce a parametrization of~$\Xss(\OQpalg)$ by~$\MQpalg$ and a group action of~$\Gss$ on~$\Xss(\OQpalg)$.
The field of definition of an element of~$\Xss(\OQpalg)$ is determined by the corresponding parameter in~$\MQpalg$, in the following sense: For every~$x$ in~$\MQpalg$, the associated deformation of~$\Fss$ can be defined over~$\cO_{\Q_{p^2}(x)}$ and this is the smallest field extension of~$\Q_{p^2}$ inside~$\Cp$ where this deformation can be defined.

Fix a completion~$\hDss$ of~$\Xss(\OQpalg)$.
Then, the parametrization of~$\Xss(\OQpalg)$ by~$\MQpalg$ extends to a parametrization of~$\hDss$ by~$\Mp$.
The following lemma implies that the action of~$\Gss$ on~$\Xss(\OQpalg)$ extends to a continuous action ${\Gss \times \hDss \to \hDss}$ that is rigid analytic in the second variable.

\begin{lemma}[\textcolor{black}{\cite[Section~14, Proposition~19.2 and Lemma~19.3]{HopGro94b}}]
  \label{l:action-regularity}
  For every~$\ss$ in~$\tSups$, the following properties hold.
  \begin{enumerate}
  \item [$(i)$]
    Each element of~$\Gss$ acts on~$\hDss$ as a rigid analytic automorphism represented by a power series with coefficients in~$\Z_{p^2}$.
    In particular, $\Gss$ acts by rigid analtyic isometries on~$\hDss$.
  \item [$(ii)$]
    For all~$N$ in~$\Z_{\ge 0}$, $r$ in~$\Z_{>0}$, element~$g$ of~$\Gss$ in~$1+p^N\Rss$, and~$x$ in~$\hDss$ satisfying~$\ord_p(x) \ge \frac{1}{r}$,
    \begin{displaymath}
      \ord_p(x - g \cdot x)
      \ge
      \frac{N + 1}{r}.
    \end{displaymath}
  \end{enumerate}
\end{lemma}

\subsection{From formal $\Zp$-modules to elliptic curves}
\label{ss:woods-hole}

Let~$R_0$, $\piem$, $\cM_0$, $\rfk$, and~$\Bbbk$ be as in Section~\ref{ss:deformation-spaces}, and let~$\ss$ be an elliptic curve defined over~$\Bbbk$.
Denote by~$\bfY(\ss, R_0)$ the space of isomorphism classes of pairs~$(E,\alpha)$ formed by an elliptic curve~$E$ given by a Weierstrass equation with coefficients in~$R_0$ and having smooth reduction, and an isomorphism ${\alpha \: \widetilde{E}\to \ss}$ defined over~$\rfk$, where two pairs~$(E,\alpha)$ and~$(E',\alpha')$ are isomorphic if there exists an isomorphism ${\psi \: E \to E'}$ defined over~$\rfk$ such that ${\alpha'\circ \widetilde{\psi}=\alpha}$.
There is a natural action of~$\Aut_{\rfk}(\ss)$ on~$\bfY(\ss, R_0)$ given, for~$\phi$ in~$\Aut_{\rfk}(\ss)$, by ${\phi\cdot (E,\alpha) = (E,\phi\circ \alpha)}$.

There is a natural map
\begin{displaymath}
  \bfY(\ss,R_0) \to \bfX(\Fss,R_0)
\end{displaymath}
mapping a class in~$\bfY(\ss,R_0)$ represented by a pair~$(E,\alpha)$ to the class in~$\bfX(\Fss,R_0)$ represented by the deformation~$(\FE, \halpha)$.
This map is a bijection by the so-called Woods Hole Theory, see \cite[Section~6]{LubSerTat64} or \cite[Theorem~4.1]{ColMcM10}.
Together with the group homomorphism ${\Aut_{\Bbbk}(\ss) \to \Aut_{\Bbbk}(\Fss)}$ given by ${\phi \mapsto \hphi}$, this bijection yields an identification
\begin{displaymath}
  \Aut_{\Bbbk}(\ss)\backslash \bfY(\ss,R_0) \xrightarrow{\sim} \Aut_{\Bbbk}(\ss) \backslash \bfX(\Fss,R_0).
\end{displaymath}
Taking direct limits as~$R_0$ runs through the ring of integers of finite field extensions of~$\Q_{p^2}$ inside~$\Cp$ yields an identification
\begin{equation}
  \label{eq:equivalence-formal}
  \{E\in \Ell_{\sups}(\Qpalg) \: \tE \text{ isomorphic to~$\ss$ over }\Fpalg\} \to \Aut(\ss)\backslash \Xss(\OQpalg).
\end{equation}
Identifying the left side above with ${\Qpalg \cap \Dss}$ and composing the natural projection from~$\Xss(\OQpalg)$ to~$\Aut(\ss)\backslash \Xss(\OQpalg)$ with the inverse of~\eqref{eq:equivalence-formal} yields a map
\begin{equation}
  \label{eq:12}
  \Piss \: \Xss(\OQpalg) \to \Ell_{\sups}(\Qpalg) \cap \Dss.
\end{equation}

Put ${\delta_{\ss} \= \# \Aut(\ss) /2}$, so ${\delta_{\ss} = 1}$ if~$j(e)$ is different from~$0$ and~$1728$ and, in all cases, ${1 \le \delta_{\ss} \le 12}$, see, \emph{e.g.}, \cite[Appendix~A, Proposition~1.2(c)]{Sil09}.

\begin{theorem}
  \label{t:Pi-properties}
  Fix~$\ss$ in~$\tSups$.
  Then, \eqref{eq:12} extends to a map
  \begin{displaymath}
    \Piss \: \hDss \to \Dss
  \end{displaymath}
  such that~$j \circ \Piss$ is represented by a power series with coefficients in~$\Z_{p^2}$ that is a ramified covering of degree~$\delta_{\ss}$.
  Moreover, for all~$x$ in~$\hDss$ and~$E$ in~$\Dss$,
  \begin{equation}
    \label{eq:13}
    \min \{ |x - x'|_p \: x' \in \Piss^{-1}(E) \}^{\delta_{\ss}}
    \le
    |j(\Piss(x)) - j(E)|_p
    \le
    \min \{ |x - x'|_p \: x' \in \Piss^{-1}(E) \}.
  \end{equation}
  In particular, $j \circ \Piss$ is an isometry if~$j(\ss)$ is different from~$0$ and~$1728$.
\end{theorem}

\begin{proof}
  For the first assertion, the ring~$\Z_{p^2}\left[\![t]\!\right]$ is a complete, local, and Noetherian~$\Zp$-algebra whose residue field is isomorphic to~$\F_{p^2}$.
  Using the parametrization of~$\bfX(\Fss, \Z_{p^2}\left[\![t]\!\right])$ by the maximal ideal of~$\Z_{p^2}\left[\![t]\!\right]$, the element~$t$ corresponds to the deformation~$(\cF^{\star}(t), \phi_{\ss})$ of~$\Fss$.
  Denote by~$E(t)$ the elliptic curve class in~$\bfY(\ss, \Z_{p^2}\left[\![t]\!\right])$ corresponding to the element ${\Aut(\ss) \cdot (\cF^{\star}(t), \phi_{\ss})}$ of ${\Aut(\ss) \backslash \bfX(\Fss, \Z_{p^2}\left[\![t]\!\right])}$.
  Since the $j$-invariant~$j(E(t))$ of~$E(t)$ is an element of~$\Z_{p^2}\left[\![t]\!\right]$, to prove the first assertion of the theorem it suffices to prove that, for every~$x$ in~$\cM_{\Qpalg}$, the equality ${\Piss (\cF^{\star}(x), \phi_{\ss}) = E(x)}$ holds.
  Let ${\psi \: \Z_{p^2}\left[\![t]\!\right] \to \OQpalg}$ by the evaluation map, defined by ${\psi(f(t)) \= f(x)}$, which is a continuous ring homomorphism.
  Moreover, denote by~$\tpsi$ the induced morphism on residue fields.
  Then, ${(\psi \cF^{\star}(t), \tpsi \phi_{\ss}) = (\cF^{\star}(x), \tpsi \phi_{\ss})}$ and the orbit ${\Aut(\ss) \cdot (\psi \cF^{\star}(t), \tpsi \phi_{\ss})}$ corresponds to the base change of~$E(t)$ under~$\psi$, which is~$E(x)$.
  This proves the equality ${\Piss(\cF^{\star}(x), \phi_{\ss}) = E(x)}$.

  To prove that~$j \circ \Piss$ is a ramified covering of degree~$\delta_{\ss}$, it suffices to prove that, for every~$E$ in ${\SupsQpalg \cap \Dss}$ such that~$j(E)$ is different from~$0$ and~$1728$, the equality ${\# \Piss^{-1}(E) = \delta_{\ss}}$ holds.
  Denote by~$\bfone_{\ss}$ the identity in~$\Aut(\ss)$, so~$\bfone_{\ss}$ and~$- \bfone_{\ss}$ act trivially on~$\Xss(\OQpalg)$.
  Thus, for each~$x$ in~$\Piss^{-1}(E)$, the stabilizer of~$x$ for the action of~$\Aut(\ss)$ on~$\Xss(\OQpalg)$ contains~$\{ \bfone_{\ss}, - \bfone_{\ss} \}$.
  Let~$\phi$ in~$\Aut(\ss)$ be in the stabilizer of~$x$ and let~$\alpha \: \tE \to \ss$ be an isomorphism such that~$(\FE, \halpha)$ represents~$x$.
  By Lemma~\ref{l:fixed}, there is~$\varphi$ in~$\Aut(E)$ such that ${\alpha \circ \tvarphi \circ \alpha^{-1} = \phi}$.
  Together with the hypothesis that~$j(E)$ is different from~$0$ and~$1728$, this implies that~$\varphi$ or~$-\varphi$ is the identity, see, \emph{e.g.}, \cite[Appendix~A, Proposition~1.2(c)]{Sil09}.
  It follows that~$\phi$ belongs to~$\{ \bfone_{\ss}, - \bfone_{\ss} \}$.
  This proves that the stabilizer of each element of~$\Piss^{-1}(E)$ equals~$\{ \bfone_{\ss}, - \bfone_{\ss} \}$.
  In particular, ${\# \Piss^{-1}(E) = \delta_{\ss}}$.

  It remains to prove~\eqref{eq:13}.
  Let~$E$ be in~$\Dss$ and~$x_1$, \ldots, $x_{\delta_{\ss}}$ the zeros of ${j \circ \Piss - j(E)}$ repeated according to multiplicity.
  Then, there is~$h(t)$ in~$\Z_{p^2} \left[\![t]\!\right]$ such that~$|h|_p$ is constant equal to~$1$ on~$\Op$ and
  \begin{displaymath}
    j \circ \Piss (t) - j(E)
    =
    h(t) \cdot \prod_{i = 1}^{\delta_{\ss}} (t - x_i),
  \end{displaymath}
  see, \emph{e.g.}, \cite[Exercise~3.2.2(1)]{FrevdP04}.
  Since for all~$i$ in~$\{1, \ldots, \delta_{\ss} \}$ and~$x$ in~$\hDss$,
  \begin{displaymath}
    \min \{ |x - x'|_p \: x' \in \Piss^{-1}(E) \}
    \le
    |x - x_i|_p
    \le
    1,
  \end{displaymath}
  \eqref{eq:13} follows.
\end{proof}

\subsection{Hecke correspondences}
\label{ss:Hecke-correspondences}
This section recalls the construction and main properties of Hecke correspondences.
For details, refer to \cite[Sections~7.2 and~7.3]{Shi71} for the general theory or the survey~\cite[Part~II]{DiaIm95}.

Let~$\K$ be an algebraically closed field of characteristic zero.
For all~$n$ in~$\Z_{>0}$ and divisor~$\fD$ in~$\Div(\Ell(\K))$,
\begin{equation}
  \label{eq:deg-Hecke-op-divisor}
  \deg(T_n(\fD))
  =
  \sigma_1(n) \deg(\fD).
\end{equation}
For ${n = 1}$, the correspondence~$T_1$ is the identity on~$\Div(\Ell(\K))$ by definition.
Extend linearly Hecke correspondences to~$\Div(\Ell(\K)) \otimes \Q$.

For each~$N$ in~$\Z_{>0}$, denote by~$\Ell_0(N)$ the \emph{modular curve of level~$N$}.
It is a quasi-projective variety defined over~$\Q$.
The points of~$\Ell_0(N)$ over~$\K$ parametrize the moduli space of equivalence classes of pairs~$(E,C)$, where~$E$ is an elliptic curve over~$\K$, $C$ a cyclic subgroup of~$E$ of order~$N$, and two such pairs~$(E,C)$ and~$(E',C')$ are equivalent if there exists an isomorphism ${E \to E'}$ over~$\K$ taking~$C$ to~$C'$.
In particular, for~$N = 1$ and every algebraically closed field~$\K$, the curve~$\Ell_0(1)(\K)$ parametrizes~$\Ell(\K)$, and~$\Ell_0(1)$ is isomorphic to the affine line $\A_{\Q}^1$.

For ${N > 1}$, denote by~$\Phi_N(X,Y)$ the \emph{modular polynomial of level~$N$}, which is a symmetric polynomial in~$\Z[X,Y]$ that is monic in both~$X$ and~$Y$, see, \emph{e.g.}, \cite[Chapter~5, Sections~2 and~3]{Lan87}.
This polynomial is characterized by the equality
\begin{equation}\label{eq:mod-polynomial}
  \Phi_N(j(E),Y)
  =
  \prod_{C\leq E \text{ cyclic of order }N}(Y-j(E/C)) \text{ for every } E \text{ in } \Ell(\K).
\end{equation}
This implies that the plane algebraic curve defined by
\begin{equation}
  \label{eq:mod-equation}
  \Phi_N(X,Y)
  =
  0
\end{equation}
provides a birational model for~$\Ell_0(N)$.
For each prime number~$q$, let ${\alpha_q \: \Ell_0(q)\to \Ell_0(1)}$ and ${\beta_q \: \Ell_0(q)\to \Ell_0(1)}$ be the rational maps defined over~$\Q$ given in terms of moduli spaces by
\begin{displaymath}
  \alpha_q(E,C) \= E
  \text{ and }
  \beta_q(E,C) \= E/C.
\end{displaymath}
In terms of the model~\eqref{eq:mod-equation} with~$N = q$, the rational maps~$\alpha_q$ and~$\beta_q$ correspond to the projections onto the~$X$ and~$Y$ coordinates, respectively.
Denote by~$(\alpha_q)_*$ and~$(\beta_q)_*$ the push-forward action of~$\alpha_q$ and~$\beta_q$ on divisors, respectively.
Denote also by~$\alpha_q^{*}$ the pull-back action of~$\alpha_q$ on divisors, defined for~$x$ in~$\Ell_0(1)(\K)$ by
\begin{displaymath}
  \alpha_q^{*}(x)
  \=
  \sum_{\substack{y\in \Ell_0(q)(\K) \\ \alpha_q(y)=x}} \deg_{\alpha_q}(y) y,
\end{displaymath}
where~$\deg_{\alpha_q}(y)$ is the local degree of~$\alpha_q$ at~$y$.
This definition extends by linearity to arbitrary divisors.
The pull-back action~$\beta_q^{*}$ of~$\beta_q$ is defined in a similar way.
Then, the Hecke correspondence ${T_q \: \Div(\Ell(\K))\to \Div(\Ell(\K))}$ is recovered as
\begin{equation}
  \label{eq:algebraic definition}
  T_q
  =
  (\alpha_q)_* \circ \beta_q^*
  =
  (\beta_q)_* \circ \alpha_q^*,
\end{equation}
where the second equality follows from the first and from the symmetry of~$T_q$.

For every integer~$n$ satisfying ${n \ge 2}$, the correspondence~$T_n$ can be recovered from different~$T_q$'s, for~$q$ running over prime divisors of~$n$, by using the identities
\begin{equation}
  \label{eq:multiplicativity-Hecke-operators}
  T_{\ell} \circ T_{m}
  =
  T_{\ell m} \, \text{ for coprime $\ell$ and~$m$ in~$\Z_{>0}$};
\end{equation}
\begin{equation}\label{eq:Hecke-Tpr}
  T_{q^r}
  =
  T_q \circ T_{q^{r-1}} - q \cdot T_{q^{r-2}} \text{ for every prime number~$q$ and integer~$r$ satisfying ${r \ge 2}$}.
\end{equation}

The action of Hecke correspondences on sets and compactly supported measures is defined as follows.
Recalling that~$\Ell(\Cp)$ is identified with~$\Cp$ with the ${j}$-invariant, put, for all~$n$ in~$\Z_{>0}$ and subset~$A$ of~$\Cp$,
\begin{displaymath}
  T_n(A)
  \=
  \bigcup_{a \in A} \supp(T_n(a)).
\end{displaymath}
This defines an action of~$T_n$ on sets that is compatible with the action on effective divisors: ${T_n(\supp(\fD)) = \supp(T_n \fD)}$.
To state further properties of this action, recall that~$T_n$ acts on the space~$C_b(\Cp)$ of continuous and bounded functions ${F \: \Cp \to \R}$ by ${T_n F(E) \= F(T_n(E))}$, see, \emph{e.g.}, \cite[Lemma~2.1]{HerMenRiv20}.
Standard approximation arguments show that the image of an open (resp.~closed, compact) set by~$T_n$ is a set of the same nature.
Moreover, the action of~$T_n$ on~$C_b(\Cp)$ is continuous.
Hence, for each Borel measure~$\mu$ on~$\Cp$ whose support is compact, the linear functional ${F \mapsto \int T_nF \dd \mu}$ is continuous, and thus defines a Borel measure on~$\Cp$ supported on the compact set~$T_n(\supp(\mu))$.
It is the push-forward of~$\mu$ by~$T_n$, denoted by~$(T_n)_*\mu$.
Then, ${\supp((T_n)_*\mu) = T_n(\supp(\mu))}$ and the total mass of~$(T_n)_* \mu$ equals~$\sigma_1(n)$ times the total mass of~$\mu$.

For all~$x$ in~$\Xss(\OQpalg)$, $n$ in ${\Z_{>0} \ssetminus p \Z}$, and isogeny~$\phi$ in~$\Hom_n(\ss, \ss')$, the isomorphism~$\hphi$ belongs to~$\Iso_{\Fpalg}(\Fss, \Fsspr)$ and
\begin{equation}
  \label{eq:14}
  T_n(\Piss(x))|_{\bfD(\ss')}
  =
  \frac{1}{\# \Aut(\ss')} \sum_{\phi \in \Hom_n(\ss, \ss')} \Pisspr(\hphi \cdot x).
\end{equation}
By continuity of~$T_n$, this holds for every~$x$ in~$\hDss$, see, \emph{e.g.}, \cite[Lemma~2.1]{HerMenRiv20}.

\section{Asymptotic distribution of integer points on $p$-adic spheres}
\label{s:p-adic-Linnik}
This section proves Theorem~\ref{t:p-adic-Linnik} in Section~\ref{s:p-adic-linnik}.
The proof relies on the following theorem.
Let~$n$, $Q$, $V_m(Q)$, $S_{\ell}(Q)$, and~$\oO_Q(\Zp)$ be as in Section~\ref{s:p-adic-linnik}.
For each~$r$ in~$\Z_{>0}$, denote by ${\red_r \: \Zp^n \to (\Z / p^r \Z)^n}$ the reduction map and~$\oO_Q(\Z / p^r \Z)$ the corresponding orthogonal group of~$Q$.
This group is finite, acts on the finite set~$(\Z / p^r \Z)^n$, and, for every~$\ell$ in~$\Zp^{\times}$, it leaves~$\red_r(S_{\ell}(Q))$ invariant.

\begin{theorem}[Modular deviation estimate]
  \label{t:finitary-deviation-estimate}
  Let~$\kappa_n$ be as in Theorem~\ref{t:p-adic-Linnik} and fix~$r$ in~$\Z_{>0}$.
  Then, for every~$\varepsilon$ in~$\R_{>0}$ if~$n \ge 4$, and all~$\varepsilon$ in~$\R_{>0}$ and~$S$ in~$\R$ satisfying ${S \ge 1}$ if~$n = 3$, there is~$C$ in~$\R_{>0}$ such that the following property holds.
  Let~$\Sigma$ be an orbit of~$\oO_Q(\Z / p^r\Z)$ in~$(\Z / p^r \Z)^n$ and~$m$ in~$\Z_{>0}$ such that
  \begin{displaymath}
    V_m(Q) \neq \emptyset
    \text{ and }
    \red_r(V_m(Q)) \subseteq \Sigma.
  \end{displaymath}
  If~$n = 3$, then suppose in addition that the largest square diving~$m$ is less than~$S$.
  Then, for every~$\sigma$ in~$\Sigma$,
  \begin{displaymath}
    \left| \frac{ \# \left\{ \bfx \in V_m(Q) \: \red_r(\bfx) = \sigma \right\} }{\#V_m(Q)} - \frac{1}{\# \Sigma} \right|
    \le
    C \frac{m^{\frac{n}{4}-\kappa_n + \varepsilon}}{\#V_m(Q)}.
  \end{displaymath}
\end{theorem}

Theorem~\ref{t:finitary-deviation-estimate} combined with Hensel's lemma yields the following corollary.
See Section~\ref{ss:proof of finitary deviation estimate} for its proof.
Endow~$\Zp^n$ with a distance compatible with its product topology and let~$\ell$ in~$\Zp \ssetminus \{ 0 \}$ be such that ${S_{\ell}(Q) \neq \emptyset}$ and the compact group~$\oO_Q(\Zp)$ acts transitively on~$S_{\ell}(Q)$.
As in the statement of Theorem~\ref{t:p-adic-Linnik}, denote by~$\mu_{\ell}$ the unique Borel probability measure on~$S_{\ell}(Q)$ that is invariant under the action of~$\oO_Q(\Zp)$, see, \emph{e.g.}, Lemma~\ref{l:existence-invariant-measure}.
This measure is uniquely determined by the property that, for all~$r$ in~$\Z_{>0}$ and~$\sigma$ in~$\red_r(S_{\ell}(Q))$,
\begin{equation}
  \label{eq:15}
  \mu_{\ell} ( S_{\ell}(Q) \cap \red_r^{-1}(\sigma) )
  =
  \frac{1}{\# \red_r(S_{\ell}(Q))}.
\end{equation}

\begin{coro}
  \label{c:functional-deviation}
  Let~$\kappa_n$ be as in Theorem~\ref{t:p-adic-Linnik} and let~$\delta$ be in~$\R_{>0}$.
  Then, for every~$\varepsilon$ in~$\R_{>0}$ if~$n \ge 4$, and all~$\varepsilon$ in~$\R_{>0}$ and~$S$ in~$\R$ satisfying ${S \ge 1}$ if~$n = 3$, there are~$C$ in~$\R_{>0}$ and~$N$ in~$\Z_{\ge 0}$, such that the following property holds.
  Let~$\ell$ in~$\Zp \ssetminus \{ 0 \}$ be such that ${S_{\ell}(Q) \neq \emptyset}$ and~$\oO_Q(\Zp)$ acts transitively on~$S_{\ell}(Q)$ and~$m$ in~$\Z_{>0}$ such that
  \begin{displaymath}
    m\equiv \ell \mod p^{N}
    \text{ and }
    V_m(Q)\neq \emptyset.
  \end{displaymath}
  If~$n = 3$, then suppose in addition that the largest square diving~$m$ is less than or equal to~$S$.
  Then, for every function ${F \: \Zp^n \to \R}$ that is constant on every ball of radius~$\delta$,
  \begin{displaymath}
    \left| \frac{1}{\#V_m(Q)} \sum_{\bfx \in V_m(Q)}F(\bfx) - \int F \dd\mu_{\ell} \right|
    \le
    C \left( \sup_{\Zp^n} |F| \right) \frac{m^{\frac{n}{4}-\kappa_n + \varepsilon}}{\#V_m(Q)}.
  \end{displaymath}
\end{coro}

The proof of Theorem~\ref{t:finitary-deviation-estimate} rephrases this result as a norm estimate on a certain finite-dimensional $L^2$-function space (Lemma~\ref{l:devformula}), as described, \emph{e.g.}, in \cite[Section~11.3]{EllMicVen13}.
The main ingredient to prove this $L^2$-norm estimate is the construction of an auxiliary modular form that has the key property of being cuspidal (Proposition~\ref{p:modular form}).
The $L^2$-norm estimate is then deduced from the bounds for Fourier coefficients of cuspidal modular forms shown by Deligne for~$n$ even~\cite{Del74}, Iwaniec~\cite{Iwa87} for~$n\ge 5$ odd, and Duke~\cite{Duk88} and Blomer \cite[Lemma~4.4]{Blo04} for ${n = 3}$.

Section~\ref{ss:modular form} defines the auxiliary modular form and proves it is cuspidal.
Section~\ref{ss:proof of finitary deviation estimate} derives Theorems~\ref{t:p-adic-Linnik} and~\ref{t:finitary-deviation-estimate} from this and proves Corollary~\ref{c:functional-deviation}.

\subsection{Auxiliary modular form}
\label{ss:modular form}
For each row vector~$\bfx$, denote by~$\bfx^\intercal$ its transpose.
Let~$n$ be an integer satisfying ${n \ge 3}$, $Q$ a positive definite quadratic form in~$\Z[X_1, \ldots, X_n]$, and~$A_Q$ the symmetric matrix in~$M_n(\Z)$ such that ${Q(\bfx) = \frac{1}{2} \bfx A_Q \bfx^\intercal}$.
Each diagonal entry of~$A_Q$ is even.
The \emph{level}~$N_Q$ of~$A_Q$ is the smallest positive integer~$N$ such that~$N A_Q^{-1}$ belongs to~$M_n(\Z)$.

Denote by~$\H$ the upper half-plane, defined by ${\H \= \{\tau \in \C \: \Im(\tau) > 0\}}$, and consider the usual action of~$\SL(2, \Z)$ on~$\H$.
For each~$\gamma$ in~$\SL(2, \Z)$ given by ${\gamma = \begin{psmallmatrix}a & b\\ c & d\end{psmallmatrix}}$, this action is defined by~$\gamma(\tau) \= \frac{a \tau + b}{c \tau + d}$.
For each~$N$ in~$\Z_{>0}$, let~$\Gamma_0(N)$ and~$\Gamma_1(N)$ be the congruence subgroups, defined by
\begin{align*}
  \Gamma_0(N)
  & \=
    \left\{ \bigl( \begin{smallmatrix} a & b \\ c & d \end{smallmatrix} \bigr) \in \SL(2, \Z) \: c \equiv 0 \mod N \right\};
  \\
  \Gamma_1(N)
  & \=
    \left\{ \bigl( \begin{smallmatrix} a & b \\ c & d \end{smallmatrix} \bigr) \in \Gamma_0(N) \: a, d \equiv 1 \mod N \right\}.
\end{align*}

For each nonempty finite set~$\Sigma$, denote by~$\cF(\Sigma)$ the vector space of complex valued functions defined on~$\Sigma$, endowed with the inner product ${\left\langle \cdot, \cdot \right\rangle_{\Sigma}}$ defined by
\begin{displaymath}
  \left\langle f,g \right\rangle_{\Sigma}
  \=
  \sum_{\sigma \in \Sigma}f(\sigma)\overline{g(\sigma)}.
\end{displaymath}
Denote by~$\| \cdot \|_{\Sigma}$ the corresponding norm and by~$\bfone_\Sigma$ the constant function in~$\cF(\Sigma)$ equal to~$1$.

The following is the main ingredient in the proof of Theorem~\ref{t:finitary-deviation-estimate}.

\begin{proposition}
  \label{p:modular form}
  Fix~$r$ in~$\Z_{>0}$ and put
  \begin{displaymath}
    N
    \=
    \lcm \left\{ 2 p^{2r} N_Q, \det(A_Q)2^{n+2} \right\}.
  \end{displaymath}
  Moreover, let~$\Sigma$ be a nonempty subset of~$(\Z / p^r \Z)^n \ssetminus \{ \bfzero \}$ and let~$f$ be a function in~$\cF(\Sigma)$.
  Then, the series given, for~$\tau$ in~$\H$, by
  \begin{displaymath}
    \vartheta_f(\tau)
    \=
    \sum_{m = 0}^{\infty} \left( \sum_{\substack{\bfx \in V_m(Q) \\ \red_r(\bfx) \in \Sigma}} f(\red_r(\bfx)) \right) \exp(2 \pi i m \tau)
  \end{displaymath}
  defines a modular form of weight~$\frac{n}{2}$ for the group~$\Gamma_1(N)$ in the sense of Shimura~\cite{Shi73}.
  If this modular form is cuspidal, then~$f$ is orthogonal to~$\bfone_\Sigma$ in~$\cF(\Sigma)$.
  If in addition~$\Sigma$ is an orbit of~$\oO_Q(\Z / p^r \Z)$, then this condition is also sufficient for~$\vartheta_f$ to be cuspidal.
\end{proposition}

After recalling basic properties of theta series in Section~\ref{ss:prel-theta-seri}, Section~\ref{ss:modularity} proves the modularity of~$\vartheta_f(\tau)$ using the transformation properties of the aforementioned theta series as found, \emph{e.g.}, in \cite[Chapter~VI]{Ogg69}, \cite[Section~2]{Shi73}, or \cite[Chapter~10]{Iwa97}, following Hecke, Pfetzer, and Schoeneberg.
Section~\ref{ss:cuspidality} completes the proof of Proposition~\ref{p:modular form} by showing the cuspidality criterion.

\subsubsection{Preliminaries on theta series}
\label{ss:prel-theta-seri}
For each odd integer~$d$, put
\begin{displaymath}
  \varepsilon_d \= 1 \text{ if~$d \equiv 1 \mod 4$}
  \text{ and }
  \varepsilon_d \= i \text{ if~$d \equiv -1 \mod 4$}.
\end{displaymath}
Moreover, for~$a$ in~$\Z$, denote by~$\left( \frac{a}{d} \right)$ the extended quadratic residue symbol as defined in~\cite[p.~442]{Shi73}, see also~\cite[p.~46]{Iwa97}.

For each~$M$ in~$\Z_{>0}$, consider the elements of~$(\Z / M \Z)^n$ as row vectors.
For all~$\xi$ in~$(\Z / M \Z)^n$ and~$N$ in~$\Z_{>0}$, denote by $N \cdot \xi$ the vector in~$(\Z / (MN) \Z)^n$ that equals~$N \bfx \mod MN$ for every~$\bfx$ in~$\Z^n$ satisfying ${\bfx \mod M = \xi}$.

For each~$\xi$ in~$(\Z / N_Q \Z)^n$ satisfying ${A_Q \xi^\intercal = \bfzero}$, define the theta series~$\Theta(\tau; Q, \xi)$ for~$\tau$ in~$\H$ by
\begin{displaymath}
  \Theta(\tau; Q, \xi)
  \=
  \sum_{\substack{\bfx \in \Z^n \\ \bfx \mod N_Q = \xi}} \exp \left( 2\pi i Q( \bfx )\tau /N_Q^2 \right).
\end{displaymath}
It satisfies
\begin{align}
  \label{eq:cuspide}
  \lim_{ \tau \rightarrow i\infty } \Theta(\tau; Q, \xi)
  & =
    \begin{cases}
      1 & \text{if } \xi = \bfzero; \\
      0 & \text{if } \xi \neq \bfzero,
    \end{cases}
  \\ \label{eq:16}
  \Theta(\tau+1; Q, \xi)
  & =
    \exp \left( 2\pi i Q(\xi)/N_Q^2 \right) \Theta(\tau; Q, \xi),
    \intertext{and, for every~$c$ in~$\Z_{>0}$,}
    \label{eq:17}
    \Theta(\tau; Q, \xi)
        & =
          \sum_{\substack{\xi' \in (\Z / cN_Q \Z)^n \\ \xi' \mod N_Q = \xi}} \Theta(c\tau; cQ, \xi').
\end{align}
For each~$z$ in~$\C\ssetminus \{0\}$, denote by~$\arg(z)$ the argument of~$z$ taking values in~$(-\pi, \pi]$ and, for each~$r$ in~$\R$, put ${z^r \= |z|^r \exp(r \arg(z) i)}$.
Then,
\begin{equation}
  \label{eq:18}
  \Theta \left( - \frac{1}{\tau}; Q, \xi \right)
  =
  \frac{(-i \tau)^{\frac{n}{2}}}{\det(A_Q)^{\frac{1}{2}}} \sum_{\substack{\xi' \in (\Z / N_Q \Z)^n \\ A_Q (\xi')^{\intercal} = \bfzero}}
  \exp \left( 2\pi i (\xi' A_Q \xi^\intercal)/N_Q^2 \right) \Theta(\tau;Q, \xi'),
\end{equation}
see \cite[Section~2]{Shi73} or \cite[Proposition~10.4]{Iwa97}.
Since for every~$\gamma$ in~$\Gamma_0(2N_Q)$, given by ${\gamma = \begin{psmallmatrix}a & b\\ c & d\end{psmallmatrix}}$, the number~$d$ is odd, the properties above imply
\begin{equation}
  \label{eq:Thetamodular}
  \Theta(\gamma(\tau); Q, \xi)
  =
  \exp \left( 2\pi i ab Q(\xi) / N_Q^2 \right) \left(\frac{\det(A_Q)}{d}\right)\left(\frac{2c}{d}\right)^n\varepsilon_d^{-n}(c\tau+d)^{\frac{n}{2}}\Theta(\tau;Q, a \xi),
\end{equation}
see \cite[Proposition~2.1 and comment~(i) below it]{Shi73} or \cite[Proposition~10.6$(ii)$]{Iwa97}.

\begin{lemma}
  \label{l:cusptheta}
  For every~$\gamma$ in~$\SL(2, \Z)$ given by ${\gamma = \begin{psmallmatrix}a & b\\ c & d\end{psmallmatrix}}$ with ${c > 0}$,
  \begin{displaymath}
    \lim_{\tau \to i\infty}\frac{\Theta(\gamma(\tau); Q, \xi)}{ (-i\tau)^{\frac{n}{2}}}
    =
    \frac{1}{\det(A_Q)^{\frac{1}{2}}} \sum_{\substack{\xi' \in (\Z / c N_Q \Z)^n \\ \xi' \mod N_Q = \xi}} \exp(2\pi i a Q(\xi') / (cN_Q^2)).
  \end{displaymath}
\end{lemma}
\begin{proof}
  By~\eqref{eq:16}, \eqref{eq:17}, and the formula ${c \gamma(\tau) = a - \frac{1}{c \tau + d}}$,
  \begin{equation}
    \label{eq:19}
    \begin{split}
      \Theta(\gamma(\tau);Q, \xi)
      & =
        \sum_{\substack{\xi' \in (\Z / c N_Q \Z)^n \\ \xi' \mod N_Q = \xi}} \Theta(c \gamma(\tau); cQ, \xi')
      \\ & =
           \sum_{\substack{\xi' \in (\Z / c N_Q \Z)^n \\ \xi' \mod N_Q = \xi}} \exp(2\pi i aQ(\xi')/(cN_Q^2)) \Theta\left(- \frac{1}{c\tau+d}; cQ, \xi' \right).
    \end{split}
  \end{equation}
  Moreover, \eqref{eq:18} implies, for every~$\xi'$ in~$(\Z / c N_Q \Z)^n$ with~$ \xi' \mod N_Q = \xi$,
  \begin{equation}
    \label{eq:20}
    \Theta\left(- \frac{1}{c\tau+d}; cQ, \xi' \right)
    =
    \frac{(-i(c\tau +d))^{\frac{n}{2}}}{(c^n\det(A_Q))^{\frac{1}{2}}}
    \sum_{\substack{\hxi \in (\Z / cN_Q \Z)^n \\ cA_Q \hxi^\intercal = \bf0}}
    \exp(2\pi i \hxi A_Q (\xi')^\intercal / (cN_Q^2)) \Theta(c \tau + d; cQ, \hxi).
  \end{equation}
  Together with~\eqref{eq:cuspide}, this implies
  \begin{displaymath}
    \lim_{\tau \to i\infty} \frac{\Theta\left(\frac{-1}{c\tau+d}; cQ, \xi' \right)}{(-i(c\tau +d))^{\frac{n}{2}}}
    =
    \frac{1}{(c^n\det(A_Q))^{\frac{1}{2}}}
  \end{displaymath}
  and thus
  \begin{displaymath}
    \lim_{\tau \to i\infty} \frac{\Theta(\gamma(\tau); Q, \xi)}{(-i(c\tau +d))^{\frac{n}{2}}}
    =
    \frac{1}{(c^n\det(A_Q))^{\frac{1}{2}}} \sum_{\substack{\xi' \in (\Z / c N_Q \Z)^n \\ \xi' \mod N_Q \equiv \xi}} \exp(2\pi i a Q(\xi') / (cN_Q^2)).
  \end{displaymath}
  Combined with $(-i(c\tau +d))^{\frac{n}{2}} \sim c^{\frac{n}{2}} (-i\tau)^{\frac{n}{2}}$ as $\tau\to i\infty$, this yields the desired result.
\end{proof}

\subsubsection{Modularity}
\label{ss:modularity}
This section proves that the series~$\vartheta_f(\tau)$ in Proposition~\ref{p:modular form} is modular.

That the theta series~$\Theta(\tau; p^rQ, N_Q \cdot \sigma)$ is well defined follows from ${A_{p^r Q} = p^r A_Q}$, ${N_{p^r Q} = p^r N_Q}$, and, for every~$\sigma$ in~$(\Z / p^r \Z)^n$, the equality $A_{p^rQ}(N_Q \cdot \sigma)^\intercal = \bfzero$ in~$(\Z / p^r N_Q \Z)^n$.

A direct computation shows that, for every~$\sigma$ in~$(\Z / p^r \Z)^n$,
\begin{displaymath}
  \Theta(p^r \tau; p^rQ, N_Q \cdot \sigma)
  =
  \sum_{\substack{ \bfx \in \Z^n \\ \red_r(\bfx) = \sigma}} \exp(2\pi i Q(\bfx)\tau)
\end{displaymath}
and hence
\begin{equation}
  \label{eq:21}
  \vartheta_f(\tau)
  =
  \sum_{ \sigma \in \Sigma } f(\sigma) \Theta(p^r\tau; p^rQ, N_Q \cdot \sigma).
\end{equation}
Thus, to prove that~$\vartheta_f(\tau)$ is modular for $\Gamma_1(N)$, it suffices to prove that, for every~$\sigma$ in~$(\Z / p^r \Z)^n$, the theta series~$\Theta(p^r\tau; p^rQ, N_Q \cdot \sigma)$ is modular for $\Gamma_1(N)$.

Let~$\begin{psmallmatrix}a & b\\ c & d\end{psmallmatrix}$ be in ${\Gamma_0(2 p^{2r} N_Q) \cap \Gamma_1(p^r)}$ and put ${\gamma \= \begin{psmallmatrix}a & b\\ c & d\end{psmallmatrix}}$.
Then, ${a(N_Q \cdot \sigma) = N_Q \cdot \sigma}$ and, by~\eqref{eq:Thetamodular} applied with~$\gamma$ replaced by~$\begin{psmallmatrix} a & p^r b\\ c/p^r & d\end{psmallmatrix}$,
\begin{equation}
  \label{eq:22}
  \begin{split}
    \Theta(p^r \gamma(\tau); p^rQ, N_Q \cdot \sigma)
    &=
      \Theta\left( \frac{a(p^r\tau) + p^r b}{(c/p^r)(p^r\tau)+d}; p^rQ, N_Q \cdot \sigma \right)
    \\ & =
         \left(\frac{\det(p^r A_Q)}{d}\right)\left(\frac{2(c/p^r)}{d}\right)^n\varepsilon_d^{-n}(c\tau+d)^{\frac{n}{2}}\Theta(p^r \tau; p^rQ, N_Q \cdot \sigma).
  \end{split}
\end{equation}
Since ${m \mapsto \left(\frac{m}{d}\right)}$ is a completely multiplicative function, see, \emph{e.g.}, \cite[3.(iii) and the last line in p.~442]{Shi73},
\begin{displaymath}
  \left(\frac{\det(p^r A_Q)}{d}\right)\left(\frac{2(c/p^r)}{d}\right)^n
  =
  \left(\frac{2^n\det(A_Q)}{d}\right)\left(\frac{c}{d}\right)^n.
\end{displaymath}
Using that~$m \mapsto \left(\frac{2^n\det(A_Q)}{m}\right)$ is a character modulo a divisor of~$2^{n+2}\det(A_Q)$, it follows that, if in addition ${d \equiv 1 \mod 2^{n + 2} \det(A_Q)}$, then ${\left(\frac{2^n \det(A_Q)}{d} \right) = 1}$.
Thus, if~$\gamma$ belongs to
\begin{equation}
  \label{eq:23}
  \Gamma_0(2p^{2r} N_Q) \cap \Gamma_1(p^r) \cap \Gamma_1(2^{n + 2} \det(A_Q)),
\end{equation}
then
\begin{displaymath}
  \Theta(p^r \gamma(\tau); p^rQ, N_Q \cdot \sigma)
  =
  \left(\frac{c}{d}\right)^n\varepsilon_d^{-n}(c\tau+d)^{\frac{n}{2}}\Theta(p^r \tau; p^rQ, N_Q \cdot \sigma).
\end{displaymath}
This implies that~$\Theta(p^r \tau;p^rQ, N_Q \cdot \sigma)$ is a modular form of weight~$\frac{n}{2}$ for $\Gamma_1(N)$ and that the same holds for $\vartheta_f(\tau)$.

\subsubsection{Cuspidality}
\label{ss:cuspidality}
This section proves the cuspidality criterion and thus completes the proof of Proposition~\ref{p:modular form}.

For the cusp $i\infty$, \eqref{eq:cuspide} implies, for every~$\sigma$ in~$(\Z / p^r \Z)^n \ssetminus \{ \bf0 \}$,
\begin{equation}
  \label{eq:24}
  \lim_{\tau \to i\infty}\Theta(\tau; p^rQ, N_Q \cdot \sigma)
  =
  0.
\end{equation}
So, \eqref{eq:21} and the hypothesis that~$\bfzero$ is outside~$\Sigma$ imply ${\lim_{\tau \to i\infty}\vartheta_f(\tau) = 0}$.

To study the behavior of~$\vartheta_f(\tau)$ at a different cusp, let~$\sigma$ be in ${(\Z / p^r \Z)^n \ssetminus \{ \bf0 \}}$ and~$\gamma$ in~$\SL(2, \Z)$ given by ${\gamma = \begin{psmallmatrix}a & b\\ c & d\end{psmallmatrix}}$ with ${c > 0}$.
Let~$s$ be the largest integer in $\{0,\ldots,r\}$ such that~$p^s$ divides~$c$.
Then, ${\gcd \left(p^{-s}c,p^{r-s}\right) = 1}$ and hence there is~$j$ in~$\Z$ such that ${j p^{-s}c\equiv d \mod p^{r-s}}$.
Put
\begin{equation}
  \label{eq:25}
  \tgamma
  \=
  \left(\begin{array}{cc}p^{r-s}a & p^sb-ja\\ p^{-s}c & p^{-r}(p^sd-cj)\end{array}\right).
\end{equation}
Then, ${\tgamma}$ belongs to~$\SL(2, \Z)$ and
\begin{displaymath}
  \left(\begin{array}{cc}p^{r} & 0\\ 0 & 1\end{array}\right)\gamma
  =
  \tgamma \left(\begin{array}{cc}p^{s} & j\\ 0 & p^{r-s}\end{array}\right).
\end{displaymath}
By Lemma~\ref{l:cusptheta},
\begin{equation}
  \label{eq:26}
  \begin{split}
    \lim_{\tau \to i\infty}\frac{\Theta(p^r \gamma(\tau); p^rQ, N_Q \cdot \sigma)}{(-i\tau)^{\frac{n}{2}}}
    & =
      \lim_{\tau \to i\infty}\frac{\Theta \left( \tgamma \left( \frac{p^s \tau + j}{p^{r - s}} \right) ; p^rQ, N_Q \cdot \sigma \right)}{(-i\tau)^{\frac{n}{2}}}
    \\ &=
         p^{(2s - r) \frac{n}{2}} \lim_{\tau \to i\infty}\frac{\Theta \left( \tgamma \left( \tau \right) ; p^rQ, N_Q \cdot \sigma \right)}{(-i\tau)^{\frac{n}{2}}}
    \\ &=
         \frac{p^{(2s - r) \frac{n}{2}}}{\det(A_{p^rQ}))^{\frac{1}{2}}}
         \sum_{\substack{\xi \in (\Z / p^{r-s}cN_Q \Z)^n \\ \xi \mod p^rN_Q = N_Q \cdot \sigma}} \exp(2\pi i a Q(\xi)/ (cN_Q^2))
    \\ &=
         \frac{1}{p^{(r - s) n} \det(A_Q)^{\frac{1}{2}}} \sum_{\substack{ \sigma' \in (\Z / p^{r-s}c \Z)^n \\ \sigma' \mod p^r =\sigma }}\exp(2\pi i a Q(\sigma')/ c).
  \end{split}
\end{equation}
Together with~\eqref{eq:21}, this implies
\begin{equation}
  \label{eq:finitecusp}
  \lim_{\tau \to i\infty}\frac{\vartheta_f(\gamma(\tau))}{(-i\tau)^{\frac{n}{2}}}
  =
  \frac{1}{p^{(r - s)n} \det(A_Q)^{\frac{1}{2}}} \sum_{\sigma \in \Sigma}f(\sigma) \sum_{\substack{\sigma' \in (\Z / p^{r-s}c \Z)^n \\ \sigma' \mod p^r = \sigma}} \exp(2\pi i a Q(\sigma')/ c).
\end{equation}
If $\vartheta_f(\tau)$ is cuspidal, then~\eqref{eq:finitecusp} with $\gamma=\begin{psmallmatrix}0 & -1\\ 1 & 0\end{psmallmatrix}$ implies~$\langle f, \bfone_{\Sigma} \rangle_{\Sigma} = 0$.
This proves the necessity condition for cuspidality in Proposition~\ref{p:modular form}.

The sufficiency condition for cuspidality in Proposition~\ref{p:modular form} is a direct consequence of~\eqref{eq:finitecusp} and the following lemma.

\begin{lemma}
  Let~$\Sigma$ be an orbit of~$\oO_Q(\Z / p^r \Z)$ in~$(\Z / p^r \Z)^n$ different from~$\{ \bfzero \}$.
  Then, for all~$a$ in~$\Z$ and~$c$, $r$, and~$t$ in~$\Z_{>0}$ satisfying~$p^r \mid p^t c$, the following function is constant,
  \begin{displaymath}
    \begin{array}{rrcl}
      \fE_{a,c} \: & \Sigma & \to & \C
      \\
                   & \sigma &\mapsto & \sum_{\substack{ \sigma' \in (\Z / p^t c \Z)^n \\ \sigma' \mod p^r = \sigma}} \exp(2 \pi i a Q(\sigma')/ c).
    \end{array}
  \end{displaymath}
\end{lemma}

\begin{proof}
  Let~$\ell$ be in~$\Z_{\ge 0}$ and~$c_0$ in~$\Z_{>0}$ such that ${c = p^{\ell}c_0}$ and ${p \nmid c_0}$.
  Choosing~$A$ and~$B$ in~$\Z$ with $Ap^{t + \ell}+Bc_0 = 1$ yields an isomorphism
  \begin{displaymath}
    \begin{array}{rcl}
      (\Z/c_0\Z)^n \times (\Z/p^{t + \ell}\Z)^n & \to & (\Z/p^{t}c\Z)^n
      \\
      (\mu, \nu) & \mapsto & A p^{t + \ell} \cdot \mu + B c_0 \cdot \nu.
    \end{array}
  \end{displaymath}
  This implies, for every~$\sigma$ in~$\Sigma$,
  \begin{equation}
    \label{eq:27}
    \begin{split}
      \fE_{a,c}(\sigma)
      &=
        \sum_{\mu \in (\Z / c_0 \Z)^n} \sum_{\substack{\nu \in (\Z / p^{t + \ell} \Z)^n \\ \nu \mod p^r = \sigma}} \exp \left ( 2\pi i a Q(A p^{t + \ell} \cdot \mu + B c_0 \cdot \nu)/(p^{\ell} c_0) \right)
      \\ &=
           \left( \sum_{\mu \in (\Z / c_0 \Z)^n} \exp \left( 2\pi i a A^2p^{2t + \ell} Q(\mu)/c_0 \right) \right)
           \left( \sum_{\substack{\nu \in (\Z / p^{t + \ell} \Z)^n \\ \nu \mod p^r = \sigma}} \exp \left( 2\pi i a B^2 c_0 Q(\nu)/ p^\ell \right) \right)
      \\ &=
           \left( \sum_{\mu \in (\Z / c_0 \Z)^n} \exp \left( 2\pi i a A^2 p^{2t + \ell} Q(\mu)/c_0 \right) \right) \fE_{aB^2c_0,p^\ell}(\sigma).
    \end{split}
  \end{equation}
  Hence, it suffices to consider the case where ${c_0 = 1}$ and ${c = p^\ell}$.

  Let~$\sigma$ and~$\hsigma$ be in~$\Sigma$.
  The hypothesis that~$\Sigma$ is an orbit of~$\oO_Q(\Z / p^r \Z)$ implies that there is~$T$ in~$\oO_Q(\Z / p^{t + \ell} \Z)$ such that $T(\sigma) = \hsigma$.
  Since
  \begin{displaymath}
    \begin{array}{rcl}
      \{\nu \in (\Z/p^{t+\ell}\Z)^n \: \nu \mod p^r = \sigma \}
      & \to &
              \{\hnu \in (\Z/p^{t+\ell}\Z)^n \: \hnu \mod p^r = \hsigma \}
      \\
      \nu & \mapsto & T(\nu)
    \end{array}
  \end{displaymath}
  is a bijection, it follows that
  \begin{multline*}
    \fE_{a, p^\ell}(\hsigma)
    =
    \sum_{\substack{\nu \in (\Z / p^{t + \ell} \Z)^n \\ \nu \mod p^r = \sigma}} \exp \left( 2\pi i a Q(T(\nu))/p^\ell \right)
    =
    \sum_{\substack{\nu \in (\Z / p^{t + \ell} \Z)^n \\ \nu \mod p^r = \sigma}} \exp \left( 2\pi i a Q(\nu)/p^\ell \right)
    =
    \fE_{a,p^\ell}(\sigma).
    \qedhere
  \end{multline*}
\end{proof}

\subsection{Proofs of Theorem~\ref{t:p-adic-Linnik} and~\ref{t:finitary-deviation-estimate} and Corollary~\ref{c:functional-deviation}}
\label{ss:proof of finitary deviation estimate}
The proofs of Theorems~\ref{t:p-adic-Linnik} and~\ref{t:finitary-deviation-estimate} rely on the following lemma.
The proof of Corollary~\ref{c:functional-deviation} is at the end of this section.

\begin{lemma}
  \label{l:devformula}
  Fix~$r$ in~$\Z_{>0}$, let~$m$ in~$\Z_{>0}$ be such that ${V_m(Q) \neq \emptyset}$ and~$\Sigma$ a subset of~$(\Z / p^r \Z)^n$ containing~$\red_r(V_{m}(Q))$.
  Then, for every orthonormal basis~$\cB_0$ of the orthogonal complement of~$\bfone_\Sigma$ in~$\cF(\Sigma)$,
  \begin{multline*}
    \Var(m, \Sigma)
    \=
    \sum_{\sigma \in \Sigma}\left(\frac{\# \left\{ \bfx \in V_m(Q) \: \red_r(\bfx) = \sigma \right\} }{\#V_m(Q) }-\frac{1}{\# \Sigma}\right)^2
    \\ =
    \frac{1}{\# V_{m}(Q)^2} \sum_{f\in \cB_0}\left| \sum_{\bfx \in V_m(Q)}f(\red_r(\bfx))\right|^2.
  \end{multline*}
\end{lemma}

\begin{proof}
  Let ${F \: \Sigma\to \C}$ be given by
  \begin{equation}
    \label{eq:28}
    F(\sigma)
    \=
    \# \left\{ \bfx \in V_m(Q) \: \red_r(\bfx) = \sigma \right\}.
  \end{equation}
  Then, ${\langle F, \bfone_\Sigma \rangle_{\Sigma} = \#V_m(Q)}$ and
  \begin{multline}
    \label{eq:l2 norm}
    \left\| F - \frac{ \langle F, \bfone_\Sigma \rangle_{\Sigma} }{\# \Sigma} \bfone_\Sigma \right\|_{\Sigma}^2
    =
    \sum_{\sigma \in \Sigma}\left( \# \left\{ \bfx \in V_m(Q) \: \red_r(\bfx) = \sigma \right\} - \frac{\#V_m(Q)}{\#\Sigma} \right)^2
    \\ =
    \# V_m(Q)^2 \cdot \Var(m, \Sigma).
  \end{multline}
  Moreover, since~$\cB_0$ is an orthonormal basis for the orthogonal complement of~$\bfone_\Sigma$ in~$\cF(\Sigma)$,
  \begin{equation}
    \label{eq:29}
    F - \frac{ \langle F, \bfone_\Sigma \rangle_{\Sigma} }{\# \Sigma} \bfone_\Sigma
    =
    \sum_{f\in \cB_0}\langle F,f\rangle_{\Sigma} f
  \end{equation}
  and thus
  \begin{equation}
    \label{eq:30}
    \begin{split}
      \left\| F - \frac{ \langle F, \bfone_\Sigma \rangle_{\Sigma} }{\# \Sigma} \bfone_\Sigma \right\|_{\Sigma}^2
      & =
        \sum_{f\in \cB_0} |\langle F,f\rangle_{\Sigma}|^2
      \\ & =
           \sum_{f\in \cB_0} \left|\sum_{\sigma \in \Sigma} \# \left\{ \bfx \in V_m(Q) \: \red_r(\bfx) = \sigma \right\} \overline{f(\sigma)}\right|^2
      \\ & =
           \sum_{f \in \cB_0} \left|\sum_{\bfx \in V_{m}(Q)}f(\red_r(\bfx))\right|^2.
    \end{split}
  \end{equation}
  Together with~\eqref{eq:l2 norm} this implies the desired identity.
\end{proof}

\begin{proof}[Proof of Theorem~\ref{t:finitary-deviation-estimate}]
  Since~$(\Z / p^r \Z)^n$ is finite, it suffices to prove the desired estimate for a given orbit~$\Sigma$ of~$\oO_Q(\Z / p^r\Z)$.
  Since the case where ${\Sigma = \{ \bfzero \}}$ is trivial, suppose ${\Sigma \subseteq (\Z / p^r \Z)^n \ssetminus \{ \bfzero \}}$.

  Let~$\cB_0$ be an orthonormal basis of the orthogonal complement of~$\bfone_\Sigma$ in~$\cF(\Sigma)$.
  By Lemma~\ref{l:devformula}, for every~$\sigma$ in~$\Sigma$,
  \begin{multline}
    \label{eq:31}
    \left| \frac{ \# \left\{ \bfx \in V_m(Q) \: \red_r(\bfx) = \sigma \right\} }{\#V_m (Q)} - \frac{1}{\# \Sigma} \right|
    \\ \le
    \sqrt{\Var(m, \Sigma)}
    =
    \frac{1}{\# V_{m}(Q)} \left( \sum_{f\in \cB_0}\left| \sum_{\bfx \in V_m(Q)}f(\red_r(\bfx))\right|^2 \right)^{\frac{1}{2}}.
  \end{multline}
  Since each~$f$ in~$\cB_0$ is orthogonal to~$\bfone_\Sigma$, by Proposition~\ref{p:modular form} the modular form~$\vartheta_f$ is cuspidal of weight~$\frac{n}{2}$ for~$\Gamma_1(N)$.
  When $n\ge 4$, for every $\varepsilon$ in~$\R_{>0}$ there exists~$C$ in~$\R_{>0}$ that only depends on~$f$ and~$\varepsilon$, such that
  \begin{displaymath}
    \left|\sum_{\bfx \in V_m(Q)}f(\red_r(\bfx))\right| = |m\text{-th Fourier coefficient of }\vartheta_f |\leq C m^{\frac{n}{4} - \kappa_n +\varepsilon},
  \end{displaymath}
  by Deligne's bound \cite[\emph{Th{\'e}or{\`e}me}~8.2]{Del74} if~$n$ is even and by Iwaniec's bound \cite[Theorem~1]{Iwa87} if~$n$ is odd.
  When ${n = 3}$ the same estimate holds for a constant~$C$ that also depends on~$S$, by Duke's \cite[Theorem~5]{Duk88} and Blomer's \cite[Lemma~4.4]{Blo04} bounds.
  Combined with~\eqref{eq:31}, this implies the desired estimate.
\end{proof}

\begin{remark}
  The literature usually states the bounds for Fourier coefficients of cuspidal modular forms for~$\Gamma_1(N)$ used in the proof above for cuspidal modular forms for~$\Gamma_0(N)$ with characters.
  These bounds also hold for cuspidal modular forms for~$\Gamma_1(N)$ because every such form is
  a finite sum of cuspidal modular forms for~$\Gamma_0(N)$ with characters, see, \emph{e.g.}, \cite[Lemma~4.3.1]{Miy89} (the proof provided there extends to the case of half-integral weight modular forms).
\end{remark}

\begin{proof}[Proof of Theorem~\ref{t:p-adic-Linnik}]
  Since the set of locally constant functions~$\Zp^n \to \R$ is dense in the space of continuous functions~$\Zp^n \to \R$, it suffices to prove that, for every locally constant function ${F \: \Zp^n \to \R}$,
  \begin{displaymath}
    \frac{1}{\# V_{m_j}(Q)} \sum_{\bfx \in V_{m_j}(Q)} F(M_{u_j}^{-1}(\bfx)) \to \int F \dd \mu_\ell
    \text{ as }
    j \to \infty.
  \end{displaymath}
  Let~$r$ in~$\Z_{>0}$ be sufficiently large so that, for every~$\sigma$ in $\Z / p^r \Z$, the function~$F$ is constant on~$\red_r^{-1}(\sigma)$, and let ${f \: (\Z / p^r \Z)^n \to \R}$ be the function determined by ${F = f \circ \red_r}$.
  Let~$\varepsilon$ be in~$\R_{>0}$ and put ${\delta \= c - (\frac{n}{4} - \kappa_n + \varepsilon) > 0}$.
  Taking~$\varepsilon$ smaller if necessary, suppose ${\delta > 0}$ and let~$C$ be given by Theorem~\ref{t:finitary-deviation-estimate}.

  The hypotheses ${S_{\ell}(Q) \neq \emptyset}$ and that~$\oO_Q(\Zp)$ acts transitively on~$S_{\ell}(Q)$ imply, for every~$\ell'$ in~$\ell (\Zp^{\times})^2$, that ${\red_r(S_{\ell'}(Q)) \neq \emptyset}$ and~$\oO_Q(\Z / p^r \Z)$ acts transitively on~$\red_r(S_{\ell'}(Q))$.
  In particular, for each~$j$ in~$\Z_{>0}$, this applies with ${\ell' = m_j}$, and~$M_{u_j}^{-1}$ maps~$S_{m_j}(Q)$ to~$S_{\ell}(Q)$ and~$\mu_{m_j}$ to~$\mu_{\ell}$.
  Denote by~$M_{u_j,r}$ the reduction modulo~$p^r$ of~$M_{u_j}$, which belongs to~$\GL_n(\Z/p^r\Z)$.
  Applying, for each sufficiently large~$j$, Theorem~\ref{t:finitary-deviation-estimate} with~$\Sigma = \red_r(S_{m_j}(Q))$ and~$m = m_j$, it follows that, for every~$\sigma$ in~$\red_r(S_{m_j}(Q))$,
  \begin{equation}
    \label{eq:32}
    \left| \frac{ \# \left\{ \bfx \in V_{m_j}(Q) \: \red_r(\bfx) = \sigma \right\} }{\#V_{m_j}(Q)} - \frac{1}{\# \red_r(S_{m_j}(Q))} \right|
    \le
    C \frac{m_j^{\frac{n}{4}-\kappa_n + \varepsilon}}{\#V_{m_j}(Q)}
    \le
    C m_j^{-\delta}.
  \end{equation}
  Moreover, the change of variables formula and~\eqref{eq:15} imply
  \begin{equation}
    \label{eq:33}
    \int F \dd \mu_\ell
    =
    \int F \dd (M_{u_j}^{-1})_* \mu_{m_j}
    =
    \int F \circ M_{u_j}^{-1} \dd \mu_{m_j}
    =
    \sum_{\sigma \in \red_r(S_{m_j}(Q))} \frac{f (M_{u_j,r}^{-1}(\sigma))}{\# \red_r(S_{m_j}(Q))}.
  \end{equation}
  Together with~\eqref{eq:32}, this implies
  \begin{displaymath}
    \left| \frac{1}{\# V_{m_j}(Q)} \sum_{\bfx \in V_{m_j}(Q)} F(M_{u_j}^{-1}(\bfx)) - \int F \dd \mu_\ell \right|
    \le
    C \left(\sum_{\sigma' \in \red_r(S_{\ell}(Q))} |f(\sigma')|\right) m_j^{-\delta},
  \end{displaymath}
  from which the desired assertion follows.
\end{proof}

The proof of Corollary~\ref{c:functional-deviation} relies on the following lemma.

\begin{lemma}
  \label{l:contention}
  Let~$r$ be in~$\Z_{>0}$ and~$\ell$ and~$m$ in~$\Zp\ssetminus \{0\}$ satisfying
  \begin{equation}
    \label{eq:34}
    | m - \ell |_p< | 2\ell |_p^2
    \text{ and }
    | m - \ell |_p \le | 2\ell |_p p^{- r}.
  \end{equation}
  Then, ${\red_r(S_{\ell}(Q))= \red_r(S_m(Q))}$.
  In particular, if in addition~$m$ belongs to~$\Z_{>0}$, then ${\red_r(V_{m}(Q)) \subseteq \red_r(S_{\ell}(Q))}$.
\end{lemma}

\begin{proof}
  Let~$\bfx$ be in~$S_{m}(Q)$, given by ${\bfx = (x_1, \ldots, x_n)}$.
  The relation
  \begin{equation}
    \label{eq:35}
    2Q(\bfx)
    =
    \sum_{i=1}^n x_i\cdot \partial_{X_i} Q(\bfx)
  \end{equation}
  implies
  \begin{equation}
    \label{eq:36}
    | \ell |_p
    =
    | m |_p
    =
    | Q(\bfx) |_p \leq | 2 |_p^{-1}\cdot \max_{i \in \{1, \ldots, n \}} \left\{ | \partial_{X_i}Q(\bfx)) |_p \right\},
  \end{equation}
  and thus
  \begin{equation}
    \label{eq:37}
    | Q(\bfx)-\ell |_p^{\frac{1}{2}}
    =
    | m - \ell |^{\frac{1}{2}}_p
    <
    | 2\ell |_p\leq \max_{i \in \{1, \ldots, n \}} \left\{ |\partial_{X_i}Q(\bfx)|_p \right\}.
  \end{equation}
  Hence, Hensel's Lemma yields~$\bfx'$ in~$S_{\ell}(Q)$, given by ${\bfx' = (x_1', \ldots, x_n')}$, such that
  \begin{equation}
    \label{eq:38}
    \max_{i \in \{1, \ldots, n \}} \left\{ | x_i' - x_i |_p \right\}
    \leq
    \frac{|Q(\bfx)-\ell|_p}{\max_{i \in \{1, \ldots, n \}} \{|\partial_{X_i}Q(\bfx)|_p\}}
    \leq
    \frac{|m-\ell|_p}{|2\ell|_p}
    \le
    p^{-r}.
  \end{equation}
  In particular, ${\bfx' \equiv \bfx \mod p^{r}}$.
  This proves that~$\bfx$ belongs to~$\red_r(S_{\ell}(Q))$ and hence ${\red_r(S_m(Q)) \subseteq \red_r(S_\ell(Q))}$.

  The reverse inclusion follows by symmetry.
\end{proof}

\begin{proof}[Proof of Corollary~\ref{c:functional-deviation}]
  Let~$r$ in~$\Z_{>0}$ be sufficiently large so that, for every~$\sigma$ in~$\Z / p^r \Z$, the set~$\red_r^{-1}(\sigma)$ is contained in a ball of radius~$\delta$, and put
  \begin{equation}
    \label{eq:39}
    \Sigma
    \=
    \red_r(S_{\ell}(Q))
    \text{ and }
    N
    \=
    \max \{ 2 \ord_p(2\ell) + 1, \ord_p(2\ell) + r \}.
  \end{equation}
  Then, $\Sigma$ is an orbit of~$\oO_Q(\Z / p^r \Z)$ and, if~$m$ is as in the statement of the corollary, then~\eqref{eq:34} holds and Lemma~\ref{l:contention} implies ${\red_r(V_m(Q)) \subseteq \Sigma}$.
  Then, the desired estimate follows from Theorem~\ref{t:finitary-deviation-estimate}.
\end{proof}

\section{\CM{} points formulae}
\label{s:CM-formulae}

This section provides several formulae for (formal) \CM{} points having supersingular reduction.
Section~\ref{ss:fixed-points-formula} provides the first formula.
It writes each \CM{} point of fundamental discriminant as projections of fixed points of certain elements of the group action described in Section~\ref{ss:from-elliptic-curves} (Theorem~\ref{t:fixed-points-formula}).
Section~\ref{ss:CM-from-canonical} provides the second formula.
For all~$r$ in~$\Z_{\ge 0}$ and discriminant~$D$ whose conductor is not divisible by~$p$, it relates~$\Lambda_{D p^{2r}}$ to~$\Lambda_D$ using the canonical branch~$\t$ of~$T_p$ (Theorem~\ref{t:CM-from-canonical}).
Section~\ref{ss:formal-CM-formulae} provides the third formula, which is the analogous relation for formal \CM{} points, and describes the relation between \CM{} and formal \CM{} points (Theorem~\ref{t:formal-CM-formulae} and Corollary~\ref{c:formalization}).

For each discriminant~$D$, this section and the rest of the paper considers~$\Lambda_D$ as a divisor.

\subsection{\CM{} points as fixed points}
\label{ss:fixed-points-formula}
Throughout this section, fix~$\ss$ in~$\tSups$ and let~$\Bss$, $\Rss$, $\Gss$, and~$\hDss$ be as in Section~\ref{ss:from-elliptic-curves}.
The \emph{Gross lattice associated with~$\ss$} is the $\Z$-lattice of dimension three~$L(\ss)$, defined by
\begin{equation}
  \label{eq:40}
  L(\ss)
  \=
  \{ \phi \in \Z + 2 \End(\ss) \: \tr(\phi) = 0 \}.
\end{equation}
It plays a central role in this section.
For each~$m$ in~$\Z_{>0}$, define
\begin{equation}
  \label{eq:41}
  V_m(\ss)
  \=
  \{ \phi \in L(\ss) \: \nr(\phi) = m \}.
\end{equation}

For each $p$-supersingular fundamental discriminant~$d$, this section writes every \CM{} point of discriminant~$d$ in~$\Dss$ as the projection of a fixed point of a certain element of the group action of~$\Gss$ on~$\hDss$.
This is done in two steps.
The first step is to define, for each~$\phi$ in~$V_{|d|}(\ss)$, a certain unit~$\Uss(\hphi)$ in the ring of integers of the subalgebra~$\Qp(\hphi)$ of~$\Bss$ (Lemma~\ref{l:unit-function}).
The second step is to prove that, as~$\phi$ varies over~$V_{|d|}(\ss)$, the projections of the fixed points of~$\Uss(\hphi)$ in~$\hDss$ run through all \CM{} points in~$\Dss$ of discriminant~$d$ (Theorem~\ref{t:fixed-points-formula}).

The statement of these results rely on the following notation.
Put
\begin{equation}
  \label{eq:42}
  \Lss
  \=
  \{ \varphi \in \Zp + 2 \Rss \: \tr(\varphi) = 0 \},
\end{equation}
so~$\Lss$ is the image of~$L(\ss) \otimes \Zp$ by the natural isomorphism ${\End(\ss) \otimes \Zp \to \Rss}$.
This set is compact because~$\Rss$ is compact and the reduced trace function is continuous.
Moreover, for every nonzero~$\varphi$ in~$\Lss$, the $p$\nobreakdash-adic number~$- \nr(\varphi)$ belongs to a $p$\nobreakdash-adic discriminant.
This motivates the definition,
\begin{displaymath}
  \Lfss
  \=
  \{ \varphi \in \Lss \: -\nr(\varphi) \text{ belongs to a fundamental $p$\nobreakdash-adic discriminant} \}.
\end{displaymath}
The set~$\Lfss$ coincides with the set of elements~$\varphi$ of~$\Bss$ such that~$\varphi^2$ belongs to a fundamental $p$\nobreakdash-adic discriminant.
Moreover, for every $p$-supersingular discriminant~$D$ whose conductor is not divisible by~$p$, the set~$V_{|D|}(\ss)$ is mapped inside~$\Lfss$ by the map~$\phi \mapsto \hphi$ (Lemma~\ref{l:p-adic-discriminants}).

\begin{lemma}[Unit function]
  \label{l:unit-function}
  Let~$\Uss \: \Lfss \to \Bss$ be defined by
  \begin{displaymath}
    \Uss(\varphi)
    \=
    \begin{cases}
      \frac{\varphi^2 + \varphi}{2}
      & \text{if~$\frac{\varphi^2 + \varphi}{2}$ belongs to~$\Gss$};
      \\
      1 + \frac{\varphi^2 + \varphi}{2}
      & \text{otherwise}.
    \end{cases}
  \end{displaymath}
  Then, $\Uss$ takes values in~$\Gss$ and, for every~$\varphi$ in~$\Lfss$, the following properties hold.
  \begin{enumerate}
  \item[$(i)$]
    The subalgebra~$\Qp(\varphi)$ of~$\Bss$ is a field extension of~$\Qp$ isomorphic to the subfield~$\Qp(\sqrt{\varphi^2})$ of~$\Cp$.
  \item[$(ii)$]
    The equality ${\cO_{\Qp(\varphi)} = \Zp[\Uss(\varphi)]}$ holds, $\Uss(\varphi)$ is a unit in~$\cO_{\Qp(\varphi)}$, and~$\disc(\Uss(\varphi))$ belongs to a fundamental $p$\nobreakdash-adic discriminant.
  \end{enumerate}
\end{lemma}

\begin{proof}
  Since ${\varphi^2 = - \nr(\varphi)}$ and~$- \nr(\varphi)$ belongs to a fundamental $p$\nobreakdash-adic discriminant, it follows that~$\varphi^2$ is outside~$\Qp^2$.
  This implies item~$(i)$.

  By~\eqref{eq:123} in Lemma~\ref{l:quadratic-extensions}$(ii)$,
  \begin{displaymath}
    \cO_{\Qp(\varphi)}
    =
    \Zp \left[ \tfrac{\varphi^2 + \varphi}{2} \right]
    =
    \Zp \left[ \Uss(\varphi) \right].
  \end{displaymath}
  In particular, $\frac{\varphi^2 + \varphi}{2}$ belongs to~$\Rss$ and thus~$\Uss(\varphi)$ belongs to~$\Gss$ and it is a unit in~$\Zp[\Uss(\varphi)]$.
  Since ${\disc(\Uss(\varphi)) = - \nr(\varphi)}$, it also follows that~$\disc(\Uss(\varphi))$ belongs to a fundamental $p$\nobreakdash-adic discriminant.
\end{proof}

For each~$\varphi$ in~$\Lfss$, define
\begin{displaymath}
  \Fixss(\varphi)
  \=
  \left\{ x \in \hDss \: \Uss(\varphi) \cdot x = x \right\}.
\end{displaymath}
For all fundamental discriminant~$d$ and~$f$ in~$\Z_{>0}$, put
\begin{displaymath}
  w_{d, f}
  \=
  \# \left( \cO_{d, f}^{\times} / \Z^{\times} \right).
\end{displaymath}
So, ${w_{-3, 1} = 3}$, ${w_{-4, 1} = 2}$, and, in all remaining cases, ${w_{d, f} = 1}$.

\begin{theorem}[Fixed points formula]
  \label{t:fixed-points-formula}
  For all~$\ss$ in~$\tSups$ and $p$-supersingular fundamental discriminant~$d$,
  \begin{equation}
    \label{eq:43}
    \Lambda_d|_{\Dss}
    =
    \frac{w_{d,1}}{\#\Aut(\ss)}\sum_{\phi\in V_{|d|}(\ss)}\sum_{x\in \Fixss(\hphi)}\Piss(x).
  \end{equation}
\end{theorem}

The proof of this theorem relies on the following version of Deuring's lifting theorem for formal $\Zp$\nobreakdash-modules, in the spirit of \cite[Proposition~2.1]{Gro86}.
For each formal group~$\cF$ over a ring~$R$, denote by ${D_{\cF} \: \End_R(\cF) \to R}$ the ring homomorphism such that, for every~$\varphi$ in~$\End_R(\cF)$, the following congruence holds in coordinates
\begin{displaymath}
  \varphi(X) \equiv D_{\cF}(\varphi) X \mod X^2.
\end{displaymath}
Moreover, for each ring homomorphism ${\delta \: R \to \OQpalg}$, denote by ${\tdelta \: R \to \Fpalg}$ the composition of~$\delta$ with the reduction morphism ${\OQpalg \to \Fpalg}$.

\begin{proposition}[Lifting formal modules]
  \label{p:formal-lifting}
  Let~$\ss$ be in~$\tSups$ and~$g_0$ in~$\Gss$ such that the subalgebra~$\Qp(g_0)$ of~$\Bss$ is a field extension of~$\Qp$ of degree two with ring of integers~$\Zp[g_0]$.
  Then, there is a bijection between the fixed points of~$g_0$ in~$\Xss(\OQpalg)$ and the continuous ring homomorphisms ${\delta \: \Zp[g_0] \to \OQpalg}$ satisfying ${D_{\Fss}|_{\Zp[g_0]} = \tdelta}$.
  For such a~$\delta$, put ${\cK \= \Q_{p^2}(\delta(g_0))}$ and let~$\cO_{\cK}$ be the ring of integers of~$\cK$.
  Then, the corresponding fixed point~$(\cF, \alpha)$ of~$g_0$ is defined over~$\cO_{\cK}$ and it is uniquely determined by the property that the unique automorphism~$\varphi_0$ in~$\Aut_{\OK}(\cF)$ such that~$ g_0=\alpha\circ \tvarphi_0 \circ \alpha^{-1}$ satisfies~$D_{\cF}(\varphi_0) = \delta(g_0)$.
\end{proposition}

\begin{proof}
  The proof begins assigning to each fixed point~$(\cF, \alpha)$ of~$g_0$ in~$\Xss(\OQpalg)$ a continuous ring homomorphism ${\delta \: \Zp[g_0] \to \OQpalg}$ as in the statement.
  Let ${\iota \: \End(\cF) \to \Rss}$ be the ring homomorphism defined by
  \begin{displaymath}
    \iota(\varphi)
    \=
    \alpha \circ \tvarphi \circ \alpha^{-1},
  \end{displaymath}
  which is continuous, see, \emph{e.g.}, \cite[Chapter~IV, Section~1, Proposition~3]{Fro68}.
  Let~$\varphi_0$ in~$\End(\cF)$ be such that~$\iota(\varphi_0) = g_0$ (Lemma~\ref{l:fixed}).
  Then, $\iota$ induces a continuous ring isomorphism ${\iota_0 \: \Zp[\varphi_0] \to \Zp[g_0]}$ and the ring homomorphism~$\delta$, defined by
  \begin{displaymath}
    \delta
    \=
    D_{\cF} \circ \iota_0^{-1} \: \Zp[g_0] \to \OQpalg,
  \end{displaymath}
  satisfies
  \begin{displaymath}
    \delta(g_0)
    =
    D_{\cF} (\iota_0^{-1}(g_0))
    =
    D_{\cF} (\varphi_0)
  \end{displaymath}
  and, for every~$g$ in~$\Zp[g_0]$,
  \begin{displaymath}
    \tdelta(g)
    =
    \widetilde{D_{\cF}(\iota_0^{-1}(g))}
    =
    D_{\tcF}(\widetilde{\iota_0^{-1}(g)})
    =
    D_{\Fss}(g).
  \end{displaymath}
  The ring homomorphism~$\delta$ is continuous because~$D_{\cF}$ is, see \cite[Chapter~IV, Section~1, Corollary~3]{Fro68}.

  Let ${\delta \: \Zp[g_0] \to \Qpalg}$ be a continuous ring homomorphism satisfying ${D_{\Fss}|_{\Zp[g_0]} = \tdelta}$ and put ${\cK \= \Q_{p^2}(\delta(g_0))}$.
  Note that ${\cK \subset \Qpalg}$.
  The next step is to prove that there is a fixed point~$(\cF_0, \alpha_0)$ in $\Xss(\OK)$ whose corresponding ring homomorphism is~$\delta$.
  Endow~$\OK$ and~$\Fpalg$ with the structure of a~$\Zp[g_0]$\nobreakdash-module with structural map~$\delta$ and~$D_{\Fss}|_{\Zp[g_0]}$, respectively.
  Then, the inclusion map ${\Zp[g_0] \to \Rss}$ gives~$\Fss$ the structure of a formal $\Zp[g_0]$\nobreakdash-module over~$\Fpalg$ in the sense of \mbox{Drinfel'd}, see \cite[Section~1]{Dri74}.
  This formal $\Zp[g_0]$\nobreakdash-module is of height one, see, \emph{e.g.}, \cite[Remark, p.~566]{Dri74}.
  Then, there is a unique deformation~$(\cF_0, \alpha_0)$ of the formal $\Zp[g_0]$\nobreakdash-module~$\Fss$ and this deformation is defined over~$\OK$, see \cite[Proposition~12.10]{HopGro94b}.
  Denote by~$\varphi_0$ the image of~$g_0$ in~$\End_{\OK}(\cF_0)$ by the structural map.
  By definition, ${D_{\cF_0}(\varphi_0) = \delta(g_0)}$.
  Moreover, since ${\alpha_0 \: \tcF_0 \to \Fss}$ is an isomorphism of formal~$\Zp[g_0]$\nobreakdash-modules, ${g_0 = \alpha_0 \circ \tvarphi_0 \circ \alpha_0^{-1}}$.
  By Lemma~\ref{l:fixed}, this proves that~$(\cF_0, \alpha_0)$, seen as a formal $\Zp$\nobreakdash-module over~$\OK$ that is a deformation of~$\Fss$, is a fixed point of~$g_0$.

  It remains to prove the uniqueness statement.
  Let~$(\cF, \alpha)$ in~$\Xss(\OQpalg)$ be another fixed point of~$g_0$, let~$\varphi$ be given by Lemma~\ref{l:fixed}, and suppose
  \begin{equation}
    \label{eq:44}
    D_{\cF}(\varphi) = \delta(g_0).
  \end{equation}
  Let~$\cK'$ be a finite field extension of~$\cK$ contained in~$\Qpalg$ such that~$(\cF, \alpha)$ belongs to~$\Xss(\cO_{\cK'})$.
  Consider~$\cO_{\cK'}$ as a~$\Zp[g_0]$\nobreakdash-module with structural map~$\delta$ and let ${\iota_0 \: \Zp[\varphi] \to \Zp[g_0]}$ be the ring isomorphism defined above.
  Then, \eqref{eq:44} ensures that the ring homomorphim ${\iota_0^{-1} \: \Zp[g_0] \to \End_{\cO_{\cK'}}(\cF)}$ endows~$\cF$ with a structure of formal $\Zp[g_0]$\nobreakdash-module over~$\cO_{\cK'}$.
  Finally, since the deformation space of the formal $\Zp[g_0]$\nobreakdash-module~$\Fss$ consists of a single point, $(\cF, \alpha)$ and~$(\cF_0, \alpha_0)$ are both isomorphic as deformations of~$\Fss$ as a formal $\Zp[g_0]$\nobreakdash-module.
  It follows that they are isomorphic as deformations of~$\Fss$ as a formal $\Zp$\nobreakdash-module.
\end{proof}

\begin{remark}\label{r:Gross-canonical-lift}
  Proposition~\ref{p:formal-lifting} relates to \cite[Proposition~2.1]{Gro86} as follows.
  Let~$g_0$ and~$\delta$ be as in the proposition and put ${\cK_0\=\Qp(\delta(g_0))}$.
  The inverse of~$\delta$ gives an embedding ${\iota_{\dag} \: \cO_{\cK_0}\to \Rss}$ that is \emph{normalized} in the sense of \cite[Section~2]{Gro86}, and the unique fixed point of~$g_0$ in~$\Xss(\OQpalg)$ attached to~$\delta$ is the \emph{canonical lifting} of the pair~$(\Fss,\iota_{\dag})$ in the sense of \cite[Section~3]{Gro86}.
\end{remark}

\begin{lemma}
  \label{l:fixed-properties}
  For every~$\ss$ in~$\tSups$, the following properties hold.
  \begin{enumerate}
  \item[$(i)$]
    For each element~$g$ of~$\Gss \ssetminus \Zp^{\times}$, every fixed point of~$g$ in~$\hDss$ belongs to~$\Xss(\OQpalg)$.
  \item[$(ii)$]
    Let~$\varphi$ be in~$\Lfss$.
    If~$\Qp(\varphi)$ is ramified (resp. unramified) over~$\Qp$, then~$\Fixss(\varphi)$ has precisely two elements (resp.~one element).
  \item[$(iii)$]
    Let~$g$ in~$\Gss \ssetminus \Zp^{\times}$ be such that ${\Zp[g]=\cO_{\Qp(g)}}$.
    Then, an element~$g'$ of ${\Gss \ssetminus \Zp^{\times}}$ has a common fixed point with~$g$ in~$\hDss$ if and only if~$g'$ belongs to~$\Qp(g)$.
  \item[$(iv)$]
    For all~$\varphi$ and~$\varphi'$ in~$\Lfss$, the sets~$\Fixss(\varphi')$ and~$\Fixss(\varphi)$ coincide if~$\varphi'$ belongs to~$\Qp(\varphi)$ and they are disjoint if~$\varphi'$ is outside~$\Qp(\varphi)$.
  \end{enumerate}
\end{lemma}

\begin{proof}
  Since~$g$ acts as a power series~$f$ with coefficients in~$\Z_{p^2}$ (Lemma~\ref{l:action-regularity}$(i)$), item~$(i)$ follows by applying, \emph{e.g.}, \cite[Exercise~3.2.2(1)]{FrevdP04} to the restriction of the power series ${f(z) - z}$ to an affinoid subdomain of~$\hDss$ containing a given fixed point of~$g$.

  For item~$(ii)$, the number of continuous ring homomorphisms ${\cO_{\Qp(\varphi)} \to \OQpalg}$ that reduce to~$D_{\Fss}|_{\cO_{\Qp(\varphi)}}$ equals two (resp.~one) if~$\Qp(\varphi)$ is ramified (resp. unramified) over~$\Qp$.
  Since ${\Zp[\Uss(\varphi)] = \cO_{\Qp(\varphi)}}$ by Lemma~\ref{l:unit-function}, Proposition~\ref{p:formal-lifting} with~$g_0 = \Uss(\varphi)$ yields the desired assertion.

  To prove item~$(iii)$, consider a fixed point of~$g$ in~$\hDss$.
  By item~$(i)$, this point belongs to~$\Xss(\OQpalg)$ and thus it is represented by a pair~$(\cF, \alpha)$.
  If~$g'$ fixes~$(\cF, \alpha)$, then by Lemma~\ref{l:fixed} both~$g$ and~$g'$ belong to the image of the map~$\Aut(\cF) \to \Gss$ given by~$\phi \mapsto \alpha \circ \tphi \circ \alpha^{-1}$.
  By the hypothesis that~$g$ is outside~$\Zp^{\times}$ and \cite[Chapter~IV, Section~1, Theorem~1$(iii)$]{Fro68}, this implies that~$g'$ belongs to~$\Qp(g)$.
  Conversely, every element of~$(\Zp[g])^{\times}$, which equals~$\Gss \cap \Qp(g)$, belongs to the image of the map ${\phi \mapsto \alpha \circ \tphi \circ \alpha^{-1}}$ and thus fixes~$(\cF, \alpha)$ by Lemma~\ref{l:fixed}.

  To prove item~$(iv)$, suppose that~$\Uss(\varphi)$ and~$\Uss(\varphi')$ have a common fixed point.
  Item~$(iii)$ implies ${\Qp(\varphi) = \Qp(\varphi')}$.
  By item~$(i)$, every element of~$\Fixss(\varphi)$ belongs to~$\Xss(\OQpalg)$ and it thus is represented by a pair~$(\cF, \alpha)$.
  By Lemma~\ref{l:fixed}, the image of the map ${\Aut(\cF) \to \Gss}$ given by~$\phi \mapsto \alpha \circ \tphi \circ \alpha^{-1}$ equals~$\cO_{\Qp(\varphi)}^{\times}$ and thus~$\cO_{\Qp(\varphi')}^{\times}$.
  Using Lemma~\ref{l:fixed} again, it follows that~$(\cF, \alpha)$ belongs to~$\Fixss(\varphi')$.
  This proves ${\Fixss(\varphi) \subseteq \Fixss(\varphi')}$.
  Reversing the roles of~$\varphi$ and~$\varphi'$, it follows that these sets are equal.
\end{proof}

Let~$d$ be a $p$-supersingular fundamental discriminant and put
\begin{displaymath}
  \epsilon_{d}
  \=
  \begin{cases}
    1
    & \text{if $p$ ramifies in~$\Q(\sqrt{d})$};
    \\
    1/2
    & \text{if $p$ is inert in~$\Q(\sqrt{d})$}.
  \end{cases}
\end{displaymath}
For all~$\ss$ in~$\tSups$ and~$f$ in~$\Z_{>0}$, denote by~$h(df^2, \ss)$ the number of conjugacy classes of optimal embeddings ${\cO_{d, f} \rightarrow \End(\ss)}$.
The work of Deuring~\cite{Deu41} yields
\begin{equation}
  \label{eq:conteo}
  \deg(\Lambda_d|_{\Dss})
  =
  \epsilon_d h(d, \ss),
\end{equation}
see \cite[Lemma~3.3]{ElkOnoYan05}.

\begin{proof}[Proof of Theorem~\ref{t:fixed-points-formula}]
  Combining \cite[(12.8) and Proposition~12.9]{Gro87}, Lemma~\ref{l:fixed-properties}$(ii)$, and~\eqref{eq:conteo}, yields
  \begin{displaymath}
    \deg(\Lambda_d|_{\Dss})
    =
    \epsilon_d h(d,e)
    =
    \frac{w_{d, 1}}{\#\Aut(\ss)}\#V_{|d|}(\ss)2\epsilon_d.
  \end{displaymath}
  Since the divisor on right side of~\eqref{eq:43} has integer coefficients and its degree equals the last term above,
  it suffices to prove that~$\supp(\Lambda_d|_{\Dss})$ is contained in the support the right side of~\eqref{eq:43}.
  To do this, let~$E$ be in~$\supp( \Lambda_d|_{\Dss})$ and ${\alpha \: \tE \to \ss}$ an isomorphism.
  Since~$E$ is a \CM{} point, it is defined over~$\Qpalg$ and thus ${E = \Piss((\FE,\halpha))}$.
  Moreover, since~$\End(E)$ is isomorphic to~$\cO_{d, 1}$, which equals~$\Z\left[\frac{d+\sqrt{d}}{2}\right]$, there exists an element~$\phi$ in~$\Z + 2\End(E)$ satisfying the equation $X^2 - d = 0$.
  This implies that the endomorphism~$\phi_0$ of~$\ss$, defined by ${\phi_0 \= \alpha\circ \tphi \circ \alpha^{-1}}$, belongs to~$L(\ss)$ and satisfies
  \begin{displaymath}
    \tr(\phi_0)=0
    \text{ and }
    \nr(\phi_0)=|d|.
  \end{displaymath}
  That is, $\phi_0$ belongs to~$V_{|d|}(\ss)$.
  Furthermore, the element~$\hphi_0$ of~$\Rss$ is the image of~$\hphi$ by the ring homomorphism
  \begin{displaymath}
    \begin{array}{rrcl}
      \iota \: & \End(\FE) & \to & \Rss
      \\ & \varphi & \mapsto & \halpha \circ \tvarphi \circ \halpha^{-1}.
    \end{array}
  \end{displaymath}
  Since~$\frac{d+\phi}{2}$ belongs to~$\End(E)$, it follows that ${\Zp\left[\frac{d+\hphi_0}{2}\right] \subseteq \iota(\End(\FE))}$.
  Since ${\Zp \left[ \Uss(\hphi_0) \right] = \Zp\left[\frac{d+\hphi_0}{2}\right]}$, Lemma~\ref{l:unit-function}$(ii)$ with ${\varphi = \hphi_0}$ implies that~$\Uss(\hphi_0)$ is a unit in~$\Zp\left[\frac{d+\hphi_0}{2}\right]$.
  Thus, $\Uss(\hphi_0)$ belongs to~$\iota(\Aut(\FE))$ and Lemma~\ref{l:fixed} implies that~$(\FE,\halpha)$ is a fixed point of $\Uss(\hphi_0)$.
  Since ${E = \Piss((\FE,\halpha))}$, this proves that~$E$ belongs to the support of the right side of~\eqref{eq:43}.
\end{proof}

\subsection{\CM{} points and the canonical branch of $T_p$}
\label{ss:CM-from-canonical}
This section proves the following formulae for~\CM{} points in~$\Sups$ whose discriminant has conductor divisible by~$p$.
The formula relies on the canonical branch of~$T_p$, whose definition follows.
Let ${v_p \: \Sups \to ]0, +\infty]}$ be Katz' valuation, as defined in \cite[Section~4A]{HerMenRiv20}, and put
\begin{displaymath}
  N_p
  \=
  \left\{ E \in \Sups \: v_p(E) < \frac{p}{p + 1} \right\}.
\end{displaymath}
For each~$E$ in~$N_p$, denote by~$H(E)$ the canonical subgroup of~$E$ \cite[Theorem~3.10.7]{Kat73}.
The \emph{canonical branch of~$T_p$} is the map~$\t \: N_p \to \Sups$ defined by ${\t(E) \= E / H(E)}$.

\begin{theorem}
  \label{t:CM-from-canonical}
  For all $p$-supersingular fundamental discriminant~$d$, ${r}$ in~$\Z_{>0}$, and~$f$ in ${\Z_{>0} \ssetminus p\Z}$,
  \begin{displaymath}
    \Lambda_{d (f p^r)^2}
    =
    \begin{cases}
      \t^*\left( \frac{\Lambda_{d f^2}}{w_{d, f}} \right) \bigm\vert_{\kval^{-1}(\frac{1}{2p})}
      & \text{if~$p$ ramifies in~$\Q(\sqrt{d})$ and ${r = 1}$};
      \\
      (\t^*)^{r - 1} (\Lambda_{d (fp)^2})
      & \text{if~$p$ ramifies in~$\Q(\sqrt{d})$ and ${r \ge 2}$};
      \\
      (\t^*)^r \left( \frac{\Lambda_{d f^2}}{w_{d, f}} \right)
      & \text{if~$p$ is inert in~$\Q(\sqrt{d})$}.
    \end{cases}
  \end{displaymath}
\end{theorem}

The proof of this theorem relies on several lemmas.

\begin{lemma}[\textcolor{black}{\cite[Theorem~B.1]{HerMenRiv20}}]
  \label{l:canonical-analyticity}
  The canonical branch~$\t$ of~$T_p$ is given by a finite sum of Laurent series, each of which converges on~$N_p$.
  Furthermore, for every~$E$ in~$\Sups$,
  \begin{equation}
    \label{eq:45}
    T_p(E)
    =
    \begin{cases}
      \t^*(E) + [\t(E)]
      & \text{if } \kval(E) \le \frac{1}{p + 1};
      \\
      \t^*(E)
      & \text{if } \kval(E) > \frac{1}{p + 1}.
    \end{cases}
  \end{equation}
\end{lemma}

The following is \cite[Lemma~4.6]{HerMenRiv20}, which is a reformulation in the current context of \cite[Theorems~3.1 and~3.10.7]{Kat73}.
See also \cite[Theorem~3.3]{Buz03}.
Let ${\kproj \: \Sups \to \left[0, \frac{p}{p + 1} \right]}$ be the map defined by
\begin{equation}
  \label{eq:46}
  \kproj
  \=
  \min\left\{ \kval, \frac{p}{p+1} \right\}.
\end{equation}

\begin{lemma}
  \label{l:sups-canonical-subgroup}
  For every~$E$ in~$N_p$,
  \begin{equation}
    \label{eq:sups-canonical-subgroup-A}
    \kproj(\t(E))
    =
    \begin{cases}
      p\kval(E)
      & \text{if $\kval(E) \in \left] 0, \frac{1}{p + 1} \right]$};
      \\
      1 - \kval(E)
      & \text{if $\kval(E) \in \left] \frac{1}{p + 1}, \frac{p}{p + 1} \right[$}
    \end{cases}
  \end{equation}
  and, for every subgroup~$C$ of~$E$ of order~$p$ that is different from~$H(E)$,
  \begin{equation}
    \label{eq:sups-canonical-subgroup-B}
    \kval(E/C) = p^{-1} \kval(E).
  \end{equation}
  Furthermore, the following properties hold.
  \begin{enumerate}
  \item[$(i)$]
    Let~$E$ be in~$\Sups$ and~$C$ a subgroup of~$E$ of order~$p$.
    When ${\kval(E) < \frac{p}{p + 1}}$, suppose in addition ${C \neq H(E)}$.
    Then,
    \begin{displaymath}
      \kval(E/C)
      =
      p^{-1}\kproj(E)
      \text{ and }
      \t(E/C)
      =
      E.
    \end{displaymath}
  \item[$(ii)$]
    For every~$E$ in~$\Sups$ satisfying ${\frac{1}{p + 1} < \kval(E) < \frac{p}{p + 1}}$, the equality ${\t^2(E) = E}$ holds.
  \end{enumerate}
\end{lemma}

The following lemma is \cite[Lemma~4.9]{HerMenRiv20}, see also \cite[Lemma~4.8]{ColMcM06} and \cite[Proposition~5.3]{Gro86}.

\begin{lemma}
  \label{l:vcm}
  Let~$D$ be a $p$-supersingular discriminant and~$m$ in~$\Z_{\ge 0}$ the largest integer such that~$p^m$ divides the conductor of~$D$.
  Then, for every~$E$ in~$\supp(\Lambda_D)$,
  \begin{displaymath}
    \kproj(E)
    =
    \begin{cases}
      \frac{1}{2} \cdot p^{-m}
      & \text{ if }p \text{ ramifies in } \Q(\sqrt{D});
      \\
      \frac{p}{p + 1} \cdot p^{-m}
      & \text{ if $p$  is inert in $\Q(\sqrt{D})$.}
    \end{cases}
    \footnote{
      When $m = 0$ and~$p$ is inert in~$\Q(\sqrt{D})$, Lemmas~\ref{l:formally-defined} and~\ref{l:katz-deformation} yield ${v_p(E) \ge 1}$.
      So, this formula does not hold with~$\kproj$ replaced by the valuation~$\kval$.
      Compare with \cite[Lemma~4.8]{ColMcM06}.
    }
  \end{displaymath}
\end{lemma}

The following lemma gathers variants of a formula of Zhang in \cite[Proposition~4.2.1]{Zha01}, see also \cite[\emph{Lemme}~2.6]{CloUll04} and \cite[Lemma~2.2]{HerMenRiv20}.
Recall that the \emph{Dirichlet convolution} of two functions~$g, \wtg \: \Z_{>0} \to \C$, is defined by
\begin{displaymath}
  (g \ast \wtg)(n)
  \=
  \sum_{d \in \Z_{>0}, d \mid n} g(d) \wtg \left(\frac{n}{d} \right).
\end{displaymath}
Denote by~$\bfone$ the constant function defined on~$\Z_{>0}$ and taking the value~1.
For each fundamental discriminant~$d$, denote by ${\psi_d \: \Z_{>0} \to \{-1, 0, 1 \}}$ the arithmetic function given by the Kronecker symbol~$\left( \frac{d}{\cdot} \right)$, put ${R_d \= \bfone \ast \psi_d}$, and denote by~$R_d^{-1}$ the inverse of~$R_d$ with respect to the Dirichlet convolution.

\begin{lemma}
  \label{l:Zhang-general}
  For all fundamental discriminant~$d$ and coprime integers~$f$ and~$\wtf$ in~$\Z_{>0}$,
  \begin{align}
    \label{eq:general-inverse-Zhang-formula}
    \frac{\Lambda_{d(f \wtf)^2}}{w_{d,f\wtf}}
    & =
      \sum_{f_0 \in \Z_{>0}, f_0 \mid f} R_d^{-1} \left( \frac{f}{f_0} \right) T_{f_0}\left(\frac{\Lambda_{d\wtf^2}}{w_{d,\wtf}}\right).
      \intertext{If in addition ${p \nmid f}$, then}
      \label{eq:p-Zhang}
      \Lambda_{d (pf)^2}
    & =
      \begin{cases}
        T_p \left( \frac{\Lambda_{d f^2}}{w_{d, f}} \right) - \frac{\Lambda_{d f^2}}{w_{d, f}}
        & \text{if $p$ ramifies in~$\Q(\sqrt{d})$};
        \\
        T_p \left( \frac{\Lambda_{d f^2}}{w_{d, f}} \right)
        & \text{if $p$ is inert in~$\Q(\sqrt{d})$}
      \end{cases}
      \intertext{and, for every integer~$m$ satisfying ${m \ge 2}$,}
      \label{eq:pm-Zhang}
      \Lambda_{d (p^m f)^2}
    & =
      \begin{cases}
        T_{p^m} \left(\frac{\Lambda_{d f^2}}{w_{d, f}} \right)
        - T_{p^{m - 1}} \left(\frac{\Lambda_{d f^2}}{w_{d, f}} \right)
        & \text{if $p$ ramifies in~$\Q(\sqrt{d})$};
        \\
        T_{p^m} \left(\frac{\Lambda_{d f^2}}{w_{d, f}} \right)
        - T_{p^{m - 2}} \left(\frac{\Lambda_{d f^2}}{w_{d, f}} \right)
        & \text{if $p$ is inert in~$\Q(\sqrt{d})$}.
      \end{cases}
  \end{align}
\end{lemma}

\begin{proof}[Proof of Theorem~\ref{t:CM-from-canonical}]
  The first step it to prove that if~$p$ ramifies (resp.~is inert) in~$\Q(\sqrt{d})$, then, for every integer~$r$ satisfying ${r \ge 1}$ (resp. ${r \ge 2}$),
  \begin{equation}
    \label{eq:47}
    T_p(\Lambda_{d (f p^r)^2})
    =
    \Lambda_{d (fp^{(r + 1)})^2} + p \frac{\Lambda_{d (fp^{(r - 1)})^2}}{w_{d, fp^{r - 1}}}.
  \end{equation}
  The proof applies several times the recursive relation~\eqref{eq:Hecke-Tpr} and the formulae~\eqref{eq:p-Zhang} and~\eqref{eq:pm-Zhang} in Lemma~\ref{l:Zhang-general}.
  If~$p$ ramifies in~$\Q(\sqrt{d})$, then
  \begin{multline*}
    w_{d, f} T_p(\Lambda_{d (fp)^2})
    =
    T_p(T_p(\Lambda_{d f^2}) - \Lambda_{d f^2})
    =
    T_{p^2}(\Lambda_{d f^2}) + p \Lambda_{d f^2} - T_p(\Lambda_{d f^2})
    \\ =
    w_{d, f} \Lambda_{d (fp^2)^2} + p \Lambda_{d f^2}
  \end{multline*}
  and, for every~$r \ge 2$,
  \begin{multline*}
    w_{d, f} T_p(\Lambda_{d (f p^r)^2})
    =
    T_p (T_{p^r}(\Lambda_{d f^2}) - T_{p^{r - 1}}(\Lambda_{d f^2}))
    \\ =
    T_{p^{r + 1}}(\Lambda_{d f^2}) + p T_{p^{r - 1}}(\Lambda_{d f^2}) - T_{p^r}(\Lambda_{d f^2}) - p T_{p^{r - 2}}(\Lambda_{d f^2}))
    \\ =
    w_{d, f} (\Lambda_{d (fp^{(r + 1)})^2} + p \Lambda_{d (fp^{(r - 1)})^2}).
  \end{multline*}
  If~$p$ is inert in~$\Q(\sqrt{d})$, then
  \begin{multline*}
    w_{d, f} T_p(\Lambda_{d (fp^2)^2})
    =
    T_p(T_{p^2}(\Lambda_{d f^2}) - \Lambda_{d f^2})
    =
    T_{p^3}(\Lambda_{d f^2}) + p T_p(\Lambda_{d f^2}) - T_p(\Lambda_{d f^2})
    \\ =
    w_{d, f} (\Lambda_{d (fp^3)^2} + p \Lambda_{d (fp)^2}),
  \end{multline*}
  and, for every~$r \ge 3$,
  \begin{multline*}
    w_{d, f} T_p(\Lambda_{d (f p^r)^2})
    =
    T_p(T_{p^r}(\Lambda_{d f^2}) - T_{p^{r - 2}}(\Lambda_{d f^2}))
    \\ =
    T_{p^{r + 1}}(\Lambda_{d f^2}) + p T_{p^{r - 1}}(\Lambda_{d f^2}) - T_{p^{r - 1}}(\Lambda_{d f^2}) - p T_{p^{r - 3}}(\Lambda_{d f^2})
    \\ =
    w_{d, f} (\Lambda_{d (fp^{(r + 1)})^2} + p \Lambda_{d (fp^{(r - 1)})^2}).
  \end{multline*}
  This completes the proof of~\eqref{eq:47}.

  By~\eqref{eq:45} in Lemma~\ref{l:canonical-analyticity} and Lemma~\ref{l:vcm}, for every~$r$ in~$\Z_{>0}$,
  \begin{equation}
    \label{eq:48}
    T_p(\Lambda_{d (f p^r)^2})
    =
    \t^*(\Lambda_{d (f p^r)^2}) + \t_*(\Lambda_{d (f p^r)^2}).
  \end{equation}
  Using Lemmas~\ref{l:sups-canonical-subgroup} and~\ref{l:vcm} to compare the support of the last divisor above with that of the right side of~\eqref{eq:47}, it follows that, if~$p$ ramifies (resp.~is inert) in~$\Q(\sqrt{d})$, then, for every integer~$r$ satisfying ${r \ge 1}$ (resp. ${r \ge 2}$),
  \begin{equation}
    \label{eq:49}
    \t^*(\Lambda_{d (f p^r)^2})
    =
    \Lambda_{d (fp^{(r + 1)})^2}.
  \end{equation}
  Suppose~$p$ ramifies in~$\Q(\sqrt{d})$.
  Then, \eqref{eq:49} implies, for every integer~$r$ satisfying ${r \ge 2}$,
  \begin{displaymath}
    \Lambda_{d (f p^r)^2}
    =
    (\t^*)^{r - 1}\Lambda_{d (fp)^2}.
  \end{displaymath}
  Moreover, by~\eqref{eq:p-Zhang} in Lemma~\ref{l:Zhang-general}, \eqref{eq:45} in Lemma~\ref{l:canonical-analyticity}, and Lemma~\ref{l:vcm},
  \begin{displaymath}
    w_{d, f} \Lambda_{d (fp)^2}
    =
    T_p(\Lambda_{d f^2}) - \Lambda_{d f^2}
    =
    \t^*(\Lambda_{d f^2}) - \Lambda_{d f^2},
  \end{displaymath}
  so Lemma~\ref{l:vcm} implies ${w_{d, f} \Lambda_{d (fp)^2} = \t^*(\Lambda_{d f^2})|_{\kval^{-1} \left( \frac{1}{2p} \right)}}$.
  Suppose~$p$ is inert in~$\Q(\sqrt{d})$.
  Then, by~\eqref{eq:p-Zhang} and~\eqref{eq:pm-Zhang} in Lemma~\ref{l:Zhang-general}, \eqref{eq:45} in Lemma~\ref{l:canonical-analyticity}, and Lemma~\ref{l:vcm},
  \begin{equation}
    \label{eq:50}
    w_{d, f} \Lambda_{d (fp)^2}
    =
    T_p(\Lambda_{d f^2})
    =
    \t^* (\Lambda_{d f^2})
  \end{equation}
  and
  \begin{displaymath}
    w_{d, f} T_p(\Lambda_{d (fp)^2})
    =
    T_p \left( T_p(\Lambda_{d f^2}) \right)
    =
    T_{p^2}(\Lambda_{d f^2}) + p \Lambda_{d f^2}
    =
    w_{d, f} \Lambda_{d (fp^2)^2} + (p + 1) \Lambda_{d f^2}.
  \end{displaymath}
  Using Lemmas~\ref{l:sups-canonical-subgroup} and~\ref{l:vcm} to compare the support of the last divisor above with that of the right side of~\eqref{eq:48} with~$r = 1$, it follows that~$\t^* (\Lambda_{d (fp)^2}) = \Lambda_{d (fp^2)^2}$.
  Combined with~\eqref{eq:49} and~\eqref{eq:50}, this implies that, for every~$r$ in~$\Z_{>0}$, the equality ${\Lambda_{d (f p^r)^2} = (\t^*)^r\left( \frac{\Lambda_{d f^2}}{w_{d, f}} \right)}$ holds.
\end{proof}

\subsection{Formal \CM{} points formulae}
\label{ss:formal-CM-formulae}
This section proves the following formulae for formal \CM{} points.
It relies on the canonical branch~$\t$ of~$T_p$ and Katz' valuation~$\kval$, recalled in Section~\ref{ss:CM-from-canonical}.
For all fundamental $p$\nobreakdash-adic discriminant~$\pfd$ and~$m$ in~$\Z_{\ge 0}$, define the affinoid~$A_{\pfd p^{2m}}$ by
\begin{equation}
  \label{eq:51}
  A_{\pfd p^{2m}}
  \=
  \begin{cases}
    \kval^{-1}(\frac{1}{2} \cdot p^{-m})
    & \text{if~$\Qpfd$ is ramified over~$\Qp$};
    \\
    \kval^{-1} ([1, +\infty])
    & \text{if~$\Qpfd$ is unramified over~$\Qp$ and~$m = 0$};
    \\
    \kval^{-1}(\frac{p}{p + 1} \cdot p^{-m})
    & \text{if~$\Qpfd$ is unramified over~$\Qp$ and~$m \ge 1$}.
  \end{cases}
\end{equation}

\begin{theorem}
  \label{t:formal-CM-formulae}
  Every formal \CM{} point has supersingular reduction.
  Furthermore, for every fundamental $p$\nobreakdash-adic discriminant~$\pfd$, the following properties hold.
  \begin{enumerate}
  \item [$(i)$]
    The inclusion ${\Lambda_{\pfd} \subseteq A_{\pfd}}$ holds.
    Moreover, ${\t(\Lambda_{\pfd}) = \Lambda_{\pfd}}$ if~$\Qpfd$ is ramified over~$\Qp$ and ${T_p(\Lambda_{\pfd}) = \Lambda_{\pfd p^2}}$ if~$\Qpfd$ is unramified over~$\Qp$.
  \item [$(ii)$]
    For every~$m$ in~$\Z_{>0}$,
    \begin{displaymath}
      \Lambda_{\pfd p^{2m}}
      =
      (\t^m \bigm\vert_{A_{\pfd p^{2m}}} )^{-1}(\Lambda_{\pfd}).
    \end{displaymath}
  \end{enumerate}
\end{theorem}

Recall that every $p$\nobreakdash-supersingular discriminant belongs to a unique $p$\nobreakdash-adic discriminant (Lemma~\ref{l:p-adic-discriminants}).

\begin{coro}
  \label{c:formalization}
  The following properties hold.
  \begin{enumerate}
  \item [$(i)$]
    A \CM{} point~$E$ is a formal \CM{} point if and only if it has supersingular reduction.
    In this case, the $p$\nobreakdash-adic discriminant of~$E$ is the unique $p$\nobreakdash-adic discriminant containing the discriminant of~$E$.
  \item [$(ii)$]
    Let~$D$ be a discriminant and~$\pd$ a $p$\nobreakdash-adic discriminant.
    If~$D$ belongs to~$\pd$, then ${\supp(\Lambda_D) \subseteq \Lambda_{\pd}}$, and if~$D$ is outside~$\pd$, then ${\supp(\Lambda_D) \cap \Lambda_{\pd} = \emptyset}$.
  \end{enumerate}
\end{coro}

The proof of Theorem~\ref{t:formal-CM-formulae} and Corollary~\ref{c:formalization} rely on several lemmas.

For every formal \CM{} point~$E$, the height of~$\FE$ is equal to two, see, \emph{e.g.}, \cite[Chapter~IV, Section~1, Theorem~1$(iii)$]{Fro68} and \cite[Chapter~IV, Corollary~7.5]{Sil09}, and thus~$E$ has supersingular reduction, see, \emph{e.g.}, \cite[Chapter~V, Theorem~3.1]{Sil09}.
In particular, for every~$m$ in~$\Z_{>0}$, the torsion subgroup~$E[p^m]$ of~$E$ is contained in the kernel of the reduction morphism ${E(\Qpalg) \to \tE(\Fpalg)}$.
In what follows, consider each endomorphism~$\varphi$ in~$\End(\FE)$ as acting on the kernel of the reduction morphism ${E(\Qpalg) \to \tE(\Fpalg)}$, see, \emph{e.g.}, \cite[Chapter~VII, Propositions~2.1 and~2.2]{Sil09}.
In particular, $\Ker(\varphi)$ is a subgroup of~$E(\Qpalg)$ and, for every~$m$ in~$\Z_{>0}$, the map~$\varphi$ is defined on~$E[p^m]$.

Let~$\pfd$ be a fundamental $p$\nobreakdash-adic discriminant.
Let~$\upsilon_{\pfd}$ be a uniformizer of~$\cO_{\Qpfd}$ if~$\Qpfd$ is ramified over~$\Qp$, and an element of~$\cO_{\Qpfd}$ whose reduction is outside~$\Fp$ if~$\Qpfd$ is unramified over~$\Qp$.
In all cases, for every~$m$ in~$\Z_{\ge 0}$,
\begin{equation}
  \label{eq:52}
  \Zp[v_{\pfd} p^m]
  =
  \Zp + p^m \cO_{\Qpfd}.
\end{equation}
For all~$m$ in~$\Z_{\ge 0}$ and~$E$ in~$\Lambda_{\pfd p^{2m}}$, let~$\varphi_E$ be an element of~$\End(\FE)$ with the same trace and norm as~$\upsilon_{\pfd} p^m$.
Then, $\varphi_E$ is outside~$p \End(\FE)$.
Conversely, if~$E$ belongs to~$\Ell(\Qpalg)$ and there is an element of ${\End(\FE) \ssetminus p \End(\FE)}$ with the same trace and norm as~$\upsilon_{\pfd} p^m$, then~$\End(\FE)$ is a $p$\nobreakdash-adic quadratic order isomorphic to~\eqref{eq:52} and thus~$E$ belongs to~$\Lambda_{\pfd p^{2m}}$.

\begin{lemma}
  \label{l:canonically-algebraic}
  Let~$\pfd$ be a fundamental $p$\nobreakdash-adic discriminant, $m$ in~$\Z_{\ge 0}$, and~$E$ in~$\Lambda_{\pfd p^{2m}}$.
  When~$\Qpfd$ is unramified over~$\Qp$, suppose in addition ${m \ge 1}$.
  Then, there is a unique subgroup~$C_E$ of~$\Ker(\varphi_E)$ of order~$p$ and the following properties hold.
  \begin{enumerate}
  \item [$(i)$]
    The quotient~$E / C_E$ belongs to~$\Lambda_{\pfd p^{2(m - 1)}}$ if ${m \ge 1}$ and to~$\Lambda_{\pfd}$ if ${m = 0}$.
  \item [$(ii)$]
    If~$C$ is a subgroup of order~$p$ of~$E$ different from~$C_E$, then~$E / C$ belongs to~$\Lambda_{\pfd p^{2(m + 1)}}$.
  \end{enumerate}
\end{lemma}

The proof of this lemma relies on a couple of lemmas.

\begin{lemma}
  \label{l:ideally-formal}
  For all formal \CM{} point~$\check{E}$ and~$m$ in~$\Z_{>0}$, an element~$\psi$ of~$\End(\cF_{\check{E}})$ belongs to~$p^m \End(\cF_{\check{E}})$ if and only if ${\check{E}[p^m] \subseteq \Ker(\psi)}$.
\end{lemma}

\begin{proof}
  Let~$[p^m]_{\check{E}}$ be the morphism of multiplication by~$p^m$ on~$\check{E}$.

  For every element~$\psi$ in~$p^m \End(\cF_{\check{E}})$, there is~$\psi'$ in~$\End(\cF_{\check{E}})$ such that ${\psi = \psi' \circ \widehat{[p^m]}_{\check{E}}}$, so
  \begin{equation}
    \label{eq:53}
    \Ker(\psi)
    \supseteq
    \Ker(\widehat{[p^m]}_{\check{E}})
    =
    \check{E}[p^m].
  \end{equation}
  Conversely, for every~$\psi$ in~$\End(\cF_{\check{E}})$ satisfying ${\check{E}[p^m] \subseteq \Ker(\psi)}$, there is~$\psi'$ in~$\End(\cF_{\check{E}})$ such that~$\psi = \psi' \circ \widehat{[p^m]}_{\check{E}}$ \cite[Theorem~1.5]{Lub67}, so~$\psi$ belongs to~$p^m \End(\cF_{\check{E}})$.
\end{proof}

Recall that, for all~$E$ and~$E'$ in~$\Ell(\Cp)$ and isogeny ${\phi \: E \to E'}$, the function ${\hphi \: \FE \to \FEpr}$ is induced by~$\phi$.

\begin{lemma}
  \label{l:canonical-generators}
  Let~$E$ be a formal \CM{} point and~$\varphi$ in ${\End(\FE) \ssetminus p \End(\FE)}$ such that~$\varphi \circ \varphi$ belongs to~$p \End(\FE)$.
  Then, the following properties hold.
  \begin{enumerate}
  \item [$(i)$]
    There is a unique subgroup~$C_0$ of~$E$ of order~$p$ contained in~$\Ker(\varphi)$.
    Moreover, $\varphi(E[p]) = C_0$.
  \item [$(ii)$]
    Let~$C$ be a subgroup of~$E$ of order~$p$, put ${E' \= E/C}$, and let ${\phi \: E \to E'}$ be an isogeny satisfying ${\Ker(\phi) = C}$.
    Then, $\hphi \circ \varphi \circ \widehat{\overline{\phi}}$ belongs to~$p \End(\FEpr)$ if and only if~$C = C_0$.
  \item [$(iii)$]
    Suppose that~$\nr(\varphi)$ belongs to~$p^2 \Zp$, put ${E_0 \= E/C_0}$, and let ${\phi \: E \to E_0}$ be an isogeny satisfying ${\Ker(\phi) = C_0}$.
    Then, there is~$\varphi_0$ in ${\End(\cF_{E_0}) \ssetminus p \End(\cF_{E_0})}$ such that ${\varphi = \widehat{\overline{\phi}} \circ \varphi_0 \circ \hphi}$.
  \end{enumerate}
\end{lemma}

\begin{proof}
  The proof applies Lemma~\ref{l:ideally-formal} repeatedly without mentioning it.

  To prove item~$(i)$, put ${C_0\= \varphi(E[p])}$.
  The group~$\Ker(\varphi)$ cannot contain two distinct subgroups of order~$p$ of~$E$, for otherwise ${E[p] \subseteq \Ker(\varphi)}$ and~$\varphi$ would belong to~$p \End(\FE)$.
  The hypothesis that~$\varphi \circ \varphi$ belongs to~$p \End(\FE)$ implies ${E[p] \subseteq \Ker(\varphi \circ \varphi)}$ and thus ${\varphi(E[p]) \subseteq \Ker(\varphi)}$.
  The group~$\varphi(E[p])$ cannot be reduced to the neutral element of~$E$ because ${E[p] \not\subseteq \Ker(\varphi)}$.
  Moreover, ${\varphi(E[p]) \neq E[p]}$ because ${E[p] \subseteq \Ker(\varphi \circ \varphi)}$.
  This implies that~$C_0$ is the unique subgroup of order~$p$ of~$\Ker(\varphi)$.

  For item~$(ii)$, ${C = C_0}$ implies
  \begin{displaymath}
    \widehat{\overline{\phi}}(E'[p])
    =
    C_0
    \subseteq
    \Ker(\varphi),
  \end{displaymath}
  so ${E'[p] \subseteq \Ker(\hphi \circ \varphi \circ \widehat{\overline{\phi}})}$, and thus~$\hphi \circ \varphi \circ \widehat{\overline{\phi}}$ belongs to~$p \End(\FEpr)$.
  If ${C \neq C_0}$, then item~$(i)$ implies
  \begin{displaymath}
    (\varphi \circ \widehat{\overline{\phi}})(E'[p])
    =
    \varphi(C)
    =
    \varphi(E[p])
    =
    C_0.
  \end{displaymath}
  The group~$C_0$ is not contained in~$\Ker(\phi)$, so ${\Ker(\hphi \circ \varphi \circ \widehat{\overline{\phi}}) \not\subseteq E'[p]}$ and thus ${\hphi \circ \varphi \circ \widehat{\overline{\phi}}}$ is outside~$p \End(\FEpr)$.

  For item~$(iii)$, the additional hypothesis implies ${E \left[ p^2 \right] \subseteq \Ker(\varphi \circ \overline{\varphi})}$, so ${\overline{\varphi}(E \left[ p^2 \right]) \subseteq \Ker(\varphi)}$.
  By item~$(i)$ applied to~$\overline{\varphi}$, the group~$\overline{\varphi}(E \left[ p \right])$ has order~$p$.
  Since
  \begin{displaymath}
    p\overline{\varphi}(E \left[ p^2 \right])
    =
    \overline{\varphi}(pE \left[ p^2 \right])
    =
    \overline{\varphi}(E \left[ p \right]),
  \end{displaymath}
  it follows that~$\overline{\varphi}(E \left[ p^2 \right])$ contains a cyclic subgroup~$\hC$ of~$E$ of order~$p^2$.
  Moreover, $C_0$ is the unique group of order~$p$ contained in~$\Ker(\varphi)$, so
  \begin{displaymath}
    p\hC
    =
    C_0
    =
    \overline{\phi}(E_0[p])
    =
    p\overline{\phi}(E_0\left[p^2\right]).
  \end{displaymath}
  Combined with
  \begin{displaymath}
    \overline{\phi} \circ \phi(E\left[p^2\right])
    =
    pE\left[p^2\right]
    =
    E[p]
    \subseteq
    \overline{\phi}(E_0[p^2]),
  \end{displaymath}
  this implies ${\hC \subseteq \overline{\phi}(E_0\left[p^2\right])}$.
  Since~$\overline{\phi}(E_0[p^2])$ also contains~$E[p]$ and is of order~$p^3$, it follows that
  \begin{equation}
    \label{eq:54}
    \overline{\phi} \left( E_0 \left[ p^2 \right] \right)
    =
    E[p] + \hC.
  \end{equation}
  Furthermore, ${\hC \subseteq \Ker(\varphi)}$ and item~$(i)$ imply
  \begin{displaymath}
    \varphi(E[p] + \hC)
    =
    \varphi(E[p])
    =
    C_0
    =
    \Ker(\phi).
  \end{displaymath}
  Together with~\eqref{eq:54}, this implies ${E_0 \left[ p^2 \right] \subseteq \Ker(\hphi \circ \varphi \circ \widehat{\overline{\phi}})}$.
  So, there is~$\varphi_0$ in~$\End(\cF_{E_0})$ such that ${p^2 \varphi_0 = \hphi \circ \varphi \circ \widehat{\overline{\phi}}}$ and thus ${\widehat{\overline{\phi}} \circ \varphi_0 \circ \hphi = \varphi}$.
  Finally, $\varphi_0$ cannot belong to~$p \End(\cF_{E_0})$ for otherwise~$\varphi$ would belong to~$p \End(\FE)$.
\end{proof}

\begin{proof}[Proof of Lemma~\ref{l:canonically-algebraic}]
  The hypotheses imply that~$\varphi_E \circ \varphi_E$ belongs to~$p \End(\FE)$, so Lemma~\ref{l:canonical-generators}$(i)$ with ${\varphi = \varphi_E}$ yields the first assertion.

  To prove item~$(i)$, put ${E_0 \= E / C_E}$ and let ${\phi \: E \to E_0}$ be an isogeny satisfying ${\Ker(\phi) = C_E}$.
  Suppose ${m \ge 1}$, so~$\nr(\varphi_E)$ belongs to~$p^2 \Zp$.
  Then, the element~$\varphi_0$ of ${\End(\cF_{E_0}) \ssetminus p \End(\cF_{E_0})}$ given by Lemma~\ref{l:canonical-generators}$(iii)$ with ${\varphi = \varphi_E}$ has the same trace and norm as~$\upsilon_{\pfd} p^{m - 1}$.
  This implies that~$E_0$ belongs to~$\Lambda_{\pfd p^{2(m - 1)}}$.
  It remains to consider the case where ${m = 0}$.
  By hypothesis, in this case~$\Qpfd$ is ramified over~$\Qp$.
  By Lemma~\ref{l:canonical-generators}$(ii)$ with~$\varphi = \varphi_E$, there is~$\varphi_0'$ in~$\End(\cF_{E_0})$ such that ${p \varphi_0' = \hphi \circ \varphi_E \circ \widehat{\overline{\phi}}}$.
  It follows that~$\varphi_0'$ has the same trace and norm as~$\upsilon_{\pfd}$ and thus that~$E_0$ belongs to~$\Lambda_{\pfd}$.

  To prove item~$(ii)$, put ${E' \= E / C}$ and let ${\phi \: E \to E'}$ be an isogeny satisfying ${\Ker(\phi) = C}$.
  Then, the endomorphism ${\hphi \circ \varphi_E \circ \widehat{\overline{\phi}}}$ of~$\FEpr$ has the same norm and trace as~$\upsilon_{\pfd} p^{m + 1}$.
  By Lemma~\ref{l:canonical-generators}$(ii)$ with ${\varphi = \varphi_E}$, this endomorphism is outside~$p \End(\FEpr)$, so~$E'$ belongs to~$\Lambda_{\pfd p^{2(m + 1)}}$.
\end{proof}

For each fundamental $p$\nobreakdash-adic discriminant~$\pfd$, denote by~$\Q_{p^2}(\sqrt{\pfd})$ the compositum of~$\Q_{p^2}$ and~$\Qpfd$.

\begin{lemma}
  \label{l:formally-defined}
  For all~$\ss$ in~$\tSups$ and fundamental $p$\nobreakdash-adic discriminant~$\pfd$,
  \begin{equation}
    \label{eq:55}
    \Piss^{-1}(\Lambda_{\pfd} \cap \Dss)
    \subseteq
    \Xss(\cO_{\Q_{p^2}(\sqrt{\pfd})}).
  \end{equation}
\end{lemma}

\begin{proof}
  Let~$(\cF, \alpha)$ be a point in ${\Piss^{-1}(\Lambda_{\pfd} \cap \Dss)}$, $\cO$ the image of~$\End(\FE)$ by ${\varphi \mapsto \alpha \circ \tvarphi \circ \alpha^{-1}}$, and~$g_0$ in~$\cO^{\times}$ such that ${\cO = \Zp[g_0]}$.
  Then, $\cO$ is isomorphic to~$\cO_{\Qpfd}$ and~$(\cF, \alpha)$ is a fixed point of~$g_0$ by Lemma~\ref{l:fixed}.
  It follows that the ring homomorphism ${\delta \: \Zp[g_0] \to \OQpalg}$ given by Proposition~\ref{p:formal-lifting} takes values in~$\Qpfd$ and thus that~$(\cF, \alpha)$ belongs to~$\Xss(\cO_{\Q_{p^2}(\sqrt{\pfd})})$.
\end{proof}

\begin{lemma}
  \label{l:katz-deformation}
  For every~$\ss$ in~$\tSups$,
  \begin{displaymath}
    \min\{ v_p \circ \Piss, 1 \}
    =
    \min\{ \ord_p, 1 \}.
  \end{displaymath}
  In particular, the map~$\kproj$ defined by ${\kproj \= \min \{ \kval, \frac{p}{p + 1} \}}$ satisfies ${\kproj \circ \Piss = \min \left\{ \ord_p, \frac{p}{p + 1} \right\}}$.
\end{lemma}

\begin{proof}
  Let~$\delta_{\ss}$ be as in Section~\ref{ss:woods-hole} and~$\sm_{\ss}$ as in \cite[Proposition~4.3]{HerMenRiv20}, so, for every~$E$ in~$\Dss$, the equality ${\kval(E) = \frac{1}{\delta_{\ss}} \ord_p(j(E) - \sm_{\ss})}$ holds.
  Using that~$\sm_{\ss}$ belongs to~$\Z_{p^2}$ \cite[Remark~4.4]{HerMenRiv20} and Theorem~\ref{t:Pi-properties}, the difference~$j \circ \Piss - \sm_{\ss}$ is represented by a power series with coefficients in~$\Z_{p^2}$ that is a ramified covering of degree~$\delta_{\ss}$ from~$\hDss$ to~$\Mp$.
  Thus, denoting by~$x_1$, \ldots, $x_{\delta_{\ss}}$ the zeros of ${j \circ \Piss - \sm_{\ss}}$ repeated according to multiplicity, there is~$h$ in~$\Z_{p^2} \left[\![t]\!\right]$ such that~$|h|_p$ is constant equal to~$1$ and
  \begin{displaymath}
    j \circ \Piss (t) - \sm_{\ss}
    =
    h(t) \cdot \prod_{i = 1}^{\delta_{\ss}} (t - x_i),
  \end{displaymath}
  see, \emph{e.g.}, \cite[Exercise~3.2.2(1)]{FrevdP04}.
  Thus, for every~$x$ in~$\hDss$,
  \begin{equation}
    \label{eq:56}
    \kval \circ \Piss(x)
    =
    \frac{1}{\delta_{\ss}} \ord_p(j \circ \Piss(x) - \sm_{\ss})
    =
    \frac{1}{\delta_{\ss}} \sum_{i = 1}^{\delta_{\ss}} \ord_p(x - x_i).
  \end{equation}
  Moreover, if~$\pfd_0$ denotes the $p$\nobreakdash-adic discriminant of~$\Z_{p^2}$, then~$\sm_{\ss}$ belongs to~$\Lambda_{\pfd_0}$ by \cite[Remark~4.4]{HerMenRiv20}.
  Thus, for each~$i$ in~$\{1, \ldots, \delta_{\ss} \}$, the point~$x_i$ belongs to~$\Xss(\Z_{p^2})$ by Lemma~\ref{l:formally-defined}.
  That is, seen as an element of~$\Mp$, the point~$x_i$ belongs to~$p \Z_{p^2}$.
  In particular, for every~$x$ in~$\hDss$,
  \begin{displaymath}
    \min \{ \ord_p(x - x_i), 1 \}
    =
    \min \{ \ord_p(x), 1 \}.
  \end{displaymath}
  Together with~\eqref{eq:56} this implies the lemma.
\end{proof}

\begin{proof}[Proof of Theorem~\ref{t:formal-CM-formulae}]
  The second paragraph after Corollary~\ref{c:formalization} proves first assertion.

  The proof of item~$(i)$ also establishes that, if~$\Qpfd$ is ramified (resp. unramified) over~$\Qp$, then every~$E$ in~$\Lambda_{\pfd}$ (resp.~$\Lambda_{\pfd p^2}$) is not too supersingular and the group~$C_E$ given by Lemma~\ref{l:canonically-algebraic} is the canonical subgroup~$H(E)$ of~$E$.
  Suppose that~$\Qpfd$ is ramified over~$\Qp$, let~$E$ be in~$\Lambda_{\pfd}$, and put ${E_0 \= E / C_E}$.
  By Lemma~\ref{l:canonically-algebraic}$(i)$, the quotient~$E_0$ belongs to~$\Lambda_{\pfd}$.
  Thus, Lemmas~\ref{l:formally-defined} and~\ref{l:katz-deformation} imply ${\kproj(E) \ge \frac{1}{2}}$ and ${\kproj(E_0) \ge \frac{1}{2}}$.
  Using Lemma~\ref{l:sups-canonical-subgroup} several times yields
  \begin{displaymath}
    \kproj(E)
    =
    \kproj(E_0)
    =
    \frac{1}{2},
    C_E = H(E),
    \t(E) = E_0,
    \text{ and }
    \t(E_0) = E.
  \end{displaymath}
  This implies item~$(i)$ when~$\Qpfd$ is ramified over~$\Qp$.
  Suppose that~$\Qpfd$ is unramified over~$\Qp$.
  Lemmas~\ref{l:formally-defined} and~\ref{l:katz-deformation} yield ${\Lambda_{\pfd} \subseteq A_{\pfd}}$.
  To prove ${\Lambda_{\pfd p^2} \subseteq T_p(\Lambda_{\pfd})}$, let~$E$ be in~$\Lambda_{\pfd p^{2}}$ and put ${E_0 \= E / C_E}$.
  Then, ${E_0}$ belongs to~$\Lambda_{\pfd}$ by Lemma~\ref{l:canonically-algebraic}$(i)$, so~$E$ belongs to~$\supp(T_p(E_0))$ and hence to~$T_p(\Lambda_{\pfd})$.
  By Lemma~\ref{l:sups-canonical-subgroup}, it also follows that~$E$ is not too supersingular and ${C_E = H(E)}$.
  It remains to prove ${T_p(\Lambda_{\pfd}) \subseteq \Lambda_{\pfd p^2}}$.
  To do this, let~$E'$ be in~$\Lambda_{\pfd}$, $E''$ in~$\supp(T_p(E'))$, and ${\phi \: E' \to E''}$ an isogeny of degree~$p$.
  Then, ${\kproj(E'') = \frac{1}{p + 1}}$ by the first assertion of item~$(i)$ and Lemma~\ref{l:sups-canonical-subgroup}, so~$E''$ is outside~$\Lambda_{\pfd}$.
  The endomorphism~$\hphi \circ \varphi_{E'} \circ \widehat{\overline{\phi}}$ of~$\cF_{E''}$ has the same trace and norm as~$\upsilon_{\pfd} p$.
  It follows that~$E''$ belongs to~$\Lambda_{\pfd}$ or~$\Lambda_{\pfd p^2}$.
  But, $E''$ is outside~$\Lambda_{\pfd}$, as shown above, so~$E''$ belongs to~$\Lambda_{\pfd p^2}$.

  The proof of item~$(ii)$ proceeds by induction, showing in addition that, for every~$E$ in~$\Lambda_{\pfd p^{2m}}$, the equality ${C_E = H(E)}$ holds.
  If~$m = 1$ and~$\Qpfd$ is unramified over~$\Qp$, then item~$(i)$ and~\eqref{eq:45} in Lemma~\ref{l:canonical-analyticity} imply that every element~$E$ of~$\Lambda_{\pfd}$ is too supersingular and ${\t^{-1}(E) = T_p(\{E\})}$.
  Using item~$(i)$ again yields item~$(ii)$.
  As shown above, for every~$E$ in~$\Lambda_{\pfd p^2}$ the equality ${C_E = H(E)}$ holds.
  To complete the proof of the base step, suppose ${m = 1}$ and that~$\Qpfd$ is ramified over~$\Qp$.
  Since for every~$\check{E}$ in~$\Lambda_{\pfd}$ the equality ${C_{\check{E}} = H(\check{E})}$ holds and item~$(i)$ implies ${\kval(\t(\check{E})) = \frac{1}{2}}$, combining~\eqref{eq:45} in Lemma~\ref{l:canonical-analyticity} and Lemma~\ref{l:canonically-algebraic}$(ii)$ yields ${(\t|_{A_{\pfd p^2}})^{-1}(\Lambda_{\pfd}) \subseteq \Lambda_{\pfd p^2}}$.
  To prove the reverse inclusion, let~$E$ be in~$\Lambda_{\pfd p^2}$ and put ${E_0 \= E / C_E}$.
  Then, ${E_0}$ belongs to~$\Lambda_{\pfd}$ by Lemma~\ref{l:canonically-algebraic}$(i)$ and item~$(i)$ implies ${\kproj(E_0) = \frac{1}{2}}$.
  If~$\kproj(E) \neq \frac{1}{2p}$, then Lemma~\ref{l:sups-canonical-subgroup} would yield ${\kproj(E) = \frac{1}{2}}$ and ${E = \t(E_0)}$.
  By item~$(i)$, this would imply that~$E$ belongs to~$\Lambda_{\pfd}$.
  This contradiction proves ${\kproj(E) = \frac{1}{2p}}$.
  Using Lemma~\ref{l:sups-canonical-subgroup} again, if follows that ${C_E = H(E)}$ and ${\t(E) = E_0}$.
  This proves ${\Lambda_{\pfd p^2} \subseteq A_{\pfd p^2}}$ and ${\t(\Lambda_{\pfd p^2}) \subseteq \Lambda_{\pfd}}$, and completes the proof of the base step.
  To prove the induction step, let~$m$ in~$\Z_{>0}$ be such that item~$(ii)$ holds and, for every~$E$ in~$\Lambda_{\pfd p^{2 m}}$, the equality ${C_E = H(E)}$ holds.
  Combined with~\eqref{eq:45} in Lemma~\ref{l:canonical-analyticity} and Lemma~\ref{l:canonically-algebraic}$(ii)$, this last property implies
  \begin{equation}
    \label{eq:57}
    \t^{-1}(\Lambda_{\pfd p^{2m}})
    \subseteq
    \Lambda_{\pfd p^{2 (m + 1)}}.
  \end{equation}
  To prove the reverse inclusion, let~$E$ be in~$\Lambda_{\pfd p^{2 (m + 1)}}$ and put ${E_0 \= E / C_E}$.
  Then, ${E_0}$ belongs to~$\Lambda_{\pfd p^{2m}}$ by Lemma~\ref{l:canonically-algebraic}$(i)$, so the induction hypothesis implies
  \begin{displaymath}
    \kproj(E_0)
    =
    \begin{cases}
      \frac{p}{p + 1} \cdot p^{- m}
      & \text{if~$\Qpfd$ is unramified over~$\Qp$};
      \\
      \frac{1}{2} \cdot p^{- m}
      & \text{if~$\Qpfd$ is ramified over~$\Qp$}.
    \end{cases}
  \end{displaymath}
  In particular, $E_0$ is not too supersingular.
  Suppose~$\kproj(E_0) \neq p \kproj(E)$.
  Then, Lemma~\ref{l:sups-canonical-subgroup} would yield ${\kproj(E) = p \kproj(E_0)}$ and ${\t(E_0) = E}$.
  By the induction hypothesis, this would imply that~$E$ belongs to~$\Lambda_{\pfd p^{2(m - 1)}}$, which is absurd.
  This contradiction proves ${\kproj(E_0) = p \kproj(E)}$.
  Using Lemma~\ref{l:sups-canonical-subgroup} again, it follows that ${C_E = H(E)}$ and ${\t(E) = E_0}$.
  This proves
  \begin{displaymath}
    \Lambda_{\pfd p^{2(m + 1)}}
    \subseteq
    A_{\pfd p^{2(m + 1)}}
    \text{ and }
    \t(\Lambda_{\pfd p^{2(m + 1)}})
    \subseteq
    \Lambda_{\pfd p^{2 m}}.
  \end{displaymath}
  Together with~\eqref{eq:57} this completes the proof of the induction step.
\end{proof}

\begin{proof}[Proof of Corollary~\ref{c:formalization}]
  For item~$(i)$, Theorem~\ref{t:formal-CM-formulae} implies that every formal \CM{} point has supersingular reduction.
  To prove the second assertion, suppose~$E$ has supersingular reduction, let~$D$ be the discriminant of~$E$, and~$\pd$ the unique $p$\nobreakdash-adic discriminant containing~$D$.
  Denote by~$d$ and~$f$ the fundamental discriminant and conductor of~$D$, respectively, so~$D = d f^2$ and~$\End(E)$ is isomorphic to~$\cO_{d, f}$.
  Moreover, let~$\pfd$ be the fundamental $p$\nobreakdash-adic discriminant and~$m$ in~$\Z_{\ge 0}$ such that ${\pd = \pfd p^{2m}}$ (Lemma~\ref{l:p-adic-discriminants-Appendix}$(i)$).
  Then, $d$ belongs to~$\pfd$, ${m = \ord_p(f)}$, and~$\End(E) \otimes \Zp$ is a $p$\nobreakdash-adic quadratic order isomorphic to~$\Zp + p^m \cO_{\Qpfd}$.
  In particular, the $p$\nobreakdash-adic discriminant of~$\End(E) \otimes \Zp$ equals~$\pd$ by Lemma~\ref{l:p-adic-discriminants-Appendix}$(ii)$.
  Consider the natural map~$\End(E) \otimes \Zp \to \End(\FE)$ induced by the ring homomorphism~$\End(E) \to \End(\FE)$.
  Its image is a $p$\nobreakdash-adic order of $p$\nobreakdash-adic discriminant~$\pd$.
  This implies that~$\End(\FE)$ is a $p$\nobreakdash-adic quadratic order and that there is~$m'$ in~$\Z_{\ge 0}$ such that the $p$\nobreakdash-adic discriminant of~$\End(\FE)$ equals~$\pfd p^{2 m'}$.
  Combining Lemma~\ref{l:vcm} and Theorem~\ref{t:formal-CM-formulae} it follows that~$m' = m$ and thus that the $p$\nobreakdash-adic discriminant of~$\End(\FE)$ equals~$\pd$.
  Thus, $E$ belongs to~$\Lambda_{\pd}$.

  Since every discriminant in~$\pd$ is $p$-supersingular (Lemma~\ref{l:p-adic-discriminants}), the first assertion of item~$(ii)$ is a direct consequence of item~$(i)$.
  To prove the second assertion, suppose~$D$ is outside~$\pd$.
  If~$D$ is not $p$-supersingular,
  then~$\supp(\Lambda_D)$ is disjoint from~$\Sups$ and thus from~$\Lambda_{\pd}$ by Theorem~\ref{t:formal-CM-formulae}.
  Suppose that~$D$ is $p$-supersingular and let~$\pd'$ be the unique $p$\nobreakdash-adic discriminant containing~$D$.
  Then, ${\supp(\Lambda_{D}) \subseteq \Lambda_{\pd'}}$ by item~$(i)$ and thus ${\supp(\Lambda_{D}) \cap \Lambda_{\pd} = \emptyset}$.

\end{proof}

\section{Asymptotic distribution of \CM{} points of fundamental discriminant}
\label{s:CM-fundamental}
This section proves the following theorem determining the asymptotic distribution of \CM{} points of fundamental discriminant.
It is one of the main ingredients in the proof of Theorem~\ref{t:CM-symmetric}.
Recall that, for all $p$\nobreakdash-adic discriminant~$\pd$ and discriminant~$D$ in~$\pd$, the inclusion ${\supp(\Lambda_D) \subseteq \Lambda_{\pd}}$ holds (Corollary~\ref{c:formalization}$(ii)$).

\begin{theorem}
  \label{t:CM-fundamental}
  For every fundamental $p$\nobreakdash-adic discriminant~$\pfd$, the subset~$\Lambda_{\pfd}$ of~$\Sups$ is compact.
  Moreover, there is a Borel probability measure~$\nu_{\pfd}$ on~$\Ell(\Cp)$ such that ${\supp(\nu_{\pfd}) = \Lambda_{\pfd}}$ and, for all~$\varepsilon$ and~$\delta$ in~$\R_{>0}$, there is~$C$ in~$\R_{>0}$ such that the following property holds.
  For all function ${F \: \Lambda_{\pfd} \to \R}$ that is constant on every ball of~$\Lambda_{\pfd}$ of radius~$\delta$ and fundamental discriminant~$d$ in~$\pfd$,
  \begin{equation}
    \label{eq:58}
    \left| \int F \dd \odelta_{d} - \int F \dd \nu_{\pfd} \right|
    \le
    C \left( \sup_{\Lambda_{\pfd}} |F| \right) |d|^{- \frac{1}{28} + \varepsilon}.
  \end{equation}
\end{theorem}

For each~$\ss$ in~$\tSups$, Section~\ref{ss:trace-zero-spheres} introduces ``zero-trace spheres'' of~$\Lfss$ and proves that each of these sets carries a natural homogeneous measure (Proposition~\ref{p:trace-zero-spheres}).
Sections~\ref{ss:CM-parametrization} and~\ref{ss:trace-zero-spheres-to-CM} prove that zero-trace spheres parametrize formal \CM{} points in~$\Dss$ of fundamental $p$\nobreakdash-adic discriminant, via fixed points of the group action described in Section~\ref{ss:from-elliptic-curves} (Propositions~\ref{p:CM-parametrization} and~\ref{p:trace-zero-spheres-to-CM}).
Section~\ref{ss:proof-CM-fundamental} deduces Theorem~\ref{t:CM-fundamental} from the results on the asymptotic distribution of integer points on $p$\nobreakdash-adic spheres in Section~\ref{s:p-adic-Linnik}, and an equidistribution result for \CM{} points in a supersingular residue disc proved in Section~\ref{ss:residual} (Theorem~\ref{t:residual}).

\subsection{Zero-trace spheres and their homogeneous measures}
\label{ss:trace-zero-spheres}
Throughout this section, fix~$\ss$ in~$\tSups$.
Let~$\Bss$, $\Rss$, and~$\Gss$ be as in Section~\ref{ss:from-elliptic-curves} and let~$L(\ss)$, $\Lss$, and~$\Lfss$ be as in Section~\ref{ss:fixed-points-formula}.
The degree function defines a ternary quadratic form~$Q_{\ss}^0$ on the lattice~$L(\ss)$, which is positive definite and defined over~$\Z$.
Using the natural map~$\End(\ss) \to \Rss$ defined by~$\phi \mapsto \hphi$ as in Section~\ref{ss:from-elliptic-curves}, the quadratic form~$Q_{\ss}^0$ extends to a quadratic form on~$\Lss$ taking values on~$\Zp$.

For each~$\ell$ in~$\Zp$ such that~$-\ell$ belongs to a fundamental $p$\nobreakdash-adic discriminant, put
\begin{displaymath}
  S_{\ell}^0(\ss)
  \=
  \{ \varphi \in \Lss \: \nr(\varphi) = \ell \}.
\end{displaymath}
The set~$S_{\ell}^0(\ss)$ is a \emph{zero-trace sphere}.
It is a metric subspace of~$\Rss$.
This section defines a natural homogeneous measure on each zero-trace sphere.

For each fundamental $p$\nobreakdash-adic discriminant~$\pfd$, put
\begin{equation}
  \label{eq:59}
  \bfL_{\ss, \pfd}
  \=
  \{ \varphi \in \Lss \: - \nr(\varphi) \in \pfd \}.
\end{equation}
Clearly, as~$\pfd$ varies these sets form a partition of~$\Lfss$.
Moreover, for each~$\pfd$,
\begin{equation}
  \label{eq:Spehre0decomposition}
  \bfL_{\ss, \pfd}
  =
  \bigsqcup_{\Delta \in \pfd} S_{- \Delta}^0(\ss)
\end{equation}
is a partition.
The action of~$\Gss$ on~$\Bss$ by conjugation preserves the reduced trace and norm, so it restricts to a left action
\begin{equation}
  \label{eq:actionBeRe}
  \begin{array}{rcl}
    \Gss \times \bfL_{\ss, \pfd} & \to & \bfL_{\ss, \pfd}
    \\
    (g, \varphi) & \mapsto & g \varphi g^{-1}.
  \end{array}
\end{equation}
Moreover, for every~$\Delta$ in~$\pfd$, this action restricts to an action of~$\Gss$ on~$S_{- \Delta}^0(\ss)$, which is the restriction to~$\Gss$ of the action of the orthogonal group~$\oO_{Q_{\ss}^0}(\Zp)$ on~$S_{- \Delta}^0(\ss)$.

\begin{proposition}
  \label{p:trace-zero-spheres}
  For all~$\ss$ in~$\tSups$ and fundamental $p$\nobreakdash-adic discriminant~$\pfd$, the following properties hold.
  \begin{enumerate}
  \item[$(i)$]
    The set~$\bfL_{\ss, \pfd}$ is nonempty and compact and, for every~$\varphi$ in~$\bfL_{\ss, \pfd}$, the subalgebra~$\Qp(\varphi)$ of~$\Bss$ is a field extension of~$\Qp$ isomorphic to~$\Qpfd$.
  \item[$(ii)$]
    For each~$\Delta$ in~$\pfd$, the set~$S_{- \Delta}^0(\ss)$ is nonempty and compact and the group~$\Gss$ acts transitively and by isometries on~$S_{- \Delta}^0(\ss)$.
    In particular, \eqref{eq:Spehre0decomposition} is the partition of~$\bfL_{\ss, \pfd}$ into orbits with respect to the action~\eqref{eq:actionBeRe}.
  \item[$(iii)$]
    For each~$\ell$ in~$\Zp$ such that~$-\ell$ belongs to~$\pfd$, there exists a unique Borel probability measure~$\nu^{\ss}_{\ell}$ on~$S_{\ell}^0(\ss)$ that is invariant under the action of~$\Gss$.
    It is also the unique Borel probability measure on~$S_{\ell}^0(\ss)$ that is invariant under the action of the orthogonal group~$\oO_{Q_{\ss}^0}(\Zp)$.
    Moreover, ${\supp(\nu^{\ss}_{\ell}) = S_{\ell}^0(\ss)}$.
  \item [$(iv)$]
    For all~$\Delta$ in~$\pfd$, $\ss'$ in~$\tSups$, and~$g$ in~$\bfG_{\ss, \ss'}$, the map ${\Rss \to \Rsspr}$ given by ${\varphi \mapsto g \varphi g^{-1}}$ maps~$S_{-\Delta}^0(\ss)$ to~$S_{-\Delta}^0(\ss')$ and~$\nu^{\ss}_{- \Delta}$ to~$\nu^{\ss'}_{- \Delta}$.
  \end{enumerate}
\end{proposition}

The proof of this proposition relies on the following general lemma.

\begin{lemma}
  \label{l:existence-invariant-measure}
  Let~$G$ be a group acting transitively and by isometries on a compact ultrametric space~$X$.
  Then, there is a unique Borel probability measure on~$X$ that is invariant by~$G$.
  Moreover, the support of this measure equals~$X$ and this measure is invariant under every isometry of~$X$.
\end{lemma}

\begin{proof}
  Denote by~$\dist_X$ the distance on~$X$.
  For each~$r$ in~$\R_{>0}$, let~$\sim_r$ be the equivalence relation on~$X$ defined by~$\dist_X(x, x') \le r$, $X_r$ the finite set of equivalence classes of~$\sim_r$, and~$\mu_r$ the probability measure on~$X_r$ assigning the same mass to each element of~$X_r$.
  For every~$r'$ in~$]0, r[$, the partition of~$X$ into equivalence classes of~$\sim_{r'}$ is finer than that of~$\sim_r$.
  The action of~$G$ on~$X$ descends to a transitive action on~$X_{r'}$, so each equivalence class of~$\sim_r$ contains the same number of equivalence classes of~$\sim_{r'}$.
  It follows that the natural projection~$X_{r'} \to X_{r}$ maps~$\mu_{r'}$ to~$\mu_{r}$.
  Since the collection of all equivalence classes of~$\sim_r$, as~$r$ varies in~$\R_{>0}$, generates the topology of~$X$, Carath{\'e}odory's theorem implies that there is a unique Borel probability measure~$\mu$ on~$X$ such that, for every~$r$ in~$\R_{>0}$, its projection to~$X_r$ equals~$\mu_r$.
  By construction, ${\supp(\mu) = X}$.

  If~$g$ is an isometry of~$X$, then, for every~$r$ in~$\R_{>0}$, the map~$g$ descends to a bijection of~$X_r$ and thus leaves~$\mu_r$ invariant.
  From the definition of~$\mu$, it follows that ${g_* \mu = \mu}$.
  In particular, $\mu$ is invariant by~$G$.
  To prove uniqueness, let~$\mu'$ be a Borel probability measure on~$X$ that is invariant under~$G$.
  Then, for each~$r$ in~$\R_{>0}$, the measure~$\mu'$ projects to a measure~$\mu_r'$ on~$X_r$ that is invariant under the induced action of~$G$.
  Since this action is transitive, ${\mu_r' = \mu_r}$.
  Since this holds for every~$r$ in~$\R_{>0}$, from the definition of~$\mu$ it follows that ${\mu' = \mu}$.
\end{proof}

\begin{proof}[Proof of Proposition~\ref{p:trace-zero-spheres}]
  To prove item~$(i)$, let~$\Delta$ be in~$\pfd$.
  The proof relies on an embedding of~$\Qp(\sqrt{\Delta})$ into~$\Bss$, see \cite[\emph{Chapitre}~II, \emph{Corollaire}~1.9]{Vig80}.
  Let~$\varphi_0$ be the image of~$\frac{\Delta + \sqrt{\Delta}}{2}$ in~$\Bss$.
  Then, ${\tr(\varphi_0) = \Delta}$ and ${\nr(\varphi_0) = \frac{\Delta^2 - \Delta}{4}}$, so~$\tr(\varphi_0)$ and~$\nr(\varphi_0)$ belong to~$\Zp$ and thus~$\varphi_0$ belongs to~$\Rss$.
  Moreover, putting ${\delta \= 2 \varphi_0 - \Delta}$, the equalities ${\tr(\delta) = 0}$ and ${\nr(\delta) = - \Delta}$ hold and thus~$\delta$ belongs to~$S_{-\Delta}^0(\ss)$.
  This proves ${S_{-\Delta}^0(\ss) \neq \emptyset}$ and thus ${\bfL_{\ss, \pfd} \neq \emptyset}$.
  The compactness of~$\bfL_{\ss, \pfd}$ follows from that of~$\Lss$ and~$\pfd$ and the continuity of the reduced norm.
  For the last assertion of item~$(i)$, ${\varphi^2 = - \nr(\varphi)}$, so~$\varphi^2$
  belongs to~$\pfd$, and thus~$\Qp(\varphi)$ is isomorphic to~$\Qpfd$.

  To prove item~$(ii)$, recall ${S_{-\Delta}^0(\ss) \neq \emptyset}$.
  Since~$\Lss$ is compact and the reduced norm is continuous, $S_{-\Delta}^0(\ss)$ is compact.
  Since the action of each element of~$\Gss$ on~$S_{-\Delta}^0(\ss)$ is the composition of a left and a right multiplication, it is an isometry.
  It remains to prove that~$\Gss$ acts transitively on~$S_{- \Delta}^0(\ss)$.
  Let~$\varphi$ and~$\varphi'$ be in~$S_{-\Delta}^0(\ss)$.
  Since~$\varphi$ and~$\varphi'$ both satisfy the equation~$X^2 - \Delta = 0$, there is an isomorphism of $\Qp$-algebras between~$\Qp(\varphi)$ and~$\Qp(\varphi')$ mapping~$\varphi$ to~$\varphi'$.
  By Skolem--Noether's theorem, this isomorphism extends to an inner automorphism of~$\Bss$, see \cite[\emph{Chapitre}~I, \emph{Th{\'e}or{\`e}me}~2.1]{Vig80}.
  In particular, there exists~$g_0$ in~$\Bss^{\times}$ such that ${g_0 \varphi g^{-1}_0 = \varphi'}$.
  Denote by~$\varpi$ a uniformizer of~$\Bss$ and put ${g \= g_0 \varpi^{- \qvalss(g_0)}}$.
  Then, $g$ belongs to~$\Gss$ and satisfies ${g \varphi g^{-1} = \varphi'}$.

  Item~$(iii)$ is a direct consequence of item~$(ii)$ and Lemma~\ref{l:existence-invariant-measure}.

  For item~$(iv)$, the map ${\varphi \mapsto g \varphi g^{-1}}$ is an isomorphism of $\Zp$\nobreakdash-algebras that extends by~$\Qp$\nobreakdash-linearity to an isomorphism of~$\Qp$-algebras ${c \: \Bss\to \bfB_{\ss'}}$.
  Since the canonical involutions of~$\Bss$ and~$\bfB_{\ss'}$ are unique, for every~$\varphi$ in~$\Bss$, the equality ${c(\overline{\varphi}) = \overline{c(\varphi)}}$ holds.
  This implies that~$c$ preserves reduced traces and norms and that it is an isometry.
  In particular, $c$ maps~$S_{- \Delta}^0(\ss)$ to~$S_{- \Delta}^0(\ss')$ isometrically.
  By item~$(iii)$, the image of~$\nu_{- \Delta}^{\ss}$ by~$c$ is a Borel probability measure on~$S_{- \Delta}^0(\ss')$ that is invariant under the action of~$\Gsspr$ and thus equals~$\nu_{- \Delta}^{\sspr}$.
\end{proof}

\subsection{Parametrizing fixed points}
\label{ss:CM-parametrization}
This section proves the following proposition that provides a natural parametrization of the fixed points associated with the elements of a given zero-trace sphere.

\begin{proposition}
  \label{p:CM-parametrization}
  For all~$\ss$ in~$\tSups$, fundamental $p$\nobreakdash-adic discriminant~$\pfd$, and~$\Delta$ in~$\pfd$, the following properties hold.
  \begin{enumerate}
  \item [$(i)$]
    If~$\Qpfd$ is unramified over~$\Qp$, then there is a continuous function
    \begin{displaymath}
      x_{\ss, \Delta} \: S_{- \Delta}^0(\ss) \to \hDss
    \end{displaymath}
    such that, for every~$\varphi$ in~$S_{- \Delta}^0(\ss)$, the equality ${\Fixss(\varphi) = \{ x_{\ss, \Delta}(\varphi) \}}$ holds.
  \item [$(ii)$]
    If~$\Qpfd$ is ramified over~$\Qp$, then there are continuous functions
    \begin{displaymath}
      x_{\ss, \Delta}^+, x_{\ss, \Delta}^- \: S_{- \Delta}^0(\ss) \to \hDss
    \end{displaymath}
    such that, for every~$\varphi$ in~$S_{- \Delta}^0(\ss)$,
    \begin{displaymath}
      x_{\ss, \Delta}^+(\varphi) \neq x_{\ss, \Delta}^-(\varphi)
      \text{ and }
      \Fixss(\varphi) = \{ x_{\ss, \Delta}^+(\varphi), x_{\ss, \Delta}^-(\varphi) \}.
    \end{displaymath}
  \end{enumerate}
\end{proposition}

The proof of this proposition relies on the following lemma.

\begin{lemma}
  \label{l:local-CM-parametrization}
  Let~$\ss$ be in~$\tSups$, $\pfd$ a fundamental $p$\nobreakdash-adic discriminant, $\Delta$ in~$\pfd$, and~$\varpi$ a uniformizer of~$\Rss$.
  For each~$\varphi_0$ in~$\bfL_{\ss, \pfd}$, put
  \begin{displaymath}
    \cC(\varphi_0)
    \=
    \{ \varphi \in S_{- \Delta}^0(\ss) \: \varphi \varphi_0^{-1} \in \bfone_{\Bss} + \varpi^3 \Rss \}.
  \end{displaymath}
  Then, there is a continuous function ${g \: \cC(\varphi_0) \to \Gss}$ such that, for every~$\varphi$ in~$\cC(\varphi_0)$, the equality ${g(\varphi) \varphi_0 g(\varphi)^{-1} = \varphi}$ holds.
\end{lemma}

\begin{proof}
  For each~$\varphi$ in~$\cC(\varphi_0)$, put ${\varsigma(\varphi) \= \varphi \varphi_0^{-1}}$.
  Then,
  \begin{equation}
    \label{eq:60}
    \varsigma(\varphi)
    \in
    \bfone_{\Bss} + \varpi^3 \Rss
    \text{ and }
    \bfone_{\Bss} + \overline{\varsigma(\varphi)}
    \neq
    0.
  \end{equation}
  The function ${g \: \cC(\varphi_0) \to \Gss}$ defined by ${g(\varphi) \= 2(\bfone_{\Bss} + \overline{\varsigma(\varphi)})^{-1}}$ is continuous.
  Moreover, the equality ${\varphi_0 \overline{\varsigma(\varphi)} = \varsigma(\varphi) \varphi_0}$ yields
  \begin{equation}
    \label{eq:61}
    g(\varphi) \varphi_0 g(\varphi)^{-1}
    =
    (\bfone_{\Bss} + \overline{\varsigma(\varphi)})^{-1} \varphi_0 (\bfone_{\Bss} + \overline{\varsigma(\varphi)})
    =
    (\bfone_{\Bss} + \overline{\varsigma(\varphi)})^{-1} (\bfone_{\Bss} + \varsigma(\varphi)) \varphi_0
    =
    \varphi.
    \qedhere
  \end{equation}
\end{proof}

\begin{proof}[Proof of Proposition~\ref{p:CM-parametrization}]
  The proof uses several times that, for each~$\varphi_0$ in~$S_{-\Delta}^0(\ss)$, the subset~$\cC(\varphi_0)$ of~$S_{-\Delta}^0(\ss)$ given by Lemma~\ref{l:local-CM-parametrization} is an open and closed.

  Suppose~$\Qpfd$ is unramified over~$\Qp$ and let~$x_{\ss, \Delta} \: S_{- \Delta}^0(\ss) \to \hDss$ be the function mapping each~$\varphi$ in~$S_{-\Delta}^0(\ss)$ to the unique element of~$\Fixss(\varphi)$ (Lemma~\ref{l:fixed-properties}$(ii)$).
  Let~$\varphi_0$ be in~$S_{- \Delta}^0(\ss)$, $x_0$ the unique element of~$\Fixss(\varphi_0)$, and ${g \: \cC(\varphi_0) \to \Gss}$ the continuous function given by Lemma~\ref{l:local-CM-parametrization}.
  Then, for each~$\varphi$ in~$\cC(\varphi_0)$, the point~$g(\varphi) \cdot x_0$ belongs to~$\Fixss(\varphi)$ and thus ${g(\varphi) \cdot x_0 = x_{\ss, \Delta}(\varphi)}$.
  In particular, the restriction of~$x_{\ss, \Delta}$ to~$\cC(\varphi_0)$ is continuous by Lemma~\ref{l:action-regularity}$(ii)$.

  Suppose~$\Qpfd$ is ramified over~$\Qp$ and denote by ${\iota \: S_{- \Delta}^0(\ss) \to S_{- \Delta}^0(\ss)}$ the involution defined by ${\iota(\varphi) \= - \varphi}$.
  For each~$\varphi_0$ in~$S_{- \Delta}^0(\ss)$, the set~$\cC(\varphi_0)$ does not contain~$- \varphi_0$.
  Moreover, for every~$\varphi_0'$ in~$S_{- \Delta}^0(\ss)$, the set~$\cC(\varphi_0')$ is either disjoint from or equal to~$\cC(\varphi_0)$.
  Since~$S_{- \Delta}^0(\ss)$ is compact, it follows that there is a finite subset~$\Phi$ of~$S_{- \Delta}^0(\ss)$ such that
  \begin{equation}
    \label{eq:62}
    \{ \cC(\varphi_0), \cC(- \varphi_0) \: \varphi_0 \in \Phi \}
  \end{equation}
  is a partition of~$S_{-\Delta}^0(\ss)$.
  For each element~$\varphi_0$ of~$\Phi$, the set~$\Fixss(\varphi_0)$ has precisely two elements by Lemma~\ref{l:fixed-properties}$(ii)$.
  Denote them by~$x_{\varphi_0}^+$ and~$x_{\varphi_0}^-$.
  Moreover, denote by~$g_{\varphi_0}$ the continuous function given by Lemma~\ref{l:local-CM-parametrization}.
  Using that~\eqref{eq:62} is a partition of~$S_{- \Delta}^0(\ss)$ by Lemma~\ref{l:action-regularity}$(ii)$ and that, for each~$\varphi_0$ in~$S_{- \Delta}^0(\ss)$, the equality ${\iota(\cC(\varphi_0)) = \cC(- \varphi_0)}$ holds, it follows that there are continuous functions ${x^+, x^- \: S_{- \Delta}^0(\ss) \to \hDss}$ such that, for each~$\varphi_0$ in~$\Phi$,
  \begin{displaymath}
    x^{\pm}|_{\cC(\varphi_0)}(\varphi)
    =
    g_{\varphi_0}(\varphi) \cdot x_{\varphi_0}^{\pm}
    \text{ and }
    x^{\pm}|_{\cC(-\varphi_0)}(\varphi)
    =
    g_{\varphi_0}(- \varphi) \cdot x_{\varphi_0}^{\pm}.
  \end{displaymath}
  Since for each~$\varphi$ in~$S_{- \Delta}^0(\ss)$, Lemma~\ref{l:fixed-properties}$(iv)$ implies ${\Fixss(\varphi) = \Fixss(-\varphi)}$, the points~$x^+(\varphi)$ and~$x^-(\varphi)$ belong to~$\Fixss(\varphi)$.
  Thus, to prove item~$(ii)$ with~${x_{\ss, \Delta}^+ = x^+}$ and~${x_{\ss, \Delta}^- = x^-}$, it suffices to prove that, for all~$\varphi_0$ in~$\Phi$ and~$\varphi$ in ${\cC(\varphi_0) \cup \cC(- \varphi_0)}$, the points~$x^+(\varphi)$ and~$x^-(\varphi)$ are different.
  Either
  \begin{displaymath}
    x^{\pm}(\varphi)
    =
    g_{\varphi_0}(\varphi) \cdot x_{\varphi_0}^{\pm}
    \text{ or }
    x^{\pm}(\varphi)
    =
    g_{\varphi_0}(- \varphi) \cdot x_{\varphi_0}^{\pm}.
  \end{displaymath}
  In both cases, ${x^+(\varphi) \neq x^-(\varphi)}$.
\end{proof}

\subsection{From zero-trace spheres to \CM{} points}
\label{ss:trace-zero-spheres-to-CM}
This section proves the following proposition relating zero-trace spheres to formal \CM{} points, and defines a natural measure on the set of formal \CM{} points of a given fundamental $p$\nobreakdash-adic discriminant and residue disc.

For all~$\ss$ in~$\tSups$ and subset~$S$ of~$\Lfss$, put
\begin{displaymath}
  \Fixss(S)
  \=
  \bigcup_{g \in S}\Fixss(g).
\end{displaymath}
The \emph{trace} of a function~$\hF \: \hDss \to \R$, is the map defined by
\begin{displaymath}
  \begin{array}{rrcl}
    \Trss(\hF) \: & \Lfss & \to & \R
    \\
                  & g & \mapsto & \frac{1}{\# \Fixss(g)} \sum_{x \in \Fixss(g)} \hF(x).
  \end{array}
\end{displaymath}
For all fundamental $p$\nobreakdash-adic discriminant~$\pfd$ such that~$\Qpfd$ is unramified (resp. ramified) over~$\Qp$ and~$\Delta$ in~$\pfd$, let~$x_{\ss, \Delta}$ (resp. ${x_{\ss, \Delta}^+}$, ${x_{\ss, \Delta}^-}$) be given by Proposition~\ref{p:CM-parametrization}.

\begin{proposition}
  \label{p:trace-zero-spheres-to-CM}
  For all~$\ss$ in~$\tSups$, fundamental $p$\nobreakdash-adic discriminant~$\pfd$, and~$\Delta$ in~$\pfd$, the following properties hold.
  \begin{enumerate}
  \item [$(i)$]
    The equality,
    \begin{equation}
      \label{eq:63}
      \Piss^{-1} \left( \Lambda_{\pfd} \cap \Dss \right)
      =
      \Fixss \left( S_{- \Delta}^0(\ss) \right)
    \end{equation}
    holds and these sets are compact.
  \item[$(ii)$]
    The measure~$\hnu_{\pfd}^{\ss}$ on~$\hDss$, defined by
    \begin{equation}
      \label{eq:64}
      \hnu_{\pfd}^{\ss}
      \=
      \begin{cases}
        (x_{\ss, \Delta})_* \nu_{- \Delta}^{\ss}
        & \text{if~$\Qpfd$ is unramified over~$\Qp$};
        \\
        \frac{1}{2} \left( (x_{\ss, \Delta}^+)_* \nu_{- \Delta}^{\ss} + (x_{\ss, \Delta}^-)_* \nu_{- \Delta}^{\ss} \right)
        & \text{if~$\Qpfd$ is ramified over~$\Qp$},
      \end{cases}
    \end{equation}
    depends on~$\ss$ and~$\pfd$ but not on~$\Delta$.
    Moreover, ${\supp(\hnu^{\ss}_{\pfd}) = \Piss^{-1} \left( \Lambda_{\pfd} \cap \Dss \right)}$ and, for every continuous ${\hF \: \hDss\to \R}$,
    \begin{equation}
      \label{eq:65}
      \int \hF \dd \hnu^{\ss}_{\pfd}
      =
      \int \Trss(\hF) \dd \nu^{\ss}_{- \Delta}.
    \end{equation}
  \end{enumerate}
\end{proposition}

\begin{proof}
  For item~$(i)$, $S_{- \Delta}^0(\ss)$ is compact by Proposition~\ref{p:trace-zero-spheres}$(ii)$.
  Since
  \begin{displaymath}
    \Fixss(S_{- \Delta}^0(\ss))
    =
    \begin{cases}
      x_{\ss, \Delta}(S_{- \Delta}^0(\ss))
      & \text{if~$\Qpfd$ is unramified over~$\Qp$};
      \\
      x_{\ss, \Delta}^+(S_{- \Delta}^0(\ss)) \cup x_{\ss, \Delta}^-(S_{- \Delta}^0(\ss))
      & \text{if~$\Qpfd$ is ramified over~$\Qp$}
    \end{cases}
  \end{displaymath}
  by Proposition~\ref{p:CM-parametrization}, it follows that~$\Fixss(S_{- \Delta}^0(\ss))$ is also compact.

  To prove that the left side of~\eqref{eq:63} is contained in the right side, let~$x$ be in ${\Piss^{-1}(\Lambda_{\pfd} \cap \Dss)}$ and put~$E \= \Piss(x)$.
  Then, $E$ is a formal \CM{} point, so it belongs to~$\SupsQpalg$ and~$x$ belongs to~$\Xss(\OQpalg)$.
  Let ${\alpha \: \tFE \to \Fss}$ be an isomorphism such that~$(\FE, \alpha)$ represents~$x$ and the ring homomorphism, defined by
  \begin{equation}
    \label{eq:66}
    \begin{array}{rrcl}
      \iota \: & \End(\FE) & \to & \Rss
      \\
               & \varphi & \mapsto & \alpha \circ \tvarphi \circ \alpha^{-1}.
    \end{array}
  \end{equation}
  Since~$\End(\FE)$ is isomorphic to~$\cO_{\Qpfd}$ and~$\cO_{\Qpfd} = \Zp \left[ \tfrac{\Delta + \sqrt{\Delta}}{2} \right]$ by~\eqref{eq:123} in Lemma~\ref{l:quadratic-extensions}$(ii)$, there is an element~$\varphi$ of~$\Zp + 2 \End(\FE)$ satisfying the equation~$X^2 - \Delta = 0$.
  Then, $\iota(\varphi)$ belongs to~$\Zp + 2 \Rss$, satisfies the equation~$X^2 - \Delta = 0$, and thus belongs to~$S_{- \Delta}^0(\ss)$.
  Lemma~\ref{l:unit-function} implies
  \begin{equation}
    \label{eq:67}
    \iota(\Aut(\FE))
    =
    \cO_{\Qp(\iota(\varphi))}^{\times}
    =
    \Zp[\Uss(\iota(\varphi))]^{\times}
    \text{ and }
    \Uss(\iota(\varphi))
    \in
    \iota(\Aut(\FE)).
  \end{equation}
  By Lemma~\ref{l:fixed}, this implies that~$x$ belongs to~$\Fixss(\iota(\varphi))$ and thus to the right side of~\eqref{eq:63}.

  To prove the reverse inclusion, recall that ${S_{- \Delta}^0(\ss) \neq \emptyset}$ by Proposition~\ref{p:trace-zero-spheres}$(ii)$, and let~$\varphi$ be in~$S_{- \Delta}^0(\ss)$ and~$x$ in~$\Fixss(\varphi)$.
  By Lemma~\ref{l:fixed-properties}$(i)$, the point~$x$ belongs to~$\Xss(\OQpalg)$.
  Put ${E \= \Piss(x)}$ and let ${\alpha \: \tFE \to \Fss}$ be an isomorphism such that~$(\FE, \alpha)$ represents~$x$.
  By Lemma~\ref{l:fixed}, the unit~$\Uss(\varphi)$ belongs to~$\iota(\End(\FE))$.
  It follows that ${\Zp[\Uss(\varphi)] \subseteq \iota(\End(\FE))}$ and thus that~$\varphi$ belongs to~$\iota(\End(\FE))$.
  This implies that~$\End(\FE)$ contains a solution of ${X^2 - \Delta = 0}$ and thus that it is a $p$\nobreakdash-adic quadratic order of $p$\nobreakdash-adic discriminant~$\pfd$.
  That is, ${E}$ belongs to~$\Lambda_{\pfd}$.
  Since~$E$ belongs to~$\Dss$ by definition, it follows that~$x$ belongs to the left side of~\eqref{eq:63}.

  The equality ${\supp(\hnu_{\pfd}^{\ss}) = \Piss^{-1}(\Lambda_{\pfd} \cap \Dss)}$ in item~$(ii)$ follows from item~$(i)$ because ${\supp(\nu_{- \Delta}^{\ss}) = S_{- \Delta}^0(\ss)}$ by Proposition~\ref{p:trace-zero-spheres}$(iii)$.
  To prove~\eqref{eq:65}, let ${\hF \: \hDss \to \R}$ be continuous.
  The change of variables formula implies
  \begin{displaymath}
    \int \hF \dd \hnu_{\pfd}^{\ss}
    =
    \int \hF \circ x_{\ss, \Delta} \dd \nu_{- \Delta}^{\ss}
    =
    \int \Trss(\hF) \dd \nu_{- \Delta}^{\ss}
  \end{displaymath}
  if~$\Qpfd$ is unramified over~$\Qp$ and, if~$\Qpfd$ is ramified over~$\Qp$, then
  \begin{displaymath}
    \int \hF \dd \hnu_{\pfd}^{\ss}
    =
    \int \frac{1}{2} (\hF \circ x_{\ss, \Delta}^+ + \hF \circ x_{\ss, \Delta}^-) \dd \nu_{- \Delta}^{\ss}
    =
    \int \Trss(\hF) \dd \nu_{- \Delta}^{\ss}.
  \end{displaymath}
  This proves~\eqref{eq:65} in all cases.
  It remains to prove that, for every~$\Delta'$ in~$\pfd$, the identity~\eqref{eq:65} holds with~$\Delta$ replaced by~$\Delta'$.
  Let~$u$ in~$\Zp^{\times}$ be such that~$\Delta' = u^2 \Delta$.
  Then, the left multiplication map ${\varphi \mapsto u\varphi}$ induces a bijective isometry ${S_{- \Delta}^0(\ss) \to S_{- \Delta'}^0(\ss)}$ and thus it maps~$\nu_{- \Delta}^{\ss}$ to~$\nu_{- \Delta'}^{\ss}$ by Proposition~\ref{p:trace-zero-spheres}$(iii)$.
  Hence, by the change of variables formula and Lemma~\ref{l:fixed-properties}$(iv)$, for every continuous ${\hF \: \hDss \to \R}$,
  \begin{displaymath}
    \int \Trss(\hF)(\varphi) \dd \nu^{\ss}_{- \Delta'}(\varphi)
    =
    \int \Trss(\hF)(u\varphi) \dd \nu^{\ss}_{- \Delta}(\varphi)
    =
    \int \Trss(\hF)(\varphi) \dd \nu^{\ss}_{- \Delta}(\varphi).
    \qedhere
  \end{displaymath}
\end{proof}

\subsection{Equidistribution of CM{} points on a supersingular residue disc}
\label{ss:residual}

This section proves the following theorem.
\begin{theorem}
  \label{t:residual}
  For every~$\varepsilon$ in~$\R_{>0}$, there is~$C$ in~$\R_{>0}$ such that, for all $p$-supersingular fundamental discriminant~$d$, $f$ in~$\Z_{>0}$, and~$\ss$ in~$\tSups$,
  \begin{displaymath}
    \left| \frac{ \deg(\Lambda_{d f^2}|_{\Dss})}{\deg(\Lambda_{d f^2})}
      - \frac{24}{(p - 1) \# \Aut(\ss)} \right|
    \le
    C |d|^{-\frac{1}{28}+\varepsilon} (f |f|_p)^{-\frac{1}{2}+\varepsilon}.
  \end{displaymath}
\end{theorem}

To state a corollary, let~$\Rsups$ by the real vector space, defined by
\begin{displaymath}
  \Rsups
  \=
  \left\{ (z_{\ss})_{\ss \in \tSups} \: z_{\ss} \in \R \right\},
\end{displaymath}
and let~$\vsups$ be the vector in~$\Rsups$, defined by
\begin{equation}
  \label{eq:68}
  \vsups_{\ss}
  \=
  \frac{24}{(p - 1)\# \Aut(\ss)}.
\end{equation}
The mass formula~\eqref{eq:Deuring-Eichler} implies that~$\vsups$ is a probability vector.
For each divisor~$\Lambda$ on~$\Sups$, denote by~$v(\Lambda)$ the vector in~$\Rsups$ defined by~$v(\Lambda)_{\ss} \= \deg(\Lambda|_{\Dss})$.

The following corollary is a direct consequence of Theorem~\ref{t:residual}.

\begin{coro}
  \label{c:residual}
  Let~$(D_n)_{n = 1}^{\infty}$ be a sequence of $p$-supersingular discriminants such that
  \begin{displaymath}
    D_n|D_n|_p \to -\infty
    \text{ as }
    n \to \infty.
  \end{displaymath}
  Then,
  \begin{displaymath}
    \lim_{n \to \infty} \frac{v(\Lambda_{D_n})}{\deg(\Lambda_{D_n})}
    =
    \vsups.
  \end{displaymath}
\end{coro}

The hypothesis that~$D_n|D_n|_p \to -\infty$ as~$n \to \infty$ cannot be weakened to ${D_n \to -\infty}$ as ${n \to \infty}$, see Remark~\ref{r:tpssing} below.

When restricted to discriminants for which~$p$ is inert in the corresponding quadratic imaginary extension of~$\Q$, Theorem~\ref{t:residual} is a particular case of \cite[Theorem~1.1]{JetKan11} and of the ``sparse equidistribution'' result of Michel \cite[Theorem~3]{Mic04} in the case of fundamental discriminants.
When restricted to the case where the $p$-supersingular discriminants in~$(D_n)_{n = 1}^{\infty}$ share the same fundamental discriminant and, for every~$n$, the conductor of~$D_n$ is the ${n}$-th power of a prime number independent of~$n$, Corollary~\ref{c:residual} follows from \cite[Theorem~1.5]{Vat02}.
See also \cite[Section~3]{Cor02} for related work on the surjectivity on CM{} points of reduction maps with respect to inert primes.

The proof of Theorem~\ref{t:residual} for fundamental discriminants constructs an auxiliary modular form of weight~$\frac{3}{2}$ that is cuspidal and derives the desired estimates from Duke's bounds for Fourier coefficients~\cite{Duk88}.
This cuspidal modular also appears in the proof of \cite[Theorem~1.4]{ElkOnoYan05}.
The proof also relies on Siegel's classical estimate: For every~$\varepsilon$ in~$\R_{>0}$, there is~$C$ in~$\R_{>0}$ such that, for every fundamental discriminant~$d$,
\begin{equation}
  \label{eq:Siegel-low-general}
  \deg(\Lambda_{d})
  \ge
  C |d|^{\frac{1}{2}-\varepsilon},
\end{equation}
see, \emph{e.g.}, \cite{Sie35b} or~\cite{Gol74}.
As in~\cite{CloUll04}, Zhang's formula (Lemma~\ref{l:Zhang-general}) allows to pass from fundamental discriminants to the general case.

The proof of Theorem~\ref{t:residual} relies on several lemmas that are only needed in the case of non-fundamental discriminants.
Lemma~\ref{l:tpssing-residual} recalls the description in~\cite{Gro87} of the action of Hecke correspondences on supersingular residue discs in terms of the Brandt matrices and treats discriminants whose conductor is divisible by~$p$.
Lemma~\ref{l:supersingular-Hecke-residual} applies Deligne's bound to estimate the norm of eigenvalues of Brandt matrices.

To statement of the next lemma relies on the following notation.
In the rest of this section, consider vectors in~$\Rsups$ as column vectors.
For all~$m$ in~$\Z_{>0}$ and~$\ss$ and $\ss'$ in~$\tSups$, denote by~$B(m)_{\ss, \ss'}$ the number of subgroup schemes~$C$ of~$\ss$ of order~$m$ such that~$\ss / C$ is isomorphic to~$\ss'$.
Put ${B(m) \= (B(m)_{\ss, \ss'})_{\ss, \ss' \in \tSups}}$.
By \cite[Proposition~2.3]{Gro87}, $B(m)$ is the Brandt matrix of degree~$m$ defined by~\cite[(1.5)]{Gro87}.
Note that~$B(1)$ is the identity matrix.
Recall that the Frobenius map~$\Frob$ maps~$\tSups$ onto itself and it induces an involution on this set, see Section~\ref{ss:ss}.
It follows that the induced linear map ${\Frob_* \: \Rsups \to \Rsups}$, defined by ${\Frob_*(v)_{\ss} \= v_{\Frob(\ss)}}$, is also an involution.
Moreover, ${\Frob_*(\vsups) = \vsups}$ because, for every~$\ss$ in~$\tSups$ that does not have a representative elliptic curve defined over~$\Fp$, the equality ${\# \Aut(\ss) = 2}$ holds.

\begin{lemma}
  \label{l:tpssing-residual}
  \hfill

  \begin{enumerate}
  \item[$(i)$]
    For all~$m$ in ${\Z_{>0} \ssetminus p\Z}$ and divisor~$\Lambda$ on~$\Sups$,
    \begin{equation}
      \label{eq:69}
      v(T_m(\Lambda))
      =
      B(m)^{\intercal} v(\Lambda).
    \end{equation}
  \item[$(ii)$]
    The equality~$\Frob_* = B(p)^{\intercal}$ holds, as linear endomorphisms of~$\Rsups$.
    Moreover, for all~$r$ in~$\Z_{>0}$ and divisor~$\Lambda$ on~$\Sups$,
    \begin{equation}
      \label{eq:70}
      v(T_{p^r}(\Lambda))
      =
      \sigma_1(p^r)\cdot \Frob_*^r (v(\Lambda)).
    \end{equation}
  \item[$(iii)$]
    For all $p$-supersingular discriminant~$D$ and~$r$ in~$\Z_{>0}$,
    \begin{equation}
      \label{eq:71}
      \frac{v(\Lambda_{D p^{2r}})}{\deg(\Lambda_{D p^{2r}})}
      =
      \frac{v(\Lambda_D)}{\deg(\Lambda_D)}.
    \end{equation}
  \end{enumerate}
\end{lemma}

\begin{proof}
  By continuity, to prove item~$(i)$ it suffices to consider the case where the divisor~$\Lambda$ is supported on~$\SupsQpalg$, see, \emph{e.g.}, \cite[Lemma~2.1]{HerMenRiv20}.
  In this case, the desired assertion follows from the fact that, for all~$E$ in~$\SupsQpalg$ and~$m$ in ${\Z_{>0} \ssetminus p\Z}$, the reduction map induces a bijection from the set of subgroups of~$E$ of order~$m$ to the set of subgroup schemes of order~$m$ of~$\tE$, see, \emph{e.g.}, \cite[Chapter~III, Corollary~6.4(b) and Chapter~VII, Proposition~3.1(b)]{Sil09}.

  The first assertion of item~$(ii)$ follows from the fact that each~$\ss$ in~$\tSups$ has a unique subgroup scheme of order~$p$ and that this subgroup scheme is the kernel of the Frobenius map from~$\ss$ to~$\Frob(\ss)$.
  The proof of~\eqref{eq:71} for ${r = 1}$ relies on the fact that the reduction~$\tPhi_p(X, Y)$ of the modular polynomial~$\Phi_p$ modulo~$p$ is given by ${\tPhi_p(X, Y) = (X - Y^p)(X^p - Y)}$, see, \emph{e.g.}, \cite[Chapter~5, Section~2, pp.~57-58]{Lan87}.
  Together with~\eqref{eq:deg-Hecke-op-divisor} with $n=p$, \eqref{eq:algebraic definition} with~$q = p$, and the definition of~$\Frob_*$, this implies~\eqref{eq:71} for~$r = 1$.
  The case where ${r \ge 2}$ follows by induction using the multiplicative property of Hecke correspondences~\eqref{eq:Hecke-Tpr} and the fact that~$\Frob$ induces an involution on~$\tSups$.

  To prove item~$(iii)$, denote by~$d$ and~$f$ the fundamental discriminant and the conductor of~$D$, respectively, so ${D = d f^2}$.
  Put
  \begin{displaymath}
    r_0 \= \ord_p(f),
    f_0 \= p^{- r_0} f,
    \text{ and }
    D_0 \= d f_0^2.
  \end{displaymath}
  By item~$(ii)$, \eqref{eq:p-Zhang}, and~\eqref{eq:pm-Zhang} in Lemma~\ref{l:Zhang-general}, to prove item~$(iii)$ it suffices to prove ${\Frob_*(v(\Lambda_{D_0})) = v(\Lambda_{D_0})}$.
  From~\eqref{eq:multiplicativity-Hecke-operators} and items~$(i)$ and~$(ii)$, for each~$m$ in ${\Z_{>0} \ssetminus p \Z}$, the maps~$\Frob_*$ and~$B(m)^{\intercal}$ commute.
  Thus, ${\Frob_*(v(\Lambda_{D_0})) = v(\Lambda_{D_0})}$ follows from~\eqref{eq:general-inverse-Zhang-formula} in Lemma~\ref{l:Zhang-general} with~$\wtf = 1$ and ${\Frob_*(v(\Lambda_d)) = v(\Lambda_d)}$, whose proof follows.
  Fix~$\ss$ in~$\tSups$.
  The endomorphism rings~$\End(\ss)$ and~$\End(\Frob(\ss))$ are isomorphic.
  Applying~\eqref{eq:conteo} first to~$\ss$ and then to~$\Frob(\ss)$ yields
  \begin{equation}
    \label{eq:72}
    v(\Lambda_d)_{\ss}
    =
    \deg(\Lambda_d|_{\Dss})
    =
    \epsilon_d h(d, e)
    =
    \epsilon_d h(d, \Frob(\ss))
    =
    \deg(\Lambda_d|_{\bfD_{\Frob(\ss)}})
    =
    \Frob_*(v(\Lambda_d))_{\ss}.
    \qedhere
  \end{equation}
\end{proof}

\begin{remark}
  \label{r:tpssing}
  For every $p$-supersingular discriminant~$D$, the sequence of vectors $\left( \frac{v(\Lambda_{D p^{2r}})}{\deg(\Lambda_{D p^{2r}})} \right)_{r = 1}^{\infty}$ is constant by Lemma~\ref{l:tpssing-residual}.
  Thus, this sequence cannot converge to~$\vsups$ unless ${\frac{v(\Lambda_D)}{\deg(\Lambda_D)} = \vsups}$.
  This proves that in Corollary~\ref{c:residual} it is insufficient to suppose that ${D_n \to -\infty}$ as ${n \to \infty}$.
\end{remark}

The statement of the next lemma relies on the following notation.
Endow~$\Rsups$ with the scalar product~$\langle \cdot, \cdot \rangle_{\sups}$ and norm~$\| \cdot \|_{\sups}$, defined by
\begin{equation}
  \label{eq:inner-product}
  \langle v, v' \rangle_{\sups}
  \=
  \sum_{\ss \in \tSups} \frac{v_\ss v_\ss'}{\vsups_{\ss}}
  \text{ and }
  \| v \|_{\sups}
  \=
  \sqrt{\langle v, v \rangle_{\sups}}.
\end{equation}

\begin{lemma}
  \label{l:supersingular-Hecke-residual}
  There is an orthonormal basis~$\cB$ of~$\Rsups$ containing the vector~$\vsups$ and such that, for every~$m$ in~$\Z_{>0}$, each vector in~$\cB$ is an eigenvector of~$B(m)^{\intercal}$.
  For each~$v$ in~$\cB$, let ${\lambda_v \: \Z_{>0} \rightarrow \C}$ be defined by ${B(m)v = \lambda_v(m) v}$.
  Then, the following properties hold.
  \begin{enumerate}
  \item[$(i)$]
    For all~$m$ in ${\Z_{>0} \ssetminus p\Z}$ and~$r$ in~$\Z_{\ge 0}$, the equality ${\lambda_{\vsups}(p^rm) = \sigma_1(m)}$ holds.
  \item[$(ii)$]
    For every~$\varepsilon$ in~$\R_{>0}$, there is~$C_1$ in~$\R_{>0}$ such that, for all~$v$ in~$\cB$ different from~$\vsups$ and~$m$ in~$\Z_{>0}$,
    \begin{displaymath}
      |\lambda_v(m)|
      \le
      C_1 m^{\frac{1}{2} + \varepsilon}.
    \end{displaymath}
  \end{enumerate}
\end{lemma}

\begin{proof}
  The proof begins recalling some facts about the space~$M_2(\Gamma_0(p))$ of holomorphic modular forms of weight~2 for~$\Gamma_0(p)$.
  This space contains the Eisenstein series~$F_p(\tau)$, defined by
  \begin{equation}
    \label{eq:73}
    F_p(\tau)
    \=
    \frac{p-1}{24}+\sum_{r=0}^{\infty}\sum_{m \in \Z_{>0}, p \nmid m}\sigma_1(m)\exp(2\pi i mp^r \tau),
  \end{equation}
  see~\cite[(5.7)]{Gro87}.
  The subspace of cuspidal modular forms~$S_2(\Gamma_0(p))$ has codimension one in~$M_2(\Gamma_0(p))$, so~$M_2(\Gamma_0(p)) = \C F_p\oplus S_2(\Gamma_0(p))$, see, \emph{e.g.}, \cite[Theorems~2.5.2 and~4.2.7]{Miy89}.
  Since the constant coefficient of~$F_p$ is nonzero, it follows that every modular form in~$M_2(\Gamma_0(p))$ whose constant coefficient is zero is cuspidal.

  For the first assertion and item~$(i)$, ${\| \vsups\|_{\sups} = 1}$ and, for all~$m$ in ${\Z_{>0} \ssetminus p \Z}$ and~$r$ in~$\Z_{\ge 0}$, the equality ${B(p^r m)^{\intercal} \vsups = \sigma_1(m) \cdot \vsups}$ holds, see \cite[Proposition~2.7(1, 6)]{Gro86}.
  Moreover, for every~$m$ in~$\Z_{>0}$, the matrix~$B(m)^{\intercal}$ is self-adjoint with respect to the inner product~\eqref{eq:inner-product} and, for every~$m'$ in~$\Z_{>0}$, the matrices~$B(m)^{\intercal}$ and~$B(m')^{\intercal}$ commute, see \cite[Proposition~2.7(5, 6)]{Gro86}.
  It follows that there is an orthonormal basis~$\cB$ of~$\Rsups$ containing~$\vsups$ and such that, for every~$m$ in~$\Z_{>0}$, each vector in~$\cB$ is an eigenvector of~$B(m)^{\intercal}$.
  This proves the first assertion and item~$(i)$.

  To prove~$(ii)$, define, for all~$v$ and $v'$ in~$\Rsups$, the function of~$\tau$ in~$\H$,
  \begin{displaymath}
    \phi(v, v')(\tau)
    \=
    \frac{p-1}{24} \langle v, \vsups \rangle_{\sups} \langle v', \vsups \rangle_{\sups}
    +\sum_{m = 1}^{\infty} \langle B(m)^{\intercal}v, v' \rangle_{\sups} \exp(2 \pi i m \tau).
  \end{displaymath}
  This function belongs to~$M_2(\Gamma_0(p))$ by \cite[Propositions~4.4 and~5.6]{Gro87}.
  In particular, for each~$v$ in~$\cB$ different from~$\vsups$, the modular form~$\phi(v, v)$ has Fourier expansion
  \begin{displaymath}
    \phi(v, v)(\tau)
    =
    \sum_{m = 1}^{\infty} \lambda_v(m) \exp(2 \pi i m \tau).
  \end{displaymath}
  Since the constant term of~$\phi(v, v)$ is zero, $\phi(v, v)$ is cuspidal and item~$(ii)$ follows from~\eqref{eq:7} and Deligne's bound \cite[\emph{Th{\'e}or{\`e}me}~8.2]{Del74}.
\end{proof}

For all fundamental discriminant~$d$ and integer~$f$ satisfying ${f \ge 2}$,
\begin{equation}
  \label{eq:74}
  \deg(\Lambda_{df^2})
  =
  \frac{\deg(\Lambda_d)}{w_{d, 1}} \left( R_d^{-1} * \sigma_1 \right)(f),
\end{equation}
by~\eqref{eq:deg-Hecke-op-divisor} and~\eqref{eq:general-inverse-Zhang-formula} in Lemma~\ref{l:Zhang-general} with ${\wtf = 1}$.

\begin{lemma}
  \label{l:R-inverse}
  For every~$\varepsilon$ in~$\R_{>0}$, there is~$C$ in~$\R_{>0}$ such that, for all~$m$ in~$\Z_{>0}$ and fundamental discriminant~$d$,
  \begin{displaymath}
    \left| R_d^{-1}(m) \right|
    \le
    C m^\varepsilon
    \text{ and }
    \left( R_d^{-1}\ast \sigma_1 \right) (m)
    \ge
    C^{-1} m^{1-\varepsilon}.
  \end{displaymath}
\end{lemma}

\begin{proof}
  Recall that ${\psi_d \: \Z_{>0} \to \{-1, 0, 1 \}}$ is the arithmetic function given by the Kronecker symbol~$\left( \frac{d}{\cdot} \right)$.
  Denote by~$\mu$ the M{\"o}bius function, so ${R_d^{-1} = \mu \ast (\mu \cdot \psi_d)}$.
  For every prime number~$q$,
  \begin{displaymath}
    R_d^{-1}(q^s)
    =
    \begin{cases}
      1
      & \text{if } s=0;
      \\
      -1 - \psi_d(q)
      & \text{if } s = 1;
      \\
      \psi_d(q)
      & \text{if } s=2;
      \\
      0 & \text{if } s\ge 3.
    \end{cases}
  \end{displaymath}
  This implies that, for every~$m$ in~$\Z_{>0}$, the inequality ${|R_d^{-1}(m)| \le d(m)}$ holds.
  So, the first inequality follows from~\eqref{eq:7}.

  To prove the second inequality, let~$N$ in~$\Z_{>0}$ be such that, for every integer~$q$ satisfying ${q \ge N}$, the inequality ${\frac{q - 1}{q} \ge q^{-\varepsilon}}$ holds and let~$C'$ in~$]0, 1[$ be such that, for every~$q$ in~$\{2, \ldots, N \}$, the inequality ${\frac{q - 1}{q} \ge C' q^{- \varepsilon}}$ holds.
  Since for every~$s$ in~$\Z_{>0}$,
  \begin{displaymath}
    (R_d^{-1} \ast \sigma_1)(q^s)
    =
    q^s - \psi_d(q) q^{s - 1}
    \ge
    q^{s - 1} (q - 1),
  \end{displaymath}
  it follows that, for every~$m$ in~$\Z_{>0}$,
  \begin{displaymath}
    \frac{(R_d^{-1} \ast \sigma_1)(m)}{m}
    \ge
    \prod_{q \mid m, \text{ prime}} \frac{q - 1}{q}
    \ge
    (C')^N m^{- \varepsilon}.
    \qedhere
  \end{displaymath}
\end{proof}

\begin{proof}[Proof of Theorem~\ref{t:residual}]
  Fix~$\varepsilon$ in~$\R_{>0}$ and put~$\varepsilon' \= \frac{\varepsilon}{3}$.
  Let~$C$ (resp.~$C_0$, $C_1$, $C_2$) be given by Siegel's estimate~\eqref{eq:Siegel-low-general} (resp.~\eqref{eq:7}, Lemma~\ref{l:supersingular-Hecke-residual}$(ii)$, Lemma~\ref{l:R-inverse}) with~$\varepsilon$ replaced by~$\varepsilon'$.
  For all fundamental discriminant~$\widehat{d}$ and~$\widehat{f}$ in~$\Z_{>0}$, let~$w_{\widehat{d}, \widehat{f}}$ be as in Section~\ref{ss:fixed-points-formula} and put ${u(\widehat{d} \widehat{f}^2) \= w_{\widehat{d}, \widehat{f}}}$.

  Suppose first ${f = 1}$, so ${D = d}$ and~$D$ is a fundamental discriminant.
  For each~$m$ in~$\Z_{>0}$, denote by~$H_p(m)$ the modified Hurwitz numbers defined by Gross~\cite[(1.8)]{Gro87} and, for each~$\ss$ in~$\tSups$, put
  \begin{displaymath}
    a_{\ss}(m)
    \=
    \frac{\# \Aut(\ss)}{2} \sum_{\substack{D' \text{ discriminant} \\ D' \mid m}} \frac{h(D', \ss)}{u(D')}
  \end{displaymath}
  if~$-m$ is a discriminant, and ${a_{\ss}(m) \= 0}$ otherwise.
  Then, the functions of~$\tau$ in~$\H$, defined by
  \begin{displaymath}
    \theta_{\ss}(\tau)
    \=
    1 + \sum_{m=1}^\infty a_{\ss}(m) \exp(2 \pi i m \tau)
    \text{ and }
    E_p(\tau)
    \=
    \frac{p-1}{12} + 2 \sum_{m=1}^\infty H_p(m) \exp(2 \pi i m \tau),
  \end{displaymath}
  are modular forms of weight~$\frac{3}{2}$ for~$\Gamma_0(4p)$ by \cite[(12.8), Proposition~12.9, and~(12.11)]{Gro87}.
  Moreover, \cite[(3.6), (3.13), and~(3.14)]{ElkOnoYan05} implies
  \begin{equation}
    \label{eq:75}
    \theta_{\ss}(\tau) - \frac{12}{p-1} E_p(\tau)
    =
    \sum_{m=1}^\infty \left( a_{\ss}(m) - \frac{24}{(p - 1)} H_p(m) \right) \exp(2 \pi i m \tau).
  \end{equation}
  In particular, the function on the left side of the equation above is a cuspidal modular form.
  Then, \eqref{eq:conteo} implies
  \begin{equation}
    \label{eq:76}
    a_{\ss}(|d|)
    =
    \frac{\# \Aut(\ss)}{2 u(d)} h(d, \ss)
    =
    \frac{\# \Aut(\ss)}{2\epsilon_d u(d)} \deg(\Lambda_d |_{\Dss})
  \end{equation}
  and \cite[(1.7) and~(1.8)]{Gro87},
  \begin{equation}
    \label{eq:77}
    H_p(|d|)
    =
    \frac{h(d)}{2\epsilon_d u(d)}
    =
    \frac{\deg(\Lambda_{d})}{2\epsilon_d u(d)}.
  \end{equation}
  Thus,
  \begin{equation*}
    \left| a_{\ss}(|d|) - \frac{24}{(p - 1)} H_p(|d|) \right|
    =
    \left(\frac{\#\Aut(\ss)}{2 \epsilon_d u(d)} \right) \deg(\Lambda_d) \left| \frac{ \deg(\Lambda_d|_{\Dss}) }{ \deg(\Lambda_d) } - \frac{24}{(p - 1) \# \Aut(\ss)} \right|.
  \end{equation*}
  Combined with Siegel's bound~\eqref{eq:Siegel-low-general} and Duke's bound \cite[Theorem~5]{Duk88} for the~$|d|$-th coefficient of the cuspidal modular form on the left side of~\eqref{eq:75}, this yields the desired estimate when ${f = 1}$.

  When ${f \ge 2}$, it suffices to consider the case where ${p \nmid f}$ by Lemma~\ref{l:tpssing-residual}$(iii)$.
  Then, ${w_{d,f} = 1}$ and~\eqref{eq:general-inverse-Zhang-formula} in Lemma~\ref{l:Zhang-general} with~$\wtf = 1$ and Lemma~\ref{l:tpssing-residual}$(i)$ imply
  \begin{equation}
    \label{eq:78}
    v(\Lambda_D)
    =
    \frac{1}{w_{d, 1}} \sum_{f_0 \in \Z_{>0}, f_0 \mid f} R_d^{-1} \left( \frac{f}{f_0} \right) B(f_0)^{\intercal} v(\Lambda_d).
  \end{equation}
  Writing~$v(\Lambda_d)$ as a linear combination of the elements in the base~$\cB$ yields
  \begin{displaymath}
    \frac{v(\Lambda_D)}{\deg(\Lambda_D)}
    =
    \frac{\deg(\Lambda_d)}{\deg(\Lambda_D)w_{d, 1}} \sum_{v \in \cB} \left( R_d^{-1} \ast \lambda_v \right) (f) \left\langle \frac{v(\Lambda_d)}{\deg(\Lambda_d)}, v \right\rangle_{\sups} v.
  \end{displaymath}
  Since
  \begin{displaymath}
    \langle v(\Lambda_d), \vsups \rangle_{\sups}
    =
    \sum_{\ss \in \tSups} \deg(\Lambda_d|_{\Dss})
    =
    \deg(\Lambda_d),
  \end{displaymath}
  Lemma~\ref{l:supersingular-Hecke-residual}$(i)$ and~\eqref{eq:74} yield
  \begin{equation}
    \label{eq:residual-deviation}
    \frac{v(\Lambda_D)}{\deg(\Lambda_D)} - \vsups
    =
    \sum_{v \in \cB, v \neq \vsups} \frac{\left( R_d^{-1} \ast \lambda_v \right) (f)}{\left( R_d^{-1} \ast \sigma_1 \right) (f)} \left\langle \frac{v(\Lambda_d)}{\deg(\Lambda_d)}-\vsups, v \right\rangle_{\sups} v.
  \end{equation}
  By the choice of~$C_0$, $C_1$, and~$C_2$, for every~$v$ in~$\cB$ different from~$\vsups$,
  \begin{displaymath}
    \left| \left(R_d^{-1} \ast \lambda_v \right)(f) \right|
    \le
    C_1 C_2 \sum_{f_0 \in \Z_{>0}, f_0 \mid f} \left(\frac{f}{f_0} \right)^{\varepsilon'} f_0^{\frac{1}{2} + {\varepsilon'} }
    \le
    C_0C_1C_2 f^{\frac{1}{2} + 2 \varepsilon'}.
  \end{displaymath}
  Combined with~\eqref{eq:residual-deviation} and the choice of~$C_2$, this implies
  \begin{displaymath}
    \left\| \frac{v(\Lambda_D)}{\deg(\Lambda_D)} - \vsups \right\|_{\sups}
    \le
    C_0C_1C_2^2 f^{- \frac{1}{2} + \varepsilon} \left\| \frac{v(\Lambda_d)}{\deg(\Lambda_d)}-\vsups \right\|_{\sups}.
  \end{displaymath}
  So, the desired estimate follows from the definition of~$\|\cdot \|_{\sups}$ and the case where ${f = 1}$ established above.
\end{proof}

\subsection{Proof of Theorem~\ref{t:CM-fundamental}}
\label{ss:proof-CM-fundamental}
Theorem~\ref{t:formal-CM-formulae} implies ${\Lambda_{\pfd} \subseteq \Sups}$.
That~$\Lambda_{\pfd}$ is compact then follows from Proposition~\ref{p:trace-zero-spheres-to-CM}$(i)$, the finiteness of~$\tSups$, and the continuity of~$\Piss$.

For each~$\ss$ in~$\tSups$, let~$\hnu_{\pfd}^{\ss}$ be as in Proposition~\ref{p:trace-zero-spheres-to-CM}$(ii)$ and put
\begin{equation}
  \label{eq:79}
  \nu_{\pfd}^{\ss}
  \=
  (\Piss)_* \hnu_{\pfd}^{\ss}.
\end{equation}
Since~$\Piss$ is continuous, ${\nu_{\pfd}^{\ss}}$ is a Borel probability measure on~$\Ell(\Cp)$ and ${\supp(\nu_{\pfd}^{\ss}) = \Lambda_{\pfd} \cap \Dss}$.
Put
\begin{equation}
  \label{eq:80}
  \nu_{\pfd}
  \=
  \frac{24}{p - 1} \sum_{\ss \in \tSups} \frac{1}{\# \Aut(\ss)} \nu_{\pfd}^{\ss}.
\end{equation}
Then, ${\nu_{\pfd}}$ is a Borel probability measure on~$\Ell(\Cp)$ by the mass formula~\eqref{eq:Deuring-Eichler} and ${\supp(\nu_{\pfd}) = \Lambda_{\pfd}}$.

The proof of~\eqref{eq:58} relies on the following consequence of Theorem~\ref{t:residual}.
Recall from~\eqref{eq:41} that, for all~$\ss$ in~$\tSups$ and~$m$ in~$\Z_{>0}$, the set~$V_m(\ss)$ is defined as~$\{ \phi \in L(\ss) \: \nr(\phi) = m \}$.

\begin{lemma}
  \label{l:residual-Siegel}
  Let~$\ss$ be in~$\tSups$ and~$\varepsilon$ in~$\R_{>0}$.
  Then, for every $p$-supersingular fundamental discriminant~$d$ such that~$|d|$ is sufficiently large,
  \begin{displaymath}
    \# V_{|d|}(\ss)
    \ge
    |d|^{\frac{1}{2}-\varepsilon}.
  \end{displaymath}
\end{lemma}

\begin{proof}
  Theorem~\ref{t:fixed-points-formula} and Lemma~\ref{l:fixed-properties}$(ii)$ imply, for every $p$-supersingular fundamental discriminant~$d$,
  \begin{displaymath}
    \# V_{|d|}(\ss)
    \ge
    \frac{1}{3} \deg(\Lambda_d|_{\Dss}).
  \end{displaymath}
  Together with Theorem~\ref{t:residual} and Siegel's estimate~\eqref{eq:Siegel-low-general}, this implies the desired assertion.
\end{proof}

The estimate~\eqref{eq:58} is a direct consequence of Theorem~\ref{t:residual}, the lemma above, and the following proposition.

\begin{proposition}
  \label{p:CM-fundamental}
  Let~$\ss$ be in~$\tSups$ and~$\pfd$ a fundamental $p$\nobreakdash-adic discriminant.
  Then, for all~$\varepsilon$ and~$\delta$ in~$\R_{>0}$, there is~$C$ in~$\R_{>0}$ such that the following property holds.
  For all function ${F \: \Lambda_{\pfd} \cap \Dss \to \R}$ that is constant on every ball of ${\Lambda_{\pfd} \cap \Dss}$ of radius~$\delta$ and fundamental discriminant~$d$ in~$\pfd$ for which ${V_{|d|}(\ss) \neq \emptyset}$,
  \begin{displaymath}
    \left| \int F \dd \odelta_{\Lambda_d|_{\Dss}}- \int F \dd \nu_{\pfd}^{\ss} \right|
    \leq
    C \left( \sup_{\Lambda_{\pfd} \cap \Dss} |F| \right) \frac{|d|^{\frac{13}{28}+\varepsilon}}{\# V_{|d|}(\ss)}.
  \end{displaymath}
\end{proposition}

\begin{proof}
  For each function ${F \: \Lambda_{\pfd} \cap \Dss \to \R}$, put ${\breve{F} \= \Trss(F \circ \Piss)}$.
  By Theorem~\ref{t:fixed-points-formula}, the definition of~$\Trss$, Proposition~\ref{p:trace-zero-spheres-to-CM}$(ii)$, and the change of variables formula, for every fundamental discriminant~$d$ in~$\pfd$ for which ${V_{|d|}(\ss) \neq \emptyset}$,
  \begin{equation}
    \label{eq:81}
    \int F \dd \odelta_{\Lambda_d|_{\Dss}}
    =
    \frac{1}{\# V_{|d|}(\ss)} \sum_{\phi \in V_{|d|}(\ss)} \breve{F}(\hphi)
    \text{ and }
    \int F \dd \nu_{\pfd}^{\ss}
    =
    \int \breve{F} \dd \nu_{|d|}^{\ss}.
  \end{equation}
  Since~$\bfL_{\ss, \pfd}$ is compact and~$\Piss$ is continuous, by Proposition~\ref{p:CM-parametrization} there is~$\breve{\delta}$ in~$\R_{>0}$ such that, if~$F$ is constant on every ball of ${\Lambda_{\pfd} \cap \Dss}$ of radius~$\delta$, then~$\breve{F}$ is constant on every ball of~$\bfL_{\ss, \pfd}$ of radius~$\breve{\delta}$.

  Given~$\varepsilon$ in~$\R_{>0}$, let~$C$ in~$\R_{>0}$ be given by Corollary~\ref{c:functional-deviation} with~$n = 3$, $\delta$ replaced by~$\breve{\delta}$, $Q = Q_{\ss}^0$, and~$S = 4$.
  Moreover, let~$d$ be a fundamental discriminant in~$\pfd$ for which ${V_{|d|}(\ss) \neq \emptyset}$ and let ${F \: \Lambda_{\pfd} \cap \Dss \to \R}$ be constant on every ball of radius~$\delta$.
  Then, by Proposition~\ref{p:trace-zero-spheres}$(ii)$ the hypotheses of Corollary~\ref{c:functional-deviation} are satisfied with~$\ell = |d|$, $m = |d|$, and with~$F$ replaced by~$\breve{F}$.
  The desired estimate is then a direct consequence of Corollary~\ref{c:functional-deviation} and~\eqref{eq:81}.
\end{proof}

\section{Equidistribution of partial Hecke orbits}
\label{s:Hecke-orbits}
This section proves the following quantitative version of Theorem~\ref{t:Hecke-orbits} in Section~\ref{ss:Hecke-orbits} and that distinct partial Hecke orbits have different limit measures (Proposition~\ref{p:orbit-measures}).
Recall from Section~\ref{ss:Hecke-orbits} that, for each~$E$ in~$\Sups$, the subgroup~$\NE$ of~$\Zp^{\times}$ equals~$(\Zp^{\times})^2$ if~$E$ is not a formal \CM{} point, and ${\left\{ \nr \left( \varphi \right) \: \varphi \in \Aut(\FE) \right\}}$ otherwise.

\begin{custtheo}{C'}
  \label{t:Hecke-orbits-pr}
  For all~$E$ in~$\Sups$ and coset~$\coset$ in~$\Qp^{\times}/ \NE$ contained in~$\Zp$, the closure~$\corbitc$ in~$\Sups$ of the partial Hecke orbit~$\corbit$ is compact and, for every~$E'$ in~$\corbitc$, the equality ${\NEpr = \NE}$ holds.
  Moreover, there is a Borel probability measure~$\mu_{\coset}^E$ on~$\Ell(\Cp)$ such that ${\supp(\mu_{\coset}^E) = \corbitc}$ and the following property holds.
  For all~$\varepsilon$ in~$\R_{>0}$ and locally constant function ${F \: \Sups \to \R}$, there is~$C$ in~$\R_{>0}$ such that, for all~$E'$ in~$\overline{\Orb_{\NE}(E)}$ and~$n$ in ${\coset \cap \Z_{>0}}$,
  \begin{equation}
    \label{eq:82}
    \left| \int F \dd \odelta_{T_{n}(E')}- \int F \dd \mu_{\coset}^{E} \right|
    \leq
    C n^{-\frac{1}{2}+\varepsilon}.
  \end{equation}
\end{custtheo}

The proof of Theorem~\ref{t:Hecke-orbits-pr} begins in Section~\ref{ss:supersingular-spheres} introducing, for all~$\ss$ and~$\ss'$ in~$\tSups$, ``supersingular spheres'' of the $p$\nobreakdash-adic space~$\Hom_{\Fpalg}(\Fss, \Fsspr)$ and proving that each of these sets carries a natural homogeneous measure (Proposition~\ref{p:supersingular-spheres}).
The proof proceeds by proving in Section~\ref{ss:supersingular-sphere-to-orbit} that each closure of a partial Hecke orbit restricted to a residue disc is the projection of a supersingular sphere by an evaluation map (Proposition~\ref{p:supersingular-sphere-to-orbit}).
Section~\ref{ss:Hecke-orbits-proof} proves Theorem~\ref{t:Hecke-orbits-pr} using the asymptotic distribution of integer points on $p$\nobreakdash-adic spheres in Section~\ref{s:p-adic-Linnik}.
The proof of Proposition~\ref{p:orbit-measures} occupies Section~\ref{ss:applications}.

Section~\ref{ss:symmetry-breaking} relies on the following corollary of Theorem~\ref{t:Hecke-orbits-pr}.
The statement relies on the action of Hecke correspondences on sets and measures, see Section~\ref{ss:Hecke-correspondences}.
For each~$E$ in~$\Sups$, denote by~$\cdot$ the multiplication in the quotient group~$\Qp^{\times} / \NE$.

\begin{coro}
  \label{c:Hecke-on-orbits}
  Let~$E$ be in~$\Sups$ and let~$\coset$ and~$\coset'$ cosets in ${\Qp^{\times} / \NE}$ contained in~$\Zp$.
  Then, for all~$E'$ in~$\overline{\Orb_{\NE}(E)}$ and~$n$ in ${\coset \cap \Z_{>0}}$,
  \begin{displaymath}
    T_n \left( \overline{\Orb_{\coset'}(E')} \right)
    =
    \overline{\Orb_{\coset \cdot \coset'}(E)}
    \text{ and }
    \frac{1}{\sigma_1(n)} (T_n)_* \mu_{\coset'}^{E'}
    =
    \mu_{\coset \cdot \coset'}^E.
  \end{displaymath}
\end{coro}

\begin{proof}
  Let~$(n_j)_{j = 1}^{\infty}$ be a sequence in ${\coset' \cap \Z_{>0}}$ tending to~$\infty$ such that, for every~$j$, the integer~$n_j$ is coprime to~$n$.
  The sequence~$(\odelta_{T_{n \cdot n_j}(E')})_{j = 1}^{\infty}$ converges to~$\mu_{\coset \cdot \coset'}^E$ as ${j \to \infty}$ by Theorem~\ref{t:Hecke-orbits-pr}.
  Since ${\NEpr = \NE}$ by Theorem~\ref{t:Hecke-orbits-pr} and, for every~$j$,
  \begin{displaymath}
    \odelta_{T_{n \cdot n_j}(E')}
    =
    \frac{1}{\sigma_1(n)} (T_n)_* \odelta_{T_{n_j}(E')}
  \end{displaymath}
  by~\eqref{eq:multiplicativity-Hecke-operators}, the sequence~$(\odelta_{T_{n \cdot n_j}(E')})_{j = 1}^{\infty}$
  converges to~$\frac{1}{\sigma_1(n)} (T_n)_* \mu_{\coset'}^{E'}$ as ${j \to \infty}$ by Theorem~\ref{t:Hecke-orbits-pr} with~$E$ replaced by~$E'$.
  This proves the equality of measures.
  The equality of sets follows by comparing the supports of these measures using Theorem~\ref{t:Hecke-orbits-pr} again.
\end{proof}

\subsection{Supersingular spheres and their homogeneous measures}
\label{ss:supersingular-spheres}
Throughout this section, fix~$\ss$ and~$\ss'$ in~$\tSups$.

The group~$\Hom(\ss, \ss')$ is a free~$\Z$\nobreakdash-module of rank~$4$.
For each an isogeny~$\phi$ in~$\Hom(\ss, \ss')$, denote by~$\overline{\phi}$ its dual isogeny in~$\Hom(\ss', \ss)$.
The ring~$\End(\ss)$ is a maximal order in the quaternion algebra $\End(\ss)\otimes \Q$ over~$\Q$ and the map $\End(\ss)\to \End(\ss)$ given by $\phi \mapsto \overline{\phi}$  extends by $\Q$-linearity to the canonical involution in $\End(\ss)\otimes \Q$.
The ring~$\End(\ss)$ has characteristic zero and the subring generated by the identity map~$\bfone_{\ss}$ on~$\ss$ equals the subset of endomorphisms~$\phi$ satisfying ${\overline{\phi}=\phi}$.
Identify this subring with~$\Z$.
Then, for every~$\phi$ in~$\Hom(\ss, \ss')$, the equality ${\overline{\phi} \phi = \deg(\phi)}$ holds.

The $\Z$-bilinear map
\begin{displaymath}
  \begin{array}{rrcl}
    \langle \ , \ \rangle \: & \Hom(\ss, \ss') \times \Hom(\ss, \ss') & \to & \End(\ss)
    \\ & (\phi_1,\phi_2) & \mapsto & \overline{\phi_1}\phi_2+\overline{\phi_2}{\phi_1}
  \end{array}
\end{displaymath}
takes values in $\Z$ and induces the quadratic form
\begin{displaymath}
  \begin{array}{rrcl}
    Q_{\ss, \ss'} \: & \Hom(\ss, \ss') & \to & \Z
    \\ & \phi & \mapsto & \frac{1}{2}\langle \phi,\phi\rangle.
  \end{array}
\end{displaymath}
This quadratic form is positive definite and defined over~$\Z$.
Furthermore, for every~$\phi$ in~$\Hom(\ss, \ss')$,
\begin{displaymath}
  Q_{\ss,\ss'}(\phi)
  =
  Q_{\ss',\ss}(\overline{\phi})
  =
  \overline{\phi} \phi
  =
  \deg(\phi)
\end{displaymath}
and, for all~$\ss''$ in~$\Sups$ and~$\psi$ in~$\Hom(\ss', \ss'')$,
\begin{equation}
  \label{eq:83}
  Q_{\ss, \ss''}(\psi \phi)
  =
  Q_{\ss',\ss''}(\psi)Q_{\ss, \ss'}(\phi).
\end{equation}

Define
\begin{displaymath}
  \bfR_{\ss, \ss'}
  \=
  \Hom_{\Fpalg}(\Fss,\Fsspr)
  \text{ and }
  \bfG_{\ss, \ss'}
  \=
  \Iso_{\Fpalg}(\Fss,\Fsspr).
\end{displaymath}
If ${\ss' = \ss}$, then ${\bfR_{\ss,\ss}=\Rss}$ and ${\bfG_{\ss, \ss} = \Gss}$.
Endow~$\bfR_{\ss, \ss'}$ with the unique distance such that, for every~$\varphi_0$ in~$\bfG_{\ss, \ss'}$, the map ${\Rss \to \bfR_{\ss, \ss'}}$ given by ${\psi \mapsto \varphi_0\circ \psi}$ is an isometry.
The natural map $\Hom(\ss, \ss') \to \bfR_{\ss, \ss'}$, denoted by~$\phi \mapsto \hphi$ as in Section~\ref{ss:from-elliptic-curves}, extends to an isomorphism of $\Zp$\nobreakdash-modules
\begin{equation}
  \label{eq:84}
  \Hom(\ss, \ss') \otimes \Zp \xrightarrow{\sim} \bfR_{\ss, \ss'}.
\end{equation}
Extend accordingly the map ${\Hom(\ss, \ss')\to \Hom(\ss', \ss)}$ given by ${\phi\mapsto \overline{\phi}}$ to a~$\Zp$-linear map ${\bfR_{\ss, \ss'}\to \bfR_{\ss', \ss}}$, the pairing ${\langle \ , \ \rangle}$ to a~$\Zp$-bilinear map ${\bfR_{\ss, \ss'}\times \bfR_{\ss, \ss'} \to \Zp}$, and~$Q_{\ss, \ss'}$ to a quadratic form on~$\bfR_{\ss, \ss'}$ taking values in~$\Zp$.
The identity~\eqref{eq:83} extends by continuity to all~$\phi$ in~$\bfR_{\ss, \ss'}$ and~$\psi$ in~$\bfR_{\ss', \ss''}$.
In particular, for all~$\varphi$ in~$\bfR_{\ss, \ss'}$ and~$\ell$ in~$\Zp$,
\begin{equation}
  \label{eq:85}
  Q_{\ss,\ss}(\varphi)
  =
  \nr(\varphi)
  \text{ and }
  Q_{\ss, \ss'}(\ell \varphi)
  =
  \ell^2 Q_{\ss, \ss'}(\varphi).
\end{equation}

For each nonzero~$\ell$ in~$\Zp$, put
\begin{displaymath}
  S_{\ell}(\ss, \ss')
  \=
  \{\varphi\in \bfR_{\ss, \ss'} \: Q_{\ss, \ss'}(\varphi)=\ell\}.
\end{displaymath}
The set~$S_{\ell}(\ss, \ss')$ is a \emph{supersingular sphere}.
It is a metric subspace of~$\bfR_{\ss, \ss'}$.
By~\eqref{eq:83} with ${\ss = \ss' = \ss''}$, the set~$S_1(\ss, \ss)$ is a subgroup of~$\Gss$ and, for every~$\varphi$ in~$S_1(\ss, \ss)$, the equality ${\varphi^{-1} = \overline{\varphi}}$ holds.

The supersingular sphere~$S_1(\ss, \ss)$ is a subgroup of~$\Gss$ and~$\Gss$ acts on~$\bfR_{\ss',\ss}$ by
\begin{equation}
  \label{eq:Gemultiplicationaction}
  \begin{array}{rcl}
    \Gss \times \bfR_{\ss',\ss}& \to & \bfR_{\ss',\ss}
    \\
    (g, \varphi) & \mapsto & g\varphi.
  \end{array}
\end{equation}
For each nonzero~$\ell$ in~$\Zp$, this action restricts to an action of~$S_1(\ss, \ss)$ on~$S_{\ell}(\ss',\ss)$.
The action~\eqref{eq:Gemultiplicationaction} is also the restriction of the natural action of the orthogonal group~$\oO_{Q_{\ss, \ss}}(\Zp)$ on~$S_{\ell}(\ss',\ss)$ to its subgroup~$S_1(\ss, \ss)$.

\begin{proposition}
  \label{p:supersingular-spheres}
  For all~$\ell$ in ${\Zp \ssetminus \{ 0 \}}$ and~$\ss$ and~$\ss'$ in~$\tSups$, the following properties hold.
  \begin{enumerate}
  \item[$(i)$]
    The supersingular sphere~$S_{\ell}(\ss, \ss')$ is nonempty and compact.
    Moreover, ${S_{\ell}(\ss, \ss') \subseteq \bfG_{\ss, \ss'}}$ if and only if~$\ell$ belongs to~$\Zp^{\times}$.
  \item[$(ii)$]
    The action of~$S_1(\ss, \ss)$ on~$S_{\ell}(\ss',\ss)$ induced by~\eqref{eq:Gemultiplicationaction} is faithful, transitive, and by isometries.
  \item[$(iii)$]
    There exists a unique Borel probability measure~$\mu^{\ss',\ss}_{\ell}$ on~$S_{\ell}(\ss',\ss)$ that is invariant under the action of~$S_1(\ss, \ss)$.
    It is also the unique Borel probability measure on~$S_{\ell}(\ss',\ss)$ that is invariant under the action of the orthogonal group~$\oO_{Q_{\ss', \ss}}(\Zp)$.
    Moreover, ${\supp(\mu^{\ss',\ss}_{\ell}) = S_{\ell}(\ss',\ss)}$.
  \item [$(iv)$]
    For all~$\ss''$ in~$\tSups$ and~$g$ in~$\bfG_{\ss, \ss'}$, the map~$\varphi \mapsto \varphi g$ maps~$S_{\ell}(\ss', \ss'')$ to~$S_{\ell \nr(g)}(\ss, \ss'')$ and~$\mu_{\ell}^{\ss', \ss''}$ to~$\mu_{\ell \nr(g)}^{\ss, \ss''}$.
  \end{enumerate}
\end{proposition}

The proof of Proposition~\ref{p:supersingular-spheres} relies on the following lemma.

\begin{lemma}
  \label{l:pseudo-group}
  For all~$\ss$ and~$\ss'$ in~$\tSups$,
  \begin{displaymath}
    \bfG_{\ss, \ss'}
    =
    \{\varphi \in \bfR_{\ss, \ss'} \: Q_{\ss, \ss'}(\varphi)\in \Zp^{\times}\}.
  \end{displaymath}
\end{lemma}

\begin{proof}
  For each~$\varphi$ in~$\bfG_{\ss, \ss'}$, \eqref{eq:83} with~$\ss'' = \ss$ implies
  \begin{displaymath}
    Q_{\ss,\ss'}(\varphi^{-1}) Q_{\ss', \ss}(\varphi)
    =
    Q_{\ss, \ss}(\varphi^{-1} \varphi)
    =
    Q_{\ss, \ss}(\widehat{\bfone_{\ss}})=1.
  \end{displaymath}
  It follows that~$Q_{\ss, \ss'}(\varphi)$ belongs to~$\Zp^{\times}$.

  Let~$\varphi$ be an element of~$\bfR_{\ss, \ss'}$ such that~$Q_{\ss, \ss'}(\varphi)$ belongs to~$\Zp^{\times}$ and put ${\ell \= Q_{\ss, \ss'}(\varphi)}$.
  Then,
  \begin{displaymath}
    \ell^{-1}\overline{\varphi}
    \in
    \bfR_{\ss',\ss},
    (\ell^{-1}\overline{\varphi})\varphi
    =
    \ell^{-1}(\overline{\varphi}\varphi)
    =
    \ell^{-1}Q_{\ss, \ss'}(\varphi)
    =
    1,
  \end{displaymath}
  and
  \begin{displaymath}
    \varphi(\ell^{-1}\overline{\varphi})
    =
    \ell^{-1}(\varphi \overline{\varphi})
    =
    \ell^{-1}Q_{\ss',\ss}(\overline{\varphi})
    =
    \ell^{-1}Q_{\ss, \ss'}(\varphi)=1.
  \end{displaymath}
  This proves that~$\ell^{-1}\overline{\varphi}$ is the inverse of~$\varphi$ and hence that~$\varphi$ belongs to~$\bfG_{\ss, \ss'}$.
\end{proof}

\begin{proof}[Proof of Proposition~\ref{p:supersingular-spheres}]
  The last assertion of item~$(i)$ is a direct consequence of Lemma~\ref{l:pseudo-group}.
  That~$S_{\ell}(\ss, \ss')$ is compact follows from the compactness of~$\bfR_{\ss, \ss'}$ and the continuity of~$Q_{\ss, \ss'}$.
  In view of~\eqref{eq:83} with~$\ss'' = \ss$, to prove ${S_{\ell}(\ss, \ss') \neq \emptyset}$ it suffices to prove that~$Q_{\ss,\ss}$ is surjective.
  Let~$n$ be in~$\Z_{\ge 0}$ and~$u$ in~$\Zp^{\times}$.
  By \cite[\emph{Chapitre}~II, \emph{Corollaire}~1.7]{Vig80}, the quaternion algebra~$\Bss$ contains an element~$\theta$ satisfying~$\theta^2=-p$ and a subalgebra~$\cK$ isomorphic to~$\Q_{p^2}$.
  By Lemma~\ref{l:integer-norms}, there is~$v$ in~$\OK$ satisfying ${\nr(v) = u}$.
  Put ${v' \= \theta^nv}$.
  Then, $v$ and~$v'$ belong to~$\Rss$ and ${Q_{\ss,\ss}(v') = \nr(v') = p^nu}$.
  This proves that~$Q_{\ss,\ss}$ is surjective and completes the proof of ${S_{\ell}(\ss, \ss') \neq \emptyset}$.

  To prove item~$(ii)$, let~$\varphi$ and~$\varphi'$ be elements of~$S_{\ell}(\ss',\ss)$ and put ${g \= \ell^{-1}\varphi'\overline{\varphi}}$.
  Then, ${g}$ belongs to~$\Rss$, ${g \varphi=\ell^{-1}\varphi'\ell=\varphi'}$, and~\eqref{eq:83} implies
  \begin{displaymath}
    Q_{\ss,\ss}(g)
    =
    \ell^{-2}Q_{\ss',\ss}(\varphi')Q_{\ss',\ss}(\varphi)
    =
    1.
  \end{displaymath}
  This proves that~$g$ belongs to~$S_1(\ss, \ss)$ and that the action of~$S_1(\ss, \ss)$ on~$S_{\ell}(\ss',\ss)$ is transitive.
  Since for all~$g$ in~$S_1(\ss, \ss)$ and~$\varphi$ in~$S_{\ell}(\ss',\ss)$ satisfying ${g\varphi = \varphi}$,
  \begin{equation}
    \label{eq:86}
    g
    =
    \ell^{-1}g(\varphi\overline{\varphi})
    =
    \ell^{-1}(g\varphi)\overline{\varphi}
    =
    \ell^{-1} \varphi\overline{\varphi}
    =
    \ell^{-1}\ell
    =
    1,
  \end{equation}
  the action of~$S_1(\ss, \ss)$ on~$S_{\ell}(\ss',\ss)$ is faithful.
  Finally, since for each~$g$ in~$\Gss$ the left multiplication map is an isometry, the action of~$S_1(\ss, \ss)$ on~$S_{\ell}(\ss',\ss)$ is by isometries.

  Item~$(iii)$ is a direct consequence of item~$(ii)$ and Lemma~\ref{l:existence-invariant-measure}.

  The first assertion of item~$(iv)$ follows from~\eqref{eq:83} and~\eqref{eq:85}, while the second from item~$(iii)$ and the fact that the left (resp.~right) multiplication map by~$g$ is an isometry.
\end{proof}

\subsection{From supersingular spheres to Hecke orbits}
\label{ss:supersingular-sphere-to-orbit}
This section proves the following proposition.
It relates supersingular spheres to partial Hecke orbits and defines a natural measure on the closure of a partial Hecke orbit inside a residue disc.

For all~$\ss$ and~$\ss'$ in~$\tSups$ and~$x$ in~$\hDss$, define the evaluation map
\begin{displaymath}
  \begin{array}{rrcl}
    \Ev^{x, \ss'} \: & \bfG_{\ss, \ss'} &\to & \hDsspr
    \\ & g & \mapsto & g \cdot x.
  \end{array}
\end{displaymath}

\begin{proposition}
  \label{p:supersingular-sphere-to-orbit}
  For all~$\ss$ and~$\ss'$ in~$\tSups$, $E$ in~$\Dss$, $x$ in~$\Piss^{-1}(E)$, coset~$\coset$ in~$\Zp^{\times} / \NE$, and~$\ell$ in~$\coset$, the following properties hold.
  \begin{enumerate}
  \item [$(i)$]
    The equality
    \begin{equation}
      \label{eq:87}
      \Pisspr^{-1} \left( \corbitc \cap \Dsspr \right)
      =
      \Ev^{x, \ss'}(S_{\ell}(\ss, \ss'))
    \end{equation}
    holds, and both sets are compact.
  \item[$(ii)$]
    The measure~$\hmu_{\coset}^{E, \ss'}$ on~$\hDsspr$, defined by
    \begin{equation}
      \label{eq:88}
      \hmu_{\coset}^{E, \ss'} \= (\Ev^{x, \ss'})_*(\mu_{\ell}^{\ss, \ss'}),
    \end{equation}
    depends on~$E$, ${\ss'}$, and~$\coset$ but on neither~$x$ nor~$\ell$.
    Moreover,
    \begin{equation}
      \label{eq:89}
      \supp \left( \hmu_{\coset}^{E, \ss'} \right)
      =
      \Pisspr^{-1} \left( \corbitc \cap \Dsspr \right).
    \end{equation}
  \end{enumerate}
\end{proposition}

The proof of this proposition relies on a couple of lemmas.
For all~$\ss$ in~$\tSups$ and~$x$ in~$\hDss$, the \emph{stabilizer}~$\bfG_{\ss,x}$ of~$x$ in~$\bfG_{\ss}$ is defined by
\begin{displaymath}
  \bfG_{\ss,x}
  \=
  \{ g\in \Gss \: g \cdot x = x \}.
\end{displaymath}

\begin{lemma}
  \label{l:norm-group}
  For all~$\ss$ in~$\tSups$, $E$ in~$\Dss$, and~$x$ in~$\Piss^{-1}(E)$, the equality
  \begin{equation}
    \label{eq:90}
    \NE
    =
    \{ \nr(g) \: g \in \bfG_{\ss, x} \}
  \end{equation}
  holds and ${(\Zp^{\times})^2 \subseteq \NE}$.
  In particular, $\NE$ is an open subgroup of~$\Zp^{\times}$ whose index is at most two if~$p$ is odd and at most four if~$p = 2$.
\end{lemma}

\begin{proof}
  If~$E$ is outside~$\SupsQpalg$, then~$\NE = (\Zp^{\times})^2$ by definition.
  Moreover, $x$ is outside~$\Xss(\OQpalg)$ and thus ${\bfG_{\ss, x} = \Zp^{\times}}$ by Lemma~\ref{l:fixed-properties}$(i)$.
  This implies that the right side of~\eqref{eq:90} equals~$(\Zp^{\times})^2$ and proves the lemma when~$E$ is outside~$\SupsQpalg$.

  Suppose~$E$ belongs to~$\SupsQpalg$ and let ${\alpha \: \tFE \to \Fss}$ be an isomorphism of formal $\Zp$\nobreakdash-modules such that~$(\FE, \alpha)$ represents~$x$.
  Consider the ring homomorphism ${\End(\FE) \to \End_{\Fpalg}(\Fss)}$ given by ${\varphi \mapsto \alpha \circ \tvarphi \circ \alpha^{-1}}$.
  Suppose~$E$ is a formal \CM{} point, let~$\cK$ be the field of fractions of~$\End(\FE)$ and~$\nr \: \cK \to \Qp$ its norm map.
  Then, ${\nr(\varphi) = \nr(\alpha \circ \tvarphi \circ \alpha^{-1})}$ and thus
  \begin{equation}
    \label{eq:91}
    \NE
    =
    \left\{ \nr(\alpha \circ \tvarphi \circ \alpha^{-1}) \: \varphi \in \Aut(\FE) \right\}.
  \end{equation}
  If~$E$ is not a formal \CM{} point, then the group~$\Aut(\FE)$ is isomorphic to~$\Zp^{\times}$ and the equality above also holds.
  In all cases, \eqref{eq:90} is a direct consequence of~\eqref{eq:91} and Lemma~\ref{l:fixed}.
  Since~$\FE$ is a formal $\Zp$\nobreakdash-module, \eqref{eq:91} implies ${(\Zp^{\times})^2 \subseteq \NE}$.
\end{proof}

Item~$(ii)$ of the following lemma is a reformulation of \cite[Theorem~1.2]{Men12}.

\begin{lemma}\label{l:count_Hom_n(ss,ss')}
  For all~$\ss$ and~$\ss'$ in~$\tSups$, the following properties hold.
  \begin{enumerate}
  \item[$(i)$]
    For every~$n$ in ${\Z_{>0} \ssetminus p\Z}$,
    \begin{displaymath}
      \#\Hom_n(\ss,\ss')
      =
      \# \Aut(\ss') \cdot \deg(T_n(\ss)|_{\{ \ss' \}}).
    \end{displaymath}
  \item[$(ii)$]
    For every~$\varepsilon$ in~$\R_{>0}$, there is~$C$ in~$\R_{>0}$ such that, for every~$n$ in ${\Z_{>0} \ssetminus p\Z}$,
    \begin{displaymath}
      \left|\frac{\#\Hom_n(\ss,\ss')}{\sigma_1(n)} - \frac{24}{(p -1)}\right|
      \le
      C n^{-\frac{1}{2}+\varepsilon}.
    \end{displaymath}
  \end{enumerate}
\end{lemma}
\begin{proof}
  Item~$(ii)$ is a direct consequence of item~$(i)$, \cite[Theorem~1.2]{Men12}, and ${\# \Aut(\ss')\leq 24}$.

  To prove item~$(i)$, put
  \begin{displaymath}
    \sC
    \=
    \left\{ C \le \ss(\Fpalg) \: \#C = n, \ss/C = \ss' \right\}.
  \end{displaymath}
  For each isogeny~$\phi$ in~$\Hom_n(\ss,\ss')$, the endomorphism~$\overline{\phi} \phi$ equals the morphism of multiplication by~$n$ on~$\ss$, which is separable, see, \emph{e.g.}, \cite[Chapter~III, Corollary~5.4]{Sil09}.
  This proves that~$\phi$ is separable and thus that its kernel~$\Ker(\phi)$ belongs to~$\sC$.
  Thus, $\phi \mapsto \Ker(\phi)$ defines a surjective map~$K \: \Hom_n(\ss,\ss') \to \sC$, see, \emph{e.g.}, \cite[Chapter~III, Proposition~4.12]{Sil09}.
  The desired identity follows from ${\# \sC = \deg(T_n(\ss)|_{\{\ss'\}})}$ and, for every~$C$ in~$\sC$, the equality ${\# K^{-1}(C) = \# \Aut(\ss')}$, see, \emph{e.g.}, \cite[Chapter~III, Corollary~4.11]{Sil09}.
\end{proof}

\begin{proof}[Proof of Proposition~\ref{p:supersingular-sphere-to-orbit}]
  To prove item~$(i)$, note first that, for each~$\ell$ in~$\coset$, the set~$\Ev^{x, \ss'}(S_{\ell}(\ss, \ss'))$ is compact because~$\Ev^{x, \ss'}$ is continuous and~$S_{\ell}(\ss, \ss')$ is compact by Proposition~\ref{p:supersingular-spheres}$(i)$.

  The proof of item~$(i)$ begins establishing~\eqref{eq:87} when~$E$ belongs to~$\SupsQpalg$.
  The first step is to prove that the left side is contained in the right side.
  Since the right side is compact and~$\Piss$ is continuous, it suffices to show that every element~$x'$ of
  \begin{displaymath}
    \Pisspr^{-1} \left( \corbit \cap \Dsspr \right)
  \end{displaymath}
  belongs to the right side of~\eqref{eq:87}.
  Put ${E' \= \Pisspr(x')}$ and let ${\alpha_0' \: \widetilde{E'} \to \sspr}$ be an isomorphism such that the induced isomorphism of formal $\Zp$\nobreakdash-modules ${\alpha' \: \tFEpr \to \Fsspr}$ is such that the deformation~$(\tFEpr, \alpha')$ of~$\Fsspr$ represents~$x'$.
  Moreover, let ${\alpha_0 \: \tE \to \ss}$ be an isomorphism of elliptic curves, so that the induced isomorphism ${\alpha \: \tFE \to \Fss}$ is such that~$(\FE, \alpha)$ represents~$x$.
  By definition of~$\corbit$, there is~$n$ in ${\coset \cap \Z_{>0}}$ such that~$E'$ belongs to~$\supp(T_n(E))$.
  That is, there is an isogeny~$\phi \: E \to E'$ such that~$\deg(\phi)$ belongs to~$\coset$.
  Denote by ${\varphi \: \FE \to \FEpr}$ the induced isomorphism and put ${g \= \alpha' \circ \tvarphi \circ \alpha^{-1}}$.
  So, ${g}$ belongs to~$\bfG_{\ss, \ss'}$ and satisfies~${Q_{\ss, \ss'}(g) = \deg(\phi)}$.
  Furthermore, the deformation~$g \cdot (\FE, \alpha)$, which equals~$(\FE, \alpha' \circ \tvarphi)$, is isomorphic to~$(\FEpr, \alpha')$ via the isomorphism~$\varphi$, so ${g \cdot x = x'}$.
  Since~$Q_{\ss, \ss'}(g)$ and~$\ell$ are both in~$\coset$, by Lemma~\ref{l:norm-group} there is~$g_0$ in~$\Gss$ such that
  \begin{displaymath}
    g_0 \cdot x = x
    \text{ and }
    \nr(g_0)
    =
    \ell Q_{\ss, \ss'}(g)^{-1}.
  \end{displaymath}
  Thus,
  \begin{displaymath}
    Q_{\ss, \ss'}(gg_0) = \ell
    \text{ and }
    (gg_0) \cdot x = g \cdot x = x'.
  \end{displaymath}
  This proves that~$x'$ belongs to the right side of~\eqref{eq:87} and completes the proof that the left side of~\eqref{eq:87} is contained in the right side when~$E$ belongs to~$\SupsQpalg$.

  To prove the reverse inclusion, recall that ${S_{\ell}(\ss, \ss') \neq \emptyset}$ by Proposition~\ref{p:supersingular-spheres}$(i)$ and let~$g$ be in~$S_{\ell}(\ss, \ss')$.
  Let~$(m_j)_{j = 1}^{\infty}$ be a sequence in~$\Z_{>0}$ tending to~$\infty$ that is contained in the coset~$\ell (\Zp^{\times})^2$ of~$\Zp^{\times} / (\Zp^{\times})^2$ and that converges to~$\ell$ in~$\Zp$.
  For each~$j$, let~$u_j$ in~$\Zp^{\times}$ be such that ${m_j = \ell u_j^2}$ and such that ${u_j \to 1}$ in~$\Zp^{\times}$ as ${j \to \infty}$.
  By~\eqref{eq:lowb-sigma1}, \eqref{eq:deg-Hecke-op-divisor}, Proposition~\ref{p:supersingular-spheres}, and Lemma~\ref{l:count_Hom_n(ss,ss')}$(ii)$, for every~$c$ in~$]\frac{1}{2}, 1[$ the hypotheses of Theorem~\ref{t:p-adic-Linnik} are satisfied for~$n = 4$ and~$Q = Q_{\ss, \ss'}$.
  Applying this theorem and using that the support of the limit measure~$\mu_\ell^{\ss, \ss'}$ equals~$S_{\ell}(\ss, \ss')$ and thus contains~$g$, it follows that, for each~$j$, there is~$\phi_j$ in~$\Hom_{m_j}(\ss, \ss')$ such that
  \begin{displaymath}
    M_{u_j}^{-1}(\hphi_j) \to g \text{ in } S_{\ell}(\ss, \ss')
    \text{ as }
    j \to \infty.
  \end{displaymath}
  Since~$u_j \to 1$ in~$\Zp$ as ${j \to \infty}$, it follows that~$\hphi_j \to g$ in~$\bfG_{\ss, \ss'}$ as ${j \to \infty}$ and thus that~$\hphi_j \cdot x \to g \cdot x$ in~$\hDsspr$ as ${j \to \infty}$.
  Since the sequence~$(m_j)_{j = 1}^{\infty}$ is contained in the coset~$\coset$, for each~$j$,
  \begin{displaymath}
    \Pisspr(\hphi_j \cdot x)
    \in
    \supp \left( T_{m_j}(E) \right) \cap \Dsspr
    \subseteq
    \corbit \cap \Dsspr.
  \end{displaymath}
  This proves that~$g \cdot x$ belongs to the closure of ${\Pisspr^{-1} \left( \corbit \cap \Dsspr \right)}$ and completes the proof of~\eqref{eq:87} when~$E$ belongs to~$\SupsQpalg$.
  In particular, ${\Pisspr^{-1} \left( \corbitc \cap \Dsspr \right)}$ is compact.

  It remains to prove~\eqref{eq:87} when~$E$ is outside~$\SupsQpalg$.
  The proof relies on the action by isometries of~$\Gss$ on~$\hDss$ (Lemma~\ref{l:action-regularity}$(i)$).
  By~\eqref{eq:14}, for every~$x_0$ in~$\hDss$, the Hausdorff distance between
  \begin{displaymath}
    \Pisspr^{-1}\left( \corbitc \cap \Dsspr \right)
    \text{ and }
    \Pisspr^{-1}\left( \overline{\Orb_{\coset}(\Piss(x_0))} \cap \Dsspr \right)
  \end{displaymath}
  and between
  \begin{displaymath}
    \Ev^{x, \ss'}(S_{\ell}(\ss, \ss'))
    \text{ and }
    \Ev^{x_0, \ss'}(S_{\ell}(\ss, \ss')),
  \end{displaymath}
  are both bounded from above by the distance between~$x$ and~$x_0$.
  Since~$\SupsQpalg$ is dense in~$\Sups$ and, for every~$x_0$ in ${\Piss^{-1}(\SupsQpalg \cap \Dss)}$, the equality~\eqref{eq:87} holds with~$x$ replaced by~$x_0$, it follows that the Hausdorff distance between the left side and the right side of~\eqref{eq:87} equals zero.
  Since both of these sets are compact, it follows that they are equal.

  To prove item~$(ii)$, for each~$g$ in~$\Gss$, denote by ${\sT_g \: \bfG_{\ss, \ss'} \to \bfG_{\ss, \ss'}}$ the right multiplication map ${h \mapsto hg}$.
  For every~$x$ in~$\hDss$,
  \begin{equation}
    \label{eq:92}
    \Ev^{x, \ss'} \circ \sT_g
    =
    \Ev^{g \cdot x, \ss'}.
  \end{equation}
  Let~$x$ and~$x'$ be in~$\Piss^{-1}(E)$ and~$\ell$ and~$\ell'$ in~$\coset$.
  Then, there exists~$\phi$ in~$\Aut(\ss)$ satisfying ${\hphi \cdot x' = x}$ and, by Lemma~\ref{l:norm-group}, there is~$g$ in~$\Gss$ satisfying ${g \cdot x' = x'}$ and ${\nr(g) \ell = \ell'}$.
  Then, ${Q_{\ss, \ss}(\hphi) = \deg(\phi) = 1}$, so~$\hphi$ belongs to~$S_1(Q_{\ss, \ss})$ and thus ${(\sT_{\hphi g})_* \mu_{\ell}^{\ss, \ss'} = \mu_{\ell'}^{\ss, \ss'}}$ by Proposition~\ref{p:supersingular-spheres}$(iv)$.
  Combined with~\eqref{eq:92}, this implies
  \begin{displaymath}
    (\Ev^{x', \ss'})_* \mu_{\ell'}^{\ss, \ss'}
    =
    (\Ev^{x', \ss'})_* ((\sT_{\hphi g})_* \mu_{\ell}^{\ss, \ss'})
    =
    (\Ev^{(\hphi g) \cdot x', \ss'})_* \mu_{\ell}^{\ss, \ss'}
    =
    (\Ev^{x, \ss'})_* \mu_{\ell}^{\ss, \ss'}.
  \end{displaymath}
  This proves the first assertion of item~$(ii)$.

  To prove the remaining assertions of item~$(ii)$, fix~$x$ in~$\Piss^{-1}(E)$ and~$\ell$ in~$\coset$.
  Then, ${\supp(\mu_{\ell}^{\ss, \ss'}) = S_{\ell}(\ss, \ss')}$ by Proposition~\ref{p:supersingular-spheres}$(iii)$ and thus
  \begin{equation}
    \label{eq:93}
    \supp \left( (\Ev^{x, \ss'})_* \mu_{\ell}^{\ss, \ss'} \right)
    =
    \Ev^{x, \ss'}(S_{\ell}(\ss, \ss')).
  \end{equation}
  Then, the desired assertion follows from item~$(i)$.
\end{proof}

\subsection{Proof of Theorem~\ref{t:Hecke-orbits-pr}}
\label{ss:Hecke-orbits-proof}
The main ingredient is the following proposition, whose proof is based on the asymptotic distribution of integer points on $p$\nobreakdash-adic spheres in Section~\ref{s:p-adic-Linnik}.

Let~$E$ be in~$\Sups$ and~$\coset$ in~$\Zp^{\times} / \NE$.
That~$\corbitc$ is compact follows from Proposition~\ref{p:supersingular-sphere-to-orbit}$(i)$, the finiteness of~$\tSups$, and the continuity of~$\Piss$.
For each~$\ss$ in~$\tSups$, let~$\hmu_{\coset}^{E, \ss}$ be the measure on~$\hDss$ given by~\eqref{eq:88} in Proposition~\ref{p:supersingular-sphere-to-orbit}$(ii)$ and put
\begin{displaymath}
  \mu_{\coset}^{E, \ss}
  \=
  (\Piss)_* \hmu_{\coset}^{E, \ss}.
\end{displaymath}
Since~$\Piss$ is continuous, ${\mu_{\coset}^{E, \ss}}$ is a Borel probability measure on~$\Ell(\Cp)$ and ${\supp(\mu_{\coset}^{E, \ss}) = \corbitc \cap \Dss}$.
Put
\begin{displaymath}
  \mu_{\coset}^{E}
  \=
  \frac{24}{p - 1} \sum_{\ss \in \tSups} \frac{1}{\# \Aut(\ss)} \mu_{\coset}^{E, \ss}.
\end{displaymath}
Then, ${\mu_{\coset}^{E}}$ is a Borel probability measure on~$\Ell(\Cp)$ by the mass formula~\eqref{eq:Deuring-Eichler} and ${\supp(\mu_{\coset}^{E}) = \corbitc}$.

\begin{proposition}
  \label{p:Hecke-orbits}
  For all~$E$ in~$\Sups$, $\ss$ in~$\tSups$, and locally constant function ${F \: \Dss \to \R}$, the following property holds.
  For all~$\varepsilon$ in~$\R_{>0}$ and coset~$\coset$ in~$\Zp^{\times} / \NE$, there is~$C$ in~$\R_{>0}$ such that, for all~$\ss'$ in~$\tSups$, $E'$ in ${\overline{\Orb_{\NE}(E)} \cap \Dsspr}$, and~$n$ in ${\coset \cap \Z_{>0}}$ for which ${\Hom_n(\ss', \ss) \neq \emptyset}$,
  \begin{displaymath}
    \left| \int F \dd \odelta_{T_{n}(E')|_{\Dss}}- \int F \dd \mu_{\coset}^{E, \ss} \right|
    \leq
    C \frac{n^{\frac{1}{2}+\varepsilon}}{\# \Hom_{n}(\ss', \ss)}.
  \end{displaymath}
\end{proposition}

The proof of this proposition relies on the following lemma.

\begin{lemma}
  \label{l:orbit-closure}
  Let~$E$ be in~$\Sups$ and~$\coset$ and~$\coset'$ cosets in~$\Zp^{\times} / \NE$.
  Then, for every~$E'$ in~$\corbitc$,
  \begin{displaymath}
    \NEpr
    =
    \NE,
    \overline{\Orb_{\coset'}(E')}
    =
    \overline{\Orb_{\coset \cdot \coset'}(E)},
    \text{ and }
    \mu_{\coset'}^{E'}
    =
    \mu_{\coset \cdot \coset'}^E.
  \end{displaymath}
\end{lemma}

\begin{proof}
  Let~$\ss$ and~$\sspr$ in~$\tSups$ be such that~$E$ and~$E'$ belong to~$\Dss$ and~$\Dsspr$, respectively.
  Moreover, fix~$\ell$ in~$\coset$, $\ell'$ in~$\coset'$, $x$ in~$\Piss^{-1}(E)$, and~$x'$ in~$\Pisspr^{-1}(E')$.
  By Proposition~\ref{p:supersingular-sphere-to-orbit}$(i)$, there is~$g_0$ in~$S_{\ell}(\ss, \sspr)$ such that ${g_0 \cdot x = x'}$.
  Then, ${\bfG_{\sspr, x'} = g_0 \bfG_{\ss, x} g_0^{-1}}$ and thus ${\NEpr = \NE}$ by Lemma~\ref{l:norm-group}.

  To prove the remaining equalities, let~$\ss''$ be in~$\tSups$ and denote by ${\sT \: \bfG_{\ss', \ss''} \to \bfG_{\ss, \ss''}}$ the right multiplication map~$h \mapsto h g_0$.
  By Proposition~\ref{p:supersingular-spheres}$(iv)$,
  \begin{displaymath}
    \sT(S_{\ell'}(\ss', \ss''))
    =
    S_{\ell \ell'}(\ss, \ss'')
    \text{ and }
    \sT_* (\mu_{\ell'}^{\ss', \ss''})
    =
    \mu_{\ell \ell'}^{\ss, \ss''}.
  \end{displaymath}
  Hence,
  \begin{multline*}
    \Ev^{x', \ss''}(S_{\ell'}(\ss', \ss''))
    =
    \{ g \cdot (g_0 \cdot x) \: g \in S_{\ell'}(\ss', \ss'') \}
    =
    \{ (g g_0) \cdot x \: g \in S_{\ell'}(\ss', \ss'') \}
    \\ =
    \{ \breve{g} \cdot x \: \breve{g} \in S_{\ell \ell'}(\ss, \ss'') \}
    =
    \Ev^{x, \ss''}(S_{\ell \ell'}(\ss, \ss''))
  \end{multline*}
  and
  \begin{displaymath}
    \Ev^{x, \ss''}_* (\mu_{\ell \ell'}^{\ss, \ss''})
    =
    \left( \Ev^{x, \ss''} \circ \sT \right)_* (\mu_{\ell'}^{\ss', \ss''})
    =
    \Ev^{x', \ss''}_* (\mu_{\ell'}^{\ss', \ss''}).
  \end{displaymath}
  Together with Proposition~\ref{p:supersingular-sphere-to-orbit} and the definition of the measures~$\mu_{\coset'}^{E'}$ and~$\mu_{\coset \cdot \coset'}^{E}$, this implies the second and third equalities.
\end{proof}

\begin{proof}[Proof of Proposition~\ref{p:Hecke-orbits}]
  Let~$\ss_0$ in~$\tSups$ be such that~$E$ belongs to~$\bfD_{\ss_0}$, $x_0$ in~$\Pi_{\ss_0}^{-1}(E)$, and put ${\breve{F} \= F \circ \Piss \circ \Ev^{x_0, \ss}}$.
  Since~$\Piss$ and~$\Ev^{x_0, \ss}$ are both continuous, ${\breve{F}}$ is locally constant.
  Let~$\delta$ in~$\R_{>0}$ be such that~$\breve{F}$ is constant on every ball of~$\bfG_{\ss_0, \ss}$ of radius~$\delta$.
  Fix~$\varepsilon$ in~$\R_{>0}$ and let~$C$ be given by Corollary~\ref{c:functional-deviation} with~$n = 4$.
  Let~$\ss'$ be in~$\tSups$ and~$E'$ in ${\overline{\Orb_{\NE}(E)} \cap \Dsspr}$.
  Then, by Proposition~\ref{p:supersingular-sphere-to-orbit}$(i)$ with~$\coset$ replaced by~$\NE$, there is~$g$ in~$S_1(\ss_0, \ss')$ such that~$g \cdot x_0$ belongs to~$\Pisspr^{-1}(E')$.
  Put ${x' \= g \cdot x_0}$.
  Denote by ${\sT \: \bfG_{\ss', \ss} \to \bfG_{\ss_0, \ss}}$ the right multiplication map given by ${g' \mapsto g' g}$ and put
  \begin{equation}
    \label{eq:94}
    \breve{F}'
    \=
    \breve{F} \circ \sT
    \text{ and }
    S_E
    \=
    \bigcup_{\coset' \in \Zp^{\times}/\NE} \overline{\Orb_{ \coset'}(E)}.
  \end{equation}
  The set~$S_E$ is compact by Lemma~\ref{l:norm-group} and~$\breve{F}'$ satisfies
  \begin{equation}
    \label{eq:95}
    \breve{F}'
    =
    F \circ \Piss \circ \Ev^{x', \ss}
    \text{ and }
    \sup_{\bfG_{\ss', \ss}} |\breve{F}'|
    =
    \sup_{S_E \cap \Dss} |F|.
  \end{equation}
  Moreover, since~$\sT$ maps~$\bfG_{\ss', \ss}$ to~$\bfG_{\ss_0, \ss}$ isometrically, the function~$\breve{F}'$ is constant on every ball of~$\bfG_{\ss', \ss}$ of radius~$\delta$.

  Let~$n$ in ${\coset \cap \Z_{>0}}$ be such that ${\Hom_n(\ss', \ss) \neq \emptyset}$.
  Then, \eqref{eq:14} and~\eqref{eq:95} yield
  \begin{equation}
    \label{eq:96}
    \int F \dd \odelta_{T_{n}(E')|_{\Dss}}
    =
    \frac{1}{\# \Hom_{n}(\ss', \ss)}\sum_{\phi \in \Hom_n(\ss', \ss)} \breve{F}'(\hphi).
  \end{equation}
  Furthermore, by Lemma~\ref{l:orbit-closure}, the definition of~$\mu_{\coset}^{E', \ss}$, and the change of variables formula,
  \begin{displaymath}
    \int F \dd \mu_{\coset}^{E, \ss}
    =
    \int F \dd \mu_{\coset}^{E', \ss}
    =
    \int \breve{F}' \dd \mu_{n}^{\ss',\ss}.
  \end{displaymath}
  Together with Proposition~\ref{p:supersingular-spheres}, \eqref{eq:96}, and Corollary~\ref{c:functional-deviation} with~$\ell$ and~$m$ equal to~$n$, this implies
  \begin{multline}
    \left| \int F \dd \odelta_{T_{n}(E')|_{\Dss}}- \int F \dd \mu_{\coset}^{E, \ss} \right|
    =
    \left| \frac{1}{\# \Hom_{n}(\sspr, \ss)}\sum_{\phi \in \Hom_n(\sspr, \ss)} \breve{F}'(\hphi) - \int \breve{F}' \dd \mu_n^{\sspr,\ss} \right|
    \\ \le
    C \left( \sup_{S_E \cap \Dss} |F| \right) \frac{n^{\frac{1}{2}+\varepsilon}}{\# \Hom_{n}(\sspr, \ss)}.
    \qedhere
  \end{multline}
\end{proof}

\begin{proof}[Proof of Theorem~\ref{t:Hecke-orbits-pr}]
  For every~$E'$ in~$\corbitc$, Lemma~\ref{l:orbit-closure} yields ${\NEpr = \NE}$.

  By the considerations at the beginning of this section, in the case where ${\coset \subseteq \Zp^{\times}}$ it only remains to prove the estimate~\eqref{eq:82}.
  In that case, \eqref{eq:82} is a direct consequence of Proposition~\ref{p:Hecke-orbits} combined with~\eqref{eq:lowb-sigma1}, \eqref{eq:deg-Hecke-op-divisor}, the definition of~$\mu_{\coset}^E$, and \cite[Theorem~1.2]{Men12} or Lemma~\ref{l:count_Hom_n(ss,ss')}$(ii)$.

  Suppose ${\coset \not\subseteq \Zp^{\times}}$ and let~$\coset_0$ be the coset in~$\Zp^{\times} / \NE$ and~$k$ in~$\Z_{>0}$ such that ${\coset = p^k \coset_0}$.
  Then, for every~$E'$ in~$\Sups$,
  \begin{displaymath}
    \Orb_{\coset}(E')
    =
    T_{p^k} \left(\Orb_{\coset_0}(E') \right)
    \text{ and }
    \overline{\Orb_{\coset}(E')}
    =
    T_{p^k} \left( \overline{\Orb_{\coset_0}(E')} \right),
  \end{displaymath}
  and this last set is compact.
  Moreover, putting
  \begin{displaymath}
    \mu_{\coset}^E
    \=
    \frac{1}{\sigma_1(p^k)} (T_{p^k})_* \mu_{\coset_0}^E,
  \end{displaymath}
  ${\supp(\mu_{\coset}^E) = \corbitc}$.
  The estimate~\eqref{eq:82} is a direct consequence of the same formula with~$\coset$ replaced by~$\coset_0$, using the change of variables formula.
\end{proof}

\subsection{On partial Hecke orbits and their limit measures}
\label{ss:applications}

This section proves the following proposition.
Section~\ref{ss:CM-orbits} relies on it.

\begin{proposition}
  \label{p:orbit-measures}
  For every~$E$ in~$\Sups$, the following properties hold.
  \begin{enumerate}
  \item [$(i)$]
    For all distinct cosets~$\coset$ and~$\coset'$ in~$\Zp^{\times} / \NE$, the sets~$\corbitc$ and~$\overline{\Orb_{\coset'}(E)}$ are disjoint.
  \item[$(ii)$]
    For all distinct cosets~$\coset$ and~$\coset'$ in~$\Qp^{\times} / \NE$ contained in~$\Zp$, the sets~$\corbitc$ and~$\overline{\Orb_{\coset'}(E)}$ are different.
    In particular, the measures~$\mu_{\coset}^E$ and~$\mu_{\coset'}^E$ are different.
  \end{enumerate}
\end{proposition}

The proof of this proposition relies on a couple of lemmas.
As in Section~\ref{ss:CM-from-canonical}, denote Katz' valuation by~$\kval$ and define

\begin{equation}
  \label{eq:97}
  \begin{array}{rcl}
    \kproj \: \Sups
    & \to & \left[0, \frac{p}{p + 1} \right]
    \\
    x & \mapsto & \min\left\{ \kval(x), \frac{p}{p+1} \right\}.
  \end{array}
\end{equation}

\begin{lemma}
  \label{l:n-th-pull-back}
  For all~$x$ in~$\left] 0, \frac{p}{p + 1} \right]$, divisor~$\fD$ on~$\kproj^{-1}(x)$, and~$n$ in~$\Z_{>0}$,
  \begin{equation}
    \label{eq:98}
    \kproj \left( \supp(T_{p^n} \fD) \right)
    \subseteq
    \left[p^{-n}x, \frac{p}{p + 1} \right]
    \text{ and }
    (T_{p^n} \fD) |_{\kproj^{-1}(p^{-n}x)}
    =
    (\t^n|_{\kproj^{-1}(p^{-n}x)})^* \fD.
  \end{equation}
\end{lemma}

The proof of this lemma relies on the following one.

\begin{lemma}[\textcolor{black}{\cite[Proposition~4.5]{HerMenRiv20}}]
  \label{l:katz kite}
  Denote by~$\uptau_0$ the identity on~$\Div \left( \left[0, \frac{p}{p + 1} \right] \right)$, $\uptau_1$ the piecewise-affine correspondence on~$\left[ 0, \frac{p}{p + 1} \right]$, defined by
  \begin{displaymath}
    \uptau_1(x)
    \=
    \begin{cases}
      [px] + p[\frac{x}{p}]
      & \text{if } x \in \left[0, \frac{1}{p + 1} \right];
      \\
      [1 - x] + p[\frac{x}{p}]
      & \text{if } x \in \left] \frac{1}{p + 1}, \frac{p}{p + 1} \right],
    \end{cases}
  \end{displaymath}
  and, for each integer~$m$ satisfying ${m \ge 2}$, $\uptau_m$ the correspondence on~$\left[ 0, \frac{p}{p + 1} \right]$, defined recursively by
  \begin{displaymath}
    \uptau_m
    \=
    \uptau_1 \circ \uptau_{m - 1} - p \uptau_{m - 2}.
  \end{displaymath}
  Then, for all~$m$ in~$\Z_{\ge 0}$ and~$n_0$ in ${\Z_{>0} \ssetminus p\Z}$,
  \begin{equation}
    \label{eq:99}
    (\kproj)_* \circ T_{p^m n_0}|_{\Sups}
    =
    \sigma_1(n_0) \cdot \uptau_m \circ (\kproj)_*.
  \end{equation}
\end{lemma}

\begin{proof}[Proof of Lemma~\ref{l:n-th-pull-back}]
  By Lemma~\ref{l:katz kite}, for all~$x'$ in~$\left] 0, \frac{p}{p + 1} \right]$ and divisor~$\fD$ on~$\kproj^{-1} \left( \left[ x', \frac{p}{p + 1} \right] \right)$, the divisor~$T_p \fD$ is supported on~$\kproj^{-1}\left(\left[ p^{-1} x', \frac{p}{p + 1} \right] \right)$.
  Together with~\eqref{eq:Hecke-Tpr} and an induction argument, this implies the first assertion.

  When ${n = 1}$, the second assertion is a direct consequence of~\eqref{eq:45} in Lemma~\ref{l:canonical-analyticity} and Lemma~\ref{l:sups-canonical-subgroup}.
  Let~$n$ in~$\Z_{>0}$ be such that the second assertion holds.
  By~\eqref{eq:Hecke-Tpr} and the first assertion,
  \begin{equation}
    \label{eq:100}
    (T_{p^{n + 1}}\fD)|_{\kproj^{-1}(p^{-(n + 1)} x)}
    =
    (T_p(T_{p^n}\fD))|_{\kproj^{-1}(p^{-(n + 1)} x)}.
  \end{equation}
  Moreover, by~\eqref{eq:45} in Lemma~\ref{l:canonical-analyticity} and Lemma~\ref{l:sups-canonical-subgroup},
  \begin{equation}
    \label{eq:101}
    \begin{split}
      (T_p(T_{p^n}\fD))|_{\kproj^{-1}(p^{-(n + 1)} x)}
      & =
        (\t|_{\kproj^{-1}(p^{-(n + 1)} x)})^* ((T_{p^n}\fD)|_{\kproj^{-1}(p^{-n}x)}))
      \\ & =
           (\t|_{\kproj^{-1}(p^{-(n + 1)} x)})^* ((\t^n|_{\kproj^{-1}(p^{-n}x)})^* \fD)
      \\ & =
           (\t^{n + 1}|_{\kproj^{-1}(p^{-(n + 1)} x)})^* \fD.
    \end{split}
  \end{equation}
  Together with~\eqref{eq:100}, this implies the second assertion with~$n$ replaced by~${n + 1}$.
  The second assertion thus follows by induction.
\end{proof}

\begin{proof}[Proof of Proposition~\ref{p:orbit-measures}]
  To prove item~$(i)$, suppose that~$\corbitc$ and~$\overline{\Orb_{\coset'}(E)}$ intersect and let~$E'$ be a common element.
  Let~$\ss$ and~$\ss'$ in~$\tSups$ be such that~$E$ and~$E'$ belong to~$\Dss$ and~$\Dsspr$, respectively, and fix~$x$ in~$\Piss^{-1}(E)$ and~$x'$ in~$\Pisspr^{-1}(E')$.
  By Proposition~\ref{p:supersingular-sphere-to-orbit}$(i)$, there are~$g$ and~$g'$ in~$\bfG_{\ss, \ss'}$ so that
  \begin{equation}
    \label{eq:102}
    Q_{\ss, \ss'}(g)
    \in
    \coset,
    Q_{\ss, \ss'}(g')
    \in
    \coset',
    \text{ and }
    g \cdot x
    =
    x'
    =
    g' \cdot x.
  \end{equation}
  This implies that~$g^{-1} g'$ belongs to~$\bfG_{\ss, x}$, and thus
  \begin{displaymath}
    Q_{\ss, \ss'}(g)^{-1} Q_{\ss, \ss'}(g')
    =
    \nr(g^{-1} g')
    \in
    \NE
    \text{ and }
    \coset
    =
    \coset'.
  \end{displaymath}

  The second assertion of item~$(ii)$ follows from the first and Theorem~\ref{t:Hecke-orbits}.
  To prove the first assertion, suppose ${\corbitc = \overline{\Orb_{\coset'}(E)}}$.
  Let~$n$ and~$n'$ be in~$\Z_{\ge 0}$ and~$\coset_0$ and~$\coset_0'$ in~$\Zp^{\times} / \NE$ such that ${\coset = p^n \coset_0}$ and ${\coset' = p^{n'} \coset_0'}$.
  Then,
  \begin{equation}
    \label{eq:103}
    \corbitc
    =
    T_{p^n} \left( \overline{\Orb_{\coset_0}(E)} \right)
    \text{ and }
    \overline{\Orb_{\coset'}(E)}
    =
    T_{p^{n'}} \left( \overline{\Orb_{\coset_0'}(E)} \right).
  \end{equation}
  Corollary~\ref{c:Hecke-on-orbits} implies
  \begin{equation}
    \label{eq:104}
    \mu_{\coset}^E
    =
    \frac{1}{\sigma_1(p^n)} (T_{p^n})_* \mu_{\coset_0}^E
    \text{ and }
    \mu_{\coset'}^E
    =
    \frac{1}{\sigma_1(p^{n'})} (T_{p^{n'}})_* \mu_{\coset_0'}^E.
  \end{equation}
  Put~$x_E \= \kproj(E)$.
  By~\eqref{eq:99} in Lemma~\ref{l:katz kite} with ${m = 0}$, the sets~$\overline{\Orb_{\coset_0}(E)}$ and~$\overline{\Orb_{\coset_0'}(E)}$ are both contained in~$\kproj^{-1}(x_E)$.
  Then, \eqref{eq:103} and Lemmas~\ref{l:n-th-pull-back} and~\ref{l:katz kite} imply that~$\corbitc$ is contained in~$\kproj^{-1} \left( \left[ p^{-n} x_E, \frac{p}{p + 1} \right] \right)$ and intersects~$\kproj^{-1}(p^{-n} x_E)$, and that~$\overline{\Orb_{\coset'}(E)}$ is contained in~$\kproj^{-1} \left( \left[ p^{-n'} x_E, \frac{p}{p + 1} \right] \right)$ and intersects~$\kproj^{-1}(p^{-n'} x_E)$.
  It follows that ${n = n'}$.
  When ${n = 0}$, the desired assertion follows from item~$(i)$.
  Suppose ${n \ge 1}$.
  By Lemma~\ref{l:sups-canonical-subgroup}, the degree~$\delta$ of~$\t^n|_{\kproj^{-1}(p^{-n}x_E)}$ equals~$p^n$ if~$x_E < \frac{p}{p + 1}$ and~$(p + 1)p^{n - 1}$ if~$x_E = \frac{p}{p + 1}$.
  In all cases,
  \begin{displaymath}
    (\t^n)_* (\t^n|_{\kproj^{-1}(p^{-n}x_E)})^*
  \end{displaymath}
  equals~$\delta$ times the identity on~$\kproj^{-1}(x_E)$.
  Thus, \eqref{eq:104} and Lemma~\ref{l:n-th-pull-back} yield
  \begin{multline*}
    \delta \mu_{\coset_0}^E
    =
    (\t^n)_* \left( (\t^n|_{\kproj^{-1}(p^{-n}x_E)})^* \mu_{\coset_0}^E \right)
    =
    (\t^n)_* \left( \sigma_1(p^n) \mu_{\coset}^E|_{\kproj^{-1}(p^{-n}x_E)} \right)
    \\ =
    (\t^n)_* \left( \sigma_1(p^n) \mu_{\coset'}^E|_{\kproj^{-1}(p^{-n}x_E)} \right)
    =
    (\t^n)_* \left( (\t^n|_{\kproj^{-1}(p^{-n}x_E)})^* \mu_{\coset_0'}^E \right)
    =
    \delta \mu_{\coset_0'}^E.
  \end{multline*}
  In particular, $\overline{\Orb_{\coset_0}(E)} = \overline{\Orb_{\coset_0'}(E)}$ by Theorem~\ref{t:Hecke-orbits} and hence ${\coset_0 = \coset_0'}$ by item~$(i)$.
  Since ${n = n'}$, this implies ${\coset = \coset'}$.
\end{proof}

\section{Equidistribution of \CM{} points along a $p$-adic discriminant}
\label{s:CM}

This section proves Theorems~\ref{t:CM-symmetric} and~\ref{t:CM-broken-symmetry}.
For each fundamental $p$\nobreakdash-adic discriminant~$\pfd$, the proof begins in Section~\ref{ss:CM-orbits} by showing how~$\Lambda_{\pfd}$ partitions into closures of partial Hecke orbits (Proposition~\ref{p:CM-orbits}).
The set~$\Lambda_{\pfd}$ coincides with a partial Hecke orbit if~$\Qpfd$ is unramified over~$\Qp$.
If~$\Qpfd$ is ramified over~$\Qp$, then~$\Lambda_{\pfd}$ partitions into precisely two closures of partial Hecke orbits.
In the latter case, Section~\ref{ss:symmetry-breaking} relies on genus theory to determine, for each discriminant~$D$ in~$\pfd$, how~$\supp(\Lambda_D)$ distributes between these parts (Proposition~\ref{p:symmetry-breaking}).
Section~\ref{s:proof-of-CM} begins deducing Theorems~\ref{t:CM-symmetric} and~\ref{t:CM-broken-symmetry} in the case of fundamental $p$\nobreakdash-adic discriminants from Theorems~\ref{t:CM-fundamental} and~\ref{t:Hecke-orbits-pr} (Propositions~\ref{p:CM-prime-to-p-I} and~\ref{p:CM-prime-to-p-II}).
It then relies on the (formal) \CM{} points formulae in Sections~\ref{ss:CM-from-canonical} and~\ref{ss:formal-CM-formulae} to deduce the general case.

\subsection{Hecke orbits of formal CM{} points}
\label{ss:CM-orbits}

For each fundamental $p$\nobreakdash-adic discriminant~$\pfd$, this section proves the following proposition describing~$\Lambda_{\pfd}$ in terms of closures of partial Hecke orbits.
Put
\begin{displaymath}
  \Npfd
  \=
  \left\{ \nr(g) \: g \in \cO_{\Qpfd}^{\times} \right\}.
\end{displaymath}

\begin{proposition}
  \label{p:CM-orbits}
  Let~$\pfd$ be a fundamental $p$\nobreakdash-adic discriminant.
  Then, for every~$E$ in~$\Lambda_{\pfd}$, the equality ${\NE = \Npfd}$ and the following properties hold.
  \begin{enumerate}
  \item [$(i)$]
    If~$\Qpfd$ is unramified over~$\Qp$, then
    \begin{equation}
      \label{eq:105}
      \Npfd
      =
      \Zp^{\times},
      \Lambda_{\pfd}
      =
      \overline{\Orb_{\Npfd}(E)},
      \text{ and }
      \nu_{\pfd}
      =
      \mu_{\Npfd}^E.
    \end{equation}
  \item [$(ii)$]
    If~$\Qpfd$ is ramified over~$\Qp$, then~$\Npfd$ has index two in~$\Zp^{\times}$,
    \begin{equation}
      \label{eq:106}
      \Lambda_{\pfd}
      =
      \overline{\Orb_{\Npfd}(E)} \sqcup \overline{\Orb_{\Zp^{\times} \ssetminus \Npfd}(E)},
      \text{ and }
      \nu_{\pfd}
      =
      \frac{1}{2} \left( \mu_{\Npfd}^E + \mu_{\Zp^{\times} \ssetminus \Npfd}^E \right).
    \end{equation}
    In particular,
    \begin{displaymath}
      \nu_{\pfd} \left ( \overline{\Orb_{\Npfd}(E)} \right)
      =
      \nu_{\pfd} \left( \overline{\Orb_{\Zp^{\times} \ssetminus \Npfd}(E)} \right)
      =
      \frac{1}{2},
    \end{displaymath}
    \begin{displaymath}
      \mu_{\Npfd}^E
      =
      2 \nu_{\pfd}|_{\overline{\Orb_{\Npfd}(E)}},
      \text{ and }
      \mu_{\Zp^{\times} \ssetminus \Npfd}^E
      =
      2 \nu_{\pfd}|_{\overline{\Orb_{\Zp^{\times} \ssetminus \Npfd}(E)}}.
    \end{displaymath}
  \end{enumerate}
\end{proposition}

The proof of this proposition relies on a couple of lemmas.

\begin{lemma}
  \label{l:semi-transitivity}
  Fix a fundamental $p$\nobreakdash-adic discriminant~$\pfd$ and~$\Delta$ in~$\pfd$.
  Then, for all~$\ss$ and~$\ss'$ in~$\tSups$, $\varphi$ in~$S_{- \Delta}^0(\ss)$, and~$\varphi'$ in~$S_{- \Delta}^0(\ss')$, there is~$g$ in ${S_1(\ss, \ss') \cup S_{-1}(\ss, \ss')}$ such that
  \begin{displaymath}
    g \varphi g^{-1} = \varphi'
    \text{ or }
    g \varphi g^{-1} = - \varphi'.
  \end{displaymath}
\end{lemma}

\begin{proof}
  Fix~$g_0$ in~$\bfG_{\ss, \ss'}$, so~$g_0^{-1} \varphi' g_0$ belongs to~$S_{- \Delta}^0(\ss)$ by Proposition~\ref{p:trace-zero-spheres}$(iv)$.
  By Proposition~\ref{p:trace-zero-spheres}$(ii)$, there is~$\rho$ in~$\Gss$ such that ${\rho^{-1} g_0^{-1} \varphi' g_0 \rho = \varphi}$.
  Suppose that~$Q_{\ss, \ss'}(g_0 \rho)$ (resp.~$- Q_{\ss, \ss'}(g_0 \rho)$) belongs to~$\Npfd$, let~$\psi$ in~$\Qp(\varphi)$ be such that~$\nr(\psi) = Q_{\ss, \ss'}(g_0 \rho)$ (resp. ${\nr(\psi) = - Q_{\ss, \ss'}(g_0 \rho)}$), and put ${g \= g_0 \rho \psi^{-1}}$.
  Then, ${g}$ belongs to~$S_1(\ss, \ss')$ (resp.~$S_{-1}(\ss, \ss')$) and
  \begin{displaymath}
    g \varphi g^{-1}
    =
    (g_0 \rho) \psi^{-1} \varphi \psi (g_0 \rho)^{-1}
    =
    (g_0 \rho) \varphi (g_0 \rho)^{-1}
    =
    \varphi'.
  \end{displaymath}
  It remains to consider the case where neither~$Q_{\ss, \ss'}(g_0\rho)$ nor~$- Q_{\ss, \ss'}(g_0\rho)$ belongs to~$\Npfd$.
  In this case, there is~$\gamma$ in~$\Gss$ such that
  \begin{displaymath}
    \gamma \varphi \gamma^{-1}
    =
    \overline{\varphi}
    \text{ and }
    \gamma^2
    =
    - Q_{\ss, \ss'}(g_0 \rho)^{-1},
  \end{displaymath}
  see Lemma~\ref{l:skew-generator}.
  Put ${g \= g_0 \rho \gamma}$.
  Then, ${g}$ belongs to~$S_{-1}(\ss, \ss')$ and
  \begin{displaymath}
    g \varphi g^{-1}
    =
    (g_0 \rho) \gamma \varphi \gamma^{-1} (g_0 \rho)^{-1}
    =
    (g_0 \rho) \overline{\varphi} (g_0 \rho)^{-1}
    =
    - (g_0 \rho) \varphi (g_0 \rho)^{-1}
    =
    - \varphi'.
    \qedhere
  \end{displaymath}
\end{proof}

\begin{lemma}
  \label{l:fixed-transposition}
  Let~$\pfd$ be a fundamental $p$\nobreakdash-adic discriminant such that~$\Qpfd$ is ramified over~$\Qp$, $\ss$ in~$\tSups$, $\varphi$ in~$\bfL_{\ss, \pfd}$, and recall that~$\Fixss(\varphi)$ has precisely two elements.
  Then, there is~$g$ in~$\Gss$ mapping~$\Fixss(\varphi)$ to itself, interchanging its elements.
  Moreover, for every such~$g$, the reduced norm~$\nr(g)$ belongs to~$\Npfd$ if and only if~$-1$ is outside~$\Npfd$.
\end{lemma}

\begin{proof}
  To prove the first assertion, recall that~$\Npfd$ has index two in~$\Zp^{\times}$ by Lemma~\ref{l:integer-norms}, so there is~$\gamma$ in~$\Gss$ satisfying
  \begin{displaymath}
    \gamma \varphi \gamma^{-1}
    =
    \overline{\varphi}
    \text{ and }
    \gamma^2 \in \Zp^{\times} \ssetminus \Npfd,
  \end{displaymath}
  see Lemma~\ref{l:skew-generator}.
  Thus, Lemma~\ref{l:fixed-properties}$(iv)$ yields
  \begin{displaymath}
    \gamma \cdot \Fixss(\varphi)
    =
    \Fixss(\overline{\varphi})
    =
    \Fixss(\varphi).
  \end{displaymath}
  Moreover, $\gamma$ is outside~$\Qp(\varphi)$, so~$\gamma$ cannot have a fixed point in~$\Fixss(\varphi)$ by Lemma~\ref{l:fixed-properties}$(iv)$.
  Since~$\Fixss(\varphi)$ has only two elements, $\gamma$ interchanges them.

  To prove the second assertion, let~$g$ in~$\Gss$ be such that~$g \cdot \Fixss(\varphi) = \Fixss(\varphi)$ and~$g$ interchanges the elements of~$\Fixss(\varphi)$.
  Then, $\gamma g$ fixes each element of~$\Fixss(\varphi)$, so it belongs to~$\Qp(\varphi)$ by Lemma~\ref{l:fixed-properties}$(iv)$.
  In particular,
  \begin{displaymath}
    \nr(\gamma g)
    =
    \nr(\gamma) \nr(g)
    =
    - \gamma^2 \nr(g)
    \in
    \Npfd.
  \end{displaymath}
  It follows that~$\nr(g)$ belongs to~$\Npfd$ if and only if~$-1$ is outside~$\Npfd$.
\end{proof}

\begin{proof}[Proof of Proposition~\ref{p:CM-orbits}]
  Let~$\ss$ in~$\tSups$ be such that~$E$ belongs to~$\Dss$ and fix~$x$ in~$\Piss^{-1}(E)$ and~$\Delta$ in~$\pfd$.
  By Proposition~\ref{p:trace-zero-spheres-to-CM}$(i)$, there is~$\varphi$ in~$S_{- \Delta}^0(\ss)$ such that~$x$ belongs to~$\Fixss(\varphi)$.

  The $p$\nobreakdash-adic quadratic orders~$\End(\FE)$ and~$\cO_{\Qpfd}$ are isomorphic by the definition of~$\Lambda_{\pfd}$ and the fact that the $p$\nobreakdash-adic discriminant is a complete isomorphism invariant for $p$\nobreakdash-adic quadratic orders (Lemma~\ref{l:p-adic-discriminants-Appendix}$(ii)$).
  Thus, $\NE = \Npfd$.

  To prove items~$(i)$ and~$(ii)$, let~$\ss'$ be in~$\tSups$ and ${\hF \: \hDsspr \to \R}$ a continuous function.
  For all~$u$ in~$\Zp^{\times}$ and~$g$ in~$\Gss$, Proposition~\ref{p:supersingular-spheres}$(iv)$ and the change of variables formula yield
  \begin{equation}
    \label{eq:int}
    \begin{split}
      \int \Trsspr(\hF)(\rho (g \varphi g^{-1})\rho^{-1}) \dd \mu_{u}^{\ss,\ss'}(\rho)
      & =
        \int \Trsspr(\hF)(\rho g \varphi (\rho g)^{-1}) \dd \mu_{u}^{\ss,\ss'}(\rho)
      \\ & =
           \int \Trsspr(\hF)(\hrho \varphi \hrho^{-1}) \dd \mu_{u \nr(g)}^{\ss,\ss'}(\hrho).
    \end{split}
  \end{equation}
  Together with Lemmas~\ref{l:fixed-properties}$(iv)$ and~\ref{l:semi-transitivity}, this implies, for every~$\varphi'$ in~$S_{- \Delta}^0(\ss')$,
  \begin{displaymath}
    \int \Trsspr(\hF)(\rho \varphi \rho^{-1}) \dd \left( \mu_{1}^{\ss,\ss'} + \mu_{-1}^{\ss,\ss'} \right) (\rho)
    =
    \int \Trsspr(\hF)(\rho \varphi' \rho^{-1}) \dd \left( \mu_{1}^{\ss,\ss'} + \mu_{-1}^{\ss,\ss'} \right) (\rho).
  \end{displaymath}
  Together with Propositions~\ref{p:trace-zero-spheres}$(iv)$ and~\ref{p:trace-zero-spheres-to-CM}$(ii)$ and the change of variables formula, this implies
  \begin{multline}
    \label{eq:107}
    \int \Trsspr(\hF)(\rho \varphi\rho^{-1}) \dd \left( \mu_{1}^{\ss,\ss'} + \mu_{-1}^{\ss,\ss'} \right) (\rho)
    \\
    \begin{aligned}
      & =
        \int \int \Trsspr(\hF)(\rho \varphi' \rho^{-1}) \dd \left( \mu_{1}^{\ss,\ss'} + \mu_{-1}^{\ss,\ss'} \right) (\rho) \dd \nu^{\ss}_{- \Delta}(\varphi')
      \\ & =
           \int \int \Trsspr(\hF)(\rho \varphi' \rho^{-1}) \dd \nu^{\ss}_{- \Delta}(\varphi') \dd \left( \mu_{1}^{\ss,\ss'} + \mu_{-1}^{\ss,\ss'} \right) (\rho)
      \\ & =
           \int \int \Trsspr(\hF)( \breve{\varphi}) \dd \nu^{\ss'}_{- \Delta}(\breve{\varphi}) \dd \left( \mu_{1}^{\ss,\ss'} + \mu_{-1}^{\ss,\ss'} \right) (\rho)
      \\ & =
           2 \int \Trsspr(\hF) \dd \nu^{\ss'}_{- \Delta}
      \\ & =
           2 \int \hF \dd \hnu^{\ss'}_{\pfd}.
    \end{aligned}
  \end{multline}

  If~$\Qpfd$ is unramified over~$\Qp$, then $\Npfd = \Zp^{\times}$ by Lemma~\ref{l:integer-norms}.
  Propositions~\ref{p:CM-parametrization}$(i)$ and~\ref{p:supersingular-sphere-to-orbit}$(ii)$ and the change of variables formula yield, for each~$u$ in~$\{ 1, -1 \}$,
  \begin{equation}
    \label{eq:108}
    \begin{split}
      \int \Trsspr ( \hF )(\rho \varphi \rho^{-1}) \dd \mu_{u}^{\ss,\ss'}(\rho)
      & =
        \int \hF(x_{\ss, \Delta}(\rho \varphi \rho^{-1})) \dd \mu_{u}^{\ss,\ss'}(\rho)
      \\ & =
           \int \hF \left( \Ev^{x, \ss'}(\rho) \right) \dd \mu_{u}^{\ss,\ss'}(\rho)
      \\ & =
           \int \hF \dd \hmu_{\Npfd}^{E, \ss'}.
    \end{split}
  \end{equation}
  Together with~\eqref{eq:107}, this implies ${\hnu^{\ss'}_{\pfd} = \hmu_{\Npfd}^{E, \ss'}}$.
  Since this holds for every~$\ss'$ in~$\tSups$, the equality ${\nu_{\pfd} = \mu_{\Npfd}^E}$ follows.
  The second equality in~\eqref{eq:105} follows from a comparison of the supports of these measures, using Theorems~\ref{t:CM-fundamental} and~\ref{t:Hecke-orbits}.

  Suppose~$\Qpfd$ is ramified over~$\Qp$, so~$\Npfd$ has index two in~$\Zp^{\times}$ by Lemma~\ref{l:integer-norms} and~$\Fixss(\varphi)$ has precisely two elements by Lemma~\ref{l:fixed-properties}$(ii)$.
  Denote by~$\breve{x}$ the element of~$\Fixss(\varphi)$ that is different from~$x$ and put~$\breve{E} \= \Piss(\breve{x})$.
  Then, Propositions~\ref{p:CM-parametrization}$(ii)$ and~\ref{p:supersingular-sphere-to-orbit}$(ii)$ and the change of variables formula yield
  \begin{multline*}
    \int \Trsspr ( \hF )(\rho \varphi \rho^{-1}) \dd \left( \mu_{1}^{\ss,\ss'} + \mu_{-1}^{\ss,\ss'} \right) (\rho)
    \\
    \begin{aligned}
      & = \frac{1}{2} \int \hF(x_{\ss, \Delta}^+(\rho \varphi \rho^{-1})) + \hF(x_{\ss, \Delta}^-(\rho \varphi \rho^{-1})) \dd \left( \mu_{1}^{\ss,\ss'} + \mu_{-1}^{\ss,\ss'} \right) (\rho)
      \\ & =
           \frac{1}{2} \int \hF \left( \Ev^{x, \ss'}(\rho) \right) + \hF \left( \Ev^{\breve{x}, \ss'}(\rho) \right) \dd \left( \mu_{1}^{\ss,\ss'} + \mu_{-1}^{\ss,\ss'} \right) (\rho)
      \\ & =
           \frac{1}{2} \int \hF \dd \left( \hmu_{\Npfd}^{E, \ss'} + \hmu_{\Npfd}^{\breve{E}, \ss'} + \hmu_{- \Npfd}^{E, \ss'} + \hmu_{- \Npfd}^{\breve{E}, \ss'} \right).
    \end{aligned}
  \end{multline*}
  Since~\eqref{eq:107} and the previous formula hold for every~$\ss'$ in~$\tSups$,
  \begin{equation}
    \label{eq:109}
    \nu_{\pfd}
    =
    \frac{1}{4} \left( \mu_{\Npfd}^{E} + \mu_{\Npfd}^{\breve{E}} + \mu_{- \Npfd}^{E} + \mu_{- \Npfd}^{\breve{E}} \right).
  \end{equation}
  Moreover, Proposition~\ref{p:supersingular-sphere-to-orbit}$(i)$ and Lemma~\ref{l:fixed-transposition} imply
  \begin{displaymath}
    \breve{E}
    \in
    \begin{cases}
      \overline{\Orb_{\Zp^{\times} \ssetminus \Npfd}(E)}
      & \text{if ${-1 \in \Npfd}$};
      \\
      \overline{\Orb_{\Npfd}(E)}
      & \text{if ${-1 \not\in \Npfd}$}.
    \end{cases}
  \end{displaymath}
  So, Lemma~\ref{l:orbit-closure} yields
  \begin{displaymath}
    \mu_{- \Npfd}^{E}
    =
    \mu_{\Npfd}^{E}
    \text{ and }
    \mu_{- \Npfd}^{\breve{E}}
    =
    \mu_{\Npfd}^{\breve{E}}
    =
    \mu_{\Zp^{\times} \ssetminus \Npfd}^{E},
  \end{displaymath}
  if~$-1$ belongs to~$\Npfd$.
  If~$-1$ is outside~$\Npfd$, then
  \begin{displaymath}
    \mu_{\Npfd}^{\breve{E}}
    =
    \mu_{\Npfd}^{E}
    \text{ and }
    \mu_{- \Npfd}^{E}
    =
    \mu_{- \Npfd}^{\breve{E}}
    =
    \mu_{\Zp^{\times} \ssetminus \Npfd}^{E}.
  \end{displaymath}
  In all cases, \eqref{eq:109} yields the equality of measures in~\eqref{eq:106}.
  Proposition~\ref{p:orbit-measures}$(i)$ implies that the closures of the partial Hecke orbits in the first equality of~\eqref{eq:106} are disjoint.
  Then, the equality of sets in~\eqref{eq:106} and the remaining assertions of item~$(ii)$ follow from a comparison of the supports of the measures~$\nu_{\pfd}$, $\mu_{\Npfd}^{E}$, and~$\mu_{\Zp^{\times} \ssetminus \Npfd}^{E}$, using Theorems~\ref{t:CM-fundamental} and~\ref{t:Hecke-orbits}.
\end{proof}

\subsection{Symmetry breaking}
\label{ss:symmetry-breaking}
Fix a fundamental $p$\nobreakdash-adic discriminant~$\pfd$ for which~$\Qpfd$ is ramified over~$\Qp$ and recall that~$\Npfd$ has index two in~$\Zp^{\times}$ (Lemma~\ref{l:integer-norms}).
Choose a point~$E_{\pfd}$ in~$\Lambda_{\pfd}$, as follows.
Suppose~$\pfd$ contains a prime discriminant~$d$ that is divisible by~$p$.
Then, $d$ is the unique fundamental discriminant in~$\pfd$ with this property and choose an arbitrary~$E_{\pfd}$ in~$\supp(\Lambda_d)$.
If~$\pfd$ does not contain a prime discriminant divisible by~$p$, then choose an arbitrary~$E_{\pfd}$ in~$\Lambda_{\pfd}$.
In all cases, put
\begin{displaymath}
  \Lambda_{\pfd}^+
  \=
  \overline{\Orb_{\Npfd}(E_{\pfd})}
  \text{ and }
  \Lambda_{\pfd}^-
  \=
  \overline{\Orb_{\Zp^{\times} \ssetminus \Npfd}(E_{\pfd})}.
\end{displaymath}
Proposition~\ref{p:CM-orbits}$(ii)$ yields the partition
\begin{equation}
  \label{eq:110}
  \Lambda_{\pfd}
  =
  \Lambda_{\pfd}^+ \sqcup \Lambda_{\pfd}^-.
\end{equation}

This section proves the following proposition describing, for each discriminant~$D$ in~$\pfd$, how~$\Lambda_D$ distributes between~$\Lambda_{\pfd}^+$ and~$\Lambda_{\pfd}^-$.
Define the divisors~$\Lambda_D^+$ and~$\Lambda_D^-$, by
\begin{displaymath}
  \Lambda_D^+
  \=
  \Lambda_D|_{\Lambda_{\pfd}^+}
  \text{ and }
  \Lambda_D^-
  \=
  \Lambda_D|_{\Lambda_{\pfd}^-},
\end{displaymath}
so ${\Lambda_D = \Lambda_D^+ + \Lambda_D^-}$.
Recall that~$\left( \frac{ \cdot }{ \cdot } \right)$ denotes the Kronecker symbol.

\begin{proposition}
  \label{p:symmetry-breaking}
  Let~$d$ be a fundamental discriminant divisible by~$p$.
  Then, for every~$f$ in ${\Z_{>0} \ssetminus p\Z}$, the following properties hold.
  \begin{enumerate}
  \item [$(i)$]
    If~$d$ is not a prime discriminant, then
    \begin{displaymath}
      \deg(\Lambda_{d f^2}^+)
      =
      \deg(\Lambda_{d f^2}^-).
    \end{displaymath}
  \item [$(ii)$]
    If~$d$ is a prime discriminant, then
    \begin{displaymath}
      \Lambda_{d f^2}^{\pm}
      =
      \begin{cases}
        \Lambda_{d f^2}
        & \text{if~$\left( \frac{d}{f} \right) = \pm 1$};
        \\
        0
        & \text{if~$\left( \frac{d}{f} \right) = \mp 1$}.
      \end{cases}
    \end{displaymath}
  \end{enumerate}
\end{proposition}

The following corollary is a direct consequence of Corollary~\ref{c:Hecke-on-orbits} and Proposition~\ref{p:CM-orbits}.
Put
\begin{displaymath}
  \nu_{\pfd}^+
  \=
  \mu_{\Npfd}^{E_{\pfd}}
  \text{ and }
  \nu_{\pfd}^-
  \=
  \mu_{\Zp^{\times} \ssetminus \Npfd}^{E_{\pfd}}.
\end{displaymath}
Since ${\supp(\nu_{\pfd}^+) = \Lambda_{\pfd}^+}$ and ${\supp(\nu_{\pfd}^-) = \Lambda_{\pfd}^-}$ by Theorem~\ref{t:Hecke-orbits} and Proposition~\ref{p:CM-orbits}$(ii)$,
\begin{displaymath}
  \nu_{\pfd}
  =
  \frac{1}{2} \left( \nu_{\pfd}^+ + \nu_{\pfd}^- \right),
  \nu_{\pfd}(\Lambda_{\pfd}^+)
  =
  \nu_{\pfd}(\Lambda_{\pfd}^-)
  =
  \frac{1}{2},
  \nu_{\pfd}^+
  =
  2 \nu_{\pfd}|_{\Lambda_{\pfd}^+},
  \text{ and }
  \nu_{\pfd}^-
  =
  2 \nu_{\pfd}|_{\Lambda_{\pfd}^-}.
\end{displaymath}

\begin{coro}
  \label{c:Hecke-on-CM}
  For all fundamental $p$\nobreakdash-adic discriminant~$\pfd$ and~$n$ in ${\Z_{>0} \ssetminus p\Z}$,
  \begin{displaymath}
    T_n \left( \Lambda_{\pfd} \right)
    =
    \Lambda_{\pfd}
    \text{ and }
    \frac{1}{\sigma_1(n)} (T_n)_*(\nu_{\pfd})
    =
    \nu_{\pfd}.
  \end{displaymath}
  If in addition~$\Qpfd$ is ramified over~$\Qp$, then
  \begin{displaymath}
    T_n(\Lambda_{\pfd}^{\pm})
    =
    \begin{cases}
      \Lambda_{\pfd}^{\pm}
      & \text{if } n \in \Npfd;
      \\
      \Lambda_{\pfd}^{\mp}
      & \text{if } n \not\in \Npfd,
    \end{cases}
    \text{ and }
    \frac{1}{\sigma_1(n)} (T_n)_*(\nu_{\pfd}^{\pm})
    =
    \begin{cases}
      \nu_{\pfd}^{\pm}
      & \text{if } n \in \Npfd;
      \\
      \nu_{\pfd}^{\mp}
      & \text{if } n \not\in \Npfd.
    \end{cases}
  \end{displaymath}
\end{coro}

The proof of Proposition~\ref{p:symmetry-breaking} relies on the following lemma.
A \emph{quadratic fundamental discriminant} is the discriminant of the ring of integers of a quadratic (real or imaginary) extension of~$\Q$.
So, a quadratic fundamental discriminant is a fundamental discriminant if and only if it is negative.
A quadratic discriminant is \emph{prime}, if it is fundamental and divisible by only one prime number.
Every quadratic fundamental discriminant can be written uniquely up to order as a product of prime quadratic discriminants that are mutually coprime, see, \emph{e.g.}, \cite[Proposition~2.2]{Lem00}.
A quadratic fundamental discriminant~$d$ divisible by~$p$ is prime if and only if~$p$ is odd and ${d =(-1)^{\frac{p-1}{2}}p}$, or if~$p = 2$ and~$d$ belongs to~$\{4,-8, 8\}$.
For all nonzero integers~$m$ and~$n$, denote by $(n, m)_p$ the Hilbert symbol over~$\Qp$, see, \emph{e.g.}, \cite[Chapter~III]{Ser73e} or \cite[Section~2.5]{Lem00}.

\begin{lemma}
  \label{l:symmetry-breaking-discriminants}
  Let~$d$ be a fundamental discriminant divisible by~$p$ and let~$p^*$ be the unique prime quadratic discriminant divisible by~$p$ in the factorization of~$d$ into prime quadratic discriminants.
  Then, the following properties hold.
  \begin{enumerate}
  \item [$(i)$]
    For every~$n$ in~$\Z_{>0}$ coprime to~$d$, the equality ${(n, d)_p = \left( \frac{p^*}{n} \right)}$ holds.
  \item[$(ii)$]
    If~$d \neq p^*$, then there is a prime number~$q$ satisfying
    \begin{displaymath}
      (q, d)_p = -1
      \text{ and }
      \left( \frac{d}{q} \right) = 1.
    \end{displaymath}
  \end{enumerate}
\end{lemma}

The proofs of Lemma~\ref{l:symmetry-breaking-discriminants} and Proposition~\ref{p:symmetry-breaking} rely on several properties of the Hilbert symbol that can be found, \emph{e.g.}, in \cite[Theorems~1 and~2, Chapter~III]{Ser73e}.
They also rely on the following notation.
For each quadratic field extension~$K$ of~$\Q$, denote by~$\Cl(K)$ the ideal class group of~$K$ and, for each fractional ideal~$\fra$ of~$K$, by~$[\fra]$ its class in~$\Cl(K)$ and by~$\Nr(\fra)$ its norm.

\begin{proof}[Proof of Lemma~\ref{l:symmetry-breaking-discriminants}]
  Put~$d' \= \frac{d}{p^*}$, so~$d'$ is a quadratic fundamental discriminant.
  Since~$(\cdot, d)_p$ and~$\left( \frac{p^*}{\cdot} \right)$ are both completely multiplicative, it suffices to prove item~$(i)$ when~$n$ is a prime number~$q$ not dividing~$d$.
  Then,
  \begin{displaymath}
    (q, d)_p
    =
    \left( q, d' \right)_p (q, p^*)_p
    =
    (q, p^*)_p
    =
    \begin{cases}
      \left(\frac{q}{p}\right)
      & \text{if $p$ is odd};
      \\
      (-1)^{\frac{q-1}{2}}
      & \text{if }p=2 \text{ and }p^*=-4;
      \\
      (-1)^{\frac{q-1}{2}+\frac{q^2-1}{8}}
      & \text{if }p=2 \text{ and }p^*=-8;
      \\
      (-1)^{\frac{q^2-1}{8}}
      & \text{if }p=2 \text{ and }p^*=8.
    \end{cases}
  \end{displaymath}
  In all cases, ${(q, d)_p = \left( \frac{p^*}{q} \right)}$ by the quadratic reciprocity law and its complementary laws, see, \emph{e.g.}, \cite[Theorems~5 and~6, Chapter~I]{Ser73e}.

  Put~$K\=\Q(\sqrt{d})$ and let ${\chi \: \Cl(K)\to \{1,-1\}}$ be the unique quadratic character such that, for every prime ideal~$\frp$ of~$\cO_K$,
  \begin{displaymath}
    \chi([\frp])
    =
    \begin{cases}
      \left( \frac{p^*}{\Nr(\frp)}\right)
      & \text{if } \gcd(\Nr(\frp),p^*)=1;
      \\
      \left(\frac{d'}{\Nr(\frp)}\right)
      & \text{if } \gcd(\Nr(\frp),d')=1,
    \end{cases}
  \end{displaymath}
  see, \emph{e.g.}, \cite[Section~2.3]{Lem00}.
  It follows from genus theory that there exists an ideal class~$[\fra]$ in~$\Cl(K)$ such that~$\chi([\fra])=-1$, see, \emph{e.g.}, \cite[Theorem~2.17]{Lem00}.
  Let~$\frb$ be an ideal of~$\cO_K$ in~$[\fra]$ whose norm is coprime to~$d$.
  Decomposing~$\frb$ into primes ideals yields a prime ideal~$\frq$ of~$\cO_K$ such that ${\chi([\frq])=-1}$.
  Then, $\Nr(\frq)$ is coprime to~$d$ and item~$(i)$ yields
  \begin{displaymath}
    (\Nr(\frq), d)_p
    =
    \left(\frac{p^*}{\Nr(\frq)}\right)
    =
    \chi([\frq])=-1.
  \end{displaymath}
  Put ${q \= \Nr(\frq)}$.
  Then, ${q}$ is a prime number.
  Since~$q$ does not divide~$d$, the equality ${\left( \frac{d}{q} \right) = 1}$ holds.
\end{proof}

\begin{proof}[Proof of Proposition~\ref{p:symmetry-breaking}]
  Let~$\pfd$ be the fundamental $p$\nobreakdash-adic discriminant containing~$d$.
  The proof uses several times that, if~$n$ in~$\Z_{>0}$ is coprime to~$d$, then ${(n, d)_p = 1}$ holds if and only if~$n$ belongs to~$\Npfd$, see, \emph{e.g.}, \cite[Proposition~1, Chapter~III]{Ser73e}.
  Put ${K \= \Q(\sqrt{d})}$, and recall from Section~\ref{ss:CM-from-canonical} that ${R_d = \bfone \ast \psi_d}$.
  Fix a field isomorphism between~$\Cp$ and~$\C$, and, for each~$E$ in~$\Ell(\Cp)$, denote by~$E \otimes \C$ the element of~$\Ell(\C)$ obtained from~$E$ by base change to~$\C$.
  Moreover, denote by~$\cE \: \Cl(K) \to \supp(\Lambda_d)$ the bijection so that, for each fractional ideal~$\fra$ of~$K$, the quotient~$\C / \fra$ is isomorphic to~$(\cE([\fra]) \otimes \C)(\C)$, see, \emph{e.g.}, \cite[Chapter~II, Section~1]{Sil94a}.

  The proof begins establishing the case where ${f = 1}$.
  To prove item~$(i)$ when~$f = 1$, let~$q$ be a prime number such that~$(q, d)_p = -1$ and~$\left( \frac{d}{q} \right) = 1$ (Lemma~\ref{l:symmetry-breaking-discriminants}(ii)).
  In particular, ${q \neq p}$, and~$q$ splits in~$K$ and is outside~$\Npfd$.
  If follows that there is an ideal~$\frq$ of~$\cO_K$ of norm~$q$ such that~$\frq\overline{\frq}=q\cO_K$.
  Thus, the map~$\fra \mapsto \fra \frq$ induces a bijection~$\iota$ of~$\supp(\Lambda_d)$ given by~$\cE([\fra]) \mapsto \cE([\fra \frq])$, whose inverse is given by~$\cE([\fra]) \mapsto \cE([\fra \overline{\frq}])$.
  Since for every fractional ideal~$\fra$ of~$K$, each of the natural maps ${\C / \fra \frq \to \C / \fra}$ and ${\C / \fra \overline{\frq} \to \C / \fra}$ is an isogeny of degree~$q$, the involution~$\iota$ interchanges~$\supp(\Lambda_d^+)$ and~$\supp(\Lambda_d^-)$ by Corollary~\ref{c:Hecke-on-CM}.
  In particular, ${\deg(\Lambda_d^+) = \deg(\Lambda_d^-)}$.

  For item~$(ii)$ with~$f = 1$, the point~$E_{\pfd}$ defining~$\Lambda_{d}^+$ and~$\Lambda_{d}^-$ at the beginning of the section belongs to~$\supp(\Lambda_{d}^+)$ by definition.
  Let~$E$ be in~$\supp(\Lambda_d)$ and ${\phi \: E_{\pfd} \to E}$ an isogeny whose degree is not divisible by~$p$ \cite[Lemma~4.8]{HerMenRiv20}.
  Let~$\fra,\fra_0$ be ideals of $\cO_K$ such that~$\cE([\fra]) = E$, $\cE([\fra_0 \fra]) = E_{\pfd}$, and such that the natural map~$\C / \fra_0 \fra \to \C / \fra$ corresponds to the isogeny~$\phi$.
  Consider the prime factorization ${\fra_0 = \frq_1^{\alpha_1} \cdots \frq_n^{\alpha_n}}$.
  Then, for each~$j$ in~$\{1, \ldots, n\}$, the norm of~$\frq_j$ is either a prime number~$q_j$ and then ${\left( \frac{d}{q_j} \right) = 1}$, or the square of a prime number~$q_j'$ and then ${\left( \frac{d}{q_j'} \right) = -1}$.
  In all cases, ${\left( \frac{d}{\Nr(\fra_0)} \right) = 1}$.
  Thus, Lemma~\ref{l:symmetry-breaking-discriminants}$(i)$ yields
  \begin{displaymath}
    (\deg(\phi), d)_p
    =
    (\Nr(\fra_0), d)_p
    =
    \left( \frac{d}{\Nr(\fra_0)} \right)
    =
    1,
  \end{displaymath}
  hence~$\deg(\phi)$ belongs to~$\Npfd$, and~$E$ belongs to~$\Lambda_d^+$ by Corollary~\ref{c:Hecke-on-CM}.
  This proves item~$(ii)$ when~$f = 1$.

  It remains to consider the case where ${f \ge 2}$.
  In this case, \eqref{eq:74} holds.
  Moreover, Corollary~\ref{c:Hecke-on-CM} yields
  \begin{multline}
    \label{eq:111}
    \deg(\Lambda_{d f^2}^{\pm})
    \\ =
    \frac{\deg(\Lambda_{d}^{\pm})}{w_{d, 1}} \sum_{\substack{f_0 \in \Z_{>0}, f_0 \mid f \\ f_0 \in \Npfd}} R_d^{-1} \left( \frac{f}{f_0} \right) \sigma_1(f_0)
    +
    \frac{\deg(\Lambda_{d}^{\mp})}{w_{d, 1}} \sum_{\substack{f_0 \in \Z_{>0}, f_0 \mid f \\ f_0 \not\in \Npfd}} R_d^{-1} \left( \frac{f}{f_0} \right) \sigma_1(f_0).
  \end{multline}
  Combined with item~$(i)$ with~$f = 1$, this implies item~$(i)$ when ${f \ge 2}$.
  By Lemma~\ref{l:symmetry-breaking-discriminants}$(i)$ and~\eqref{eq:111}, to deduce item~$(ii)$ for~$f \ge 2$ from the case where ${f = 1}$ it suffices to show the following: For every~$r$ in~$\Z_{>0}$ that belongs to~$\Zp^{\times} \ssetminus \Npfd$, the equality ${R_d^{-1}(r) = 0}$ holds.
  Since the function~$R_d^{-1}$ is multiplicative and~$(\cdot, d)_p$ is completely multiplicative, it suffices to show that, for all prime number~$q_0$ different from~$p$ such that~$(q_0, d)_p = -1$ and odd integer~$s$ in~$\Z_{>0}$, the equality ${R_d^{-1}(q_0^s) = 0}$ holds.
  Noting that~$\psi_d(q_0) = \left( \frac{d}{q_0} \right) = -1$ by Lemma~\ref{l:symmetry-breaking-discriminants}$(i)$ and denoting the M{\"o}bius function by~$\mu$, this follows from a direct computation using the formula~$R_d^{-1} = \mu \ast (\mu \cdot \psi_d)$.
\end{proof}

\subsection{Proof of Theorems~\ref{t:CM-symmetric} and~\ref{t:CM-broken-symmetry}}
\label{s:proof-of-CM}
The proof relies on a couple of propositions.
For each fundamental $p$\nobreakdash-adic discriminant~$\pfd$ for which~$\Qpfd$ is ramified over~$\Qp$, let~$\Lambda_{\pfd}^+$, $\Lambda_{\pfd}^-$, $\nu_{\pfd}^+$, and~$\nu_{\pfd}^-$ be as in Section~\ref{ss:symmetry-breaking}.

\begin{proposition}
  \label{p:CM-prime-to-p-I}
  For all~$\varepsilon$ in~$\R_{>0}$ and locally constant function ${F \: \Sups \to \R}$, there is~$C$ in~$\R_{>0}$ such that the following property holds.
  Let~$\pfd$ be a fundamental $p$\nobreakdash-adic discriminant, $d$ a fundamental discriminant in~$\pfd$, and~$f$ in ${\Z_{>0} \ssetminus p\Z}$.
  Then,
  \begin{equation}
    \label{eq:112}
    \left| \int F \dd \odelta_{df^2} - \int F \dd \nu_{\pfd} \right|
    \le
    C f^{- \frac{1}{2} + \varepsilon}
  \end{equation}
  if~$\Qpfd$ is unramified over~$\Qp$, and if~$\Qpfd$ is ramified over~$\Qp$, then
  \begin{equation}
    \label{eq:113}
    \left| \int F \dd \odelta_{df^2} - \frac{\deg(\Lambda_{d f^2}^+)}{\deg(\Lambda_{d f^2})} \int F \dd \nu_{\pfd}^+ - \frac{\deg(\Lambda_{d f^2}^-)}{\deg(\Lambda_{d f^2})} \int F \dd \nu_{\pfd}^- \right|
    \le
    C f^{- \frac{1}{2} + \varepsilon}.
  \end{equation}
\end{proposition}

\begin{proof}
  Put~$\varepsilon' \= \frac{\varepsilon}{4}$ and let~$C_0$ and~$C_1$ be given by~\eqref{eq:7} and Lemma~\ref{l:R-inverse}, respectively, with~$\varepsilon$ replaced by~$\varepsilon'$.

  Suppose that~$\Qpfd$ is unramified over~$\Qp$, so~$\Npfd = \Zp^{\times}$.
  Fix~$E$ in~$\Lambda_{\pfd}$ and let~$C_2$ be given by Theorem~\ref{t:Hecke-orbits-pr} in Section~\ref{s:Hecke-orbits} with~$\varepsilon$ replaced by~$\varepsilon'$ and~$\coset = \Npfd$.
  Proposition~\ref{p:CM-orbits}$(i)$ yields ${\mu_{\Npfd}^E = \nu_{\pfd}}$.
  Thus, applying~\eqref{eq:deg-Hecke-op-divisor}, \eqref{eq:general-inverse-Zhang-formula} in Lemma~\ref{l:Zhang-general} with ${\wtf = 1}$, Theorem~\ref{t:Hecke-orbits-pr} to each element~$E'$ of~$\supp(\Lambda_d)$ and each divisor~$f_0$ of~$f$, and~\eqref{eq:74} yields
  \begin{multline}
    \label{eq:114}
    \left| \int F \dd \odelta_{df^2} - \int F \dd \nu_{\pfd} \right|
    \\
    \begin{aligned}
      & =
        \frac{w_{d, f}}{w_{d, 1} \deg(\Lambda_{df^2})} \left| \sum_{f_0 \in \Z_{>0}, f_0 \mid f} R_d^{-1} \left( \frac{f}{f_0} \right) \sigma_1(f_0) \sum_{E' \in \supp(\Lambda_{d})} \left( \int F \dd \odelta_{T_{f_0}(E')} - \int F \dd \mu_{\Npfd}^E \right) \right|
      \\ & \le
           C_2 \frac{w_{d, f} \deg(\Lambda_d)}{w_{d, 1} \deg(\Lambda_{df^2})} \sum_{f_0 \in \Z_{>0}, f_0 \mid f} \left| R_d^{-1} \left( \frac{f}{f_0} \right) \sigma_1(f_0) \right| f_0^{- \frac{1}{2} + \varepsilon'}.
    \end{aligned}
  \end{multline}
  Using~\eqref{eq:74} again and the choice of~$C_0$ and~$C_1$ yields
  \begin{multline}
    \label{eq:115}
    \frac{w_{d, f} \deg(\Lambda_d)}{w_{d, 1} \deg(\Lambda_{df^2})} \sum_{f_0 \in \Z_{>0}, f_0 \mid f} \left| R_d^{-1} \left( \frac{f}{f_0} \right) \sigma_1(f_0) \right| f_0^{- \frac{1}{2} + \varepsilon'}
    \\ \le
    C_0 C_1^2 f^{-1 + \varepsilon'} \sum_{f_0 \in \Z_{>0}, f_0 \mid f} \left( \frac{f}{f_0} \right)^{\varepsilon'} f_0^{\frac{1}{2} + 2\varepsilon'}
    \le
    C_0^2 C_1^2 f^{- \frac{1}{2} + 4\varepsilon'}.
  \end{multline}
  Together with~\eqref{eq:114}, this yields~\eqref{eq:112} with ${C = C_0^2 C_1^2 C_2}$.

  Suppose that~$\Qpfd$ is ramified over~$\Qp$ and recall that~$\Npfd$ has index two in~$\Zp^{\times}$.
  Fix~$E^+$ in~$\Lambda_{\pfd}^+$ and~$E^-$ in~$\Lambda_{\pfd}^-$ and let~$C_2'$ be the maximum value of the constant given by Theorem~\ref{t:Hecke-orbits-pr} with~$\varepsilon$ replaced by~$\varepsilon'$, $E$ in~$\{E^+, E^-\}$, and~$\coset$ in~$\{\Npfd, \Zp^{\times} \ssetminus \Npfd\}$.
  Applying~\eqref{eq:general-inverse-Zhang-formula} in Lemma~\ref{l:Zhang-general} with ${\wtf = 1}$ and Corollary~\ref{c:Hecke-on-CM} yields
  \begin{displaymath}
    \frac{\Lambda_{df^2}^+}{w_{d, f}}
    =
    \sum_{\substack{f_0 \in \Z_{>0}, f_0 \mid f \\ f_0 \in \Npfd}} R_d^{-1} \left( \frac{f}{f_0} \right) T_{f_0} \left( \frac{\Lambda_d^+}{w_{d, 1}} \right)
    + \sum_{\substack{f_0 \in \Z_{>0}, f_0 \mid f \\ f_0 \not \in \Npfd}} R_d^{-1} \left( \frac{f}{f_0} \right) T_{f_0} \left( \frac{\Lambda_d^-}{w_{d, 1}} \right).
  \end{displaymath}
  Moreover, Lemma~\ref{l:orbit-closure} implies ${\mu_{\Npfd}^{E^{\pm}} = \nu_{\pfd}^{\pm}}$ and ${\mu_{\Zp^{\times} \ssetminus \Npfd}^{E^{\pm}} = \nu_{\pfd}^{\mp}}$.
  Thus, applying~\eqref{eq:deg-Hecke-op-divisor}, Theorem~\ref{t:Hecke-orbits-pr} to each element~$E'$ of~$\supp(\Lambda_d)$ and each divisor~$f_0$ of~$f$, \eqref{eq:111}, and~\eqref{eq:115} yields
  \begin{multline*}
    \left| \int_{\Lambda_{\pfd}^+} F \dd \odelta_{df^2} - \frac{\deg(\Lambda_{d f^2}^+)}{\deg(\Lambda_{d f^2})} \int F \dd \nu_{\pfd}^+ \right|
    \\
    \begin{aligned}
      & =
        \frac{w_{d, f}}{w_{d, 1} \deg(\Lambda_{df^2})} \left| \sum_{\substack{f_0 \in \Z_{>0}, f_0 \mid f \\ f_0 \in \Npfd}} R_d^{-1} \left( \frac{f}{f_0} \right) \sigma_1(f_0) \sum_{E' \in \supp(\Lambda_{d}^+)} \left( \int F \dd \odelta_{T_{f_0}(E')} - \int F \dd \mu_{\Npfd}^{E^+} \right)
      \right. \\ & \quad + \left.
                   \sum_{\substack{f_0 \in \Z_{>0}, f_0 \mid f \\ f_0 \not\in \Npfd}} R_d^{-1} \left( \frac{f}{f_0} \right) \sigma_1(f_0) \sum_{E' \in \supp(\Lambda_{d}^-)} \left( \int F \dd \odelta_{T_{f_0}(E')} - \int F \dd \mu_{\Zp^{\times} \ssetminus \Npfd}^{E^-} \right) \right|
      \\ & \le
           C_2' \frac{w_{d, f} \deg(\Lambda_d)}{w_{d, 1} \deg(\Lambda_{df^2})} \sum_{f_0 \in \Z_{>0}, f_0 \mid f} \left| R_d^{-1} \left( \frac{f}{f_0} \right) \sigma_1(f_0) \right| f_0^{- \frac{1}{2} + \varepsilon'}
      \\ & \le
           C_0^2 C_1^2 C_2' f^{- \frac{1}{2} + 4\varepsilon'}.
    \end{aligned}
  \end{multline*}
  A similar argument shows that the same estimate holds with~$\Lambda_{\pfd}^+$, $\Lambda_d^+$, and~$\nu_{\pfd}^+$ replaced by~$\Lambda_{\pfd}^-$, $\Lambda_d^-$, and~$\nu_{\pfd}^-$, respectively.
  Combined, these estimates yield~\eqref{eq:113} with ${C = 2C_0^2 C_1^2 C_2'}$.
\end{proof}

\begin{proposition}
  \label{p:CM-prime-to-p-II}
  Let~$\pfd$ be a fundamental $p$\nobreakdash-adic discriminant.
  Then, for all~$\varepsilon$ and~$\delta$ in~$\R_{>0}$, there is~$C'$ in~$\R_{>0}$ such that the following property holds.
  For all function~$F \: \Lambda_{\pfd} \to \R$ that is constant on every ball of~$\Lambda_{\pfd}$ of radius~$\delta$, fundamental discriminant~$d$ in~$\pfd$, and~$f$ in ${\Z_{>0} \ssetminus p\Z}$,
  \begin{displaymath}
    \left| \int F \dd \odelta_{d f^2} - \int F \dd \nu_{\pfd} \right|
    \le
    C' \left( \sup_{\Lambda_\pfd} |F| \right) |d|^{- \frac{1}{28} + \varepsilon} f^{\varepsilon}.
  \end{displaymath}
\end{proposition}

The proof of this proposition relies on the following lemma.

\begin{lemma}
  \label{l:Hecke-uniformity}
  For all~$\delta$ in~$]0,1]$ and function ${F \: \Sups \to \R}$ that is constant on every ball of radius~$\delta$, the following property holds.
  For every~$n$ in ${\Z_{>0} \ssetminus p\Z}$, the function~$T_nF$ is constant on every ball of radius~$\delta^{12}$.
\end{lemma}

\begin{proof}
  Let~$\ss$ be in~$\tSups$ and recall from Section~\ref{ss:woods-hole} that ${\delta_{\ss} = \# \Aut(\ss) / 2 \le 12}$.
  By~\eqref{eq:14}, for each~$x$ in~$\hDss$,
  \begin{displaymath}
    T_n F \circ \Piss (x)
    =
    \sum_{\ss' \in \tSups} \frac{1}{\# \Aut(\ss')} \sum_{\phi \in \Hom_n(\ss, \ss')} F \circ \Pisspr (\hphi \cdot x).
  \end{displaymath}
  Since for each~$\ss'$ in~$\tSups$ the action of~$\Gsspr$ on~$\hDsspr$ is by isometries (Lemma~\ref{l:action-regularity}$(i)$), by~\eqref{eq:13} in Theorem~\ref{t:Pi-properties} the function~$T_n F \circ \Piss$ is constant on every ball of~$\hDss$ of radius~$\delta$.
  Using~$\delta_{\ss} \le 12$ and~\eqref{eq:13} in Theorem~\ref{t:Pi-properties} again, it follows that the function~$T_n F$ is constant on every ball of~$\Dss$ of radius~$\delta^{12}$.
  Since~$\ss$ in~$\tSups$ is arbitrary, this implies the lemma.
\end{proof}

\begin{proof}[Proof of Proposition~\ref{p:CM-prime-to-p-II}]
  Put ${\varepsilon' \= \frac{\varepsilon}{3}}$ and let~$C_0$ (resp.~$C_1$) be given by~\eqref{eq:7} (resp. Lemma~\ref{l:R-inverse}) with~$\varepsilon$ replaced by~$\varepsilon'$.
  Moreover, put ${\delta'\=\min\{1,\delta\}^{12}}$ and let~$C$ in~$\R_{>0}$ be given by Theorem~\ref{t:CM-fundamental} with~$\delta$ replaced by~$\delta'$.

  Let~$d$ be a fundamental discriminant in~$\pfd$ and~$f$ in ${\Z_{>0} \ssetminus p \Z}$.
  Lemma~\ref{l:Hecke-uniformity} implies that the function~$G$, defined by
  \begin{displaymath}
    G
    \=
    \frac{1}{\left( R_d^{-1} \ast \sigma_1 \right)(f)} \sum_{f_0 \in \Z_{>0}, f_0 \mid f} R_d^{-1} \left( \frac{f}{f_0} \right) T_{f_0} F,
  \end{displaymath}
  is constant on every ball of radius~$\delta'$.
  Moreover, combining~\eqref{eq:general-inverse-Zhang-formula} in Lemma~\ref{l:Zhang-general} with~$\wtf = 1$, \eqref{eq:deg-Hecke-op-divisor}, Theorem~\ref{t:CM-fundamental}, Corollary~\ref{c:Hecke-on-CM}, and the change of variables formula yields
  \begin{equation}
    \label{eq:116}
    \left| \int F \dd \odelta_{df^2} - \int F \dd \nu_{\pfd} \right|
    =
    \left| \int G \dd \odelta_{d} - \int G \dd \nu_{\pfd} \right|
    \le
    C \left( \sup_{\Lambda_{\pfd}} |G| \right) |d|^{- \frac{1}{28} + \varepsilon}.
  \end{equation}
  For every~$E$ in~$\Lambda_{\pfd}$, Corollary~\ref{c:Hecke-on-CM} and the choice of~$C_0$ and~$C_1$ yield
  \begin{multline*}
    |G(E)|
    \le
    \frac{1}{\left( R_d^{-1} \ast \sigma_1 \right) (f)} \sum_{f_0 \in \Z_{>0}, f_0 \mid f} \left| R_d^{-1} \left( \frac{f}{f_0} \right) \sigma_1(f_0) \right| \left( \sup_{\Lambda_\pfd} |F| \right)
    \\ \le
    C_0 C_1^2 \left( \sup_{\Lambda_\pfd} |F| \right) f^{-1 + \varepsilon'} \sum_{f_0 \in \Z_{>0}, f_0 \mid f} \left( \frac{f}{f_0} \right)^{\varepsilon'} f_0^{1 + \varepsilon'}
    \le
    C_0^2 C_1^2 \left( \sup_{\Lambda_\pfd} |F| \right) f^{3\varepsilon'}.
  \end{multline*}
  Together with~\eqref{eq:116} this implies the proposition with ${C' = C C_0^2 C_1^2}$.
\end{proof}

Let~$\pfd$ be a fundamental $p$\nobreakdash-adic discriminant and~$m$ in~$\Z_{>0}$, and put ${\pd \= \pfd p^{2m}}$.
So, ${\pd}$ is a $p$\nobreakdash-adic discriminant.
Define the Borel measure~$\nu_{\pd}$ on~$\Ell(\Cp)$, by
\begin{equation}
  \label{eq:117}
  \nu_{\pd}
  \=
  \begin{cases}
    \frac{1}{p^m} ( \t^m \bigm\vert_{A_{\pd}} )^* \nu_{\pfd}
    & \text{if~$\Qpfd$ is ramified over~$\Qp$};
    \\
    \frac{1}{p^{m - 1}(p + 1)} ( \t^m \bigm\vert_{A_{\pd}} )^* \nu_{\pfd}
    & \text{if~$\Qpfd$ is unramified over~$\Qp$}.
  \end{cases}
\end{equation}
It is a probability measure because~$\nu_{\pfd}$ is a probability measure and, for every~$m$ in~$\Z_{>0}$, the map~$\t|_{A_{\pfd p^{2m}}}$ is of degree~$p$, unless~$m = 1$ and~$\Qpfd$ is unramified over~$\Qp$ in which case the degree is~$p + 1$, see Lemma~\ref{l:sups-canonical-subgroup}.

\begin{proof}[Proof of Theorems~\ref{t:CM-symmetric} and~\ref{t:CM-broken-symmetry}]
  Denote by~$\pfd$ the fundamental $p$\nobreakdash-adic discriminant and~$m$ in~$\Z_{\ge 0}$ such that ${\pd = \pfd p^{2m}}$ (Lemma~\ref{l:p-adic-discriminants-Appendix}$(i)$).

  Suppose ${m = 0}$, so~$\pd = \pfd$.
  Theorem~\ref{t:CM-fundamental} yields the first assertion of Theorem~\ref{t:CM-symmetric}.
  The second assertion of Theorem~\ref{t:CM-symmetric} and Theorem~\ref{t:CM-broken-symmetry} are a direct consequence of Propositions~\ref{p:symmetry-breaking}, \ref{p:CM-prime-to-p-I}, and~\ref{p:CM-prime-to-p-II}.

  Suppose ${m \ge 1}$.
  The first assertion of Theorem~\ref{t:CM-symmetric} follows
  from the compactness of~$\Lambda_{\pfd}$, Theorem~\ref{t:formal-CM-formulae}$(ii)$, from the fact that~$A_{\pd}$ is an affinoid, and the analyticity of~$\t$ given by Lemma~\ref{l:canonical-analyticity}.
  Using Theorem~\ref{t:formal-CM-formulae}$(ii)$ again and ${\supp(\nu_{\pfd}) = \Lambda_{\pfd}}$, it follows that ${\supp(\nu_{\pd}) = \Lambda_{\pd}}$.
  The equidistribution statement in Theorem~\ref{t:CM-symmetric} for~$\pd$ follows from that for~$\pfd$, using Theorem~\ref{t:CM-from-canonical}, the change of variables formula, and the fact that the fundamental discriminant of every discriminant in~$\pd$ belongs to~$\pfd$ (Lemma~\ref{l:p-adic-discriminants}).
  To complete the proof of Theorem~\ref{t:CM-broken-symmetry}, note that the compact sets~$\Lambda_{\pd}^+$ and~$\Lambda_{\pd}^-$, defined by
  \begin{displaymath}
    \Lambda_{\pd}^+
    \=
    (\t^m \bigm\vert_{A_{\pd}})^{-1}(\Lambda_{\pfd}^+)
    \text{ and }
    \Lambda_{\pd}^-
    \=
    (\t^m \bigm\vert_{A_{\pd}})^{-1}(\Lambda_{\pfd}^-),
  \end{displaymath}
  form a partition of~$\Lambda_{\pd}$.
  Define the Borel probability measure~$\nu_{\pd}^+$ (resp.~$\nu_{\pd}^-$) by~\eqref{eq:117} with~$\nu_{\pfd}$ replaced by~$\nu_{\pfd}^+$ (resp.~$\nu_{\pfd}^-$).
  Then, the remaining assertions of Theorem~\ref{t:CM-broken-symmetry} for~$\pd$ follow from those for~$\pfd$, using Theorem~\ref{t:CM-from-canonical} and the change of variables formula.
\end{proof}

\appendix
\section{Quadratic field extensions of~$\Qp$ and $p$-adic discriminants}
\label{ss:Apendix-A}

This appendix proves Lemma~\ref{l:p-adic-discriminants} and gathers basic facts about quadratic field extensions of~$\Qp$ and $p$\nobreakdash-adic discriminants.
The proof of Lemma~\ref{l:p-adic-discriminants} relies on the notation and terminology in Section~\ref{ss:discriminants} and a couple of lemmas.

\begin{lemma}
  \label{l:p-adic-discriminants-Appendix}
  \hfill

  \begin{enumerate}
  \item
    For every $p$\nobreakdash-adic discriminant~$\pd$, there are unique fundamental $p$\nobreakdash-adic discriminant~$\pfd$ and~$m$ in~$\Z_{\ge 0}$ satisfying ${\pd = \pfd p^{2m}}$.
    Conversely, every set of this form is a $p$\nobreakdash-adic discriminant.
  \item
    For all fundamental $p$\nobreakdash-adic discriminant~$\pfd$ and~$m$ in~$\Z_{\ge 0}$, every $p$\nobreakdash-adic quadratic order of $p$\nobreakdash-adic discriminant~$\pfd p^{2m}$ is isomorphic to the $\Zp$\nobreakdash-order ${\Zp + p^m \cO_{\Qpfd}}$ in~$\Qpfd$.
    In particular, the $p$\nobreakdash-adic discriminant is a complete isomorphism invariant of $p$\nobreakdash-adic quadratic orders.
  \item
    The set of fundamental $p$\nobreakdash-adic discriminants equals
    \begin{equation}
      \label{eq:118}
      \left\{ \Zp^{\times} \ssetminus \Zp^2, p (\Zp^{\times})^2, p (\Zp^{\times} \ssetminus \Zp^2) \right\}
    \end{equation}
    if~$p$ is odd, and if~$p = 2$, then it equals
    \begin{equation}
      \label{eq:119}
      \left\{ -3 + 8 \Z_2, -4 + 32 \Z_2, 12 + 32 \Z_2, 8 + 64 \Z_2, -8 + 64 \Z_2, 24 + 64 \Z_2, -24 + 64 \Z_2 \right\}.
    \end{equation}
  \end{enumerate}
\end{lemma}

The proof of this lemma relies on the following one.
Denote by~$\Quap$ the set of quadratic field extensions of~$\Qp$ inside~$\Cp$.
Recall that~$\Q_{p^2}$ denotes the unique unramified extension of~$\Qp$ in~$\Quap$ and that, for each~$\Delta$ in~$\Qp$, the field~$\Qp(\sqrt{\Delta})$ is the unique element of~$\Quap$ containing a root of ${X^2 - \Delta}$.
For all quadratic field extension~$\cK$ of~$\Qp$ and~$x$ in~$\cK$, let~$\tr(x)$, $\nr(x)$, and~$\disc(x)$ be as in Section~\ref{s:preliminaries}.

\begin{lemma}
  \label{l:quadratic-extensions}
  If~$p$ is odd, then let~$A$ be an integer that is not a square modulo~$p$.
  \begin{enumerate}
  \item [$(i)$]
    Every quadratic field extension of~$\Qp$ is isomorphic to a unique element of~$\Quap$.
    Moreover,
    \begin{align}
      \Q_{p^2}
      & =
        \begin{cases}
          \Qp(\sqrt{A})
          & \text{if~$p$ is odd};
          \\
          \Q_2(\sqrt{-3})
          & \text{if~$p = 2$},
        \end{cases}
        \intertext{and}
        \Quap
      & =
        \begin{cases}
          \left\{ \Qp(\sqrt{A}), \Qp(\sqrt{p}), \Qp(\sqrt{Ap}) \right\}
          & \text{if~$p$ is odd};
          \\
          \left\{ \Q_2(\sqrt{d_0}) \: d_0 \in \{ -1, -2, -3, -5, -6, -10, -14 \} \right\}
          & \text{if~$p = 2$}.
        \end{cases}
    \end{align}
  \item[$(ii)$]
    Let~$d_0$ be in~$\{ A, p, Ap \}$ if~$p$ is odd and in~$\{ -1, -2, -3, -5, -6, -10, -14 \}$ if ${p = 2}$, and put ${\cQ \= \Qp(\sqrt{d_0})}$.
    Then,
    \begin{equation}
      \label{eq:120}
      \OQ
      =
      \begin{cases}
        \Z_2 \left[ \tfrac{1 + \sqrt{-3}}{2} \right]
        & \text{if~$p = 2$ and~$d_0 = -3$};
        \\
        \Zp \left[ \sqrt{d_0} \right]
        & \text{otherwise},
      \end{cases}
    \end{equation}
    the $p$\nobreakdash-adic discriminant of~$\OQ$ equals
    \begin{equation}
      \label{eq:121}
      \begin{cases}
        d_0 (\Zp^{\times})^2
        & \text{if~$p$ is odd, or~$p = 2$ and~$d_0 = -3$};
        \\
        4d_0 (\Zp^{\times})^2
        & \text{if~$p = 2$ and~$d_0 \neq -3$},
      \end{cases}
    \end{equation}
    and the subgroup~$\nr(\OQ^{\times})$ of~$\Zp^{\times}$ is given by
    \begin{equation}
      \label{eq:122}
      \nr(\OQ^{\times})
      =
      \begin{cases}
        \Zp^{\times}
        & \text{if~$\cQ = \Q_{p^2}$};
        \\
        (\Zp^{\times})^2
        & \text{if~$p$ is odd and~$d_0 = p$ or~$Ap$};
        \\
        1 + 4 \Z_2
        & \text{if~$p = 2$ and~$d_0 = -1$ or~$-5$};
        \\
        (1 + 8 \Z_2) \cup (3 + 8 \Z_2)
        & \text{if~$p = 2$ and~$d_0 = -2$ or~$-10$};
        \\
        (1 + 8 \Z_2) \cup (-1 + 8 \Z_2)
        & \text{if~$p = 2$ and~$d_0 = -6$ or~$-14$}.
      \end{cases}
    \end{equation}
    In particular, for every~$\Delta$ in the $p$\nobreakdash-adic discriminant of~$\OQ$,
    \begin{equation}
      \label{eq:123}
      \cQ
      =
      \Qp(\sqrt{\Delta})
      \text{ and }
      \OQ
      =
      \Zp \left[ \tfrac{\Delta + \sqrt{\Delta}}{2} \right],
    \end{equation}
    and the index of~$\nr(\OQ^{\times})$ in~$\Zp^{\times}$ is one if~$\cQ = \Q_{p^2}$ and two otherwise.
  \end{enumerate}
\end{lemma}

\begin{proof}
  Since~$\Cp$ contains an algebraic closure of~$\Qp$, every quadratic field extension of~$\Qp$ is isomorphic to one in~$\Quap$.
  Two distinct elements of~$\Quap$ can not be isomorphic since every quadratic extension of fields is normal.
  This proves that every quadratic field extension of~$\Qp$ is isomorphic to a unique element of~$\Quap$.
  The description of~$\Quap$ in item~$(i)$ follows from straightforward computations with the representatives of cosets in~$\Qp^{\times}/(\Qp^{\times})^2$, as found, \emph{e.g.}, in \cite[Section~3.3, Chapter~II]{Ser73e}.
  To prove the formula for~$\Q_{p^2}$, put $\cQ \= \Qp(\sqrt{A})$ (resp~$\Q_2(\sqrt{-3})$) if~$p$ is odd (resp. ${p = 2}$).
  Then, ${\cQ}$ is the splitting field of ${X^2 - A}$ (resp. ${X^2 + X + 1}$) over~$\Qp$.
  Since the reduction of this polynomial is irreducible over~$\Fp$, it follows~$\cQ$ is unramified over~$\Qp$.

  To prove~\eqref{eq:120} in item~$(ii)$, suppose first ${p = 2}$ and ${d_0 = -3}$, and put ${u \= \frac{1 + \sqrt{-3}}{2}}$.
  Then, $\tr(u) = \nr(u) = 1$, thus~$u$ belongs to~$\OQ$, and ${\Z_2[u] \subseteq \OQ}$.
  To prove the reverse inclusion, let~$\alpha$ and~$\beta$ in~$\Q_2$ be such that
  ${\alpha + \beta u}$ belongs to~$\OQ$ and put ${h \= \alpha + \beta u}$.
  Since~$\disc(h)$ belongs to~$\Z_2$ and ${\disc(h) = - 3 \beta^2}$, it follows that~$\beta$ belongs to~$\Z_2$.
  Combined with ${\alpha = h - \beta u}$, this implies that~$\alpha$ belongs to~$\Z_2$.
  This proves ${\OQ = \Z_2[u]}$.
  Suppose~$p$ is odd, or that~$p = 2$ and~$d_0 \neq -3$, and let~$\alpha$ and~$\beta$ in~$\Qp$ be such that ${h \= \alpha + \beta \sqrt{d_0}}$ belongs to~$\OQ$.
  Then
  \begin{displaymath}
    \tr(h) = 2 \alpha
    \text{ and }
    \disc(h) = d_0 (2 \beta)^2,
  \end{displaymath}
  so~$2 \alpha$ and~$d_0 (2 \beta)^2$ belong to~$\Zp$.
  Since~$d_0$ is an integer and ${p^2 \nmid d_0}$, this implies that~$2\beta$ belongs to~$\Zp$.
  If~$p$ is odd, then it follows that~$\alpha$ and~$\beta$ belong to~$\Zp$, and ${\OQ = \Zp \left[ \sqrt{d_0} \right]}$.
  If ${p = 2}$, then ${d_0 \neq -3}$ by hypotheses, so
  \begin{displaymath}
    - d_0 \equiv 1, 2 \mod 4
    \text{ and }
    (2\alpha)^2 - d_0 (2\beta)^2
    =
    4 \nr(h)
    \equiv 0 \mod 4,
  \end{displaymath}
  and thus~$\alpha$ and~$\beta$ belong to~$\Z_2$.
  This proves ${\OQ = \Z_2 \left[ \sqrt{d_0} \right]}$ and completes the proof of~\eqref{eq:120}.

  A direct computation using~\eqref{eq:120} yields that the $p$\nobreakdash-adic discriminant of~$\OQ$ equals
  \begin{displaymath}
    \begin{cases}
      d_0 (\Zp^{\times})^2
      & \text{if~$p = 2$ and~$d_0 = -3$};
      \\
      4 d_0 (\Zp^{\times})^2
      & \text{otherwise}.
    \end{cases}
  \end{displaymath}
  This equals~\eqref{eq:121} when ${p = 2}$.
  When~$p$ is odd, this also equals~\eqref{eq:121} because~$4$ belongs to~$(\Zp^{\times})^2$.

  To prove~\eqref{eq:122}, suppose first~$p$ is odd and~$d_0 = A$.
  Since ${\Zp^{\times} \ssetminus (\Zp^{\times})^2 = A (\Zp^{\times})^2}$ and the norm map from the residue field of~$\cQ$ to~$\Fp$ is surjective, ${\nr(\OQ^{\times}) = \Zp^{\times}}$.
  Suppose~$p$ is odd and ${d_0 \neq A}$, or ${p = 2}$.
  If~$p$ divides~$d_0$, then~$\OQ^{\times} = \Zp^{\times} + \sqrt{d_0} \Zp$ by~\eqref{eq:120} and
  \begin{displaymath}
    \nr(\OQ^{\times})
    =
    \left\{ 1 - d_0 \ell^2 \: \ell \in \Zp \right\} (\Zp^{\times})^2
    =
    \begin{cases}
      (\Zp^{\times})^2
      & \text{if~$p$ is odd};
      \\
      (1 + 8 \Z_2) \cup (1 - d_0 + 8 \Z_2)
      & \text{if~$p = 2$}.
    \end{cases}
  \end{displaymath}
  It remains to consider the case where ${p = 2}$ and~$d_0 = -1$, $-3$, or~$-5$.
  Since ${(\Z_2^{\times})^2 = 1 + 8 \Z_2}$, when ${d_0 = -3}$,
  \begin{displaymath}
    3
    =
    \nr(\sqrt{-3}),
    7
    =
    \nr(2 + \sqrt{-3}),
    \text{ and }
    13
    =
    \nr(1 + 2 \sqrt{-3}).
  \end{displaymath}
  This implies ${\nr(\OQ^{\times}) = \Z_2^{\times}}$.
  If~$p = 2$ and~$d_0 = -1$ or~$-5$, then ${\OQ^{\times} = \Z_2^{\times} + (1 + \sqrt{d_0}) \Z_2}$ and
  \begin{displaymath}
    \nr(\OQ^{\times})
    =
    \left\{ 1 + 2 \ell + (1 - d_0) \ell^2 \: \ell \in \Z_2 \right\} (\Z_2^{\times})^2
    =
    1 + 4 \Z_2.
    \qedhere
  \end{displaymath}
\end{proof}

\begin{proof}[Proof of Lemma~\ref{l:p-adic-discriminants-Appendix}]
  To prove the uniqueness statement in item~$(i)$, let~$\pfd$ and~$\pfd'$ be fundamental $p$\nobreakdash-adic discriminants and~$m$ and~$m'$ in~$\Z_{\ge 0}$ such that ${\pfd p^{2m} = \pfd' p^{2m'}}$.
  Then, ${\Qpfd = \Qp(\sqrt{\pfd'})}$ and~$\pfd$ and~$\pfd'$ are both equal to the $p$\nobreakdash-adic discriminant of~$\cO_{\Qpfd}$.
  It follows that ${m = m'}$.
  To prove existence, let~$\pd$ be a $p$\nobreakdash-adic discriminant and~$\cO$ a $p$\nobreakdash-adic quadratic order of $p$\nobreakdash-adic discriminant~$\pd$.
  Then, the field of fractions of~$\cO$ has the same discriminant as~$\Qpd$, and it is therefore isomorphic to it.
  So, there is~$m$ in~$\Z_{\ge 0}$ such that~$\cO$ is isomorphic to the $\Zp$-order ${\Zp + p^m \cO_{\Qpd}}$ in~$\Qpd$.
  Thus, if~$\pfd$ denotes the $p$\nobreakdash-adic discriminant of~$\cO_{\Qpd}$, then $\pfd$ is a fundamental $p$\nobreakdash-adic discriminant and ${\pd = \pfd p^{2m}}$.

  To prove item~$(ii)$, let~$\pfd$ be a fundamental $p$\nobreakdash-adic discriminant, $m$ in~$\Z_{\ge 0}$, and~$\cO$ a $p$\nobreakdash-adic quadratic order of $p$\nobreakdash-adic discriminant~$\pfd p^{2m}$.
  Then, ${\Qp(\sqrt{\pfd p^{2m}}) = \Qpfd}$ and, as shown above with ${\pd = \pfd p^{2m}}$, there is~$m'$ in~$\Z_{\ge 0}$ such that~$\cO$ is isomorphic to the $\Zp$-order ${\Zp + p^{m'} \cO_{\Qpfd}}$ in~$\Qpfd$.
  Since the $p$\nobreakdash-adic discriminant of the latter is equal to~$\pfd p^{2m'}$, ${\pfd p^{2m'} = \pfd p^{2m}}$ and ${m' = m}$ by item~$(i)$.

  When~$p$ is odd, item~$(iii)$ follows from the combination of~\eqref{eq:121} in Lemma~\ref{l:quadratic-extensions} and, for every integer~$A$ that is not a square modulo~$p$, the equality ${A (\Zp^{\times})^2 = \Zp^{\times} \ssetminus \Zp^2}$.
  When ${p = 2}$, it is a direct consequence of~\eqref{eq:121} in Lemma~\ref{l:quadratic-extensions} and ${(\Z_2^{\times})^2 = 1 + 8 \Z_2}$.
\end{proof}

\begin{proof}[Proof of Lemma~\ref{l:p-adic-discriminants}]
  In view of Lemma~\ref{l:p-adic-discriminants-Appendix}$(i)$, to prove the first assertion it suffices to prove that a fundamental discriminant~$d$ belongs to a fundamental $p$\nobreakdash-adic discriminant if and only if it is $p$-supersingular.
  If~$p$ is odd, then by Lemma~\ref{l:p-adic-discriminants-Appendix}$(iii)$ the union of all fundamental $p$\nobreakdash-adic discriminants equals ${\Zp\ssetminus \left( p^2 \Zp \cup \Zp^2 \right)}$.
  So, $d$ belongs to a fundamental $p$\nobreakdash-adic discriminant if and only if ${\left( \frac{d}{p} \right) \neq 1}$.
  As remarked in Section~\ref{ss:discriminants}, the latter is equivalent to~$d$ being $p$\nobreakdash-supersingular.
  If~$p = 2$, then~\eqref{eq:8} yields
  \begin{displaymath}
    d \equiv 1 \mod 4
    \text{ or }
    d \equiv - 4, 8 \mod 16.
  \end{displaymath}
  Together with Lemma~\ref{l:p-adic-discriminants-Appendix}$(iii)$, this implies that~$d$ belongs to a fundamental $2$\nobreakdash-adic discriminant if and only if ${d \not \equiv 1 \mod 8}$.
  As remarked in Section~\ref{ss:discriminants}, the latter is equivalent to~$d$ being $2$\nobreakdash-supersingular.
  This completes the proof the first assertion.

  In view of Lemma~\ref{l:p-adic-discriminants-Appendix}$(i)$, to prove the second assertion it suffices to restrict to the case where the $p$\nobreakdash-adic discriminant~$\pd$ is fundamental.
  Put ${\pfd \= \pd}$, and let~$\Delta$ be in~$\pfd$ and~$r$ an integer satisfying ${r \ge 6}$.
  By Lemma~\ref{l:p-adic-discriminants-Appendix}$(iii)$, every~$\Delta'$ in~$\Zp$ satisfying ${\ord_p(\Delta - \Delta') \ge r}$ belongs to~$\pfd$.
  Suppose first~$p$ is odd, so~$\Delta$ belongs to~$\Zp^{\times}$ or~$p \Zp^{\times}$ by Lemma~\ref{l:p-adic-discriminants-Appendix}$(iii)$.
  By Dirichlet's theorem on prime numbers in arithmetic progressions, there is a prime number~$p'$ satisfying
  \begin{displaymath}
    p' \equiv -1 \mod 4
    \text{ and }
    \ord_p(\Delta + p') \ge r
  \end{displaymath}
  in the former case, and
  \begin{displaymath}
    p' \equiv -p \mod 4
    \text{ and }
    \ord_p(\Delta/p + p') \ge r
  \end{displaymath}
  in the latter.
  Putting~$d \= - p'$ in the former case and~$d \= - p p'$ in the latter, it follows that~$d$ is a fundamental discriminant and ${\ord_p(\Delta - d) \ge r}$.
  In particular, $d$ belongs to~$\pfd$.
  It remains to consider the case where ${p = 2}$.
  Lemma~\ref{l:p-adic-discriminants-Appendix}$(iii)$ implies that~$\Delta$ belongs to either
  \begin{displaymath}
    -3 + 8 \Z_2,
    -4 + 16 \Z_2,
    \text{ or }
    8 + 16 \Z_2.
  \end{displaymath}
  Let~$p'$ be a prime number satisfying
  \begin{displaymath}
    \ord_2(\Delta + p') \ge r, \ord_2(\Delta/4 + p') \ge r, \text{ or } \ord_2(\Delta/8 + p') \ge r,
  \end{displaymath}
  and put~$d \= -p', -4p'$, or $-8p'$, respectively.
  Then, $d$ is a fundamental discriminant satisfying~$\ord_2(\Delta - d) \ge r$ and it thus belongs to~$\pfd$.
\end{proof}

\bibliographystyle{alpha}

\end{document}